\documentclass{amsart}
\usepackage[margin=1in]{geometry}
\usepackage{todonotes}
\usepackage{xspace}
\usepackage{url}
\usepackage{amsthm}
\usepackage{amssymb}
\usepackage[normalem]{ulem}
\usepackage{pgfplots}
\usepackage{hyperref}

\newtheorem{thm}{Theorem}
\newtheorem{cor}{Corollary}

\newtheorem{definition}{Definition}

\newtheorem{remark}{Remark}

\usepackage{caption}
\usepackage{subcaption}

\newcommand{\R}{\mathbb{R}}
\newcommand{\D}{\mathcal{D}}
\newcommand{\F}{\mathcal{F}}
\newcommand{\PkG}{\mathcal{P}_{k,\epsilon}(G)}

\newcommand{\abar}{\bar{a}}
\newcommand{\bbar}{\bar{b}}
\newcommand{\sdot}{\! \cdot \!}

\newcommand{\Fre}{Fr\'{e}chet\xspace}
\newcommand{\edge}[2]{\overline{#1#2}}
\newcommand{\dom}{\mathop{\rm dom}}
\newcommand{\im}{\mathop{\rm im}}

\newcommand{\pdsubfig}[1]{%
\resizebox{\columnwidth}{!}{%
\begin{tikzpicture}
\begin{axis}[
xmin=0, xmax=1, ymin=0, ymax=1.1, clip=false,
xtick={0.25,0.5,0.75}, ytick={0.25,0.5,0.75},
   enlargelimits=false, axis equal image, axis on top
]
\addplot graphics [xmin=0, xmax=1,ymin=0, ymax=1.1] {#1}; 
\addplot[gray, dashed] coordinates{(0.5,0) (0.5,1)};
\addplot[gray, dashed] coordinates{(0,0) (1,1)};
\node at (axis cs: 0.03,1.02) {$\infty$};
\end{axis}
\end{tikzpicture}%
}
}
\newcommand{\pdtwosubfig}[2]{%
\resizebox{\columnwidth}{!}{%
\begin{tikzpicture}
\begin{axis}[
xmin=0, xmax=1, ymin=0, ymax=1.1, clip=false,
xtick={0.25,0.5,0.75}, ytick={0.25,0.5,0.75},
   enlargelimits=false, axis equal image, axis on top
]
\addplot[fill opacity=0.85, draw opacity=0.85] graphics [xmin=0, xmax=1,ymin=0, ymax=1.1] {#1}; 
\addplot[fill opacity=0.5, draw opacity=0.5] graphics [xmin=0, xmax=1,ymin=0, ymax=1.1] {#2}; 
\addplot[gray, dashed] coordinates{(0.5,0) (0.5,1)};
\addplot[gray, dashed] coordinates{(0,0) (1,1)};
\node at (axis cs: 0.03,1.02) {$\infty$};
\end{axis}
\end{tikzpicture}%
}%
}

\title{The (homological) persistence of gerrymandering}
\date{}
\author{Moon Duchin, Tom Needham and Thomas Weighill}
\begin{document}

\begin{abstract}
We apply persistent homology, the dominant tool from the field of topological data analysis, to study electoral redistricting. Our method combines the geographic information from a political districting plan with election data to produce a persistence diagram. We are then able to visualize and analyze large ensembles of computer-generated districting plans of the type commonly used in modern redistricting research (and court challenges).  We set out three applications:  zoning a state at each scale of districting, comparing elections, and seeking signals of gerrymandering.  Our case studies focus on redistricting in Pennsylvania and North Carolina, two states whose legal challenges to enacted plans have raised considerable public interest in the last few years.

To address the question of robustness of the persistence diagrams to perturbations in vote data and in district boundaries, we translate the classical stability theorem of Cohen--Steiner et al.\ into our setting and find that it can be phrased in a manner that is easy to interpret.  We accompany the theoretical bound with an empirical demonstration to illustrate diagram stability in practice.
\end{abstract}

\maketitle


\section{Introduction}

In this paper, we bring the techniques of topological data analysis to bear on electoral redistricting.  We begin by introducing the key notions we will need below.

\subsection{Redistricting}

In electoral politics, {\em redistricting} is the process of drawing new boundary lines for electoral regions called {\em districts}, in which an election will be conducted to select one or more people as representatives for a governing body.  For instance, the United States Congress has a House of Representatives that has had 435 seats for over a hundred years.  These seats are apportioned to the states after every decennial Census---for instance, the 2010 Census left Pennsylvania with 18 seats and North Carolina with 13, in rough proportion to their population---and then it is left to the various states to partition their territory into geographically-delimited districts that will elect a single member each.  In most cases, it is the state legislature that controls the process of drawing not only Congressional boundaries, but also the boundaries for their own legislative districts.  In the U.S. system, which is dominated by two major parties, this opens the door to {\em partisan gerrymandering}, where the lines are carefully arranged to maximize the seats secured by their side, given the anticipated pattern of voting.  For instance, the last few years have seen successful legal challenges in both Pennsylvania and North Carolina, where enacted plans  produced by Republican-controlled legislatures were found to be Republican-favoring ``gerrymanders"---in other words, to be so favorable to the party in charge that they impermissibly subordinate neutral criteria to a partisan agenda.  There are many other kinds of gerrymandering, including agendas with racial components, or to favor or target particular incumbents, but the present paper focuses on developing mathematical tools to understand the interplay of district lines and party outcomes.

\subsection{Persistent homology}

Persistent homology is the most widely recognized tool used in the burgeoning field of topological data analysis (TDA), where  concepts from algebraic topology are used to simplify, summarize and compare complex datasets. Methods from topological data analysis have proven successful in neuroscience \cite{dabaghian2012topological,bendich2016persistent}, medical diagnostics \cite{nicolau2011topology,Crawford_2019}, and machine learning \cite{adams2017persistence,carriere2019perslay,gabrielsson2019topology,hofer2019learning}, among many other applications---these references are just a selection of many dozens of examples in each of these areas.
However, applications to demographic or electoral geospatial data seem to be comparatively limited; we only know of a few examples.  In  \cite{bajardi2015unveiling}, Bajardi et al.\ use mobile phone data to study patterns of immigration.  
In \cite{stolz2016topological}, Stoltz et al.\ use TDA to study Brexit voting patterns, finding for instance that Switzerland is surrounded by countries with a different voting tendency. In \cite{banman2018mind}, Banman and Ziegelmeier apply persistent homology to global health and wealth statistics to uncover subtle developmental disparities between geographically close countries. Finally, in the example closest to our own approach, Feng and Porter introduce in \cite{feng2019persistent} a method for producing a topological signature from districting and vote data which is similar to the one we introduce below. The applications considered in \cite{feng2019persistent} are quite different from those considered here: Feng and Porter were primarily interested in studying the vote patterns at the precinct level for fixed regions using their topological signatures, whereas our goal is to compare districting plans.
We will do that by studying topological features across a collections or ``ensembles" containing thousands of alternative  plans.

\begin{figure}[ht]
\centering
\includegraphics[width=4in]{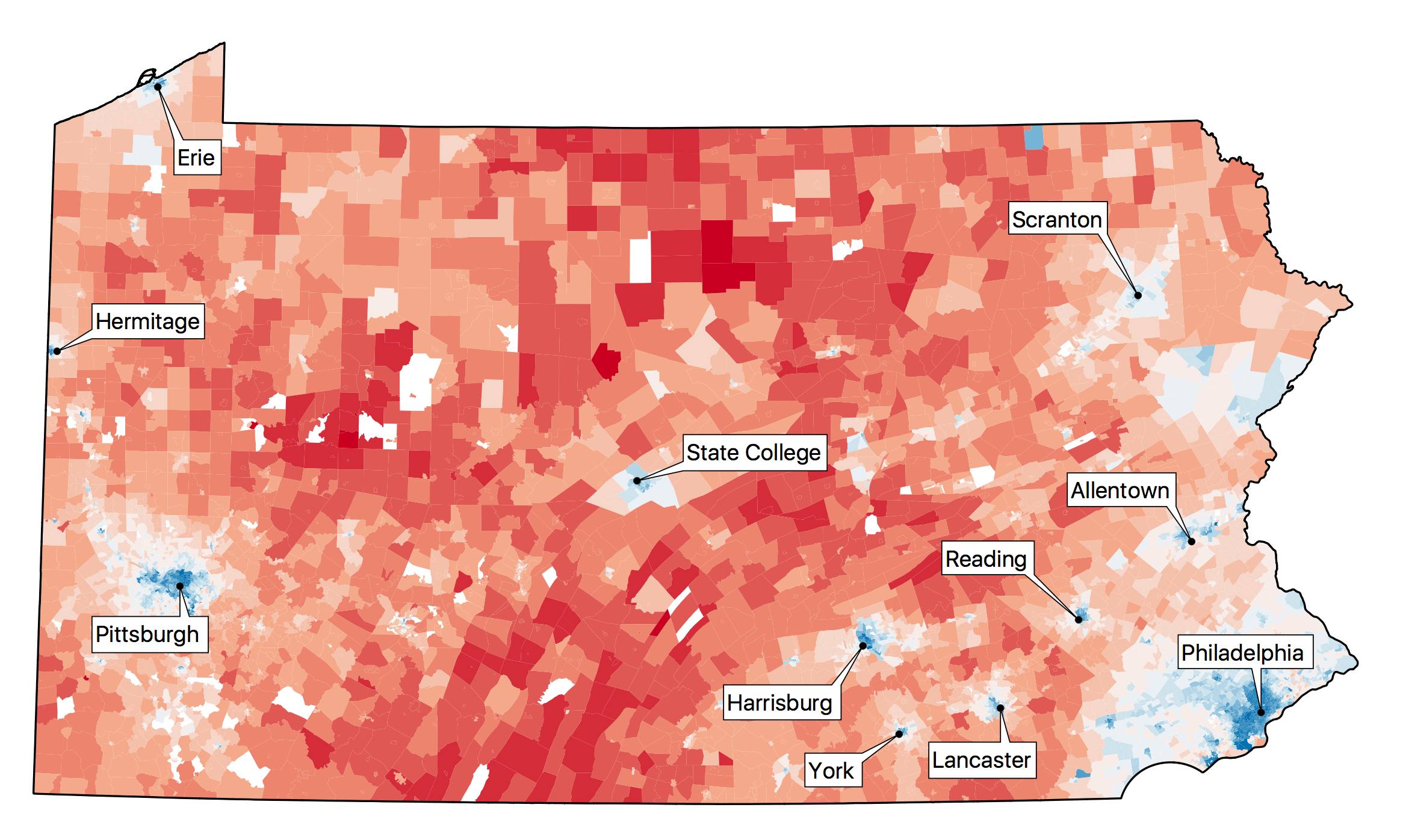}

\begin{subfigure}{0.48\textwidth}
\centering
\includegraphics[width=\textwidth]{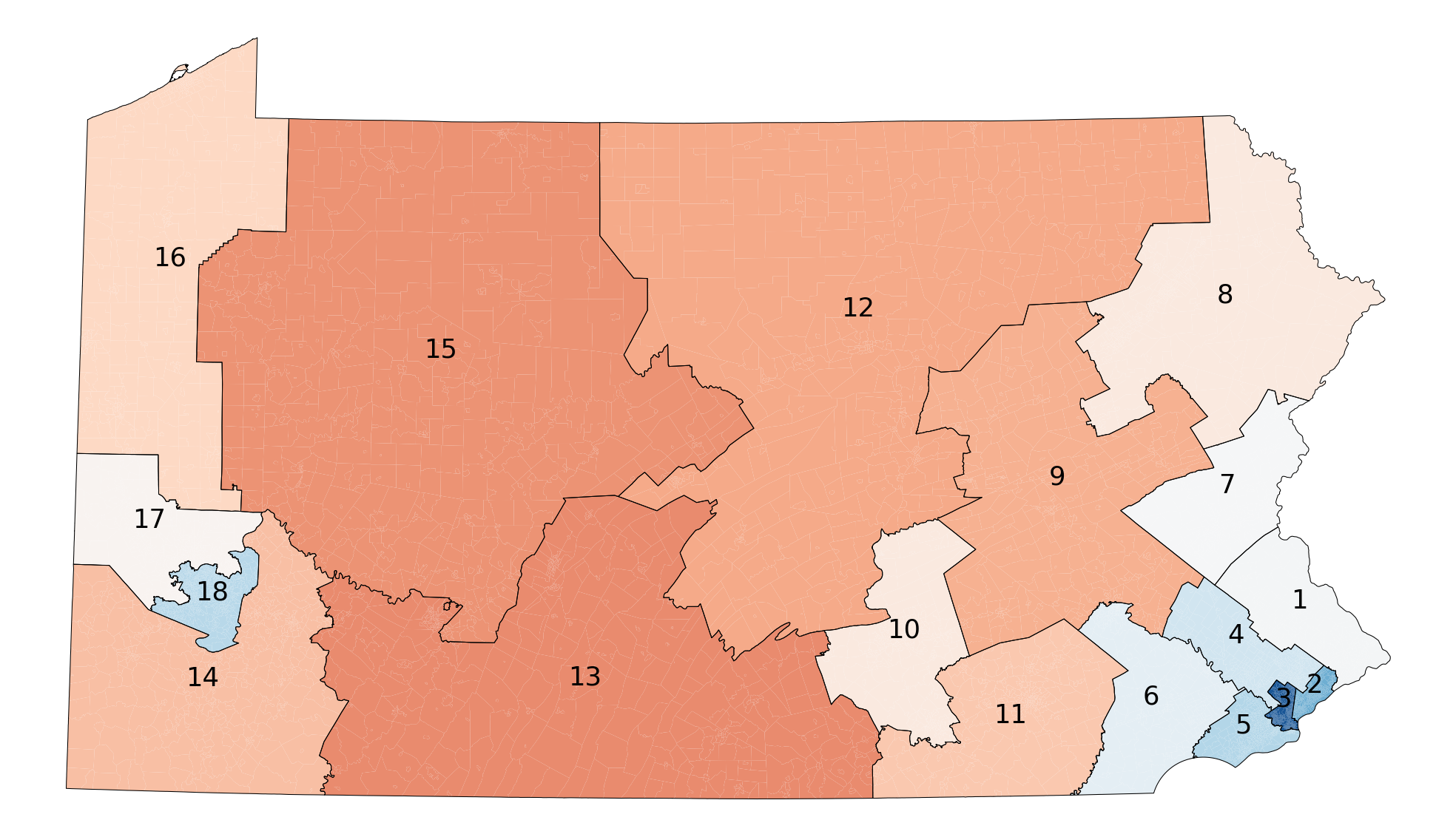}
\caption*{Congressional plan (enacted 2018)}
\end{subfigure}
\begin{subfigure}{0.48\textwidth}
\centering
\includegraphics[width=\textwidth]{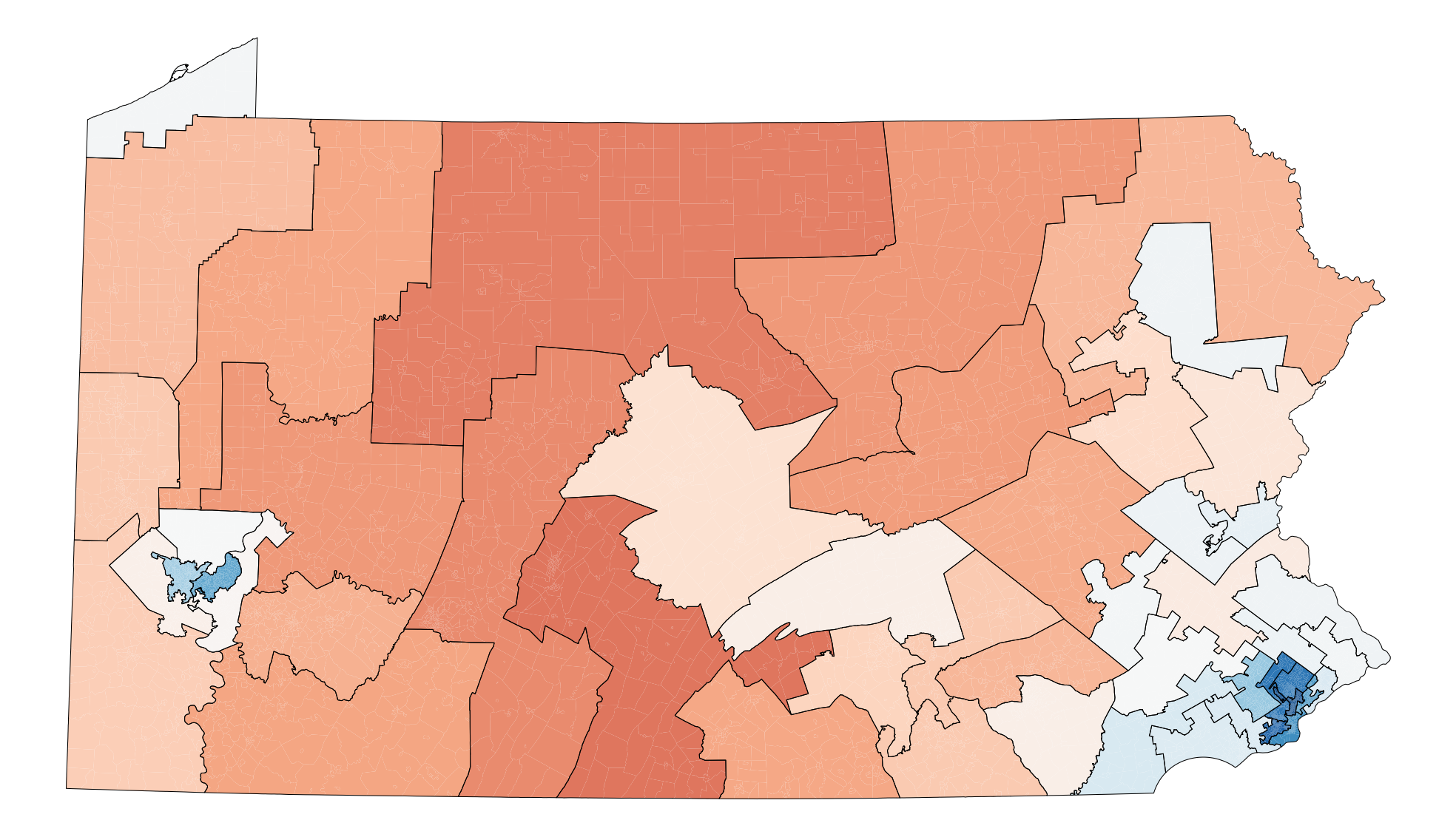}
\caption*{State Senate plan (enacted 2011)}
\end{subfigure}

\caption{This paper studies the  aggregation of votes from geographic units to districts, which is intractable to study by complete enumeration of possibilities.  We develop persistent homology techniques to summarize the districting problem.  
This figure shows a precinct map of Pennsylvania with key cities identified.  The precincts are colored by 
voting results from the 2016 Presidential election (blue for Democratic, red for Republican).
Current districts in Pennsylvania are shown with the same coloring.  Note that the Congressional districts are balanced to within  one-person population deviation (by 2010 Census count).  }
\label{fig:PAbasic}
\end{figure}

\subsection{Studying districting ensembles with persistent homology}

Suppose we want to evaluate a districting plan, either proposed or enacted, to label  it as  either a partisan gerrymander or not. Some authors have proposed universal metrics for this purpose, which do not take the specific state, time period,  or level of redistricting into account: the same score, denominated in the same units, is supposed to flag gerrymanders for Iowa's Congressional seats in 2020 or for Mississippi state House in 1960.  
 Examples which are widely cited include the efficiency gap \cite{stephanopoulos2015partisan}, and the mean-median score \cite{mcdonald2015unfair}, which belongs to a family of ``partisan symmetry" metrics  \cite{grofman2007future}. 
 
On the other hand, there has been growing consensus around the idea that evaluating a map requires comparing it to a range of alternatives that hold the geography, demography, and redistricting rules constant. In the last several decades, the use of ``demonstration plans" has been rising in court cases:  in addition to the map under consideration, a court also views plans submitted by various advocacy groups and political actors. While traditionally all of these benchmark plans are drawn ``by hand,'' the last ten years has seen the adoption of computer algorithms in court cases, where large generated ensembles of possible districting plans are used as a backdrop of comparison for a  proposed or enacted plan.  The use of algorithmic {\em ensemble methods} has already had a major impact on a number of recent court cases, notably Pennsylvania and North Carolina, where ensemble methods were cited by courts in the successful invalidation of the state's enacted maps as partisan gerrymanders \cite{mattingly-report, pegden2017, duchinpa}. 

Algorithmic ensembles can easily boost the number of alternatives to consider into the hundreds of thousands, far beyond the point where each map can be carefully analyzed with the human eye. For this reason, {\em summary statistics} are computed for the maps in the ensemble, then compared to the plan under evaluation. For example, one might choose one or several  recent elections in order to fix a background pattern of voting, then count the number of districts won by each party as the district lines vary. If the map under evaluation gives the controlling party the majority in an outlying number of districts---in other words, if the plan secures more seats for Party A than the vast majority of alternatives---this could reveal that the map was made with impermissible partisan intent and/or has unreasonably skewed effects.  

We can therefore say that modern research into redistricting forces us to be able to understand the essential qualities of large collections of territorial partitions, preferably in a way which captures both geography and relevant summary statistics, such as those connected to party or race. 
The method we propose in this paper is illustrated with respect to partisan statistics, but can be applied more generally.  We combine and summarize the adjacency and vote share information from a districting plan into a {\em dual graph} (see \S\ref{sec:setup}), decorated by vote shares under some given electoral pattern.  We can now ask questions like: how many connected components of support for Party A are there? 
A more sophisticated version of this question, which considers multiple thresholds at once, is answered by persistent homology, the mathematical theory underlying our method.

In this paper we apply methods from TDA to summarize ensembles of districting plans generated via Markov Chain Monte Carlo (MCMC) sampling. There has been some prior work using TDA to study similar ensembles of generated data, such as simulated polymer conformations \cite{emmett2016multiscale} or simulated neuronal spike trains \cite{dabaghian2012topological, chowdhury2018importance}. To our knowledge, this is the first geospatial application.

\subsection{Broad goals and novel contributions}
Before we proceed, we highlight an important guiding principle in this paper:  namely, we will measure the success of our topological data analysis techniques  by their ability to uncover information which cannot be easily observed from maps ``with one's eyes." Consider Figure \ref{fig:PAbasic}, for example, which shows the distribution of votes  in the 2016 Presidential election in Pennsylvania (following the contemporary U.S. color convention of red for Republican, blue for Democratic). No advanced data analysis techniques are required to observe that Philadelphia and Pittsburgh are far bluer than the rural north of the state.  What is not clear, however, is where the blue {\em districts} will fall in a typical redistricting plan. This is partly because the blue clusters in the map have widely different population sizes. A small Democratic-leaning city like State College 
($\approx$50,000 people) will have little impact on a 700,000-person Congressional district, but may dominate the voting in a smaller 63,000-person state House district if the boundaries fall just so. 

The main contributions of this paper are
\begin{itemize}
\item {\bf Scale and zoning.} We propose  TDA  as a tool for studying nuanced scale effects in partition problems.  In particular, persistence diagrams let us read off a {\em zoning}  of the state into redistricting-relevant regions at each districting scale.
\item {\bf Comparing elections.} Fixing a scale, we can study the difference between two vote patterns in natural units.
\item {\bf Signals of gerrymandering.} Finally, we consider the use of TDA to study differences in partisan-biased ensembles of plans.
\end{itemize}

These techniques are applied to Pennsylvania and North Carolina, where in each case we have several examples of human-made plans for which we have a ground-truth label of {\em gerrymandered} or {\em not gerrymandered}.  We find that each of our proposals is able to detect 
phenomena in those states that can be corroborated by secondary analysis.  These findings pass the ``can't-be-seen-with-your-eyes" test, 
and we hope that they will generate interest in the wider community of TDA practitioners.

\subsection{Outline of the paper}  In Section \ref{sec:mathematical_framework}, we define districting plans in the language of graph theory and briefly describe the methods that generate ensembles of plans. The persistent homology framework for studying districting plans is laid out in detail; to keep the exposition self-contained, we describe the relevant tools from TDA with a narrow focus on the application at hand. Section \ref{sec:applications} presents the three principal applications described above. We use the persistent homology formalism to compare scales of redistricting, to compare data from different elections on the same districting plan ensembles, and to determine distinguishing features from ensembles of partisan-tilted districting plans.
In Sections~\ref{sec:casestudyPA}-\ref{sec:casestudyNC}, we apply our methods to ensembles of plans in Pennsylvania and North Carolina. 
Finally, the stability and robustness of the TDA signatures is explored in Section \ref{sec:stability}. Drawing on a standard result from  \cite{cohen2007stability} on the stability of persistence diagrams to perturbations in their input data, we derive a guaranteed upper bound on  our districting plan representations when perturbing electoral and geographical data. These theoretical results are accompanied by data that raises our confidence in the robustness of our analysis.

\section{Background and mathematical framework}\label{sec:mathematical_framework}
In this section we describe the details of our topology-based approach to studying districting plans. We begin by precisely defining the notion of a districting plan and briefly describing our method for generating ensembles of plans. Next we apply the \emph{persistent homology} pipeline to certain filtered graphs associated to the plans. To keep the description accessible, we opt to restrict our exposition of persistent homology to our specific application rather than appealing to any general theory. We use \cite{edelsbrunner2010computational,carlsson2014topological} as general references for persistent homology.

\subsection{Redistricting as a graph partition problem}\label{subsec:redistricting_as_graph_partition}\label{sec:setup}

A \emph{districting plan} is a partition of a state into geographical regions meeting certain criteria, which vary by state. The building blocks are geographic units such as census blocks or election precincts. 
A very brief summary of key redistricting rules and principles is needed to set up the problem.  The number of Congressional districts in a state is deterministically related to the Census population by a formula that has been in place for nearly sixty years; the new apportionments will be announced in 2021.  The number of districts in the state legislative houses is typically fixed in  state constitutions and is rarely changed.
A criterion for a districting plan which is common to most states (and essentially active in all) is that each district must be ``contiguous," i.e., the induced subgraphs on the pieces of the partition should be connected.  And it is a universal requirement that the districts be 
population-balanced:  each district must have equal population, up to a small tolerance for deviation.  
For our purposes we will fix a deviation limit $\epsilon$, so that each of the districts should be no less than $(1-\epsilon)$ 
and at most $(1+\epsilon)$ times ideal size, which is  total Census count divided by $k$.

After fixing the choice of units, we are able to define a districting plan more formally as a graph partition, with notation following \cite{duchinAAAS, abrishami2019geometry, deford2019recombination}.
Given a tiling of a state by geographic units that overlap only along boundaries, let $G = (V,E)$ be the \emph{unit dual graph}---for instance, we can form a precinct dual graph or a census block dual graph. This is defined by having each vertex in $V$ correspond to a geographic unit of the state, then connecting two vertices by an edge if and only if their units are geographically adjacent (i.e., share a boundary of positive length). 
We denote vertices in $V$ by $v$ or $w$ and an edge in $E$ with endpoints $v$ and $w$ is denoted $\edge vw$. 
For a subset $V'\subseteq V$, recall that the {\em induced subgraph} on $V'$ is the graph $G' = (V',E')$ containing all the edges of $E$ that have both endpoints in $V'$.  We adopt the slight abuse of notation to refer by $G'$ to either the induced subgraph or its vertex set $V'$, as needed.  We will frequently make use of functions $p:V\to \R$, such as a population count for every unit.
We will write $p(G')$ to refer to the sum $\sum_{v\in V'} p(v)$.

\begin{definition}[Districting plans as balanced connected partitions]
Given  a graph $G=(V,E)$, a number of districts $k\ge 2$, a weighting function $p:V\to \R$,
and a balance threshold $\epsilon$, let $\PkG$ be the set of 
partitions $P = \{P_1,\dots,P_k\}$ where $P_i = (V_i,E_i)$ are induced graphs on a vertex partition $V=V_1\sqcup \dots\sqcup V_k$ with each $P_i$ connected and satisfying the balance condition
$$
(1-\epsilon) \cdot \frac{p(G)}k \le p(P_i) \le (1+\epsilon)\cdot  \frac{p(G)}k.
$$
Each $P\in \PkG$ is called a {\em districting plan}, and each $P_i$ is called a {\em district}.
\end{definition}

As alluded to above, districting plans are typically subjected to other tests of quality. One prominent example is ``compactness," which requires that the shapes of districts look reasonable on a map. Other criteria prefer plans that hew to county and municipal boundaries, that respect ``communities of interest," and so on.  And finally, but by no means insignificantly, the whole country is subject to the Voting Rights Act of 1965 (VRA), a rule that requires adequate access to representation for certain kinds of demographic minorities.  
Generating a districting plan meeting all relevant validity criteria is a complicated and delicate process, and in this paper we will only 
operationalize population balance, contiguity, and compactness.\footnote{For work that layers in more complicated districting criteria,
see for instance \cite{mgggvacriteria}.}

Perhaps the most concise summary of the geographic layout of a redistricting plan is the data of which districts are next to which. This information is naturally captured in a graph structure that is formed by aggregation from the dual graph of the state's geographic units: each vertex of the graph now corresponds to a district, and edges connect nodes from adjacent districts (see Figures \ref{fig:plan_to_diagram} and \ref{fig:NCperturb}). This construction is referred to as a \emph{district dual graph} for the plan.

\begin{figure}[ht]
\centering
\includegraphics[width = \textwidth, height=1.3in]{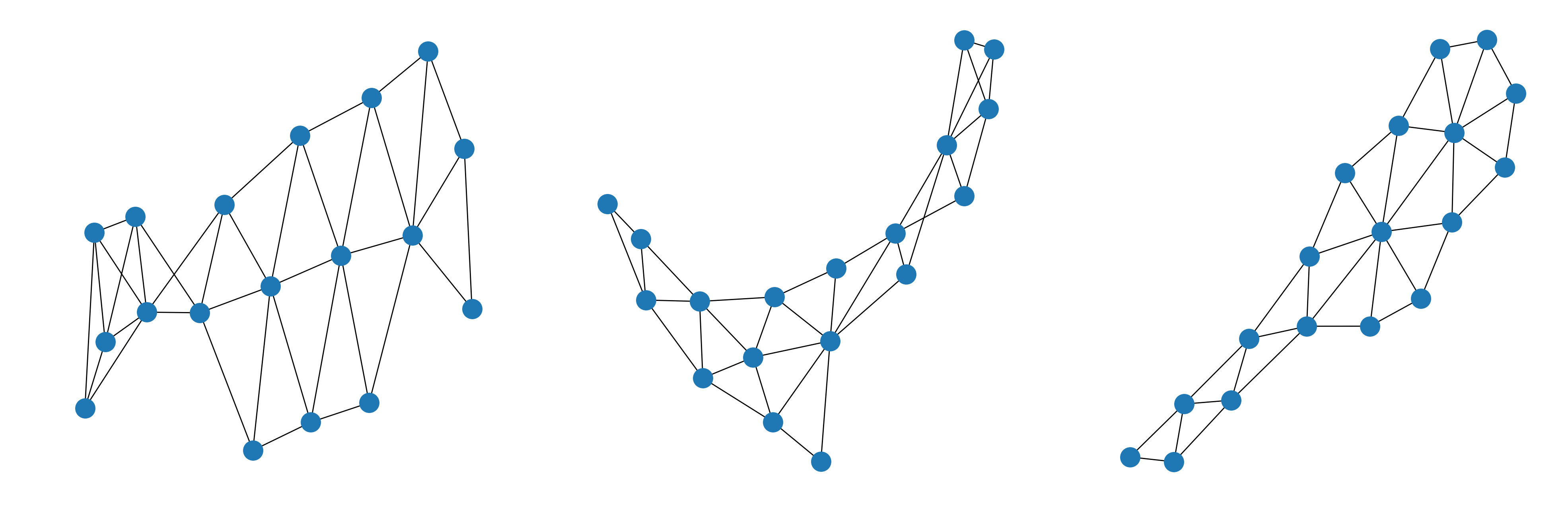}
\caption{Abstract graph representations of three 18-district  plans for Pennsylvania.} 
\label{fig:18DistrictGraphExamples}
\end{figure}

One might hope that even within a large ensemble of plans, the types of graphs one obtains are relatively restricted, thereby uncovering a common feature of the plans in the ensemble. It turns out in practice that there is a huge variety in the types of graphs in even a modest ensemble of 1000 plans. In a generated ensemble of 18-district plans for Pennsylvania (details of ensemble generation are provided in Section \ref{sec:mathematical_framework}), 991 distinct graph isomorphism classes are represented among the 1000 dual graphs. Examples of some of these graphs are shown in Figure \ref{fig:18DistrictGraphExamples} and some basic graph statistics of the ensemble are shown in Figure \ref{fig:graph_statistics}.  As the number of districts in the plan increases, so does the variety of graphs: the dual graphs constructed from our generated ensembles of 1000 50-district plans are typically all unique up to isomorphism. This suggests that passing from districting plans to dual graphs is not enough of a reduction in representation complexity to make direct analysis of the space of a plans a tractable option.

\subsection{Markov chains and districting ensembles}
In this paper we will generate plans via a Markov chain algorithm called ReCom (short for Recombination), but our TDA methods can be used to study any ensemble of districting plans regardless of how it was generated. We briefly describe the ReCom algorithm; more details can be found in \cite{deford2019recombination} (see also \cite{mgggva} for an early application of this algorithm to redistricting in Virginia). The GerryChain package, which implements the ReCom algorithm, is publicly available at \url{github.com/mggg/GerryChain}.

Let $G = (V,E)$ denote a unit dual graph and suppose that ${P} = \{P_1,\ldots,P_k\}$ is a valid districting plan (meaning that it satisfies all of the geographic and demographic criteria required by the particular state). 
ReCom is a Markov chain whose state space is $\PkG$.
A step in the chain is taken as follows. First, two distinct adjacent districts $P_{i}$ and $P_{j}$ are chosen uniformly at random. Next, a random spanning tree of the connected subgraph $P_{i} \cup P_{j} \subset G$ is constructed. Deleting an edge from this spanning tree repartitions $P_{i} \cup P_{j}$ into $P_i'\cup P_j'$, which together with the other districts defines a new partition of $G$. 
After a validity check, the chain can move there as its next state.
If no edge is found whose deletion produces a valid partition, then a new spanning tree is constructed, and the process
iterates. 

Details of the ReCom algorithm, as well as theoretical and experimental evidence of its  performance, are provided in \cite{deford2019recombination}. The main takeaways are that ReCom runs quickly enough to generate large ensembles, with  convergence heuristics indicating approximate independence of starting point.
The ReCom stationary distribution is well-approximated by a distribution on $\PkG$ that weights a plan $P$ in proportion to its compactness, given by an appropriate isoperimetric score on the districts \cite{CDRR}.

\subsection{Filtrations on dual graphs} 

Let $P = \{P_1,\ldots,P_k\} \in \mathcal{P}_{k,\epsilon}(G)$ be a districting plan for a state with unit dual graph $G$. As described in Section \ref{subsec:redistricting_as_graph_partition}, we form the \emph{district dual graph} associated to $P$, denoted $G_{P} = (V_{P},E_{P})$, as the quotient graph of $G$ associated to the partition ${P}$, with vertices $w_i$ corresponding to districts $P_i$.

Since the vertices of the dual graph $G$ represent geographic units,  numerical data assigned to the units (such as population, demographic counts, or vote counts) induce functions 
$f:V \rightarrow \R$. The main data of interest for the applications in this paper will be the number of votes for a Republican candidate $r: V \rightarrow \R$ and the number of votes for a Democratic candidate $d: V \rightarrow \R$ in a particular election. Given this data, we can compute the Republican vote share by district $R: \PkG\times\{1,\dots,k\} \rightarrow \R$ as follows:
$$R_P(i)= R(P,i)=\frac{r(P_i)}{r(P_i)+d(P_i)}.$$
Note that this can be thought of as an average of a binary Republican voting function over the major-party voters of the district.

Although the analysis in this paper will focus on the Republican vote share function $R$, we remark that our method can be used to analyze any function of interest $f:V \rightarrow \R$. From such a function, one can always induce a function  $F$ summing or  averaging  $f$ over the units that make up each district. We regard such an $F$ as a \emph{filtration function} on the vertices of $G_{P}$. For the sake of simplicity, we sometimes make the (mild) assumption that the values of any filtration function $F$ are distinct over the districts of a plan, and we will also assume that the range of $F$ is in the interval $[0,1]$.
Then the output of any filtration function $F$ on a plan with $k$ districts is $k$ values, which 
can be organized in increasing order, $0\le t_1 < t_2 < \cdots < t_k\le 1$.

\begin{definition}[Filtration]\label{def:filtered_dual_graph}
Let $G_{P}$ be a district dual graph for a plan with $k$ districts and a fixed filtration function $F$, and suppose that the   vertices $\{w_1,\dots,w_k\}$ are indexed in increasing order $F(w_1)<\dots<F(w_k)$.  Let $t_i=F(w_i)$.  
Then we define the  \emph{filtered dual graph} associated to $(G,P,F)$ to be the sequence of graphs
$$G_{P}^{t_1} \subset G_{P}^{t_2} \subset \cdots \subset G_{P}^{t_k},$$
where each $G_{P}^{t_j}$ is the induced graph on the vertices $\{w_1,\dots,w_j\}$.
\end{definition}

\begin{figure}[ht]
\centering
\begin{tikzpicture}
\node at (0,9.5) {\includegraphics[height=4cm]{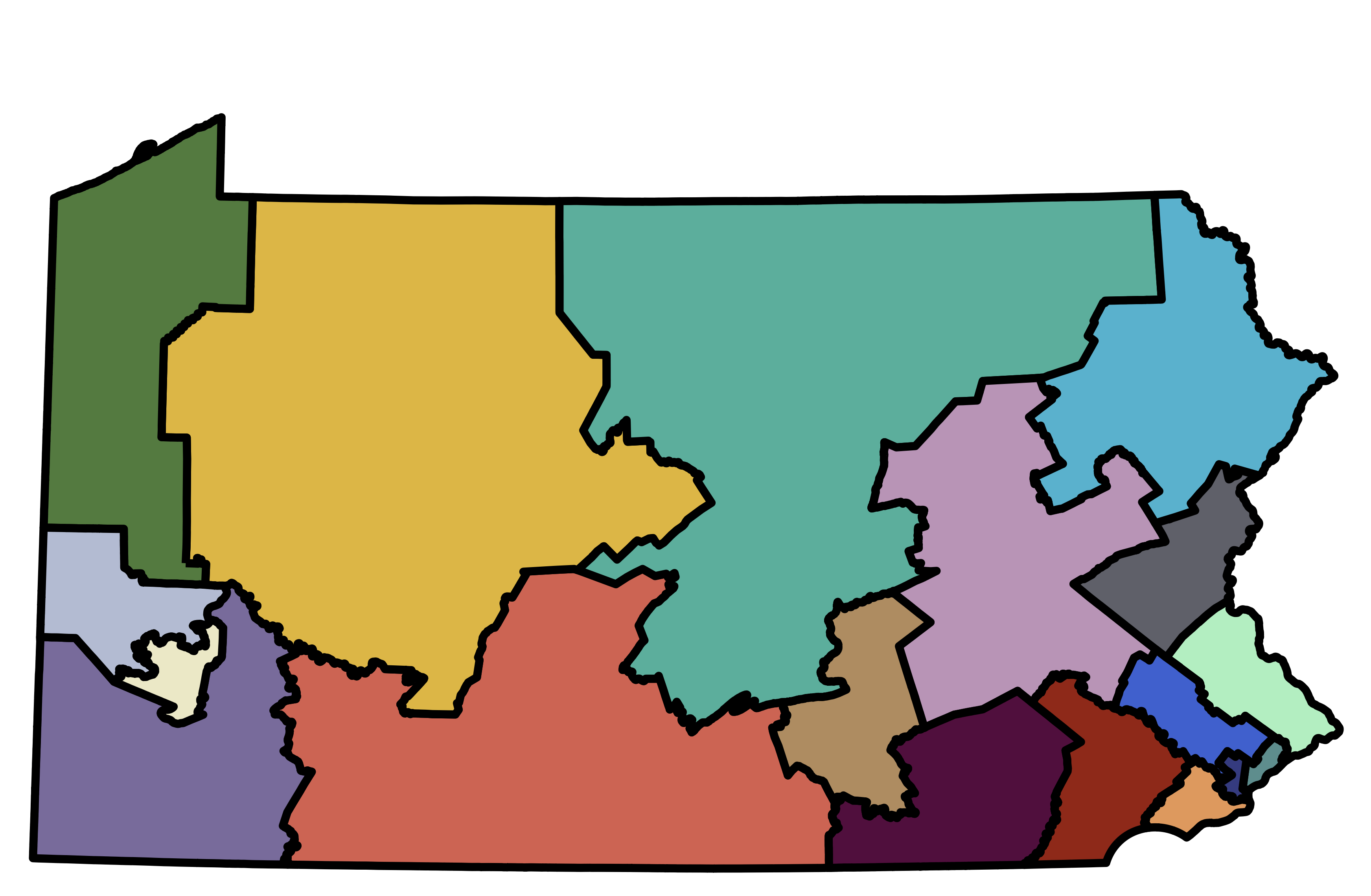}};
\node at (0,5.5) {\includegraphics[height=4cm]{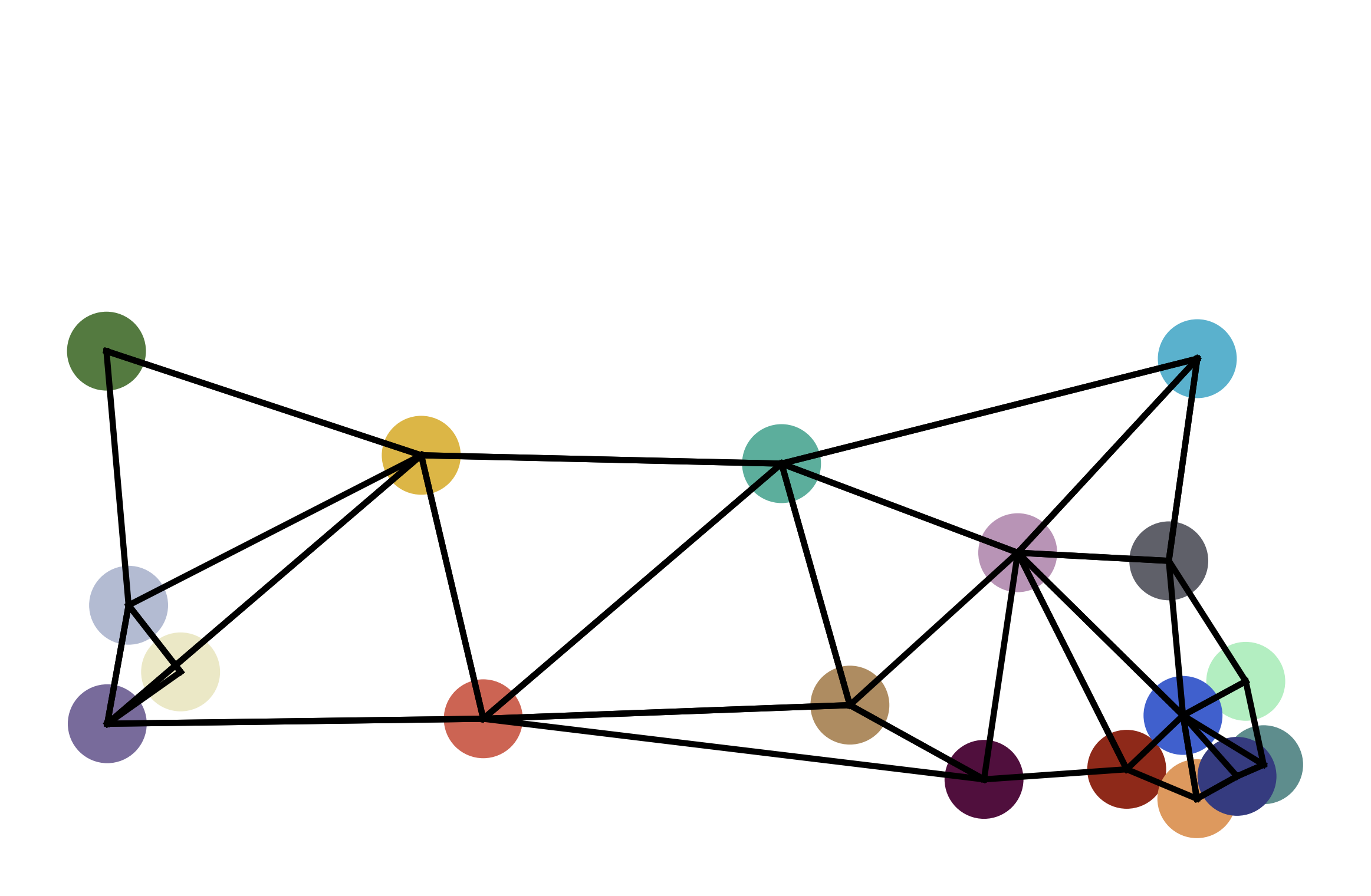}};
\node at (0,0) {\includegraphics[height=5cm]{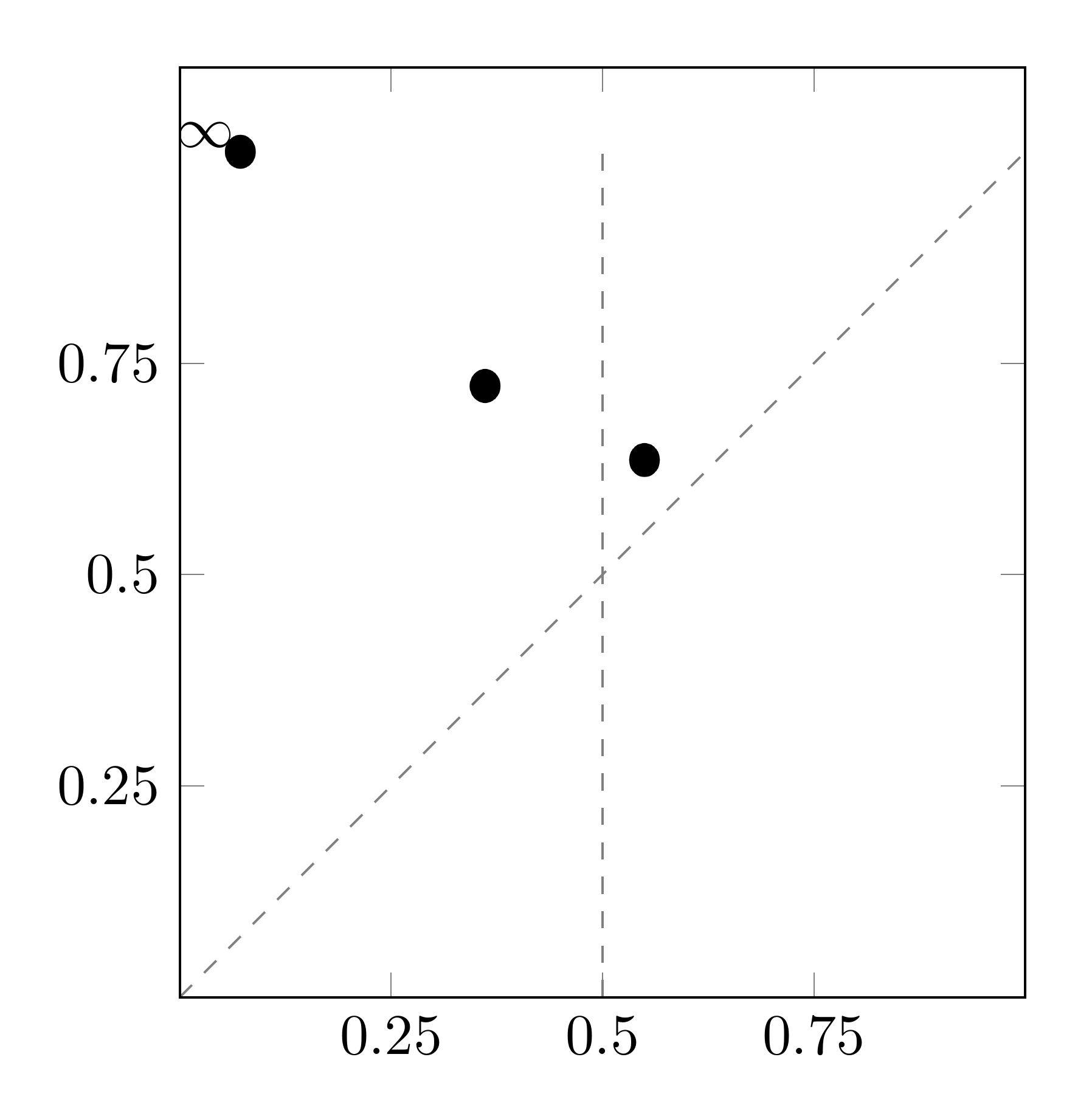}};

\foreach \x/\y/\z in {0/8/8,1/6/37,2/4/55,3/2/66,4/0/100}
{\begin{scope}[xshift=8cm,yshift=1.3*\y cm]
\draw[clip] (-2.2,-1.5) rectangle (2.2,1);
\node at (0,0) {\includegraphics[height=3.3cm]{images/REMEDIAL_filter\x.png}};
\node at (0,.6) {\footnotesize $\le$ \z \% Repub};
\end{scope}
}

\draw [->, line width=6] (0,7.3)--(0,6.5);
\draw [->,line width=6] (3.7,4.8)--(5,4.8);
\draw [->,line width=6] (5,0)--(3,0);

\node at (8,-2) {Filtration by Repub share};
\node at (0,11) {Districting plan};
\node at (0,6.1) {Dual graph};
\node at (0.3,2.5) {Persistence diagram};
\end{tikzpicture}
\caption{The main method: turning a districting plan into a persistence diagram. The persistence diagram captures both the adjacency information and the vote shares for the districts in a plan.}
\label{fig:plan_to_diagram}
\end{figure}

This is illustrated in Figure~\ref{fig:plan_to_diagram}.
The definition of $E_{P}^{t_j}$ simply means that edges from $G_{P}$ are added to the filtered dual graph as soon as possible (i.e., as soon as both of their endpoint vertices appear).
In the standard language of persistent homology, the filtered dual graph is a \emph{one-dimensional filtered simplicial complex} obtained via a \emph{sublevel set filtration} \cite{cohen2007stability}.

\subsection{Persistence diagrams}\label{sec:persistence_diagrams}

We would like to track the evolution of the shape of a filtered dual graph as we run through the values $t_j$. We focus on the evolution of the simplest topological features of the graphs: their number of connected components. This evolution is encoded in a \emph{persistence diagram} $\D=\{(b_1,d_1),\dots,(b_k,d_k)\}$, thought of as recording ``birth" and ``death" times of each connected component in the filtration.  
For each $i$, set  $b_i=t_i$  and set $d_i=t_{i'}$, where $i'\ge i$ is the minimum index such that the component of $w_i$ in  $G_{P}^{t_i}$ contains a  node whose index is less than $i$.  (Since no such index exists for $i=1$, we put $d_1=\infty$.)

The persistence diagram is a set of points $(b,d)$ with $b \in \R$, $d \in \R \cup \{\infty\}$ and  $0 \leq b \leq d$, so can be  drawn in a unit square above the diagonal, with $\infty$ represented at $d=1.1$ (say), as in Figure~\ref{fig:plan_to_diagram}.
In general, a persistence diagram could be a \emph{multiset}, containing repeated points. In our application, we have assumed that the $t_i$ values are distinct, so the points in $\D$ have multiplicity one.

Persistence diagrams are a ubiquitous summary statistic used in TDA. The precise, general definition of a persistence diagram relies on a decomposition theorem from representation theory known as Gabriel's theorem, which was proved in \cite{gabriel1972unzerlegbare}, with the connection to topological data analysis first observed in \cite{carlsson2009zigzag}. 
When representing a persistence diagram one typically omits points along the diagonal, as these represent ``trivial" topological features.

\begin{remark}
\

\begin{enumerate}
\item The persistence diagram defined here describes the evolution of the zeroth-homology of the filtered dual graph; this is the so-called \emph{degree-$0$ persistent homology} of the filtered dual graph. Likewise, persistent homology in higher degrees can be computed, tracking appearances and disappearances of higher-dimensional topological features such as loops and voids. We hope to explore first homology in future work.
\item The definition above entails a 
 convention sometimes referred to as the \emph{Elder Rule} in the persistent homology literature:  when two components merge, the component with the more recent birth time is annihilated.  
See \cite{edelsbrunner2010computational} for a discussion of this convention, as well as the more recent theoretical justification provided in \cite{curry2018fiber}.
\end{enumerate}
\end{remark}

The entire pipeline described so far is summarized in Figure \ref{fig:plan_to_diagram}. In the end, this pipeline reduces a very complex object (a geographical partition of a state together with a function describing, e.g., election data) to a simple, if somewhat abstract, summary statistic (a persistence diagram encoding the evolution of connectivity of districts according to the function). It is the thesis of this paper that this large reduction in complexity preserves, and in fact illuminates, interesting geographical and political properties of the plans. Our study of persistence diagrams for redistricting ensembles will require some quantitative tools described in the following subsections. 

\subsection{Wasserstein distances}\label{sec:wasserstein_distance}

Let $\D_1$ and $\D_2$ be a pair of persistence diagrams. We introduce  a standard family of distance metrics used to compare persistence diagrams.

\begin{definition}[Wasserstein $p$-distance]
Let $p \in [1,\infty]$ and let $(b,d)$ be a point in a persistence diagram. Define the \emph{diagonal $p$-distance} $\Delta_p(b,d)$ to be the $\ell_p$ distance from $(b,d)$ to the nearest point on the \emph{diagonal} $\{(x,x) \mid x \geq 0\}$.  A \emph{partial bijection} between $\D_1$ and $\D_2$ is a bijection $\phi$ from a subset of $\D_1$ onto a subset of $\D_2$. We use the notation $\phi:\D_1 \nrightarrow \D_2$ for a partial bijection, $\dom(\phi)$ for the domain of $\phi$ and $\im(\phi)$ for its image. For $p < \infty$, define the \emph{$p$-cost} of a partial bijection $\phi$ by
$$\left(\mathrm{cost}_p(\phi)\right)^p = \sum_{(b,d) \in \dom(\phi)} {\|(b,d) - \phi(b,d)\|_p}^p 
+ \sum_{(b,d) \not \in \dom(\phi)} \Delta_p(b,d)^p  + \sum_{(b',d') \not \in \im(\phi)} \Delta_p(b',d')^p.$$

This extends to a $p = \infty$ cost in the usual way, as a supremum:
$$\mathrm{cost}_\infty(\phi) = \max \left\{ \max_{(b,d) \in \dom(\phi)} \|(b,d) - \phi(b,d)\|_\infty,\right. 
 \left. \max_{(b,d) \not\in \dom(\phi)} \Delta_\infty(b,d), \max_{(b',d') \not \in \im(\phi)} \Delta_\infty(b',d')\right\}.$$

Finally, define the \emph{Wasserstein $p$-distance} between $\D_1$ and $\D_2$ by
\begin{equation}
d_p(\D_1,\D_2) = \inf_{\phi:\D_1 \nrightarrow \D_2} \mathrm{cost}_p(\phi). \label{eq:costmin}
\end{equation}
\end{definition}

The Wasserstein $p$-distances are naturally described in the language of matching features. The goal is to match points in the diagram, which we now recall represent persistent topological features. A partial bijection $\phi$ defines such a matching of features; those points not in the domain or range of $\phi$ are considered to be matched with a ``trivial feature" on the diagonal (i.e., a feature with zero lifespan). The $p$-cost $\mathrm{cost}_p(\phi)$ of matching $\phi$ measures how far points have to be moved in $L^p$ distance, and the goal is to find a matching which minimizes this cost. For a given pair of persistence diagrams, we will refer to the partial bijection $\phi$ that achieves the infimum in (\ref{eq:costmin}) above as an \emph{optimal $L^p$ matching}. This will be referenced later when comparing Fr\'{e}chet means for different ensembles or vote data. 

Traditionally, the $p=\infty$ Wasserstein distance, also known as \emph{bottleneck distance}, has been the standard choice for comparing persistence diagrams in topological data analysis. It was one of the first metrics of this type to be introduced \cite{cohen2007stability} and its matching cost is the simplest to interpret: the cost of matching $(b,d)$ to $(b',d')$ is whichever is larger, the difference in birth times or the difference in death times. It was moreover shown in that paper that the bottleneck distance enjoys stability under perturbations of filtrations generating persistence diagrams (this result will be stated and used in Section~\ref{sec:stability} below). The use of bottleneck distance was later justified at a more theoretical level by the famous \emph{Isometry Theorem}. This theorem, established in \cite{chazal2009proximity} and \cite{lesnick2015theory}, states that the bottleneck distance between persistence diagrams coincides with the so-called \emph{interleaving distance}, a distance which is defined purely in terms of morphisms between persistence modules using the language of category theory. This category-theoretic interpretation of bottleneck distance was extended to the other Wasserstein distances in the recent paper \cite{bubenik2018wasserstein}. On the computational side, the optimal partial bijections required by the Wasserstein distances can be computed using the Hungarian algorithm (see, e.g., \cite{chazal2009gromov}), a polynomial-time algorithm for solving the assignment problem required to find an optimal matching.

\subsection{Fr\'{e}chet means}
Let $(X,d)$ denote an arbitrary metric space.  For a finite  sample $x_1,\ldots,x_n$ from $X$, define a
discrete probability measure $\frac{1}{n} \sum_{j =1}^n \delta_{x_j}$, where $\delta_{x_j}$ the Dirac measure at $x_j$.
The associated \emph{Fr\'{e}chet functional} $F: X \rightarrow \R$ is defined by
$$F(x) = \frac{1}{n}\sum_{j=1}^n d(x,x_j)^2.$$
A \emph{Fr\'{e}chet mean} of the sample is a minimizer of this Fr\'{e}chet functional, i.e., a {\em barycenter} for the collection of points.  Below, we will make heavy use of barycenters to summarize ensembles of diagrams. 

Metric barycenters are guaranteed to exist uniquely in spaces of nonpositive curvature.  We are not in that setting, but are nonetheless able to leverage nice metric features to get a well-behaved barycenter construction.
With the $p=2$ Wasserstein distance, Turner et al.\ show that the space of all persistence diagrams is an Alexandrov space of nonnegative curvature: a complete, separable geodesic space with curvature $\ge 0$ in the sense of comparison triangles \cite{turner2014frechet}. 
Despite the positive curvature, these metric spaces have helpful structure:  they admit a  notion of a tangent space at each point, which can be endowed with an inner product  (see \cite{lytchak2006open} for general theory). Moreover, semiconcave functions  have well-defined gradients.  It is shown in \cite{turner2014frechet} that the Fr\'{e}chet functional is semiconcave, motivating its optimization via gradient descent.

\begin{figure}[ht]
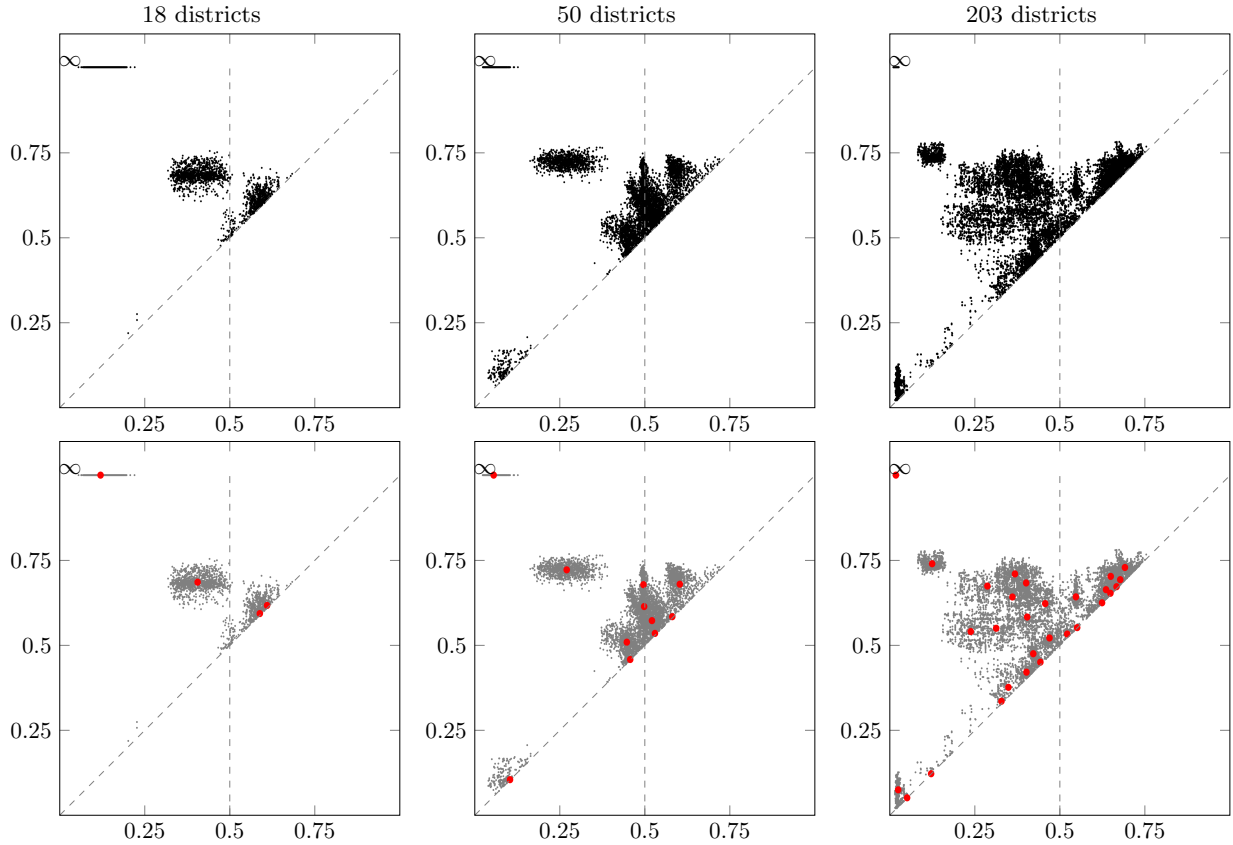

\centering
\foreach \k in {18, 50, 203}
{
\begin{subfigure}{0.32\textwidth}
\caption*{ \k \ districts}
\pdsubfig{PA_noaxes/\k districts_PRES16_PD_overlaid.png}%
\end{subfigure}
}

\foreach \k in {18, 50, 203}
{
\begin{subfigure}{0.32\textwidth}
\pdsubfig{PA_noaxes/\k districts_PRES16_FrechetPD_overlaid.png}%
\end{subfigure}
}
\caption{Pairing each plan in a districting ensemble with PRES16 vote data gives an ensemble of persistence diagrams.
These plots show overlays of the diagram ensembles for $k=18,50,203$ respectively.
In the second row, we add the Fr\'echet mean $\mathcal{F}$ for each ensemble, shown in red, representing an ``average'' persistence diagram for the ensemble.}
\label{PAFrechetoverlaid}
\end{figure}

\begin{figure}[ht]
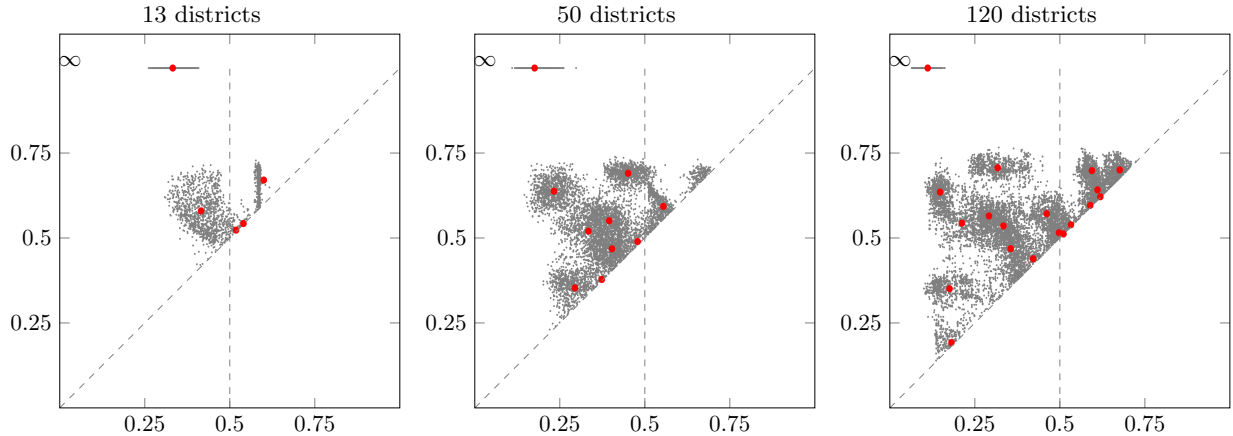

\centering
\foreach \k in {13, 50, 120}
{
\begin{subfigure}{0.32\textwidth}
\caption*{ \k \ districts}
\pdsubfig{NC_noaxes/\k districts_PRES16_FrechetPD_overlaid.png}%
\end{subfigure}
}
\caption{Overlaid persistence diagrams with Fr\'{e}chet mean (in red) for each of the NC ensembles (with PRES16 vote data). }
\label{NCFrechetoverlaid}
\end{figure}

The \cite{turner2014frechet}  authors then propose the following algorithm. For persistence diagrams $\{\D_1,\ldots,\D_m\}$, suppose that $\D^{(i)}$ is the current Fr\'{e}chet mean candidate. An update is performed as follows:
\begin{enumerate}
\item Compute optimal partial bijections $\phi^{(i)}_j:\D^{(i)} \nrightarrow \D_j$;
\item Extend $\phi^{(i)}_j$ to a function $\phi_j^{(i)} \D^{(i)} \rightarrow \D_j \cup \{(x,y) \in \R^2 \mid x = y\}$ by declaring $\phi^{(i)}_j(b,d)$ to be the point on the diagonal nearest to $(b,d)$ for each $(b,d) \not \in \dom(\phi^{(i)}_j)$;
\item For each $(b,d) \in \D^{(i)}$, update to the arithmetic mean of the set of images $\phi^{(i)}_j(b,d)$;
\item Let the new Fr\'{e}chet mean candidate $\D^{(i+1)}$ be the  updated point collection.
\end{enumerate}
This  process is shown to converge, at least to a local minimum of the Fr\'{e}chet functional. As this is an iterative algorithm, the local minimum found depends on the seed $\D^{(1)}$, that is, the initial Fr\'{e}chet mean candidate. To mitigate this dependence in our analysis, we run the algorithm multiple times, using each persistence diagram in our collection as a seed, and take the local minimum with the lowest Fr\'{e}chet functional value as our final Fr\'{e}chet mean.

\subsection{Computational details}\label{subsec:computational_details}
The main code base used for this paper is available on GitHub \footnote{\url{github.com/mggg/TDA-redistricting}}. The ensembles used in this paper were generated using the open-source package GerryChain developed by the Metric Geometry and Gerrymandering Group \cite{gerrychain}. GerryChain implements numerous Markov chains on $\PkG$ for generating ensembles of districting plans, including the recombination and flip chains discussed in this paper. For the main ensembles used in Sections \ref{app:scales} and \ref{app:elections}, we sampled every 50th plan from an ensemble of 50,000 plans generated using a ReCom chain with 
$\epsilon=.02$ as the maximum population deviation. The Congressional and state Senate ensembles were generated from dual graphs constructed from the precinct shapefiles available from MGGG \cite{mggg-states}. For the state House ensembles, more granularity is required in order to get satisfactory population balance, so we used census block dual graphs. In the case of North Carolina, the state releases electoral data at the block level for redistricting purposes. For Pennsylvania, election results were disaggregated from the precinct level down to blocks using the geospatial data package called MAUP \cite{MAUP}. The persistence diagrams were generated using the library GUDHI (\url{https://gudhi.inria.fr}) \cite{maria2014gudhi}.

\section{Methods}\label{sec:applications}\label{sec:principal}

We propose three types of analysis enabled by  persistent homology on ensembles of districting plans.

\subsection{Scale and zoning} \label{app:scales}
\begin{figure}[ht]
\centering
\includegraphics[width=.55\textwidth]{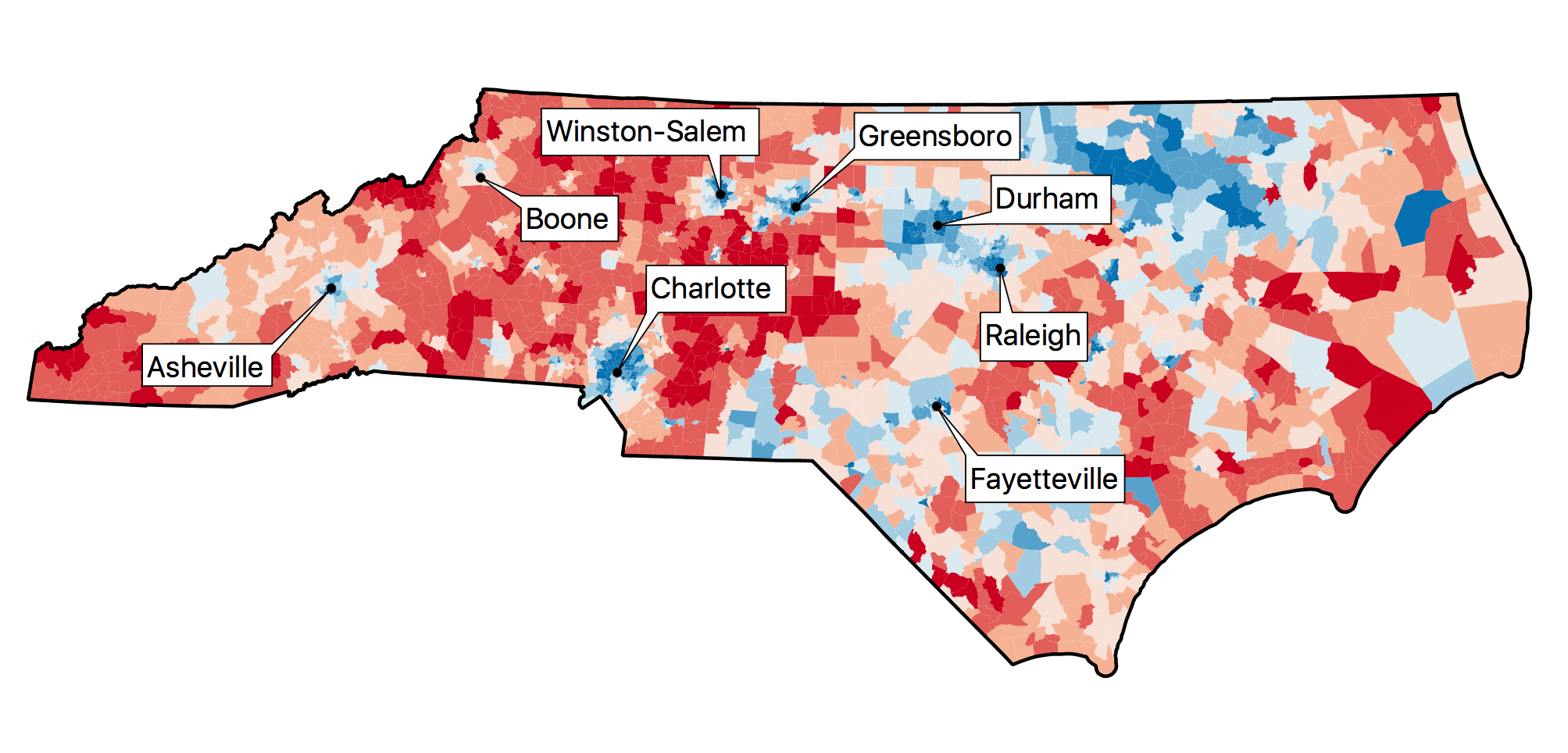}

\includegraphics[width=.55\textwidth]{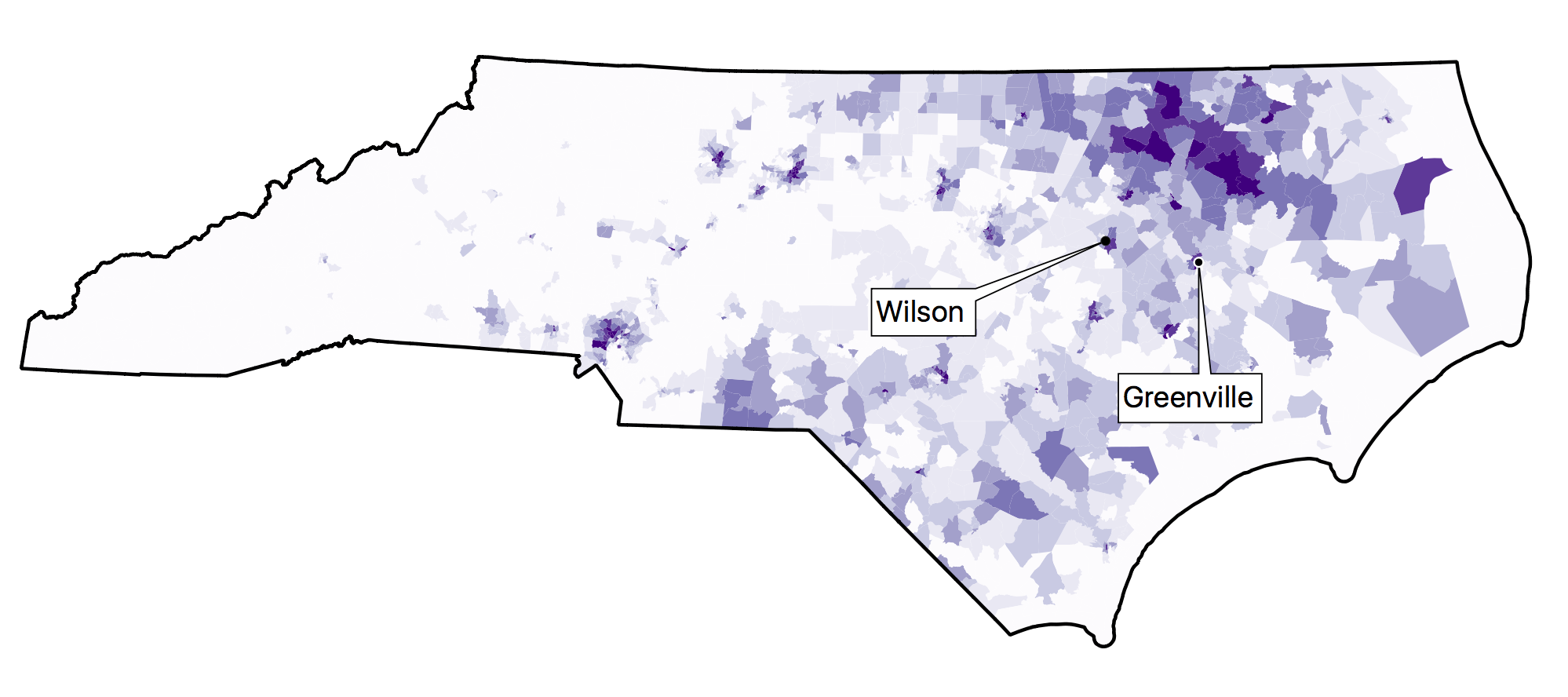}
\caption{Maps of North Carolina.  Top: results from the 2016 Presidential election (blue for Democratic, red for Republican). Bottom: North Carolina has a significant rural Black population, unlike Pennsylvania. Greenville and Wilson are the largest towns in the Northeast, but they are not the hubs of Black population.}
\label{labelledNC}
\end{figure}

The winner-take-all plurality structure of the American districting system introduces significant non-linearity to the relationship between votes and seats, whose patterns we seek to capture.  
The key question here asks how  the fine-grained vote data by precinct (see Figures \ref{fig:PAbasic} and \ref{labelledNC}) translates into an outcome at the district level.  
We would like to isolate the effects of district {\em size} on the geography and partisan properties of redistricting plans.
A scale study was carried out in \cite{RoddenWeighill2020} by a direct comparison of the district size to the range of observed partisan outcomes in an ensemble, using recent Pennsylvania voting patterns.    
In this section, we use topological signatures to probe scale effects further.  The main finding is that at each scale, the state decomposes into zones of partisan support whose anchor points can be detected by persistent homology techniques.

Before we discuss the application of our  method to ensembles of maps, we first illustrate the information captured  on individual plans. Figure \ref{fig:PAbasic} shows the currently enacted Congressional and Senate districts for Pennsylvania, colored by Trump vs. Clinton share (Pres 2016). 
Comparing to the choropleth at the top of Figure \ref{fig:PAbasic}, we can see that at the Congressional scale, the Democratic voting in Erie and State College is outweighed by the surrounding Republican votes, while the Senate scale allows those cities to pull their districts more strongly.  
Thinking in terms of peaks of Democratic vote share, Pittsburgh, Harrisburg, and Philadelphia are all pronounced peaks on the precinct level, with the corresponding Congressional districts (18, 10, and 3, respectively) also showing up as peaks on the district level.  

Homological persistence will capture a {\em relative height} of these peaks in the following way:
 Each district which is more Democratic than all its neighbors shows up as a point $(b,d)$ in the persistence diagram, where $b$ is precisely its Republican vote share. A value of $d$ on the `death' axis means that any path from that district to a more Democratic one must pass through a district whose Republican share is at least $d$. 
Thus the persistence $p=d-b$ is a kind of relative intensity of Democratic voting.   
For example, District 10 is only $9$ percentage points more Democratic than  neighboring District 11, which in turn is adjacent to District 6, a more Democratic district than District 10.   So the persistence of its feature is 9.  
Contrast this with District 18: any path from District 18 to a more Democratic district must pass through either District 13 or 15, each of which is more than $36$ points less Democratic than District 18.   This gives Pittsburgh a persistence of 36.  
It is clear that a Senate plan, with a larger number of smaller districts, will tend to have more peaks, such as Erie in the  northwest corner; the greater number of features shows up in the Senate plan's persistence diagram, together with the solutions to the minimax problems that describe their relative height.

So far we have discussed the interpretation of persistence for specific plans, but of course different plans may have different structure.  
The use of Fr\'echet means will allow us to examine how much common structure there is across the plans of a districting ensemble.

\begin{enumerate}
\item[Step 1.] {\bf Overlay.}
Start with a state, a number of districts $k$, and some fixed vote data such as from a recent election.  We generate an ensemble of one thousand plans (see Section \ref{subsec:computational_details}) and form diagrams $\D_1,\dots,\D_{1000}$.
As a first step,  we can overlay the thousand diagrams (taking the unions of the sets $\{(b,d)\}$) as the top row of Figure \ref{PAFrechetoverlaid} does for Pennsylvania ensembles. 
To tease apart the structure of these plots, we compute the Fr\'{e}chet mean $\F$ of the ensemble as in the next row of Figure \ref{PAFrechetoverlaid}. This gives us one summarized view of the ensemble as an ``average'' persistence diagram.

\item[Step 2.] {\bf Marking.}
Next, we want to extract a {\em marking} of the individual persistence diagrams.  Where $k$ is the number of districts, we can choose $\ell \le k$ points in $\F$ that are deemed sufficiently far from the diagonal to be significant.  
Then for every  diagram $\D_i$ in our ensemble, we match the points in $\D_i$ to the chosen points in $\F$ using an optimal $L^2$ matching (see Section \ref{sec:wasserstein_distance}), inducing a partial labeling of the individual persistence diagrams by elements of the Fr\'{e}chet mean. 
Returning to the overlaid diagram, we can now create individual Fr\'echet point plots corresponding to each significant feature of $\F$, with the property that each diagram $\D_i$ contributes at most one point to this plot.

\item[Step 3.] {\bf Localization.}
For a significant feature of $\F$, many of the $\D_i$ have a corresponding feature, picking out a district that initiated the connected component in the filtration of the district dual graph.  We can then measure how often in the ensemble each geographical unit belonged to that district.  This gives us a {\em heat map} that shows us where the feature is located in the state.  Figure~\ref{PAclasses18} shows Fr\'echet point plots and corresponding precinct heat maps for PA Congressional and Senate districts.  (Supplemental Figure~\ref{PAclasses203} uses a census block heat map for the House ensemble.)

\item[Step 4.] {\bf Zoning.}  If the filtration is in terms of one party's share (here, Republican), we can then zone the state into clusters of districts won by the other party (here, Democratic), as follows.  There is a one-to-one correspondence between persistence features in the northwest quadrant of the diagram $\D$  (i.e., $b<.5,d>.5$)
and 
clusters of Democratic-won districts.  The birth time ensures that the cluster is anchored by a Democratic-won district and the death time ensures that it is separated from other clusters by Republican-won districts.  For each localized point in   
$\F$, we consider northwest-quadrant points marked by it in the $\D_i$, representing Democartic clusters of districts.  We can map these clusters and report their average number of districts.  
\end{enumerate}

When studying scale effects, we are interested in which areas  represent peaks at each scale (usually cities from the Democratic point of view). Crucially, cities only show up as individuated peaks in Democratic vote strength when they are separated from other cities by redder {\em districts}---a property that it would be hard to see with one's eyes.

Note that there is no a priori reason that the localization step should produce a geographically specific position for each Fr\'echet  feature.  Picking out an identifiable location at this step is a signal that there is consistency across the ensemble in the correspondence of geography and aggregated voting.  
The method has succeeded in zoning the state at a given scale if the heat maps are geographically coherent.
If localization succeeds, the workflow lets us refer to the Fr\'echet features by geographical names and expected number of Democratic-won districts for a given election:  the Philadelphia zone (5.6 Congressional districts), the Pittsburgh zone (1.2 Congressional disricts), and so on.

\subsection{Comparing elections} \label{app:elections}
In the previous section, we studied ensembles of plans relative to some fixed vote pattern. Our second type of analysis will involve varying the vote data. We ask several questions.  
Do the Fr\'echet means themselves tell us something interesting about the relationship in voter preferences between elections?  
How robust are the patterns learned above to a shift in the chosen electoral baseline?

For a given state and district size, we again use an ensemble of one thousand districting plans. We then vary the vote data (but not the ensemble of plans) to produce multiple ensembles of diagrams,  with an associated Fr\'{e}chet mean plot for each vote pattern. Suppose that $\F$ and $\F'$ are Fr\'echet means for the same plans and two different elections.
We can pair points from $\F$ to $\F'$
(by either a geographic pairing or an optimal $L^2$ matching) to provide a summary of electoral difference.

For a pair of features that is matched, the direction of displacement carries information.  Suppose the mean point is 
$(b,d)$ in $\F$ and $(b',d')$ in $\F'$.  If $b>b'$, this zone was anchored by a district that is more Republican in the first election than the second.  If $d>d'$, then the area surrounding the peak was more Republican in the first election than the second.   
Putting these together can give a summary of the voting performance of a city relative to its exurban ring at district scale.  

If a feature in $\F$ is not matched (i.e., is matched to the diagonal) in $\F'$, this gives us a different kind of information.
This represents a substantial difference in the basic geography of representation under these vote data, picking out a district cluster that represents a peak in one election but not the other.

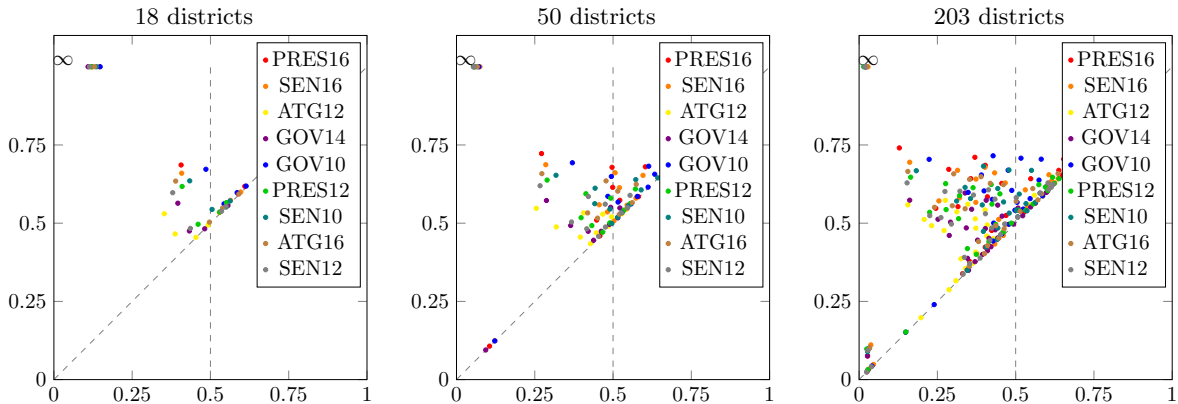
\begin{figure}[ht] 
\centering
\foreach \k in {18, 50, 203}
{
\begin{subfigure}{0.31\textwidth}
\caption*{\k \ districts}
\resizebox{\columnwidth}{!}{%
\begin{tikzpicture}
\begin{axis}[xmin=0, xmax=1, ymin=0, ymax=1.1, clip=false,
xtick={0,0.25,0.5,0.75,1}, ytick={0,0.25,0.5,0.75},
   enlargelimits=false, axis equal image, axis on top
]
\addplot[scatter, 
scatter/classes={
PRES16={red},
SEN16={orange},
ATG12={yellow},
GOV14={violet},
GOV10={blue},
PRES12={green!80!black},
SEN10={teal},
ATG16={brown},
SEN12={gray}
},
only marks, mark=*, mark size=1, scatter src=explicit symbolic]
 table[x=X, y=Y, col sep=comma, meta=label]{PA_csv/\k Frechet_means.csv};
\addplot[gray, dashed] coordinates{(0.5,0) (0.5,1)};
\addplot[gray, dashed] coordinates{(0,0) (1,1)};
\node at (axis cs: 0.03,1.02) {$\infty$};
\legend {PRES16, SEN16, ATG12, GOV14, GOV10, PRES12, SEN10, ATG16, SEN12};
\end{axis}
\end{tikzpicture}
}
\end{subfigure}
}
\caption{Overlaid Fr\'echet means for various elections in PA.}
\label{allFrechetPA}
\end{figure}

\subsection{Signals of gerrymandering}
The examination of zoning and electoral shifts described above is not targeted at identifying gerrymandering, but at understanding the tendencies of partisan-neutral districting.  
Since the detection of gerrymandering is a difficult and high-profile problem, we now turn attention to whether persistent homology can flag properties of plans that are designed to extract advantage for one major party over the other.  

Given a state, number of districts $k$, and vote pattern, our first step is to generate two ensembles:  one creating an unusually large number of Republican-majority districts and  a second favoring Democrats.  In order to avoid flagging extremely competitive plans, we will identify a {\em safe seat} for Party A when the party's vote share exceeds $.53$.  (Noting, of course, that for practical gerrymandering purposes, there may be a more complicated calculus for trading off safety with win maximizaiton.)

\begin{figure}[ht]
\centering
\begin{tikzpicture}[scale=.66]
\begin{scope}[yshift=1cm,scale=.75] 
\begin{axis}[clip=false,
   enlargelimits=false, axis on top,
      xtick={16,18,20,22,24}
]
\addplot graphics [xmin=14, xmax=26,ymin=0, ymax=350] {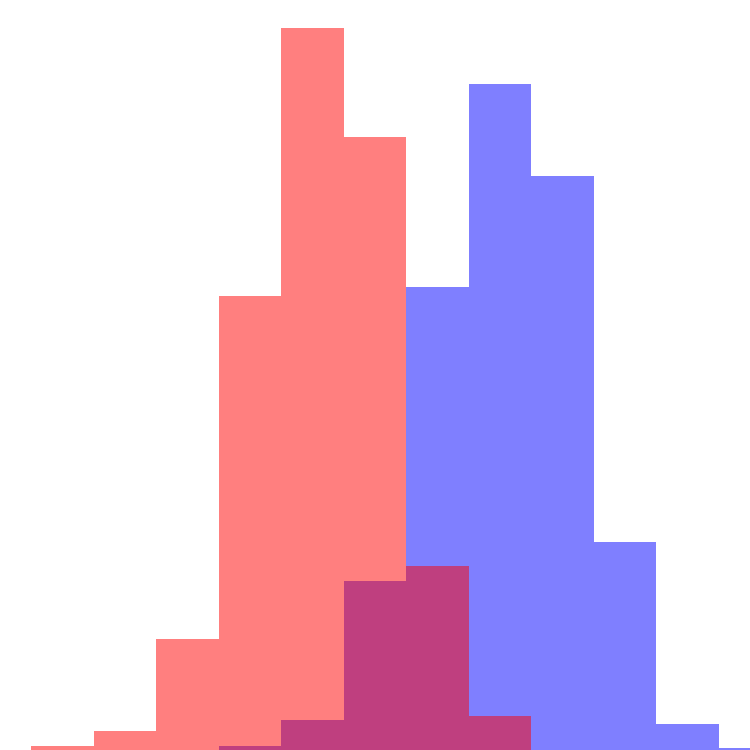}; 
\end{axis}
\end{scope}%

\begin{scope}[xshift=6.5cm] 
\begin{axis}[
xmin=0, xmax=1, ymin=0, ymax=1.1, clip=false,
xtick={0.25,0.5,0.75}, ytick={0.25,0.5,0.75},
   enlargelimits=false, axis equal image, axis on top
]
\addplot graphics [xmin=0, xmax=1,ymin=0, ymax=1.1] {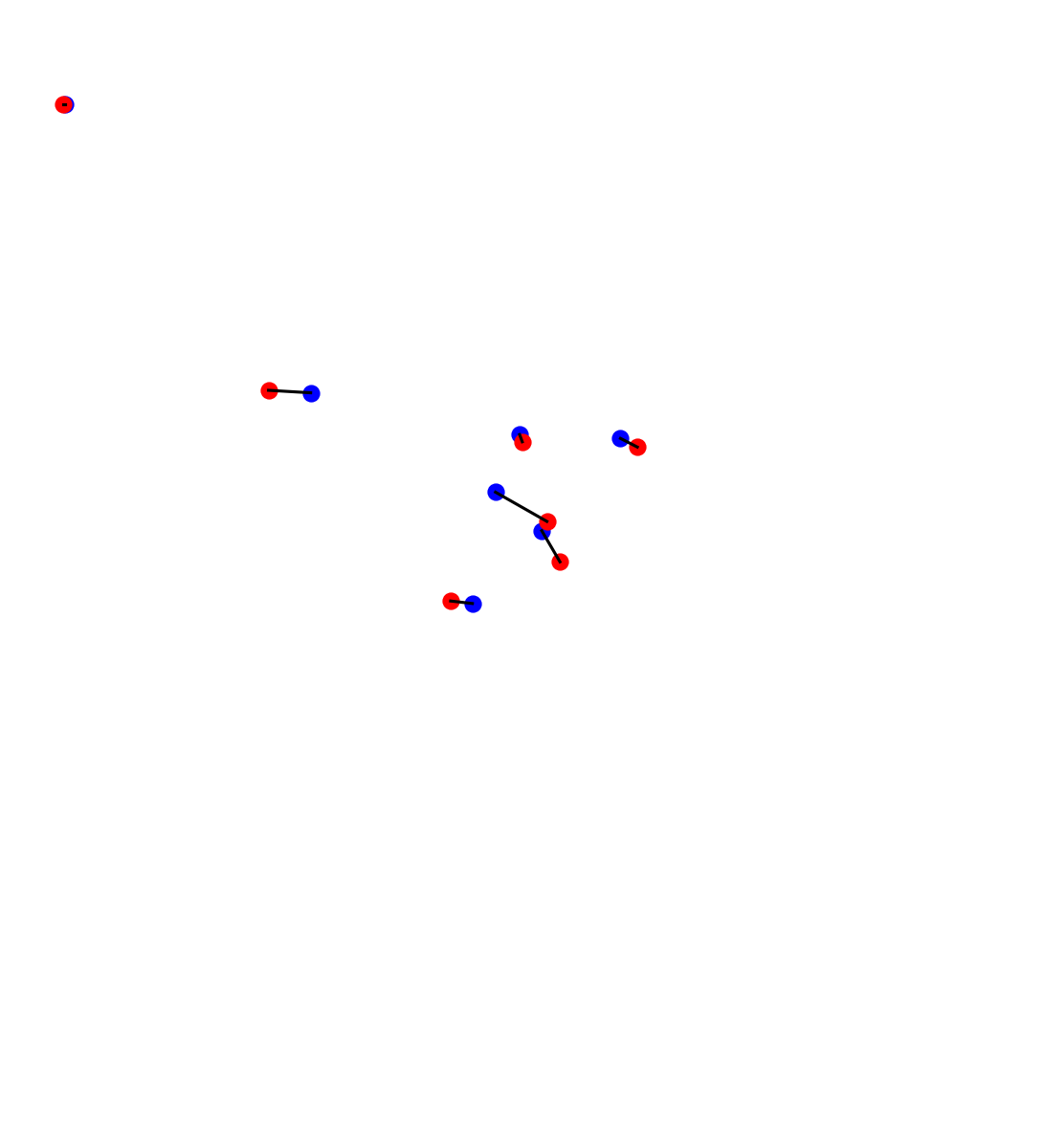}; 
\addplot[gray, dashed] coordinates{(0.5,0) (0.5,1)};
\addplot[gray, dashed] coordinates{(0,0) (1,1)};
\node at (axis cs: 0.03,1.02) {$\infty$};
\end{axis}
\end{scope}%

\begin{scope}[xshift=13cm] 
\begin{axis}[
xmin=0, xmax=1, ymin=0, ymax=1.1, clip=false,
xtick={0.25,0.5,0.75}, ytick={0.25,0.5,0.75},
   enlargelimits=false, axis equal image, axis on top
]
\addplot graphics [xmin=0, xmax=1,ymin=0, ymax=1.1] {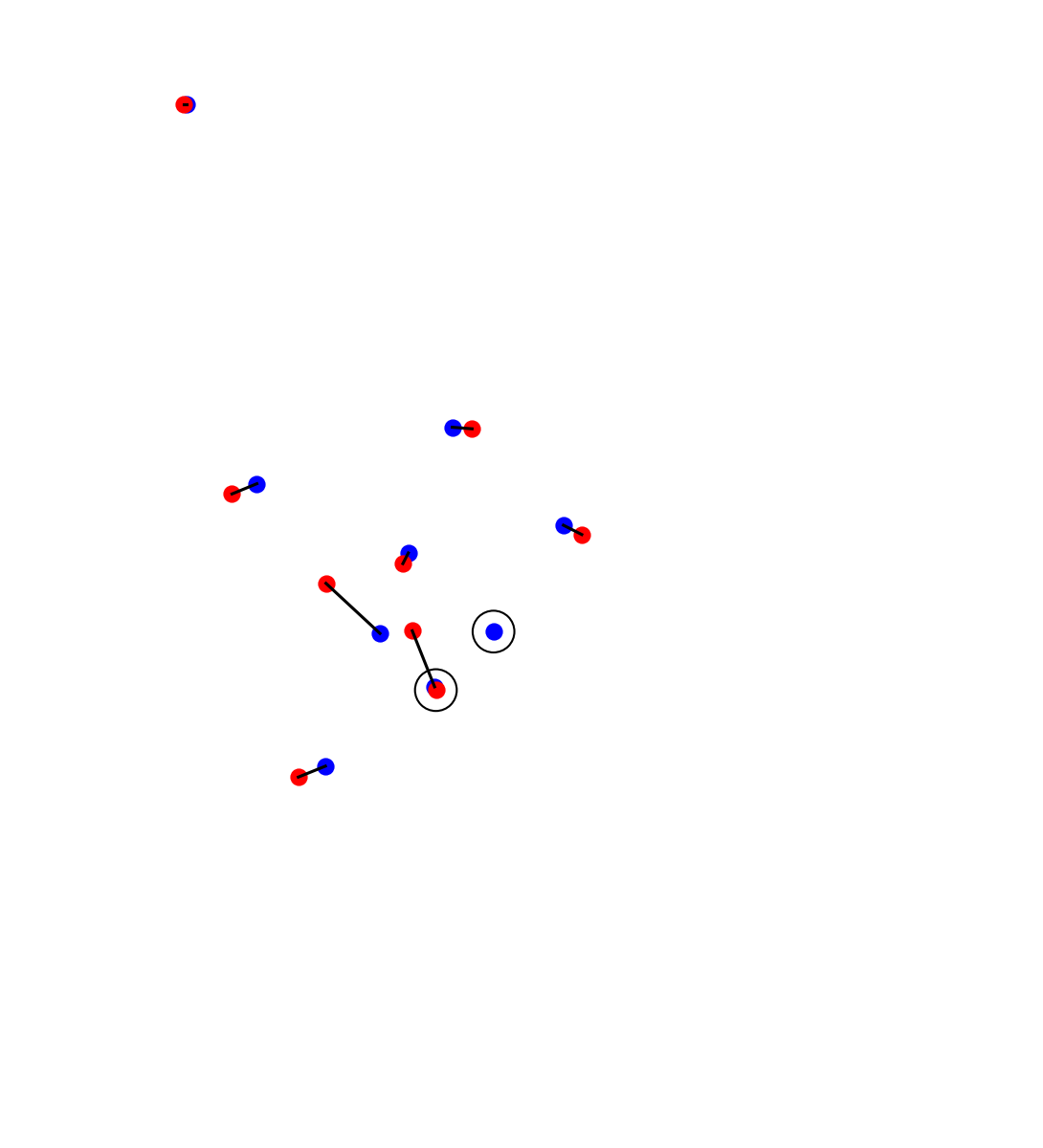}; 
\addplot[gray, dashed] coordinates{(0.5,0) (0.5,1)};
\addplot[gray, dashed] coordinates{(0,0) (1,1)};
\node at (axis cs: 0.03,1.02) {$\infty$};
\end{axis}
\end{scope}%

\begin{scope}[xshift=19.5cm,yshift=1cm,scale=.75]
\begin{axis}[clip=false,
   enlargelimits=false, axis on top,
   xtick={16,18,20,22,24}
]
\addplot graphics [xmin=14, xmax=26,ymin=0, ymax=350] {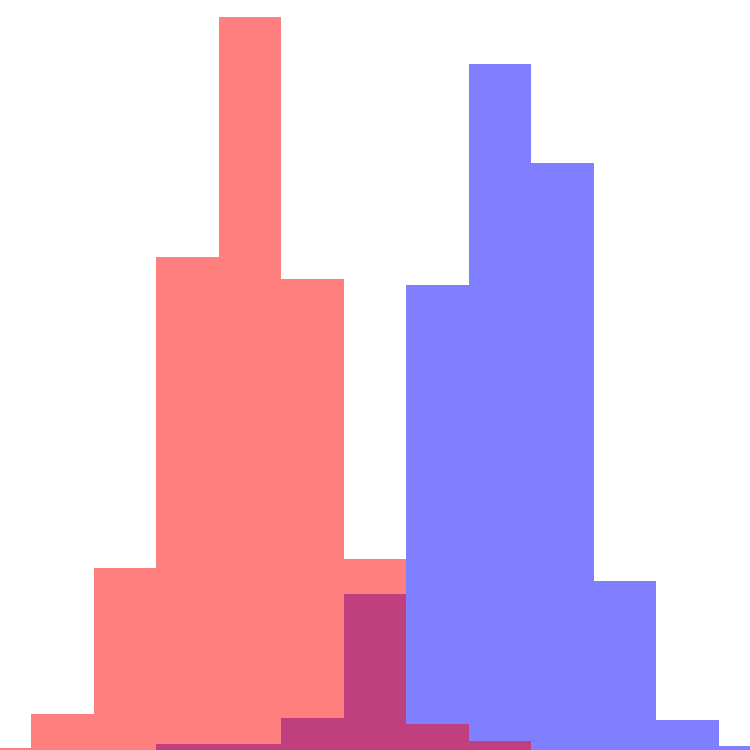}; 
\end{axis}
\end{scope}
\node at (9.2,-1) {PA state Senate};
\node at (15.7,-1) {NC state Senate};
\node at (2.65,0) {\# Dem seats};
\node at (22.15,0) {\# Dem seats};
\end{tikzpicture}
\caption{Comparison for the weighted ensembles of Senate plans, with Democratic-favoring in blue and Republican-favoring in red.  The simple weighting adjustment to the Markov chain has successfully created a difference of 3-4 seats between the ensemble means.  Despite the significant change in the partisan outcome, we see only small movement in the Fr\'echet means.}
\label{DRhist}
\end{figure}
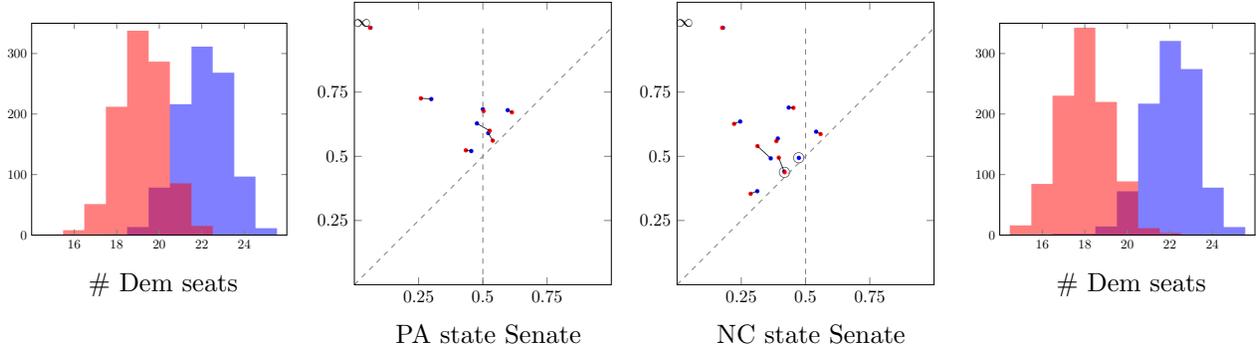

To generate an ensemble of plans biased towards Party A, we use the ReCom algorithm with a Metropolis-style weighting variant.  When a new plan $P$ is proposed by the algorithm, it is immediately accepted if the number of safe seats does not decrease. A standard Metropolis weighting with respect to a score $\sigma$ is to choose a coefficient $\beta$ and accept a plan scoring worse with probability $e^{-\beta\sdot\Delta \sigma}$.
We will do that here:  if the safe seat count decreases by $\Delta s$, then we accept $P$ with probability $e^{-2\Delta s}$.%
\footnote{We note that for standard recombination moves, this can only occur with $\Delta s=1$, corresponding to a roughly 
13.5\% probability of acceptance.  The standard protocol is to adjust $\beta$ until you have a run with a decent acceptance
rate and an ensemble that passes heuristic tests of convergence and quality.}
The purpose of the Metropolis rule is to avoid getting stuck at local optima, while favoring better plans long-term.
We use this algorithm to generate two ensembles of one thousand plans each:  a Democratic-favoring ensemble and a Republican-favoring one.  To illustrate the effectiveness of the weighted Markov chains, see Figure  \ref{DRhist}, showing a difference in expectation of 3-4 seats out of 50.

From the two biased ensembles we then generate Fr\'{e}chet point plots and associated heat maps, then compare the Fr\'echet means with either a geographic or optimal $L^2$ matching (see Section \ref{sec:wasserstein_distance} for the definition of an optimal $L^2$ matching).  As before, we have no guarantee that these two matching methods agree, but in our case studies below we find that they do.  

The aim of this analysis is to find and interpret geographical correlates of partisan advantage---not to classify entire plans as gerrymandered or not.  For that, we are unlikely to do better than simply comparing the number of seats won under various electoral assumptions to a neutral ensemble.  

In our case studies below, we can identify a few Fr\'{e}chet features which show a significant difference between the two ensembles, and we can make inferences about the reasons for these differences. However, on the whole there is surprisingly little difference between the persistence data for the two ensembles despite their very different partisan seat shares. 
This further emphasizes that persistent homology is not primarily useful for classifying plans, while reinforcing the narrative of robustly identifying relationships between geography and aggregated vote data.

\section{Case study: Pennsylvania}
\label{sec:casestudyPA}

Pennsylvania has 18 Congressional districts, 50 state Senate districts, and 203 state House districts.  
We have cleaned electoral data from nine recent statewide elections:  
two Presidential elections (2012, 2016), three U.S. Senate elections (2010, 2012, 2016), two Gubernatorial elections (2010, 2014), and two Attorney General races (2012, 2016).  The statewide party performance in these elections is shown here in Table~\ref{statewidetablePA}.

\begin{table}[ht]
\begin{tabular}{|c|c|c|c|c|c|c|c|c|c|}
\hline
 & PRES12 & PRES16 & SEN10 & SEN12 & SEN16 & GOV10 & GOV14 & ATG12 & ATG16  \\ \hline
Repub \% & 47.29 &50.35  & 51.05  &  45.44  &  50.72 &  54.52  & 45.22  & 42.58 & 48.57  \\ \hline
\end{tabular}
\caption{Overall Republican vote shares (with respect to two-party vote) for a range of recent statewide elections in Pennsylvania. }
\label{statewidetablePA}
\end{table}

\vspace{-.1in}

\subsection{Scale and zoning in PA}
\begin{figure}
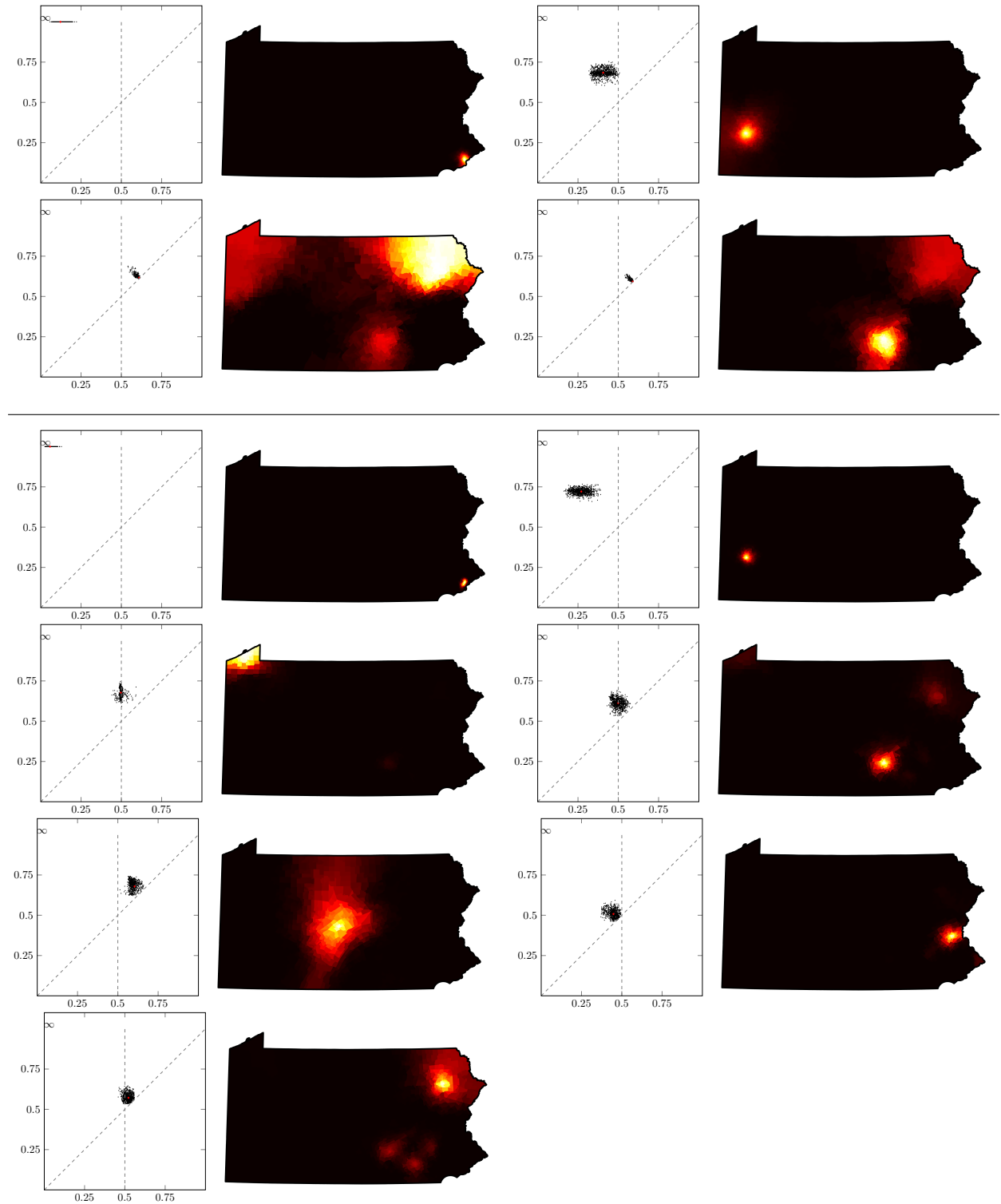

\foreach[count=\i] \n in {0,1,2,3}
{
\begin{subfigure}{0.19\textwidth}
\pdsubfig{PA_noaxes/m2m_18districts_PRES16_PD\n.png}%
\end{subfigure}
\begin{subfigure}{0.29\textwidth}
\includegraphics[width=\textwidth]{PA_plots/18mappedFrechetPRES16_\n.png}
\end{subfigure}
}
\begin{subfigure}{0.48\textwidth}
\hfill
\end{subfigure}

\vspace{.1in}

\hrule

\vspace{.1in}

\foreach[count=\i] \n in {0,1,2,3,4}
{
\begin{subfigure}{0.19\textwidth}
\pdsubfig{PA_noaxes/m2m_50districts_PRES16_PD\n.png}%
\end{subfigure}
\begin{subfigure}{0.29\textwidth}
\includegraphics[width=\textwidth]{PA_plots/50mappedFrechetPRES16_\n.png}
\end{subfigure}
}
\foreach[count=\i] \n in {5,6}
{
\begin{subfigure}{0.19\textwidth}
\pdsubfig{PA_noaxes/m2m_50districts_PRES16_PD\n.png}%
\end{subfigure}
\begin{subfigure}{0.29\textwidth}
\includegraphics[width=\textwidth]{PA_plots/50mappedFrechetPRES16_\n.png}
\end{subfigure}
}
\begin{subfigure}{0.48\textwidth}
\hfill
\end{subfigure}

\caption{Geographical localization of  Fr\'echet features in Pennsylvania Congressional and Senate plans 
$(k=18,k=50)$ with respect to PRES16 voting. The Senate features fairly clearly capture the medium-sized cities of Pennsylvania.}
\label{PAclasses18}\label{PAclasses50}
\end{figure}

We begin with partisan-neutral ensembles of maps for each choice of $k=18,50,203$.
We carry out the four steps of the scale analysis for each:  overlay, marking, localization, and zoning.
We find four zones of relative Democratic strength at the Congressional level (Figure~\ref{PAclasses18}, top), corresponding to Philadelphia, Pittsburgh, 
Scranton, and Harrisburg, in that order.  Philadelphia and Pittsburgh are highly persistent while other two zones are weaker; only those two correspond to Democratic district clusters, comprising 5.6 and 1.2 districts, respectively
(Supplemental Figure~\ref{fig:PAcluster18}).

In the state Senate ensemble (Figure~\ref{PAclasses50}, bottom), we find seven zones:  Philadelphia, Pittsburgh, Erie, Harrisburg, State College, Lancaster, and Scranton, in sequence.  It is quite valuable to note that the order as well as the membership is different between Congressional and Senate zoning.  Erie, Harrisburg, and Scranton anchor zones that will be competitive overall, while the State College zone favors Republicans and the Lancaster zone favors Democrats.  The zoning analysis indicates 14.7 Senate districts in the Philadelphia zone and 3.4 in Pittsburgh; even at this scale, the others each contribute less than one district on average.  

For the state House ensemble (Supplemental Figure~\ref{PAclasses203}), the degree of geographic coherence (localization) is weaker, but we can read  off up to eleven zones in persistence order as follows:  Philadelphia, Pittsburgh, Harrisburg, State College, Reading, Lancaster, Erie, Allentown,
Scranton, York, Hermitage.   And now all but Hermitage should expect to anchor clusters of Democratic-won districts (not just Democratic relative to their neighbors).  
We find 55.8 House districts in the Philadelphia zone---a massive bloc of the 203 seats.  

\newpage
\subsection{Comparing elections in PA}
The \Fre means for the nine elections under consideration were plotted for each scale of redistricting in Figure \ref{allFrechetPA}. 
For the Congressional ensemble, every election has two highly persistent features in its \Fre mean $\F$ (corresponding to Philadelphia and Pittsburgh). Five out of nine elections (GOV14, SEN12, PRES12, SEN10, and most prominently ATG12) have a third feature away from the diagonal. By looking at the geographic data (not shown), we can tag all of these third features as belonging to a Scranton zone.  The presence of a Scranton feature means that Scranton is separated from Philadelphia by a district redder than itself.  Since the Scranton mean points always lie below $d=.5$, the intervening districts between Scranton and Philaldelphia are still Democratic-favoring.

\begin{figure}[ht]
\centering
\foreach \k in {18, 50, 203}
{
\begin{subfigure}{0.31\textwidth}
\caption*{\k \ districts}
\resizebox{\columnwidth}{!}{%
\begin{tikzpicture}
\begin{axis}[xmin=0, xmax=1, ymin=0, ymax=1.1, clip=false,
xtick={0,0.25,0.5,0.75,1}, ytick={0,0.25,0.5,0.75}, axis equal image, axis on top
]
\addplot[scatter, 
scatter/classes={
PRES16={red},
PRES12={green!80!black}
},
only marks, mark=*, mark size=1, scatter src=explicit symbolic]
 table[x=X, y=Y, col sep=comma, meta=label]{PA_csv/\k Frechet_means_some.csv};
\addplot graphics [xmin=0, xmax=1,ymin=0, ymax=1.1] {PA_csv/\k some_means_overlay.png}; 
\addplot[gray, dashed] coordinates{(0.5,0) (0.5,1)};
\addplot[gray, dashed] coordinates{(0,0) (1,1)};
\node at (axis cs: 0.03,1.02) {$\infty$};
\legend {PRES16, PRES12};
\end{axis}
\end{tikzpicture}
}
\end{subfigure}
}
\caption{We isolate the Fr\'{e}chet means for two successive Presidential elections in PA. The mean for PRES12 has a third off-diagonal feature at Scranton, while PRES16 does not. 
A uniform partisan swing would show each red point displaced by $(.03,.03)$ relative to its paired green point, so the diagram summarizes the non-uniformity of the electoral shift.
}
\label{PRES12PRES16Frechet}
\end{figure}
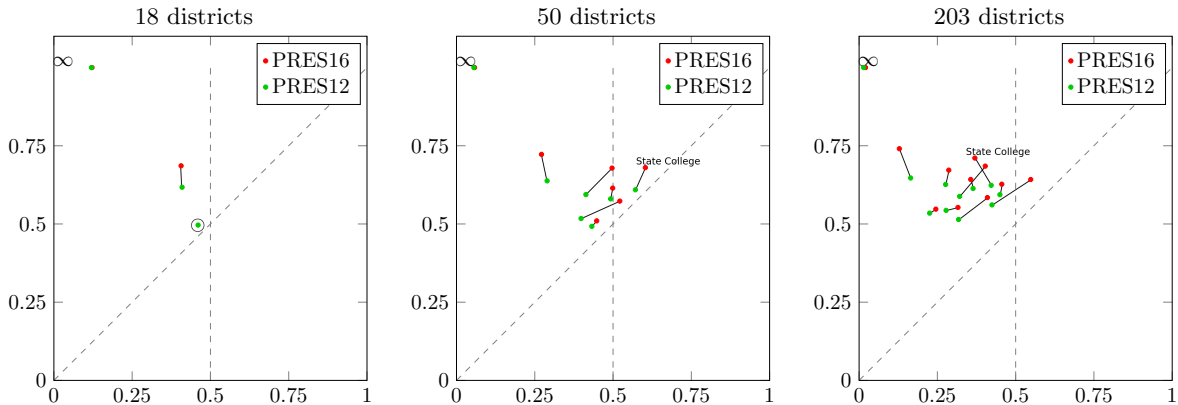

We isolate  two elections of interest, PRES12 and PRES16, for a head-to-head comparison in Figure \ref{PRES12PRES16Frechet}. We draw lines to indicate the geographic pairing between the top eleven most persistent features (this pairing differs from the optimal $L^2$ matching between the Fr\'{e}chet means).
Recall that the Republican statewide share is $.4729$ in PRES12 versus $.5035$ in PRES16, a difference of about 3 percentage points.  If the elections were related by a {\em uniform partisan swing} (in which each geographic unit has its vote share additively perturbed by the same amount), then each $b_i$ and each $d_i$ would increase by the same amount, shifting the whole diagram to the northeast with a $(.03,.03)$ translation. 
This is not at all what we observe.  We note that the Pittsburgh feature actually shifts in a Democratic direction at all three scales, most notably at the smallest (House) scale.  The difference is in the middle of the state, which separates Pittsburgh from Philadelphia.  There, the Trump election has a significantly higher $d$ value, indicating the reddening of the mid-state districts at all three scales. The State College feature is especially interesting, and is marked in the state Senate and state House plots. In the state Senate ensemble, the State College feature more or less adheres to the uniform partisan swing direction: $b$ and $d$ are both lower for PRES16 than PRES12. However, for the state House election, the trend reverses. This is because state House districts are small enough to separate the city of State College (which was more Democratic in 2016 than 2012) from the surrounding area (which was less Democratic in 2016 than 2012). Indeed, this effect can be seen in the geographic heat maps for this feature in Figures \ref{PAclasses50} and \ref{PAclasses203}---the State College feature is much more diffuse in the state Senate ensemble than in the state House ensemble. In the Congressional ensemble, the Scranton feature is present in PRES12, but not PRES16, telling us that the level of Trump preference in the Scranton district is not significantly distinct from its surroundings (in contrast with Romney support).

\subsection{Signals of gerrymandering in PA}
We now turn to the party-biased ensembles. We focus on the state Senate and the PRES16 vote for an illustrative discussion.

The matching between the two biased ensembles is shown with line segments in  Figure \ref{DRhist}, helping organize the overlaid plots in  Figure \ref{biasedmeansPA}. All the off-diagonal points were successfully paired, and the geographic pairing and optimal $L^2$ matchings were identical, indicating that the persistence structure of the two ensembles was very similar despite the significant difference in the overall number of Republican districts. 

It is instructive to compare the Pittsburgh and the Harrisburg features in the two ensembles, which show opposite effects.
The Harrisburg district is quite competitive, and the Republican ensemble shows the effect you would expect:  the district is several percentage points more Republican than in the Democratic ensemble---enough to secure a Republican win in most cases.  On the other hand, an overall Republican advantage is correlated with a more {\em Democratic} Pittsburgh district; this is the well-known phenomenon called ``packing," where wastefully high vote totals in some districts lead to overall disadvantage.
In interesting contrast to Harrisburg, the Erie  cluster was very rarely anchored by a ``safe"  district for {\em either} party, despite the weighting of the algorithm towards safe seats.   Finally, the Republican-favoring ensemble clearly ``cracks" Democratic support near  Scranton  (6th feature), while the Democratic-favoring plot resembles the neutral ensemble, indicating that in 1000 largely independent attempts, careful line-drawing near Scranton was not a frequent property of maps selected for local improvements in Democratic seats overall.  

\begin{figure}[ht]
\begin{subfigure}{0.24\textwidth}
\centering
\caption*{Philadelphia}
\pdtwosubfig{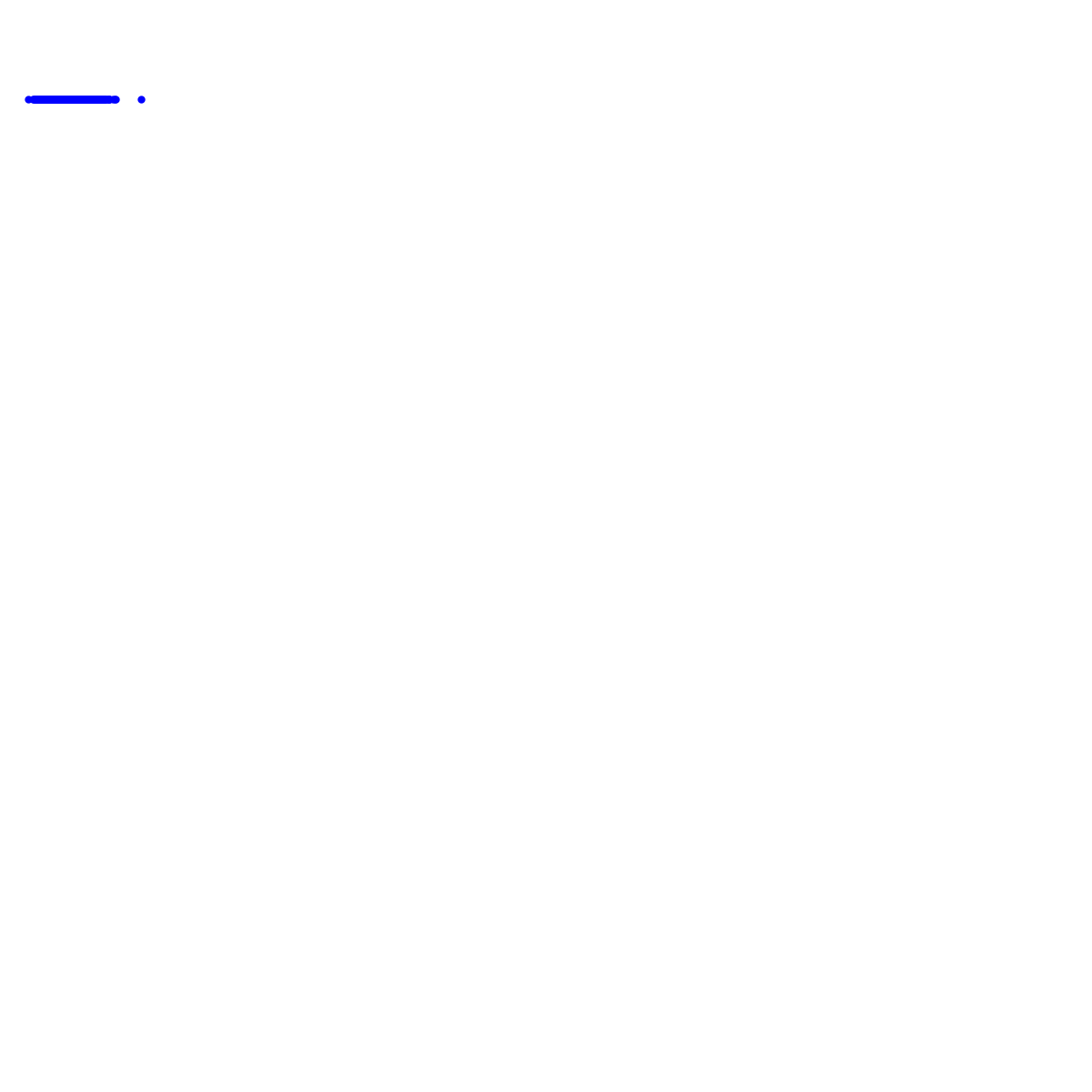}{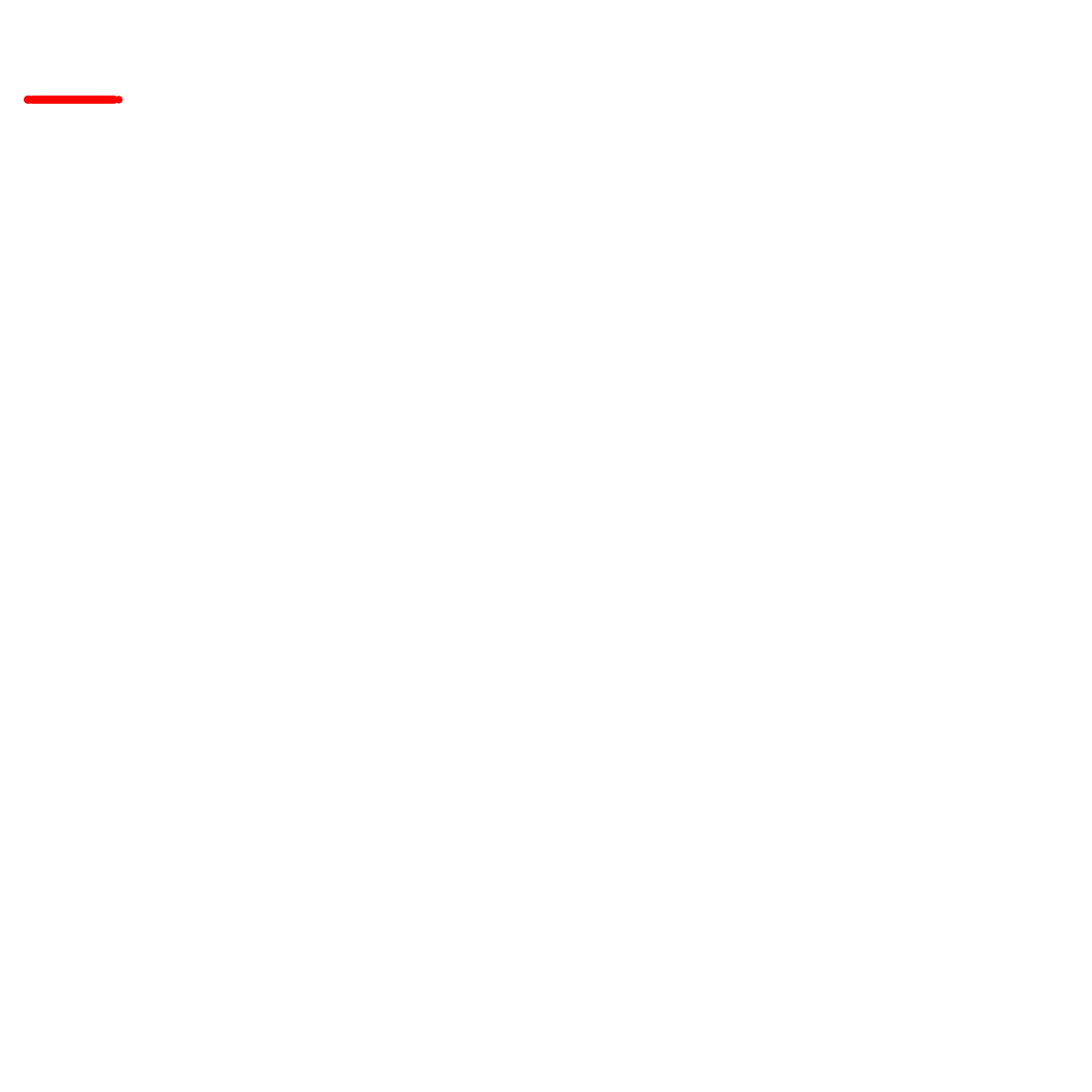}%
\end{subfigure}
\begin{subfigure}{0.24\textwidth}
\centering
\caption*{Pittsburgh}
\pdtwosubfig{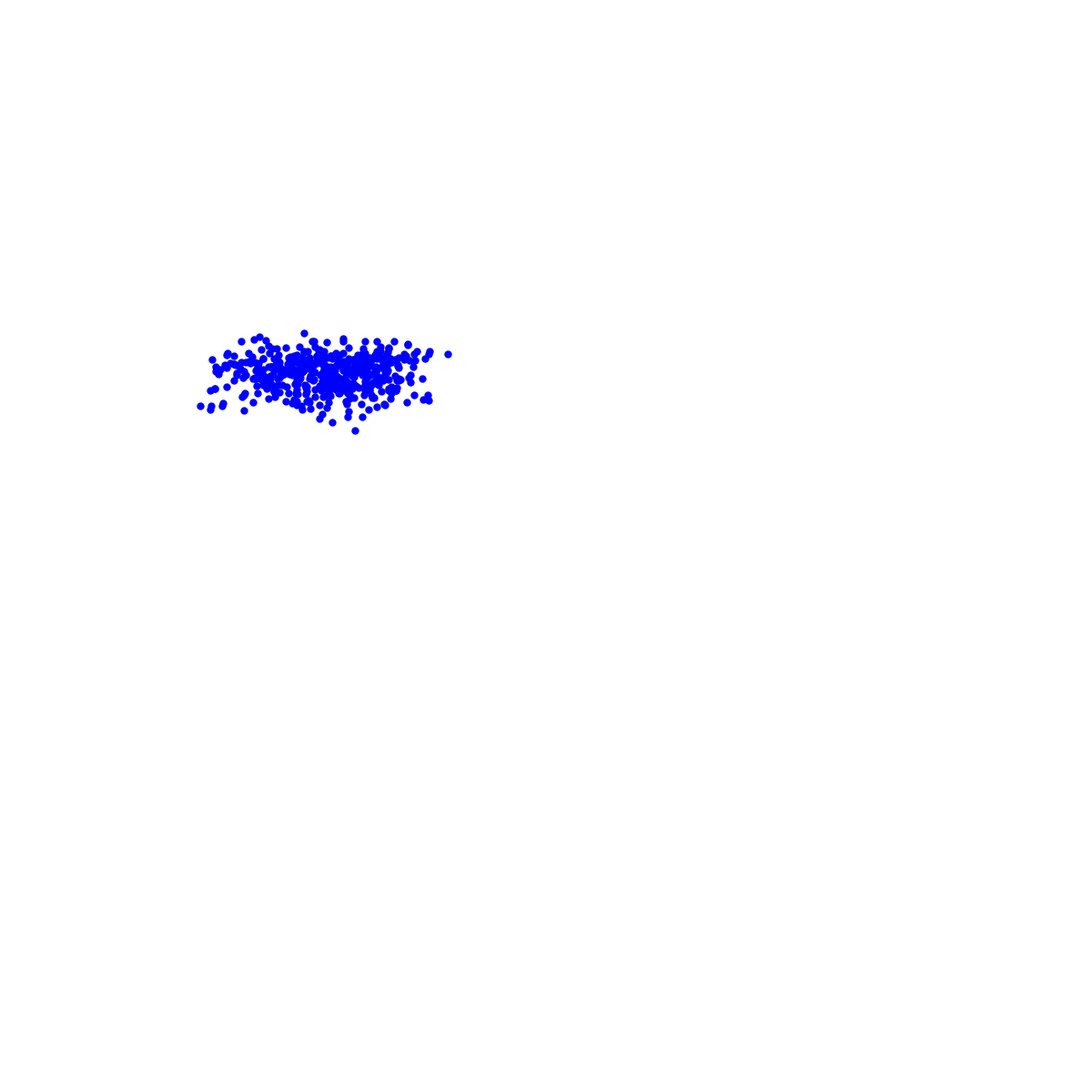}{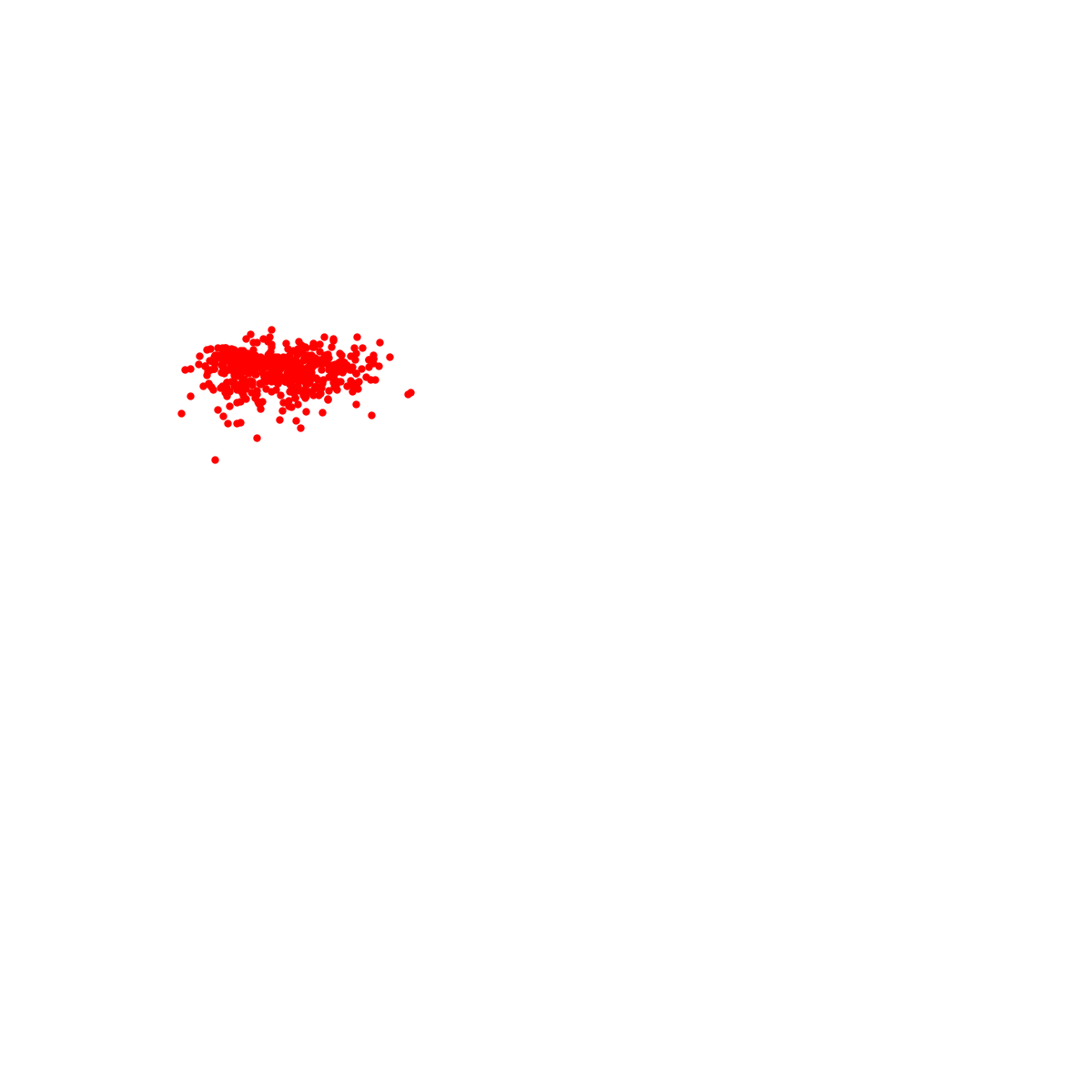}%
\end{subfigure}
\begin{subfigure}{0.24\textwidth}
\centering
\caption*{Erie}
\pdtwosubfig{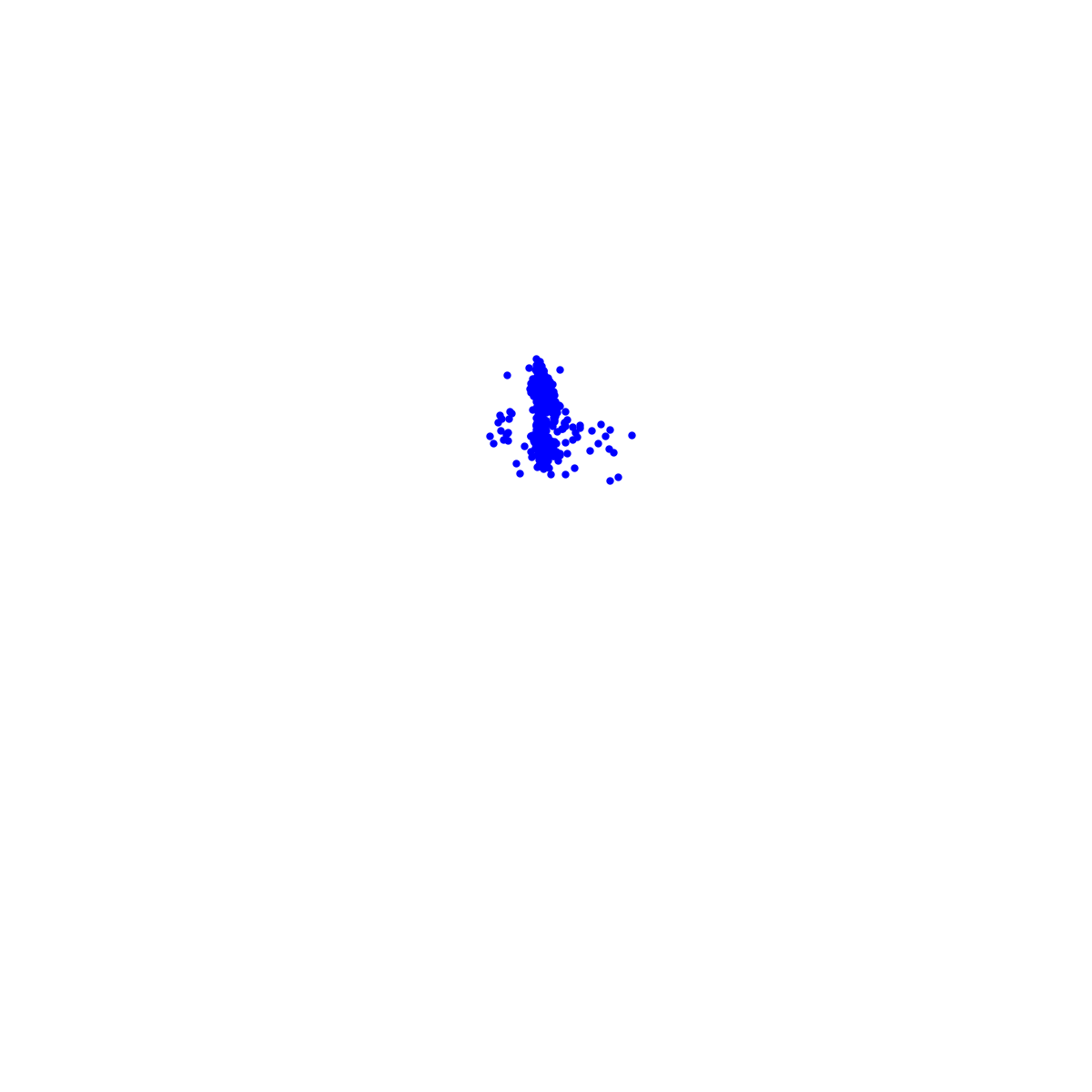}{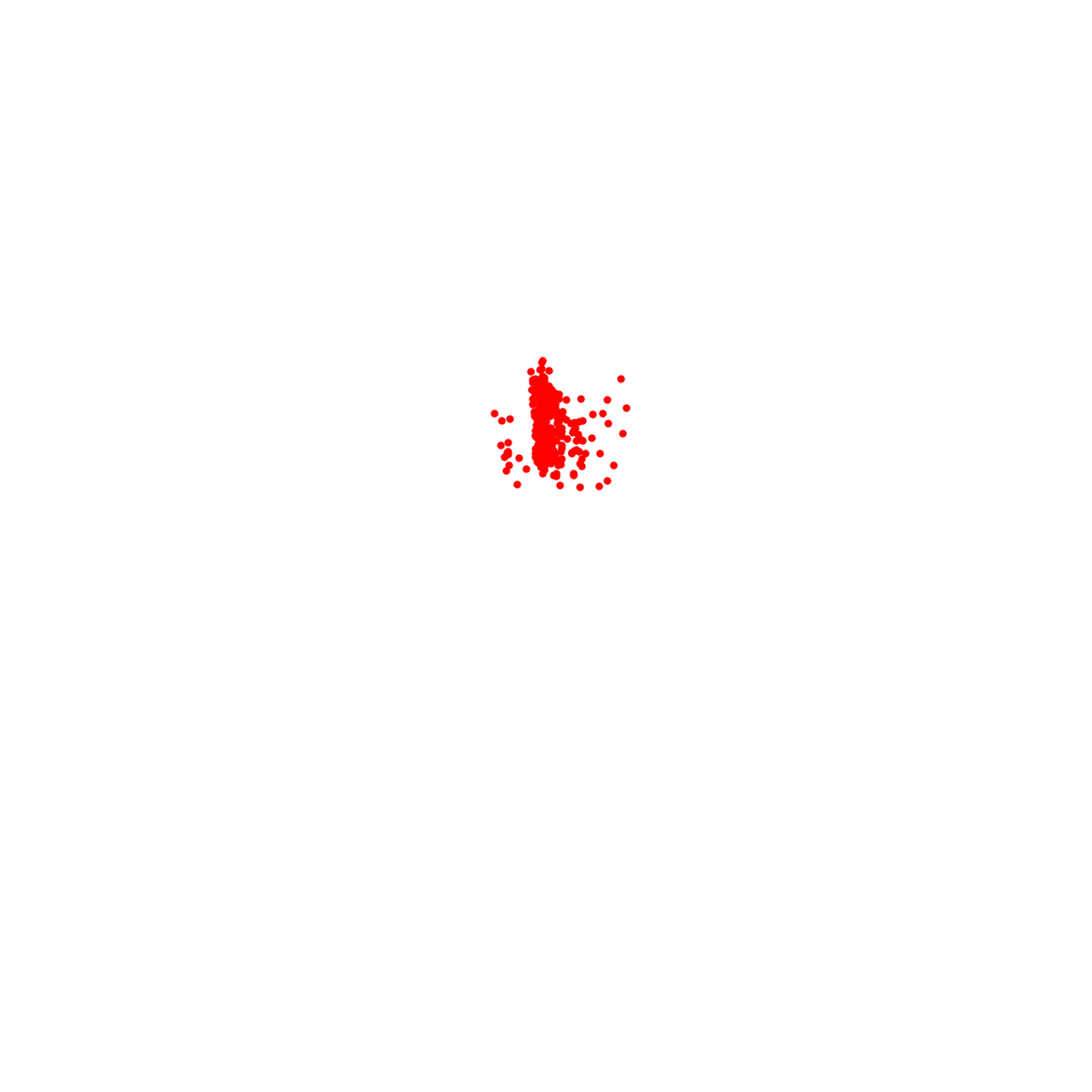}%
\end{subfigure}
\begin{subfigure}{0.24\textwidth}
\centering
\caption*{Harrisburg}
\pdtwosubfig{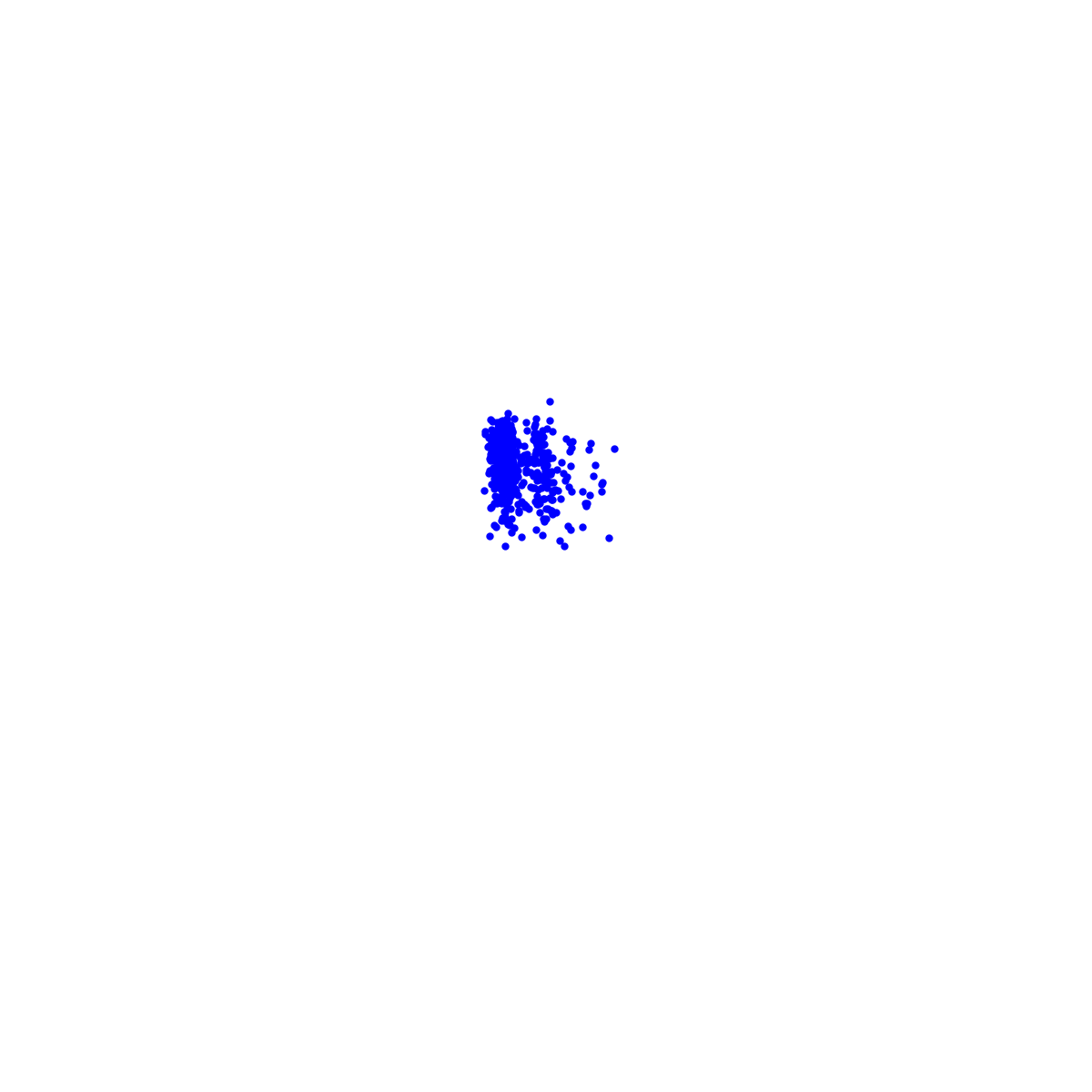}{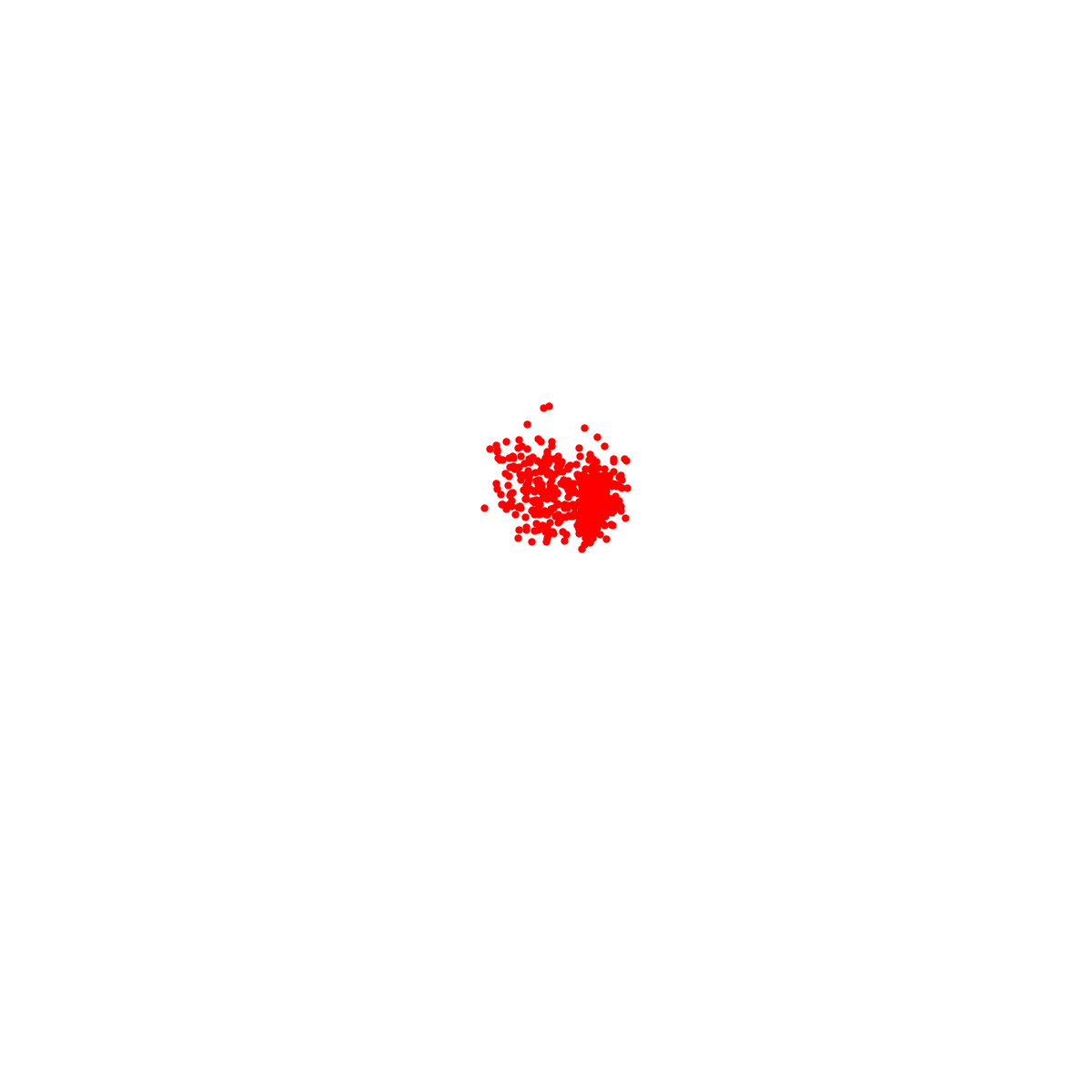}%
\end{subfigure}
\begin{subfigure}{0.24\textwidth}
\centering
\caption*{State College}
\pdtwosubfig{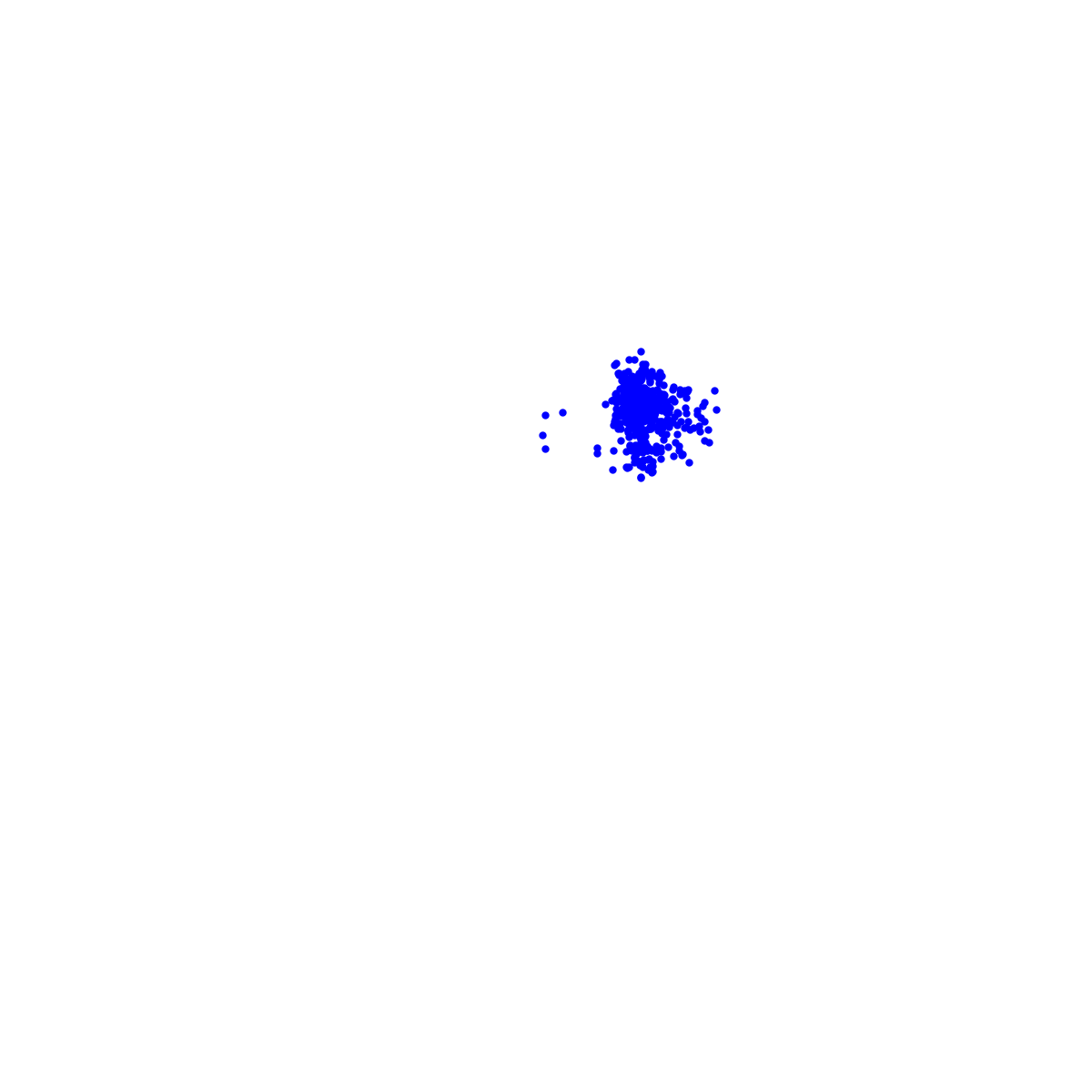}{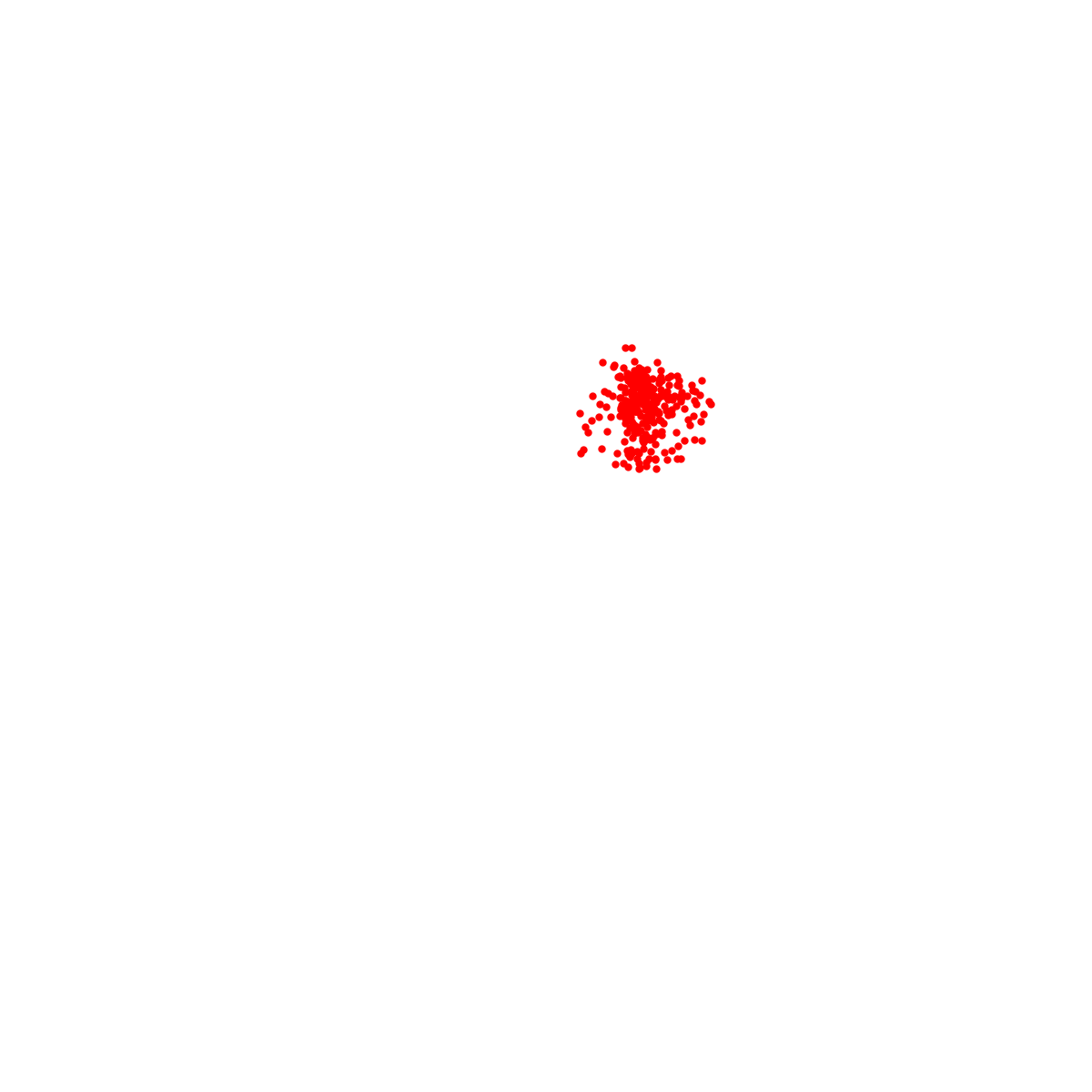}%
\end{subfigure}
\begin{subfigure}{0.24\textwidth}
\centering
\caption*{Scranton}
\pdtwosubfig{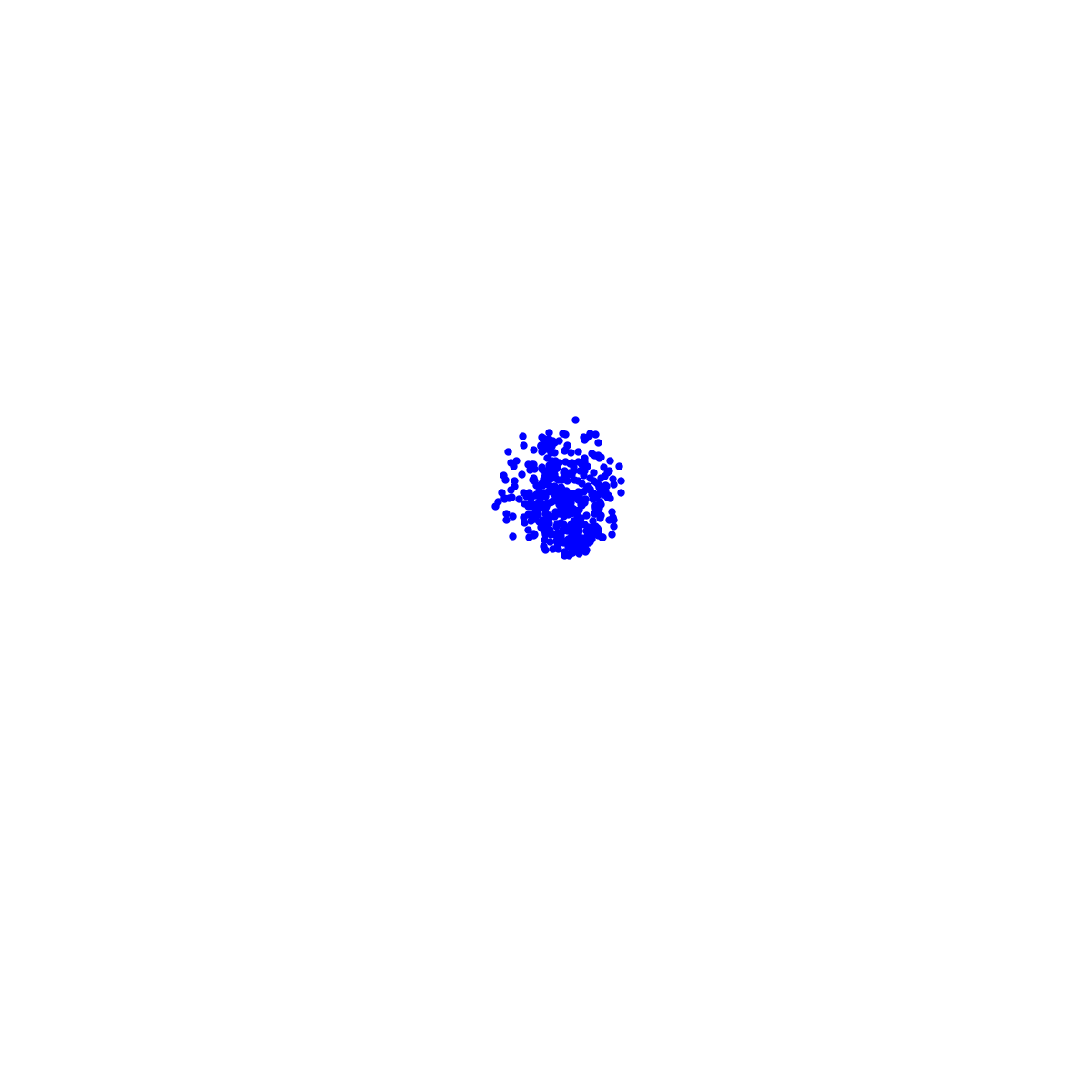}{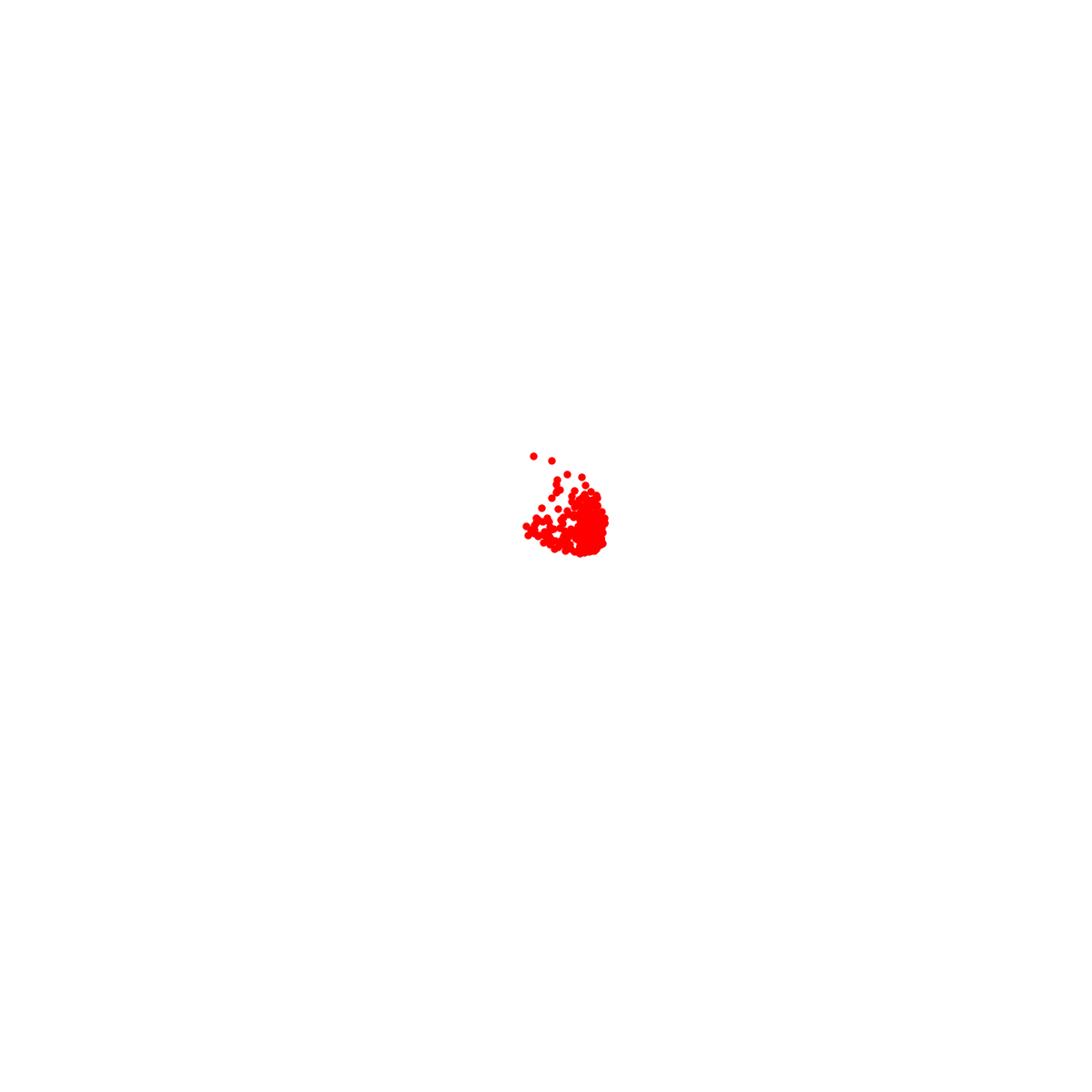}%
\end{subfigure}
\begin{subfigure}{0.24\textwidth}
\centering
\caption*{Allentown}
\pdtwosubfig{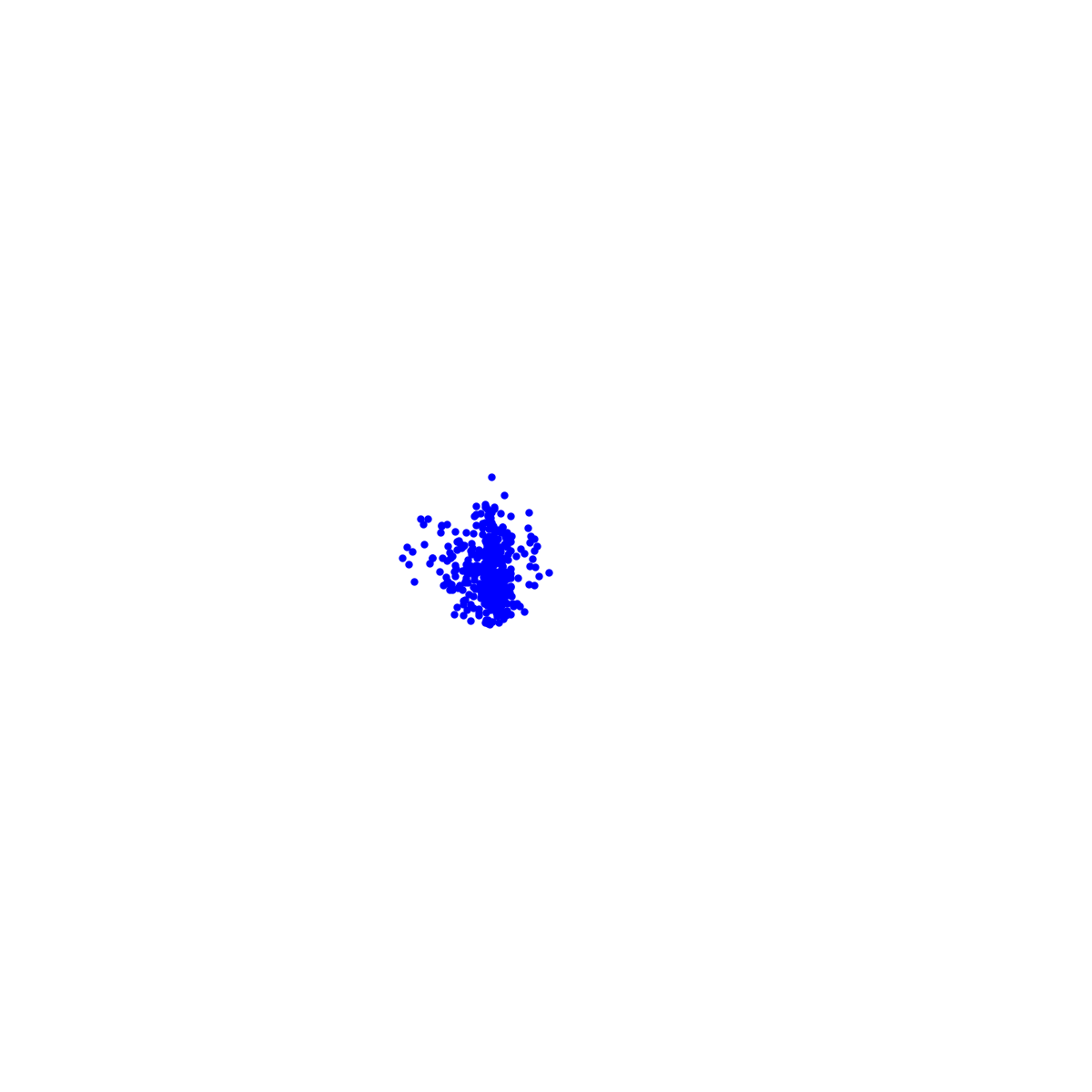}{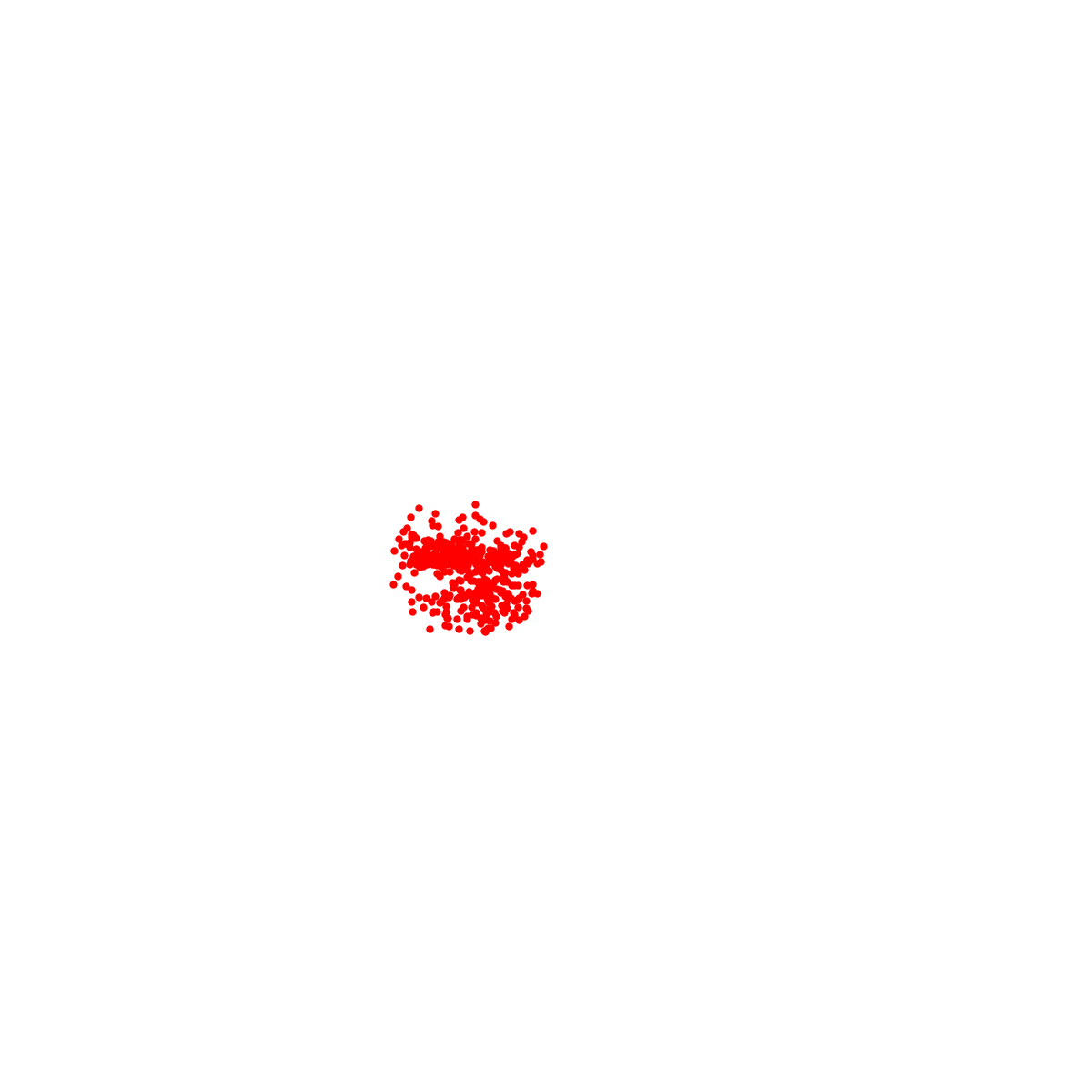}%
\end{subfigure}
\caption{Comparison of point plots for the successive Fr\'echet features in the PA Senate ensembles that are biased for Democratic and Republican safe seats, respectively. Note the contrast between the rigidity of the Erie feature (similar placement for both ensembles) and the manipulability of the Harrisburg feature, though both anchor competitive districts.}
\label{biasedmeansPAshort}
\end{figure}

\section{Case study: North Carolina}\label{sec:casestudyNC}

We turn to North Carolina, which has 13 Congressional districts, 50 state Senate districts, and 120  state House districts.  We have cleaned data for seven statewide elections in North Carolina:  Presidential elections from 2012 and 2016, the US Senate elections from 2010, 2014 and 2016, and the Gubernatorial elections from 2012 and 2016. The statewide two-way vote shares for these elections can be found in Table \ref{statewidetableNC}. 

\begin{table}[ht]
\begin{tabular}{|c|c|c|c|c|c|c|c|}
\hline
 & PRES12 & PRES16 & SEN10 & SEN14 & SEN16 & GOV12 & GOV16  \\ \hline
Repub \% & 51.08 & 51.98 & 56.02 & 49.17 & 53.02 & 55.87 & 49.95 \\ \hline
\end{tabular}
\caption{Overall Republican vote shares (with respect to two-party vote) for a range of recent statewide elections in North Carolina. }
\label{statewidetableNC}
\end{table}

\subsection{Scale and zoning in NC}
\begin{figure}[ht]
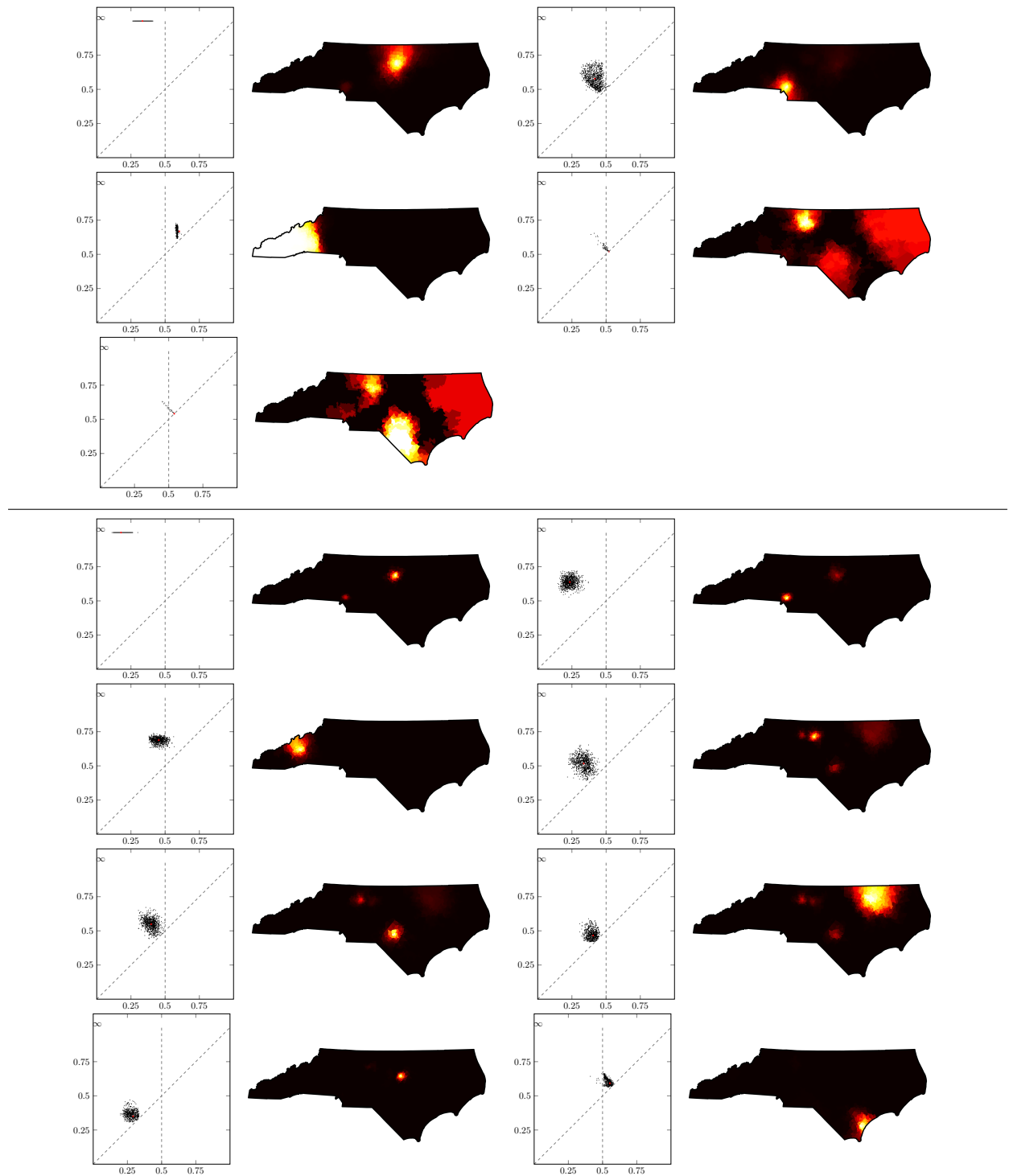
 
\foreach[count=\i] \n in {0,1,2,3,4}
{
\begin{subfigure}{0.16\textwidth}
\pdsubfig{NC_noaxes/m2m_13districts_PRES16_PD\n.png}%
\end{subfigure}
\begin{subfigure}{0.26\textwidth}
\includegraphics[width=\textwidth]{NC_plots/13mappedFrechetPRES16_\n.png}
\end{subfigure}
}
\begin{subfigure}{0.42\textwidth}
\hfill
\end{subfigure}

\medskip

\hrule

\medskip

\foreach[count=\i] \n in {0,1,2,3,4,5,6,7}
{
\begin{subfigure}{0.16\textwidth}
\pdsubfig{NC_noaxes/m2m_50districts_PRES16_PD\n.png}%
\end{subfigure}
\begin{subfigure}{0.26\textwidth}
\includegraphics[width=\textwidth]{NC_plots/50mappedFrechetPRES16_\n.png}
\end{subfigure}
}
\caption{Geographical localization of Fr\'echet features in North Carolina Congressional and Senate plans 
$(k=13,k=50)$ with respect to PRES16 voting.}
\label{NCclasses13}\label{NCclasses50}
\end{figure}

The \Fre point plots and accompanying heat maps using vote data from the 2016 Presidential election are shown in Figure \ref{NCclasses13} and Supplemental Figure \ref{NCclasses120}. 
Geographic localization shows clear success in all three cases.

At the Congressional level, North Carolina admits three clear zones for this vote data:  Durham, Charlotte, and Asheville---interestingly, Congressional districts are large enough that the Asheville district contains everything in the western tip of the state.  A fourth zone anchored in Winston-Salem is sometimes present, and there are indications of a possible fifth zone stretching south from 
Fayetteville.  These geographically more diffuse zones are often quite a bit more persistent than the corresponding zones in Pennsylvania.

In the state Senate ensemble, as with the Congressional ensemble, two highly persistent features appear near Durham and Charlotte. In contrast with the Congressional ensemble, the Asheville zone is now more geographically localized and is anchored by a Democratic-won district ($b<.5$)  in all maps.
It is notable that the fourth feature is localized near the fourth from the Congressional ensemble, but a close inspection shows that it has shifted from Winston-Salem to the nearby city of Greensboro, which is bluer but smaller.  For the first time, a rural area appears as a feature:  the heavily African-American northeast of the state anchors the sixth \Fre feature, and not the medium-sized cities of Rocky Mount, Greenville, and Wilson that are nearby.

Raleigh is represented as a  seventh Senate zone: a feature with $b<.5$ by a substantial amount, but low persistence.  This  represents a cluster anchored by very Democratic Senate districts, but that quickly merges with the Durham cluster. Its low persistence is an indicator  that it is is not a completely independent zone from Durham. 

Moving from the Senate to the  House scale, we get new (low-persistence) zones in the 9th and 10th positions corresponding to 
Wilmington and Boone.  This stands in contrast to the PA case, where the House added four  medium-persistence features compared to the  Senate. In particular, none of the small towns in the Northeast of North Carolina such as Greenville or Wilson are picked out in these plots.  Rather, they are embedded in the large Northeast zone.  This is not solely attributable to size:  each of Greenville, NC and Reading, PA has the population of 1--1.5 House districts in its respective state.  Rather, it reflects a difference in the degree of urban concentration in the Democratic vote.

 Since TDA is comparing vote totals aggregated to districts, it is not automatic that it would detect this rural Democratic strength in North Carolina, but its signature is unmistakable at the Senate and House scale, and already visible in Congressional ensembles.

\subsection{Comparing elections in NC}

This time, we will compare a larger set of elections:  all seven statewide elections in our dataset, ranging from 
2010 to 2016.  Since the association of regions with Fr\'echet mean points was successful in the last section, we can again use the plots to summarize whether  peaks of Democratic support swing in sync with their periphery between elections.  

 For the most persistent feature (the Durham zone, at  $d=\infty$), the spacing shows you the difference in Democratic
support between elections.  If voting preferences swung linearly and uniformly around the state, the displacement of the colored dots would be diagonal, with both coordinates displacing by the same amount as the first feature.  
So we learn something nontrivial from this summary diagram:  vote patterns in Charlotte and surroundings swing much more linearly than in Asheville.  In particular, Donald Trump overperformed both in heavily Democratic Asheville and (even more) in its Republican periphery, relative to a uniform swing prediction.  More generally, Asheville's partisan preference is unpredictably related to the rest of Western North Carolina.  

\begin{figure}[ht] 
\centering
\foreach \k in {13, 50, 120}
{
\begin{subfigure}{0.31\textwidth}
\caption*{\k \ districts}
\centering
\includegraphics[width=\textwidth]{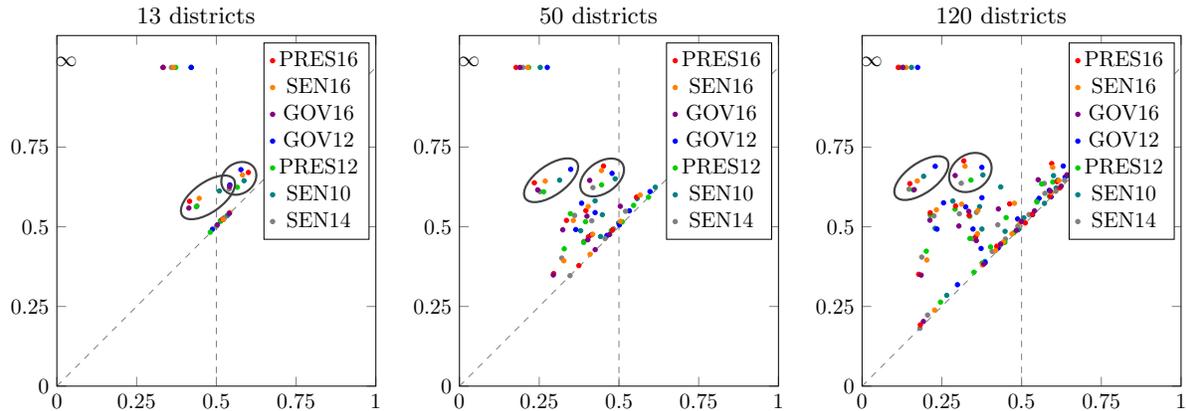}
\end{subfigure}
}
\caption{Overlaid Fr\'echet means for various elections in NC. In all three levels of redistricting, the geography-to-voting signature is strong enough to label the second and third features, shown circled, as Charlotte (left) and Asheville (right).}
\label{allFrechetNC}
\end{figure}

\subsection{Signals of gerrymandering in NC}

\begin{figure}[ht]
\centering
\begin{minipage}{.25\textwidth}
\foreach \e in {PRES16}
{
\begin{subfigure}{\textwidth}
\resizebox{\columnwidth}{!}{%
\begin{tikzpicture}
\begin{axis}[xmin=0, xmax=1, ymin=0, ymax=1.1, clip=false,
xtick={0,0.25,0.5,0.75,1}, ytick={0,0.25,0.5,0.75},
axis equal image, enlargelimits=false
]
\addplot[scatter, 
scatter/classes={
2016={red},
2012={green!80!black},
mean={black}
},
only marks, mark=*, mark size=1.5, scatter src=explicit symbolic]
 table[x=X, y=Y, col sep=comma, meta=label]{NC_csv/NC_human_Frechet\e .csv};
\addplot[gray, dashed] coordinates{(0.5,0) (0.5,1)};
\addplot[gray, dashed] coordinates{(0,0) (1,1)};
\node at (axis cs: 0.03,1.02) {$\infty$};
\legend {2016, 2012, mean};
\end{axis}
\end{tikzpicture}
}
\end{subfigure}}
\end{minipage}
\begin{minipage}{.7\textwidth}

\begin{subfigure}{0.5\textwidth}
\centering
\includegraphics[width=\textwidth]{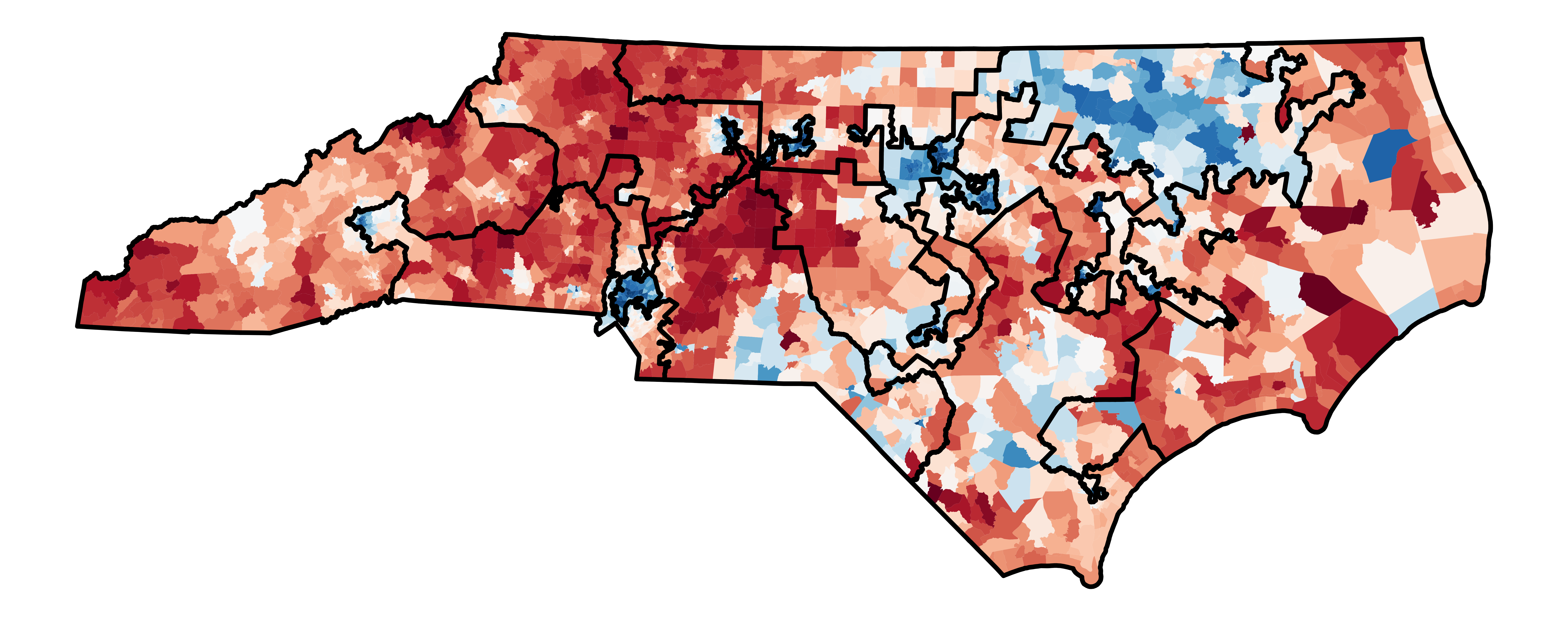}
\caption{2012 boundaries}
\end{subfigure}
\begin{subfigure}{0.5\textwidth}
\centering
\includegraphics[width=\textwidth]{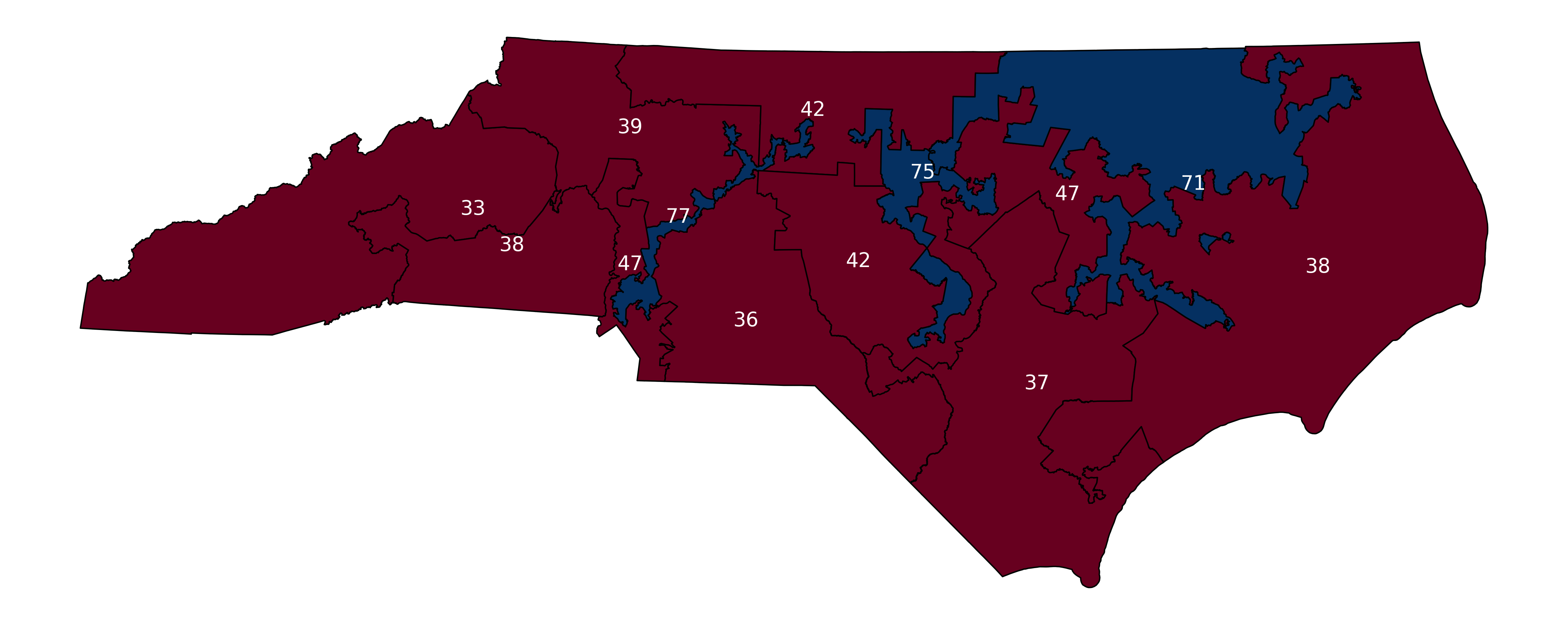}
\caption{2012 vote shares}
\end{subfigure}

\begin{subfigure}{0.5\textwidth}
\centering
\includegraphics[width=\textwidth]{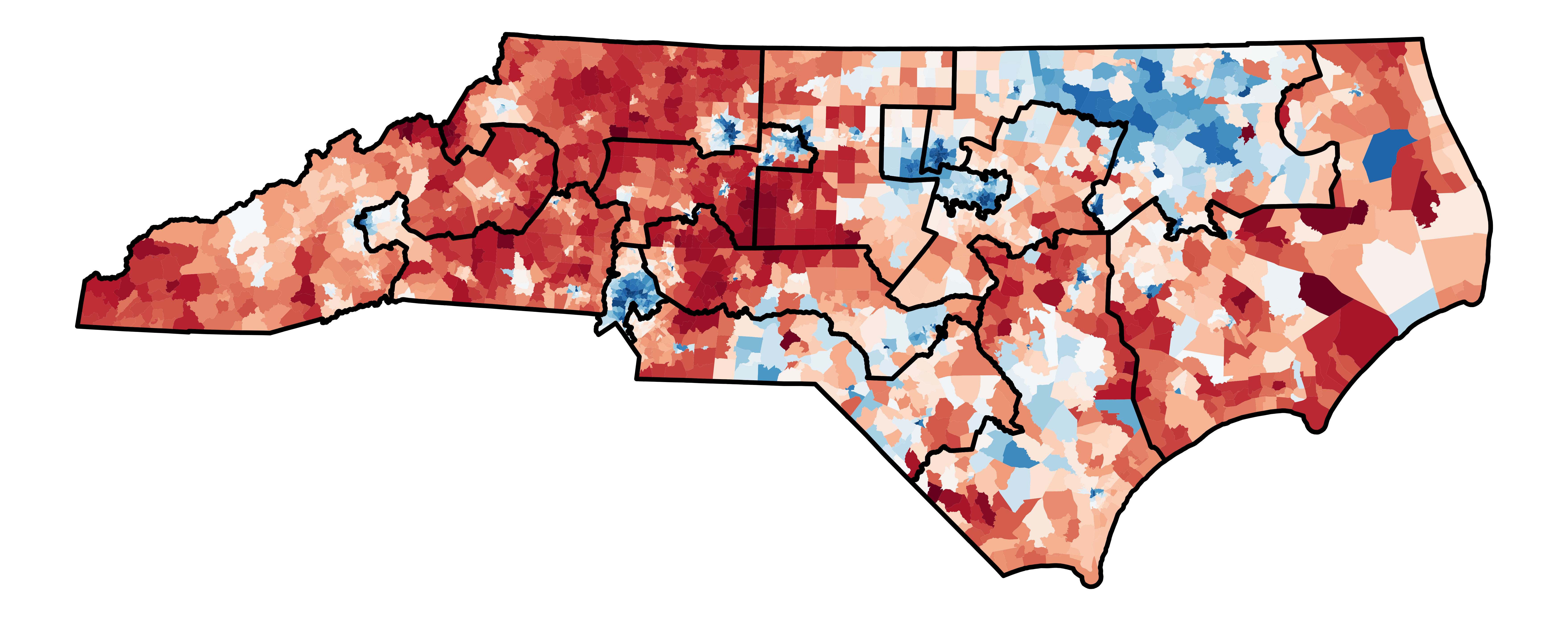}
\caption{2016 boundaries}
\end{subfigure}
\begin{subfigure}{0.5\textwidth}
\centering
\includegraphics[width=\textwidth]{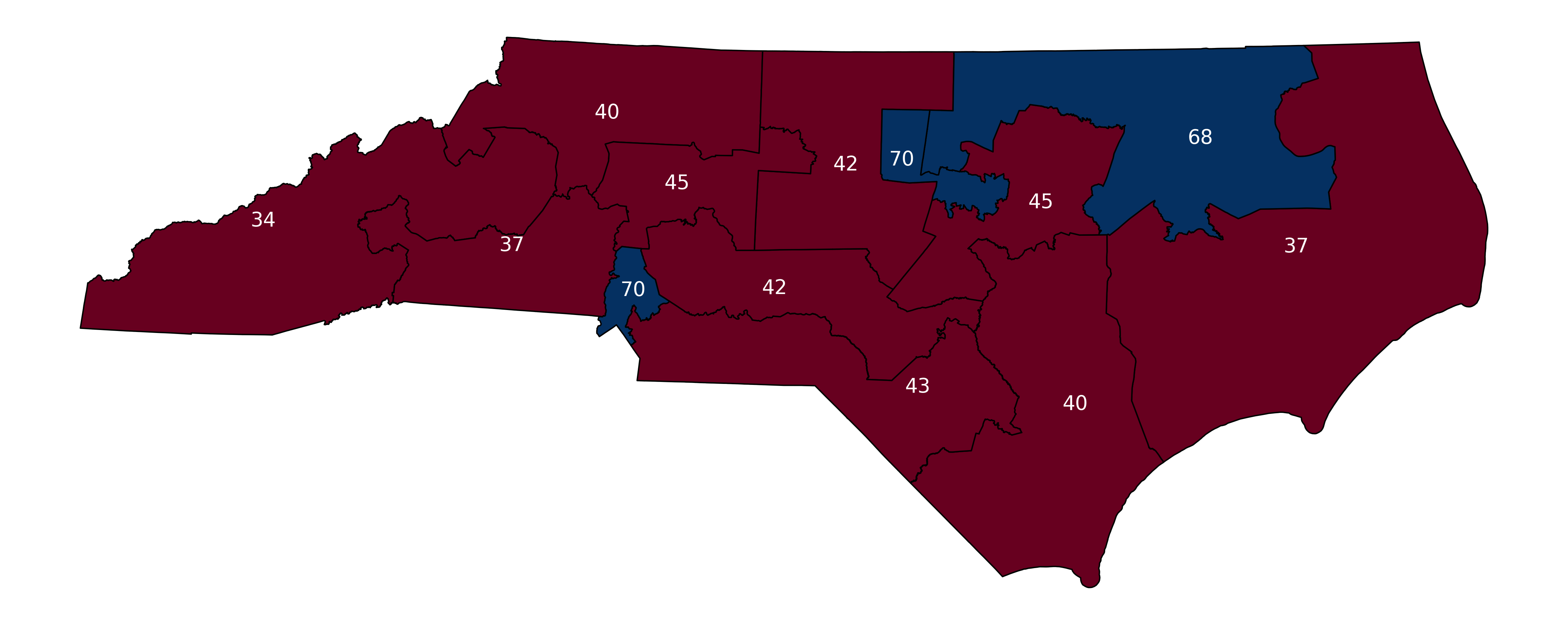}
\caption{2016 vote shares}
\end{subfigure}
\end{minipage}

\caption{Comparing two enacted Congressional plans for North Carolina with the ensemble mean, against the PRES16 vote pattern in all three cases.}
\label{NCmeans}\label{NChumanplans}
\end{figure}
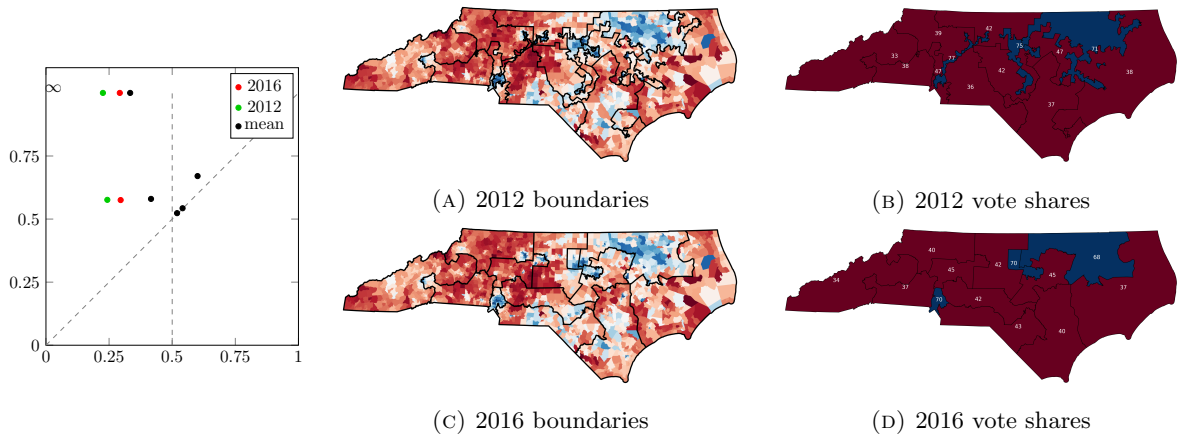

We begin by evaluating the two enacted Congressional maps for North Carolina (2012 , 2016) against the PRES16 vote pattern (Figure \ref{NCmeans}). 
Since both plans were identified by courts as Republican gerrymanders, we treat that as a reliable label.  We overlay the persistence diagrams for the two enacted plans with the \Fre means for the partisan-neutral ensemble of Congressional plans.  

The features for the enacted plans have significantly lower $b$ values than the ensemble, indicating especially Democratic districts anchoring two distinct clusters. The enacted plans each have only these two features, while the mean---and over 80\% of the maps in the neutral ensemble---includes a third feature corresponding to Asheville. This is because both enacted plans slice a large part of Asheville into a district that avoids the moderate partisan areas nearby; this district touches the Charlotte district and so becomes part of the Charlotte zone.
Note the suggestively carved boundaries in Figure~\ref{NCmeans}A,C.
This finding squares with public perception:  Asheville has been at the center of gerrymandering protests and activism for some time.

As in Pennsylvania, we close by examining  the party-biased ensembles with respect to state Senate districts and  PRES16 voting.
Figure \ref{DRhist}D shows the matching between the \Fre means used to organize the \Fre point plots in Figure \ref{biasedmeansNC} (again the geographic and optimal $L^2$ matchings agreed).

\begin{figure}[ht]
\begin{subfigure}{0.19\textwidth}
\centering
\caption*{Durham}
\pdtwosubfig{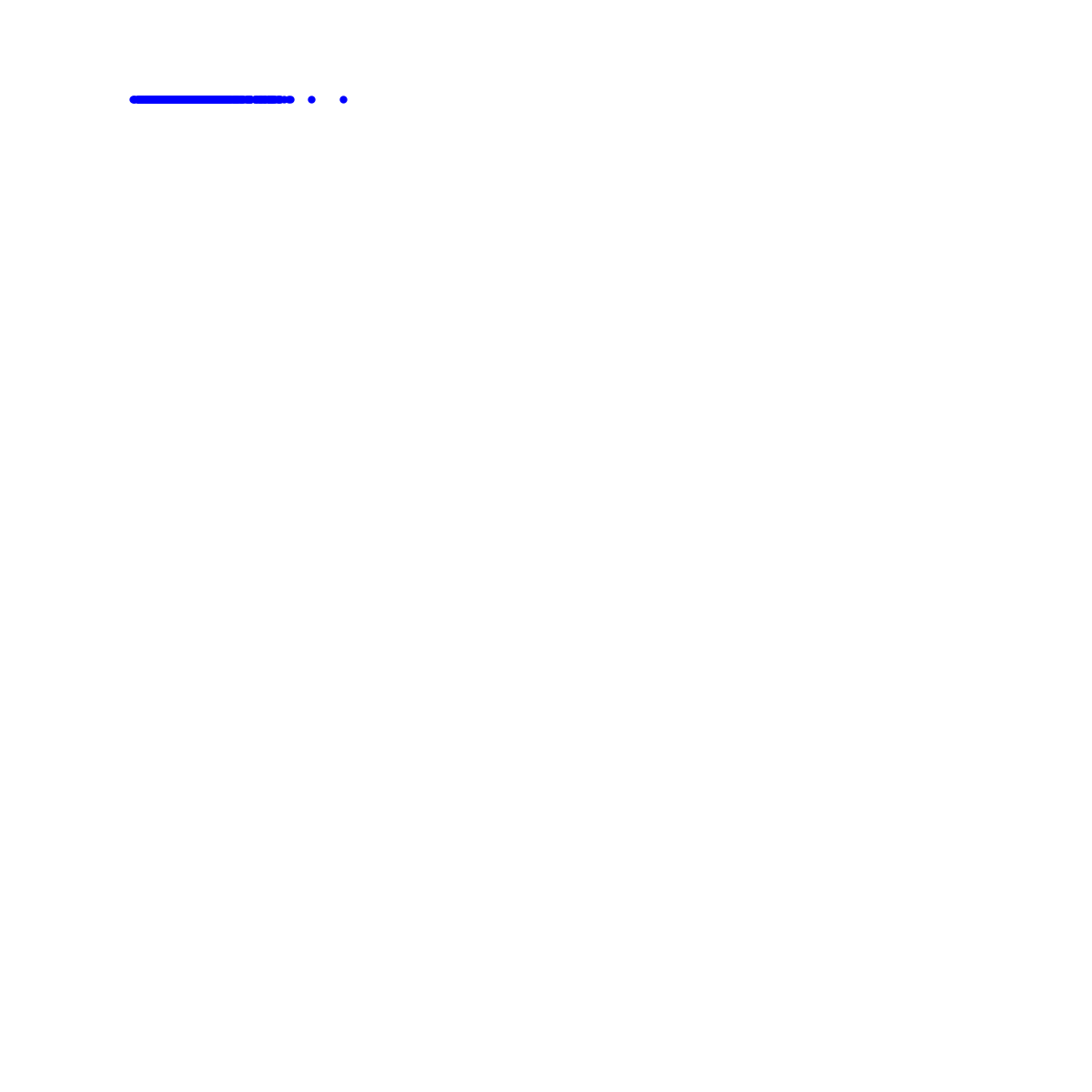}{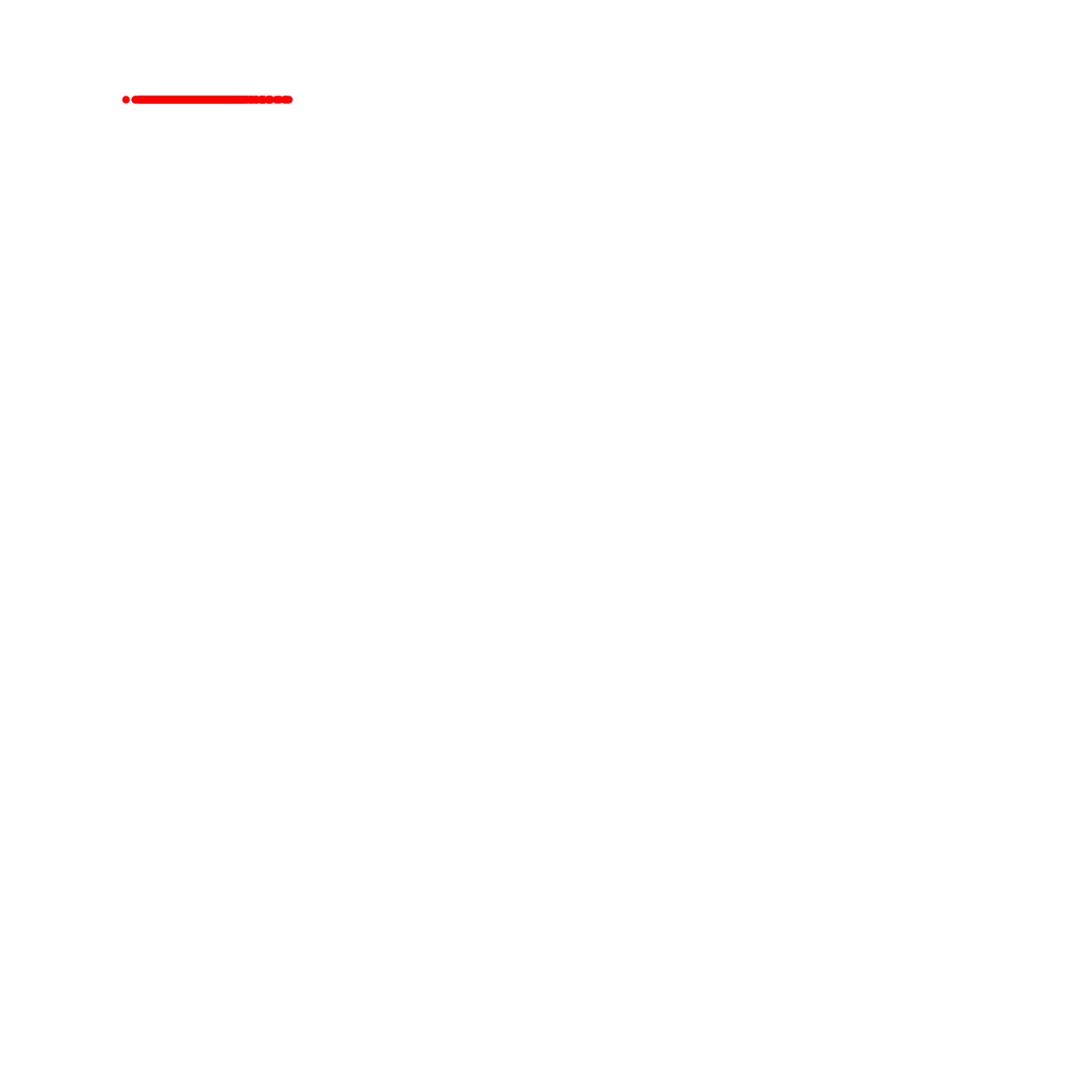}%
\end{subfigure}
\begin{subfigure}{0.19\textwidth}
\centering
\caption*{Charlotte}
\pdtwosubfig{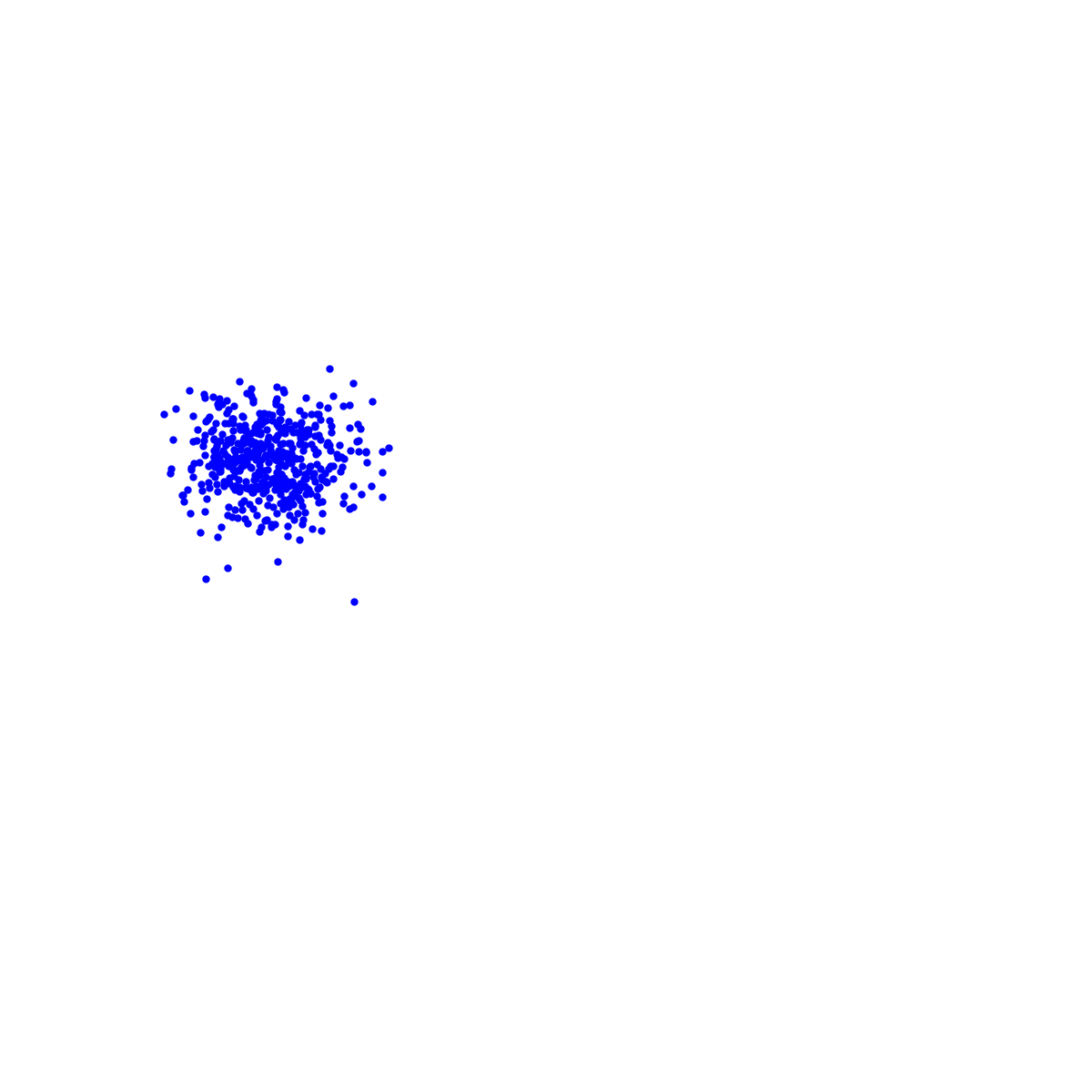}{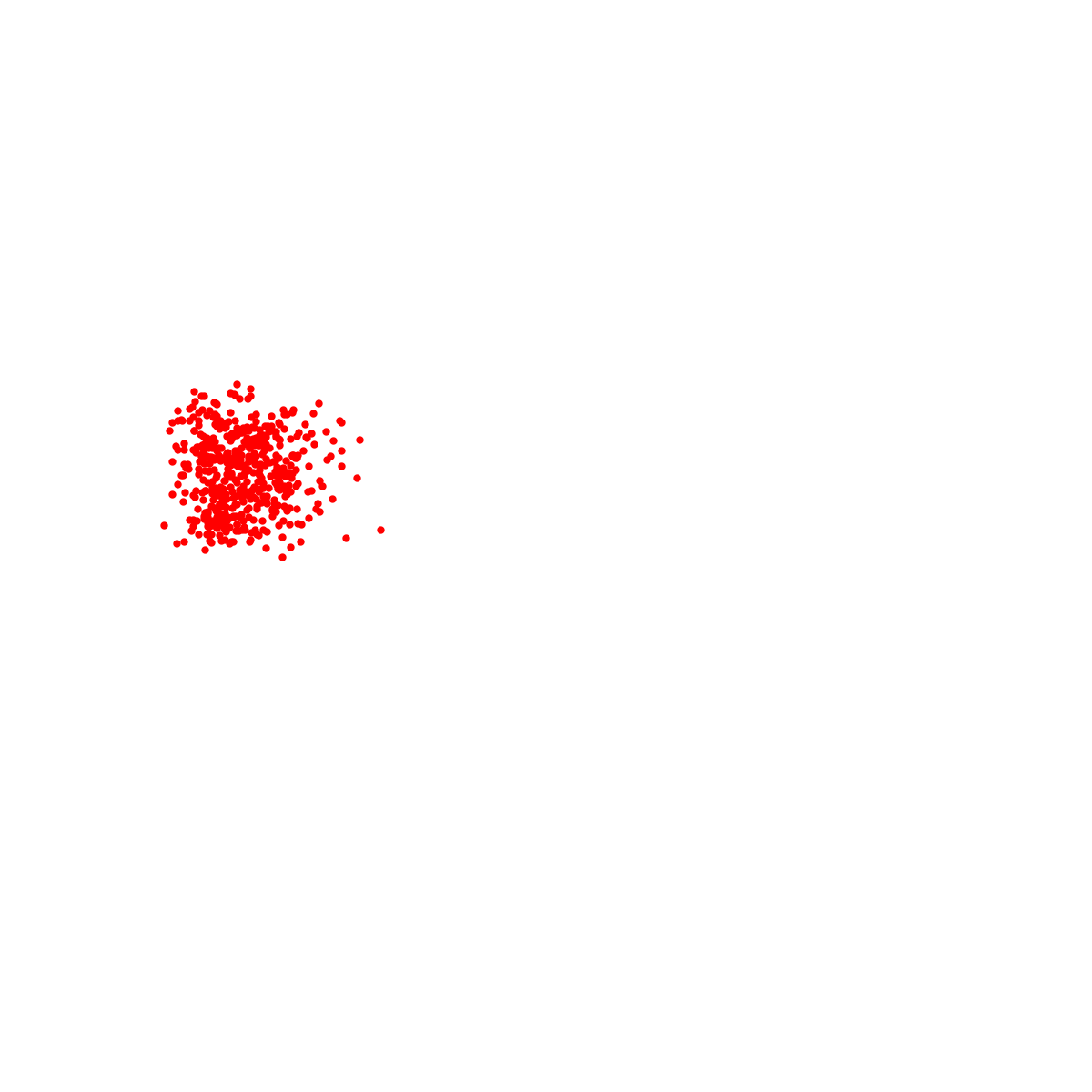}%
\end{subfigure}
\begin{subfigure}{0.19\textwidth}
\centering
\caption*{Asheville}
\pdtwosubfig{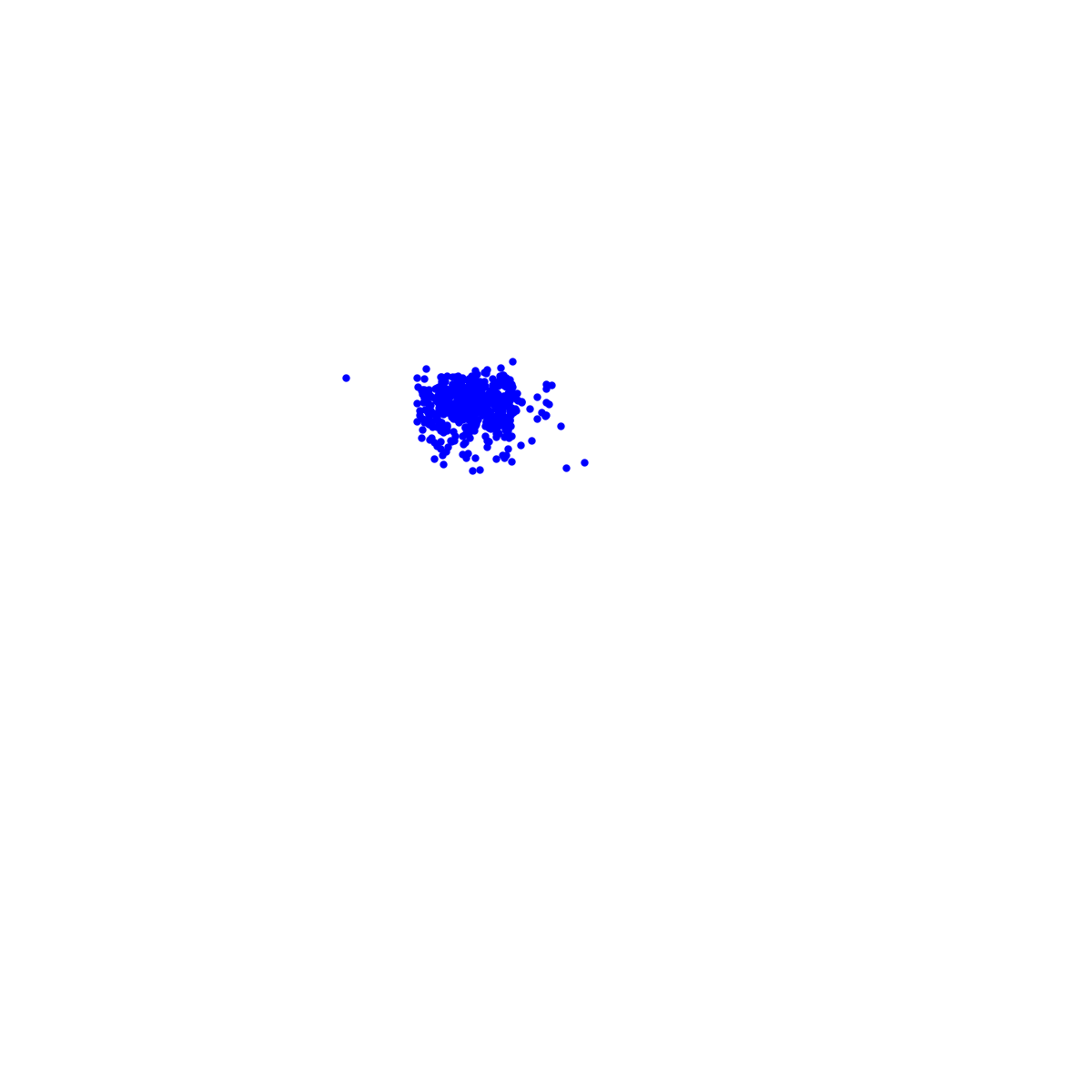}{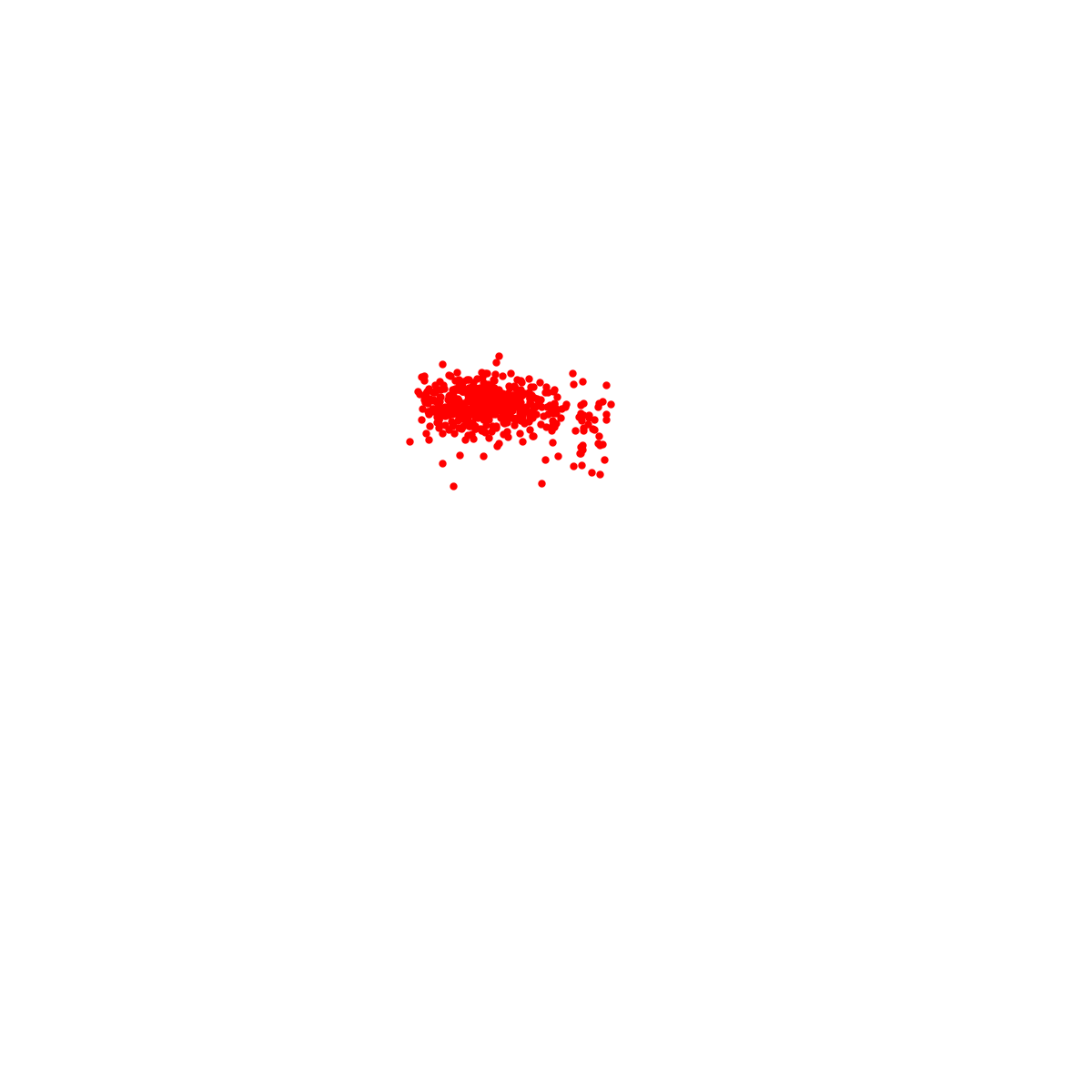}%
\end{subfigure}
\begin{subfigure}{0.19\textwidth}
\centering
\caption*{Fayetteville}
\pdtwosubfig{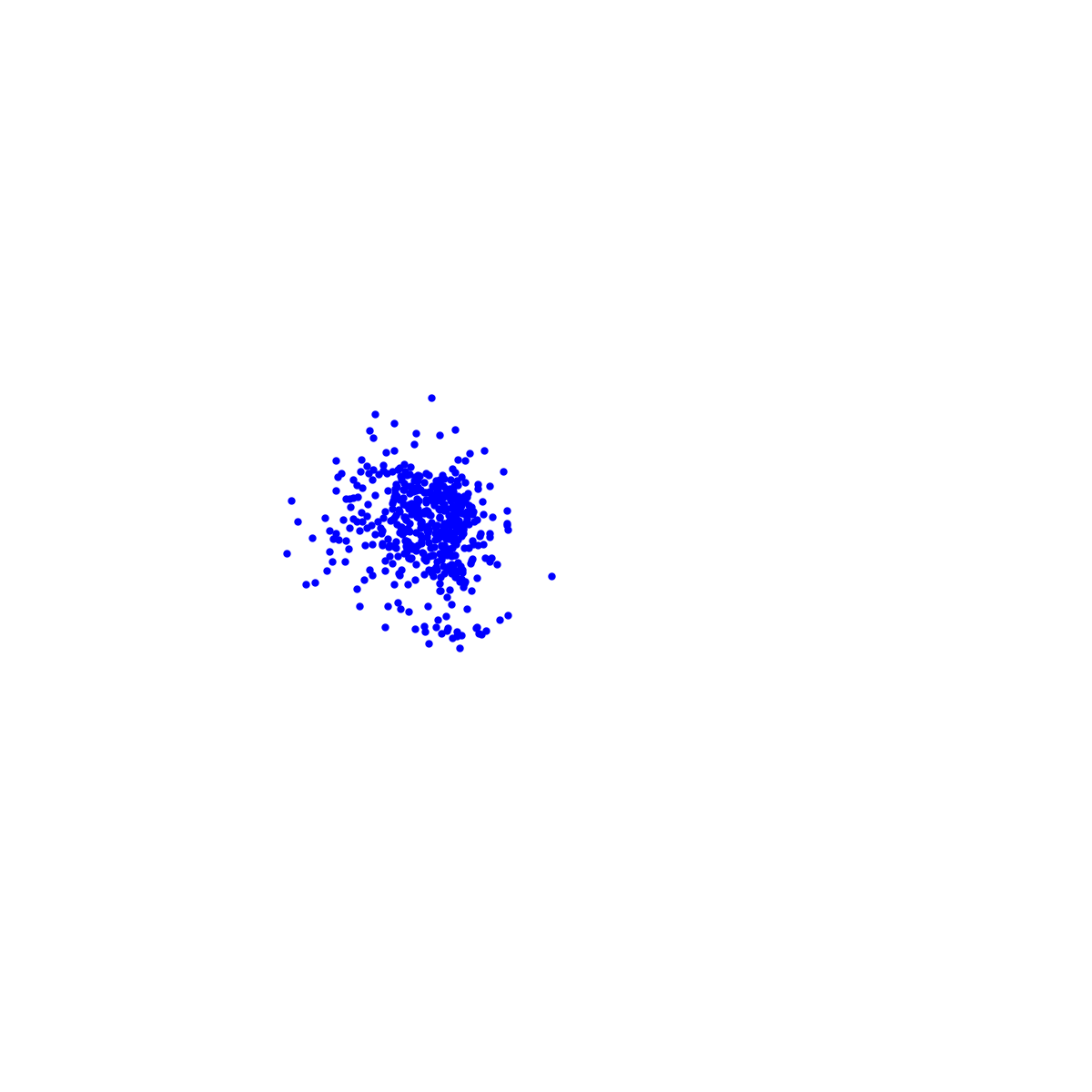}{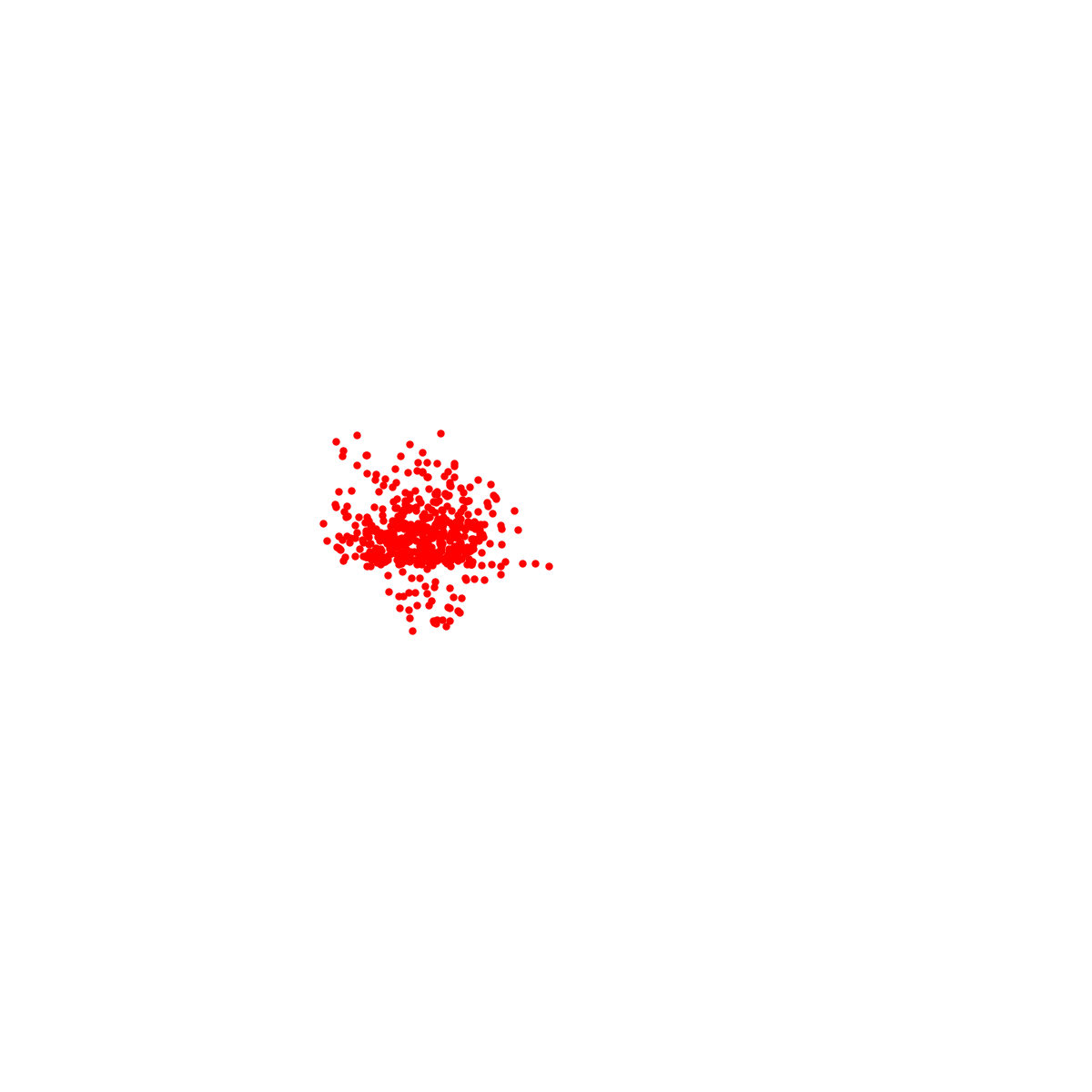}%
\end{subfigure}
\begin{subfigure}{0.19\textwidth}
\centering
\caption*{Greensboro}
\pdtwosubfig{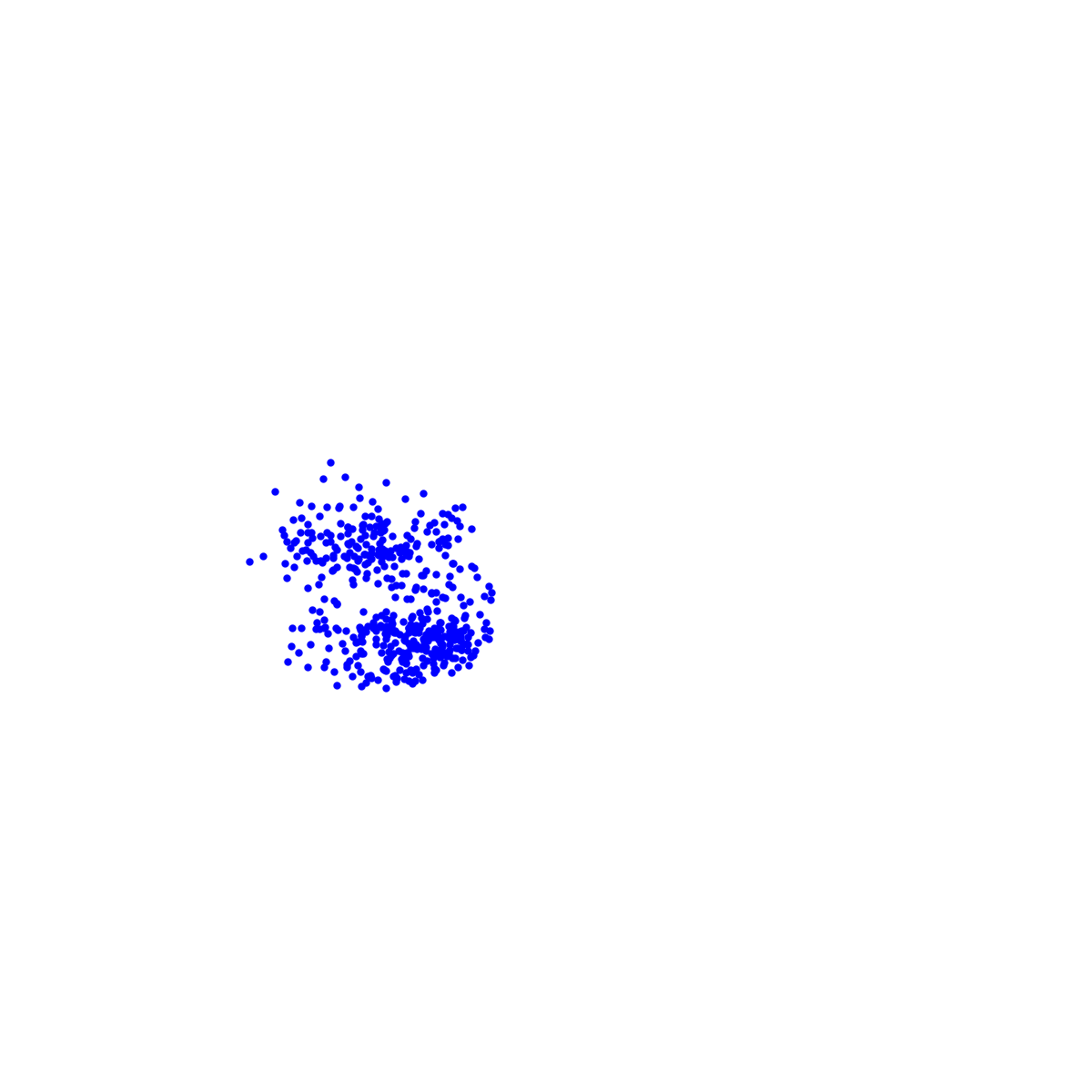}{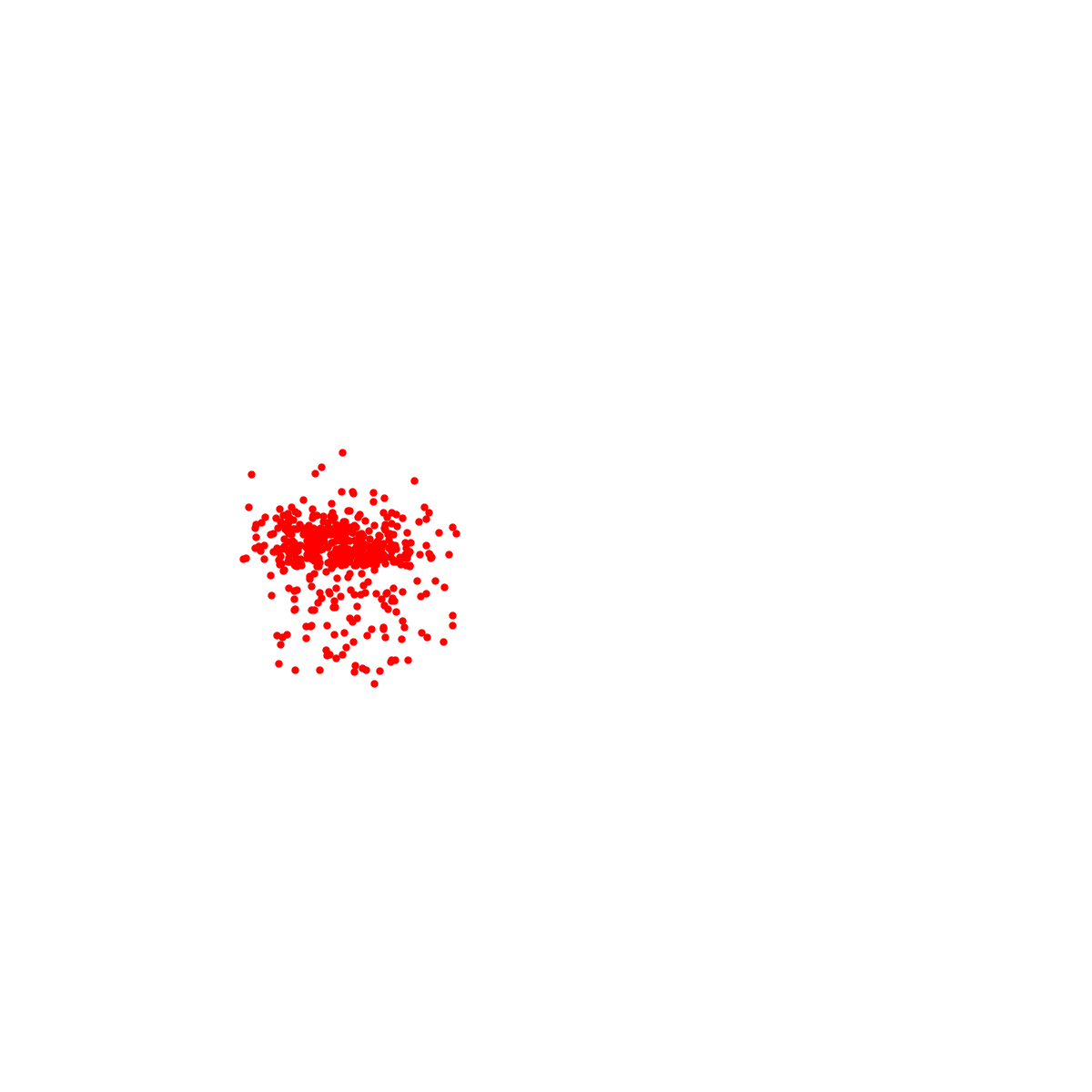}%
\end{subfigure}
\begin{subfigure}{0.19\textwidth}
\centering
\caption*{Wilmington}
\pdtwosubfig{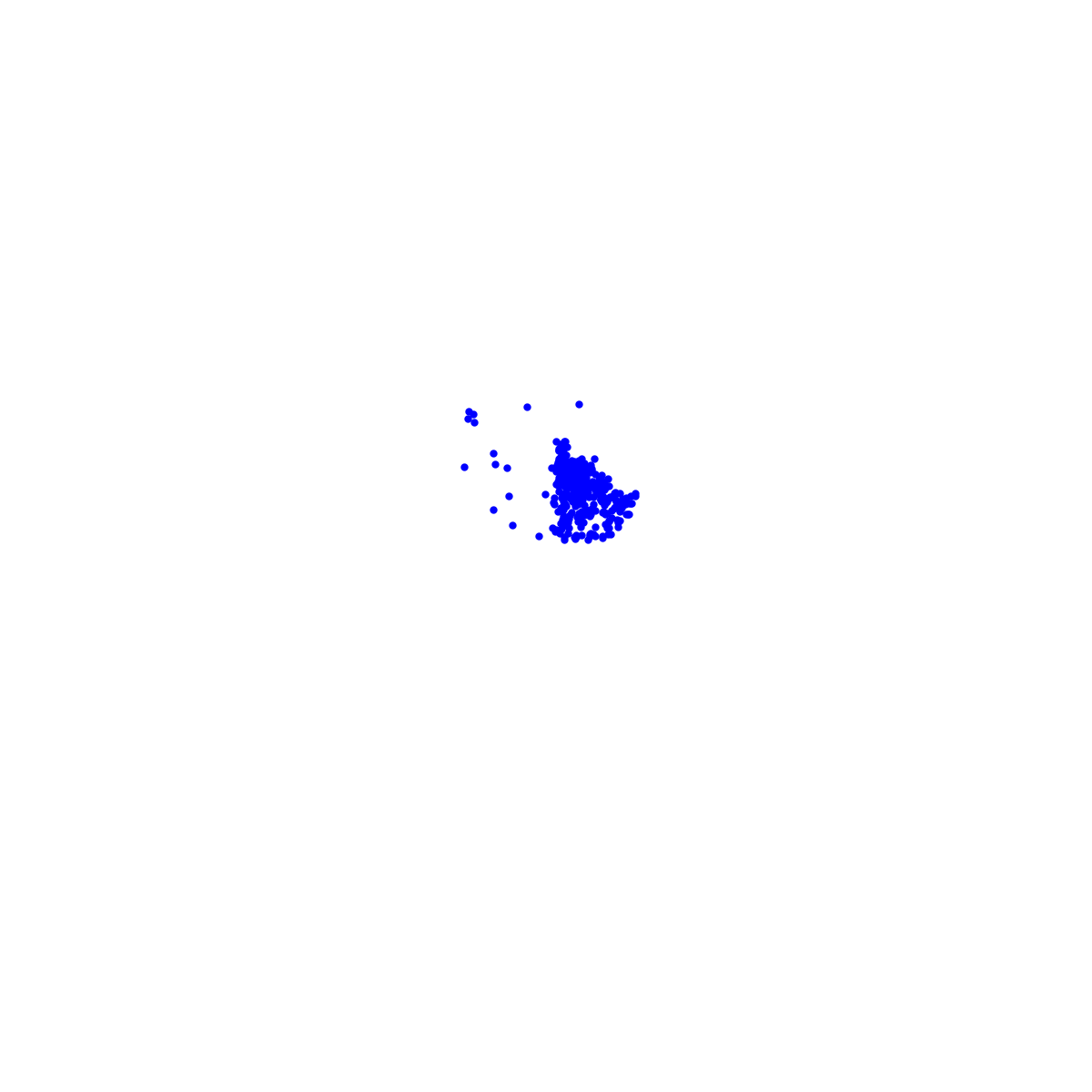}{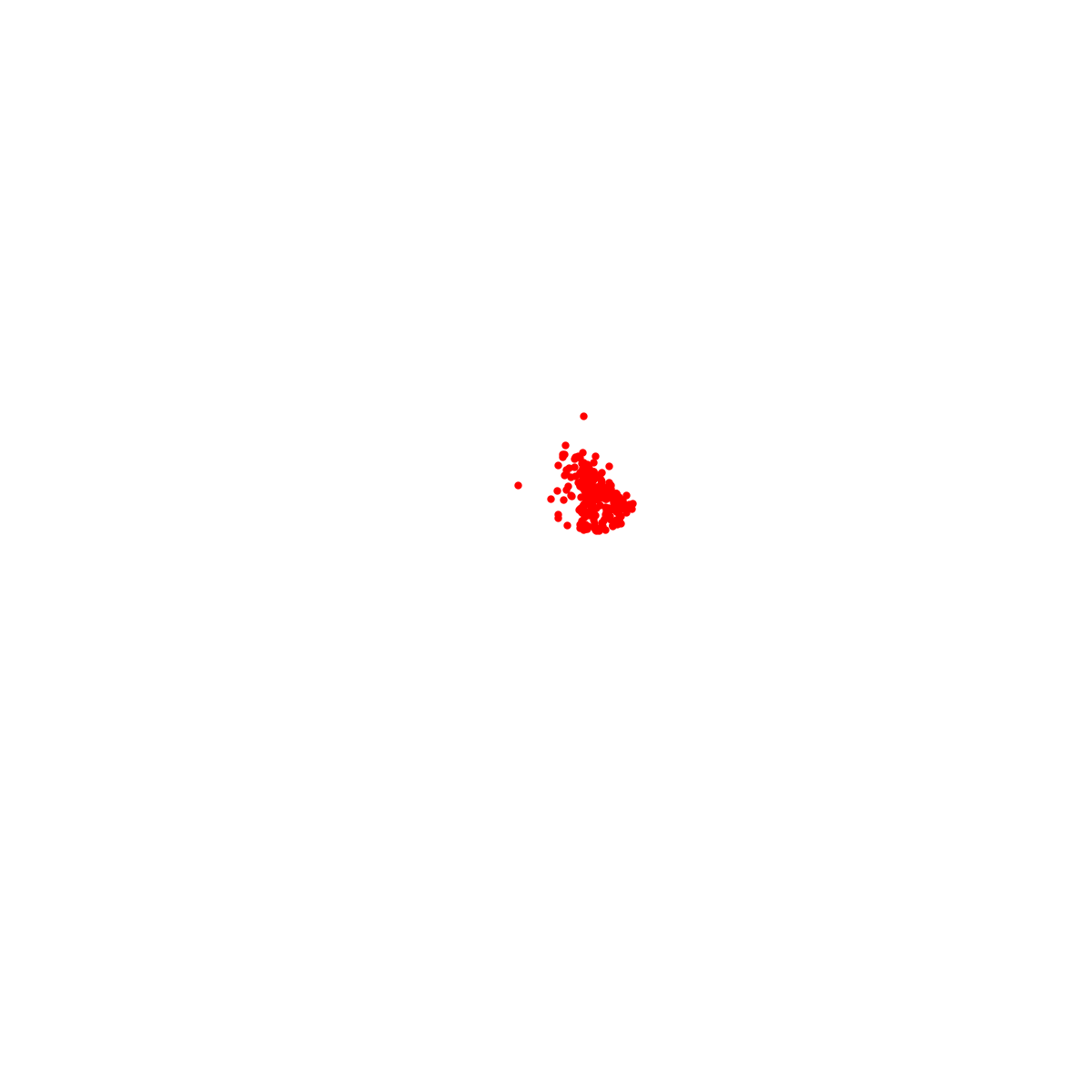}%
\end{subfigure}
\begin{subfigure}{0.19\textwidth}
\centering
\caption*{Raleigh}
\pdtwosubfig{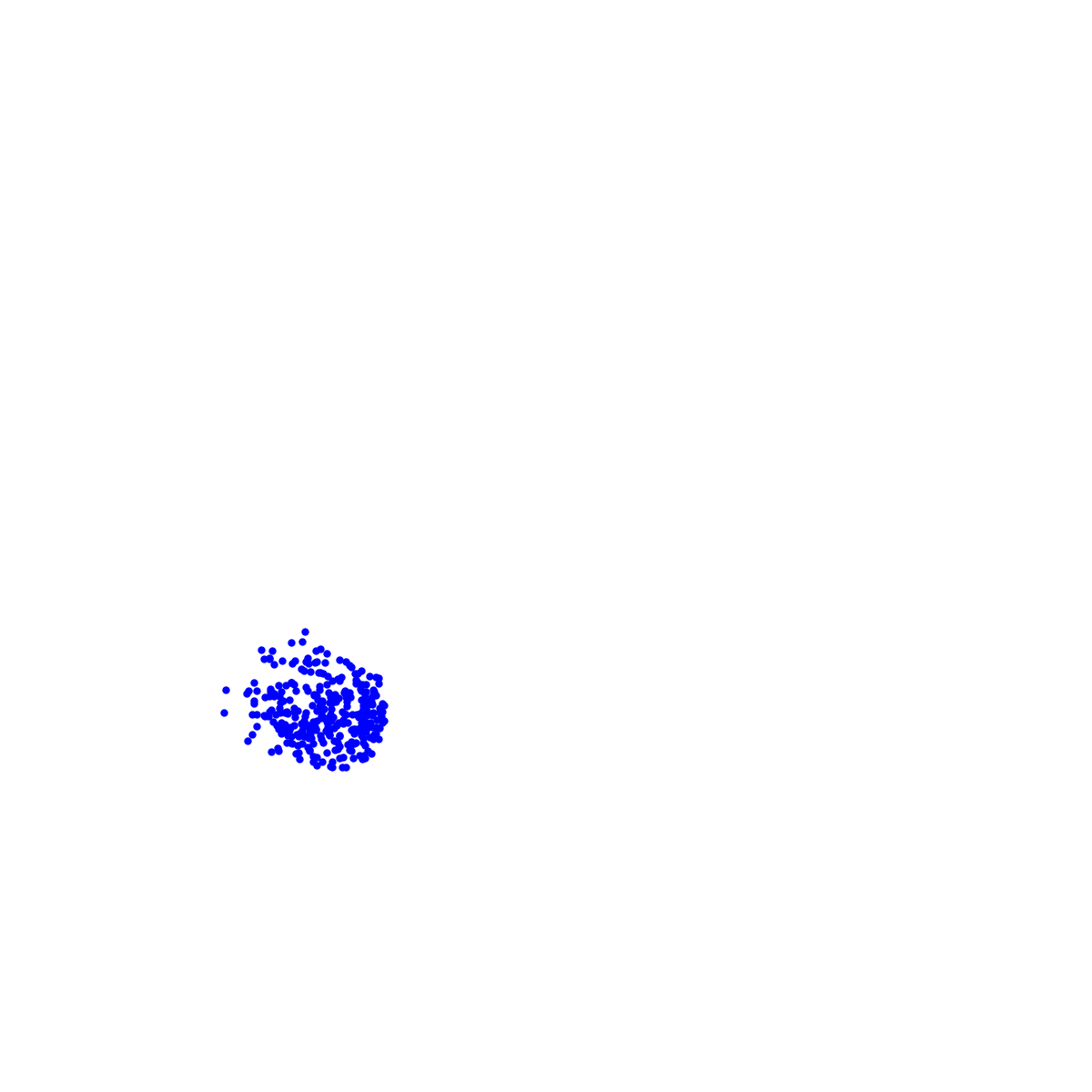}{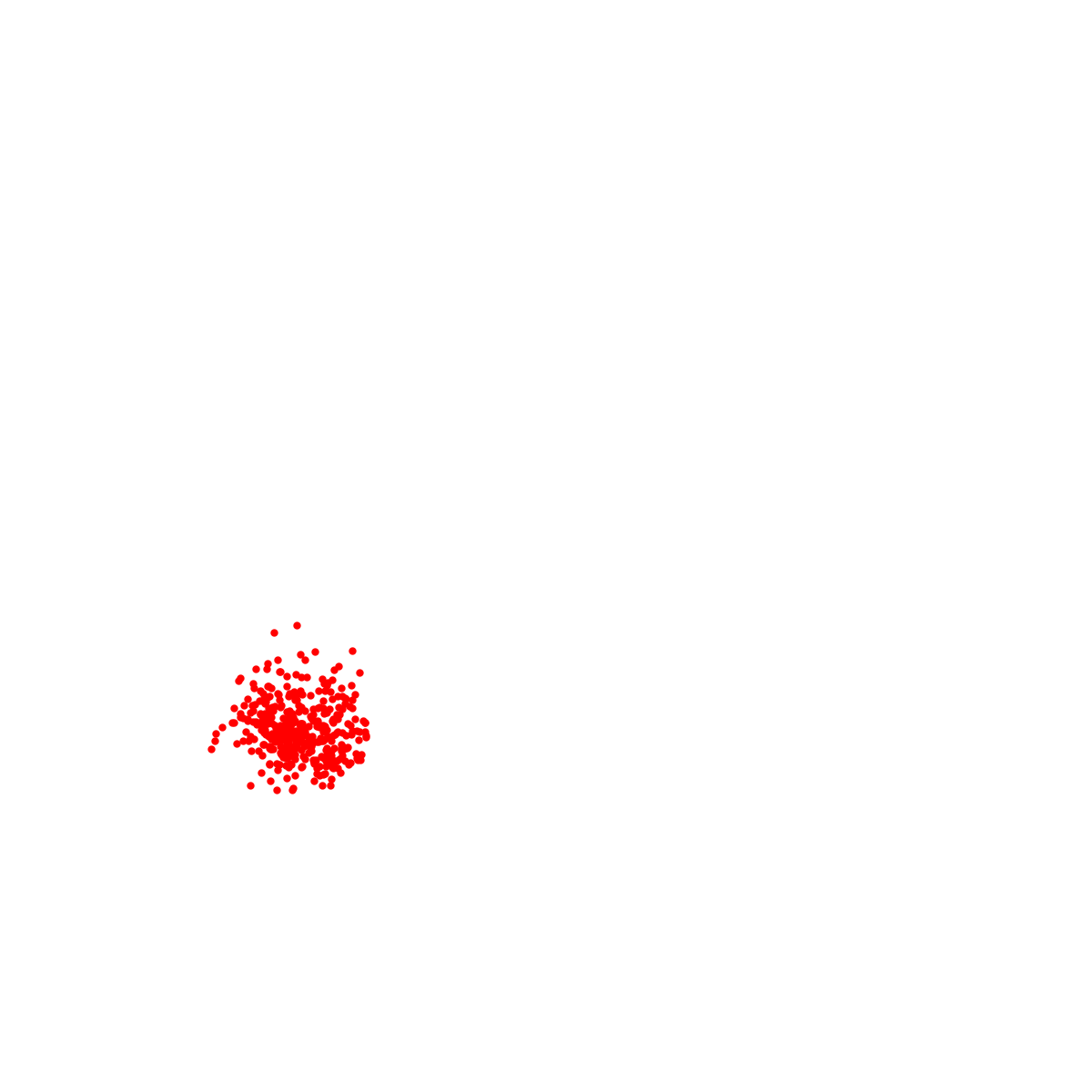}%
\end{subfigure}
\begin{subfigure}{0.19\textwidth}
\centering
\caption*{Northeast}
\pdtwosubfig{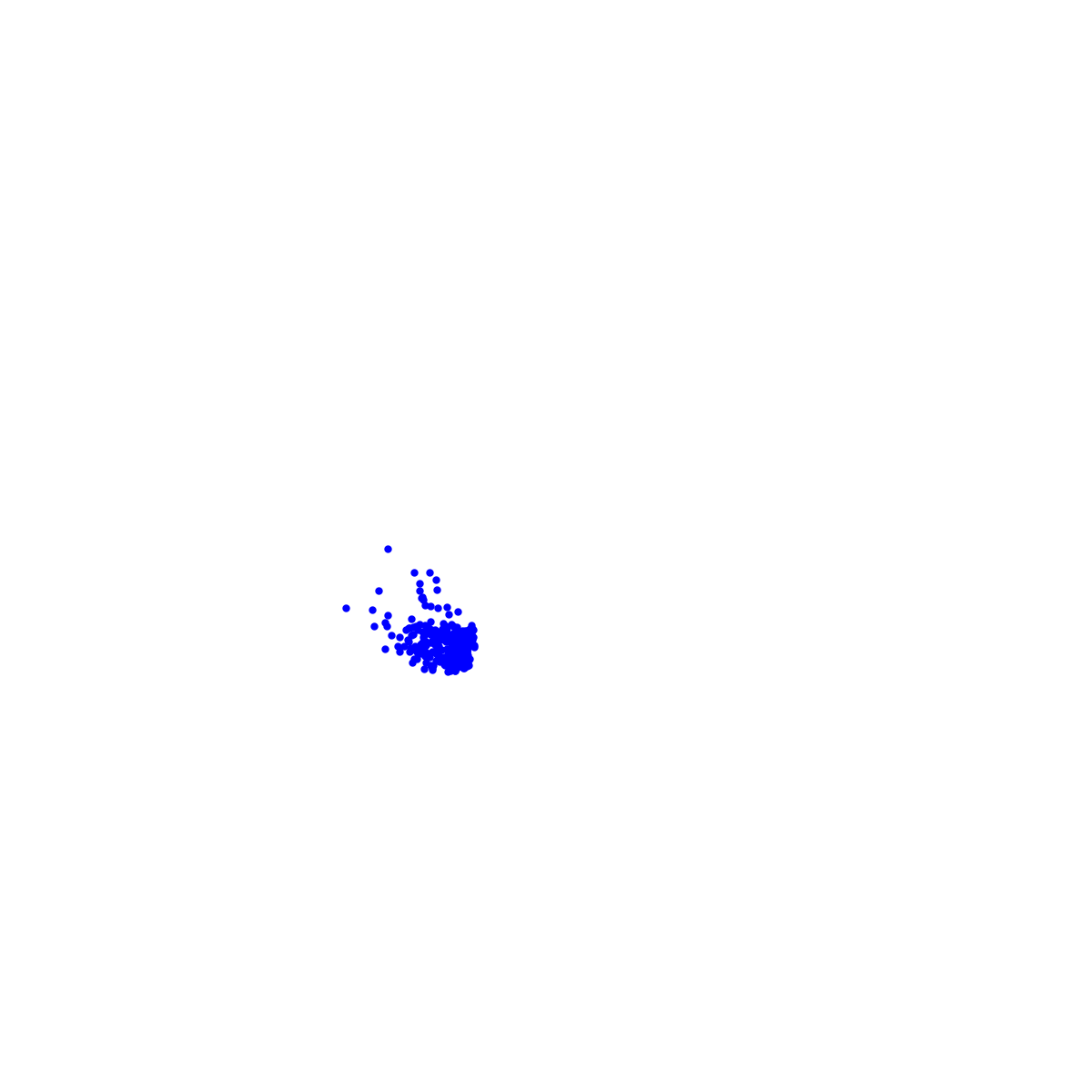}{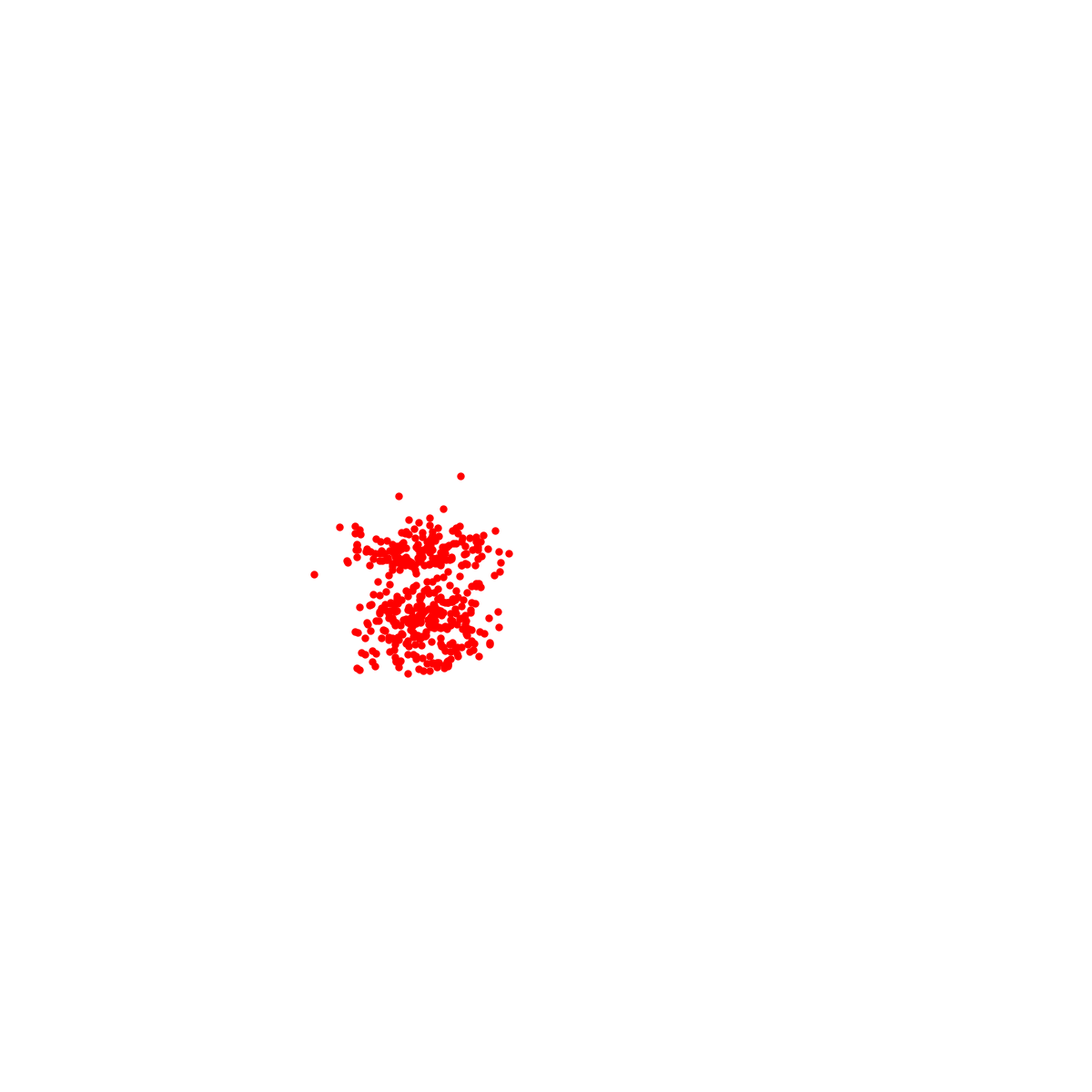}%
\end{subfigure}
\begin{subfigure}{0.38\textwidth}
\centering
\caption*{Winston-Salem (D)/Northeast (R)}
\resizebox{0.5\columnwidth}{!}{%
\begin{tikzpicture}%
\begin{axis}[
xmin=0, xmax=1, ymin=0, ymax=1.1, clip=false,
xtick={0.25,0.5,0.75}, ytick={0.25,0.5,0.75},
   enlargelimits=false, axis equal image, axis on top
]
\addplot[fill opacity=0.85, draw opacity=0.85] graphics [xmin=0, xmax=1,ymin=0, ymax=1.1] {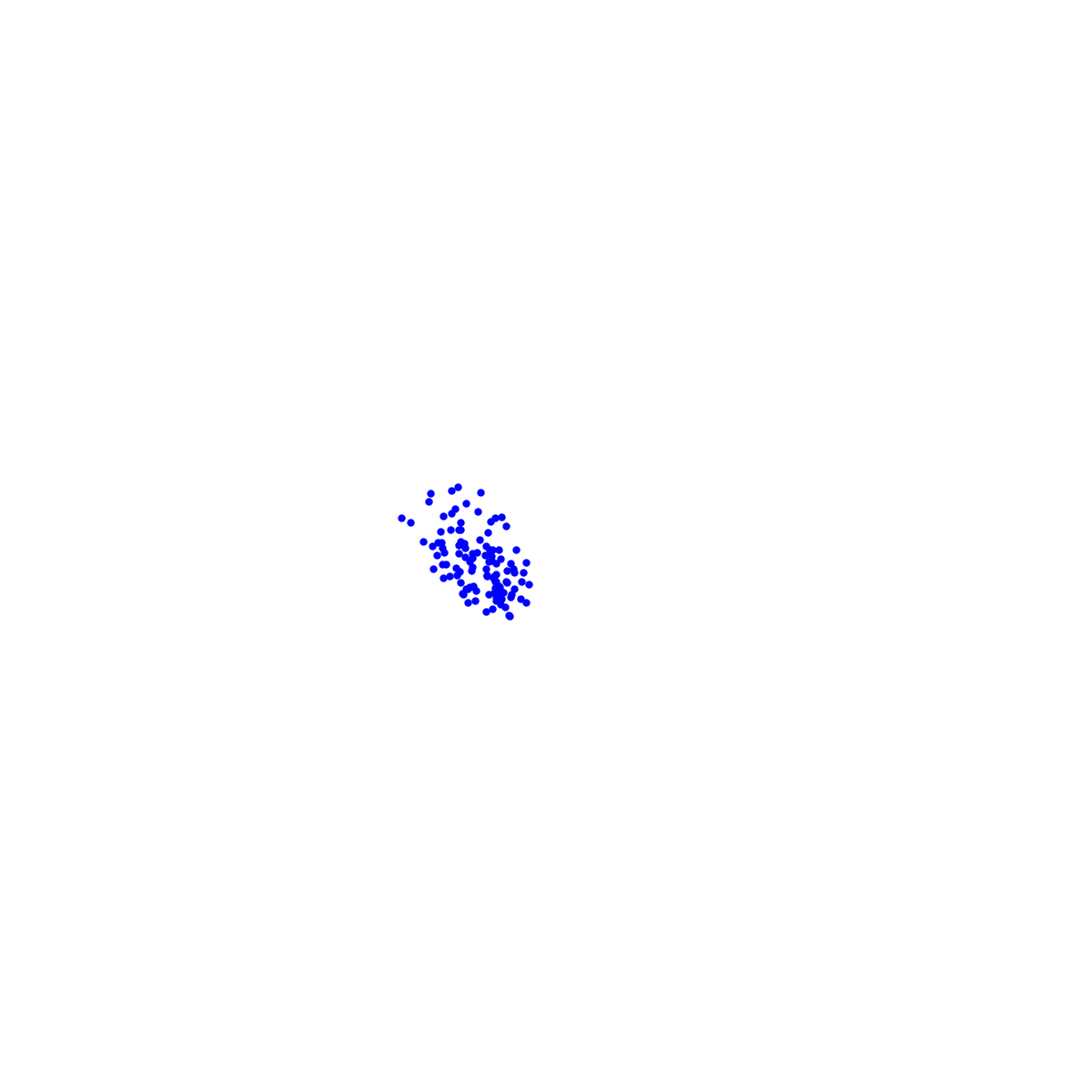}; 
\addplot[fill opacity=0.5, draw opacity=0.5] graphics [xmin=0, xmax=1,ymin=0, ymax=1.1] {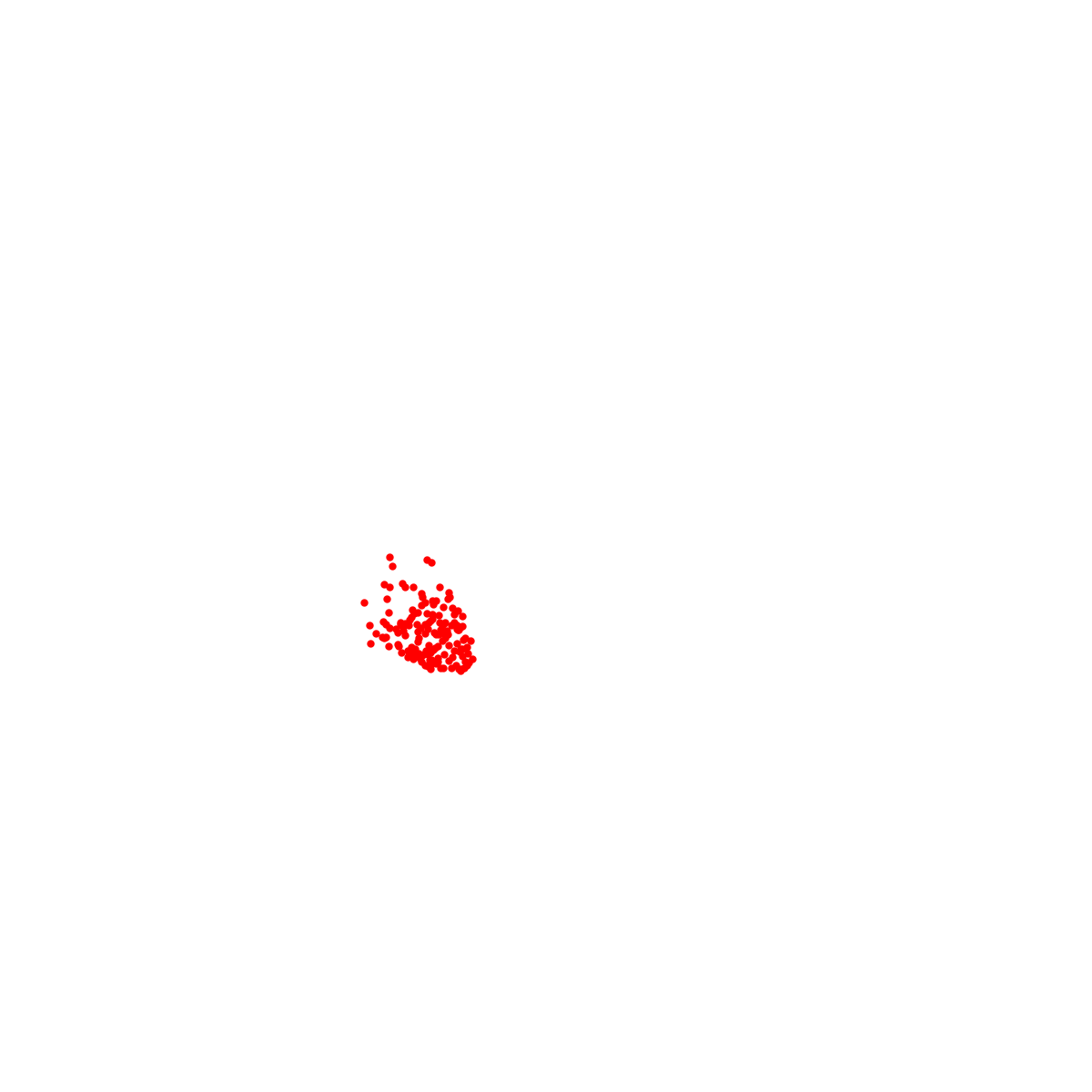}; 
\addplot[gray, dashed] coordinates{(0.5,0) (0.5,1)};
\addplot[gray, dashed] coordinates{(0,0) (1,1)};
\node at (axis cs: 0.03,1.02) {$\infty$};
\end{axis}
\end{tikzpicture}%
}%
\end{subfigure}
\caption{Comparison of point plots for the successive Fr\'echet features in the NC Senate ensembles that are biased for Democratic and Republican safe seats, respectively.
}
\label{biasedmeansNCshort}
\end{figure}
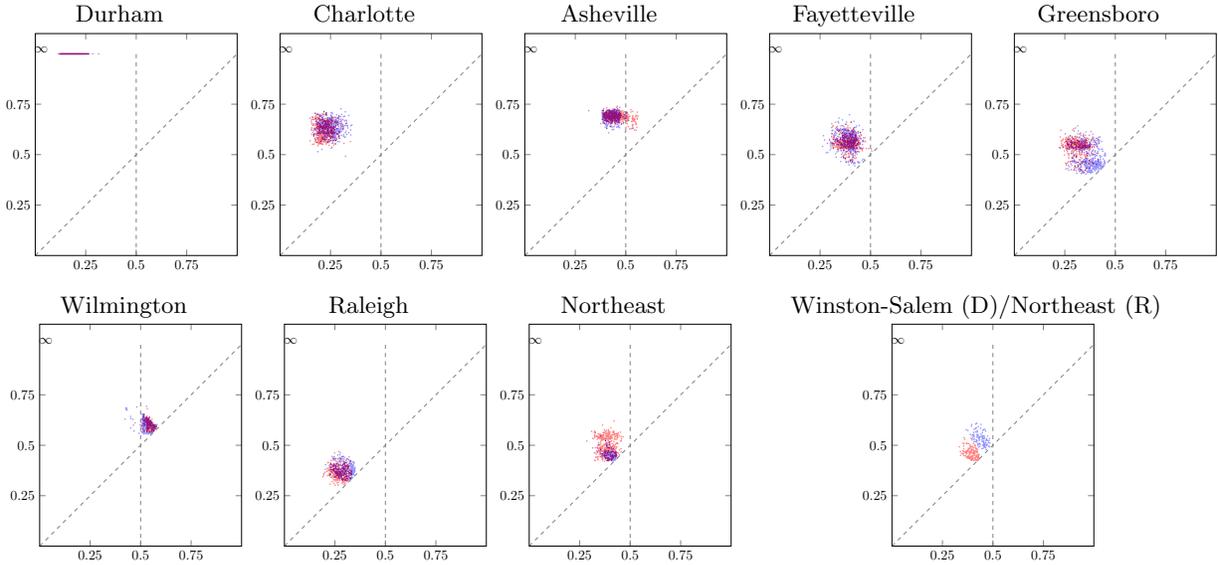

\section{Stability}\label{sec:stability}

We now address the question of stability of our persistent homology signatures of vote and districting data. In order to do so, we appeal to a stability result of \cite{cohen2007stability} which is well known in the field. In order to keep our exposition self-contained, we state the result in a special case which is relevant to our application.

\begin{thm}[\cite{cohen2007stability}]\label{thm:stability_theorem}
Let $G = (V,E)$ be a graph with $f$ and $g$ real-valued functions on $V$ and corresponding persistence diagrams $\D_f$ and $\D_g$. Then 
$$
d_\infty(\D_f,\D_g) \leq \left\|f - g\right\|_\infty,
$$
where $\|\cdot\|_\infty$ is the $\ell_\infty$ norm on functions $V \rightarrow \R$. 
\end{thm}

We then immediately obtain a stability result in our application, for the district dual graphs and Republican share filtration used above.

\begin{cor}[Stability under perturbations in vote data]\label{cor:stability_vote_data}
The bottleneck distance between persistence diagrams associated to two elections and the same districting plan is bounded above by the largest change in Republican vote share an any district.
\end{cor}

\subsection{Theoretical stability under geographic perturbations}

\begin{figure}[ht]
\centering
\begin{subfigure}{0.45\textwidth}
\includegraphics[width=\textwidth]{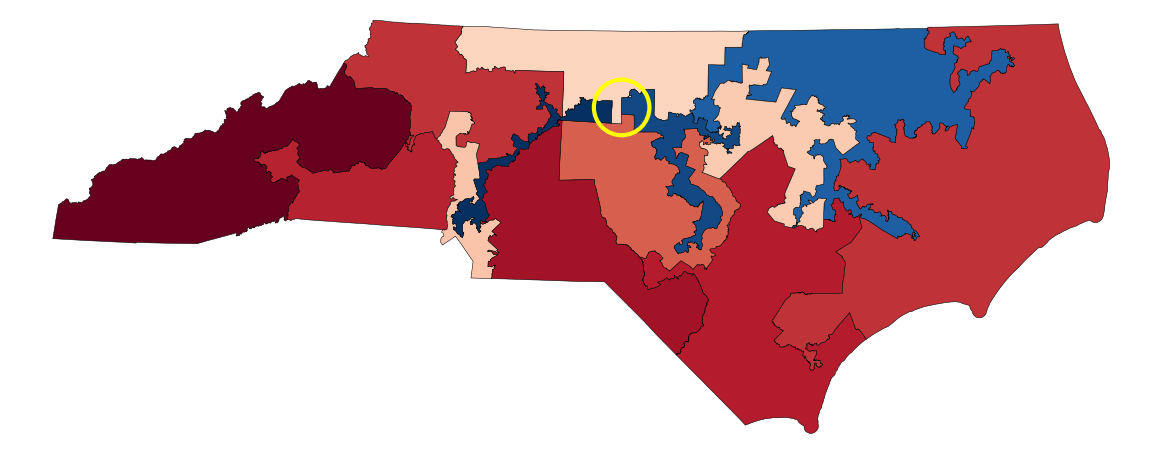}
\end{subfigure}
\begin{subfigure}{0.45\textwidth}
\includegraphics[width=\textwidth]{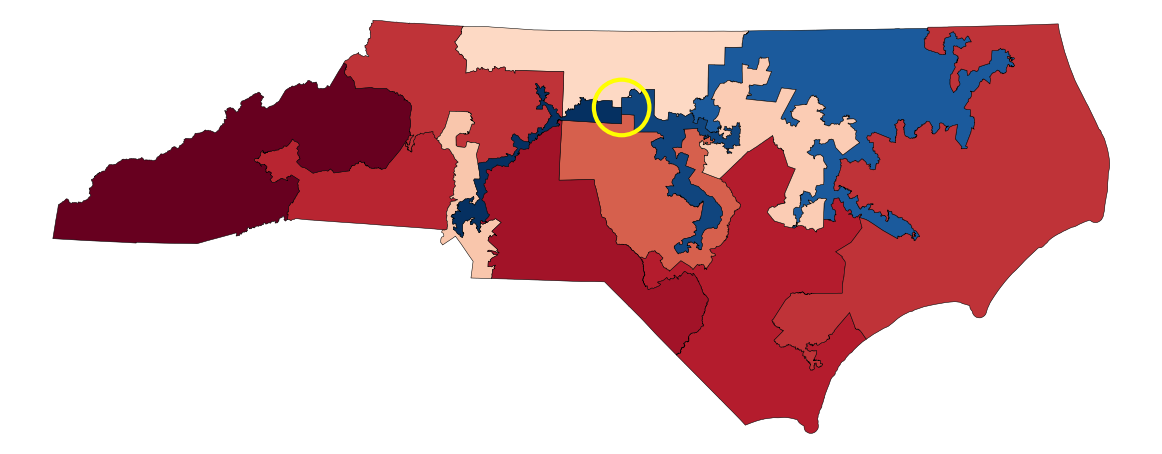}
\end{subfigure}

\begin{subfigure}{0.24\textwidth}
\resizebox{\columnwidth}{!}{%
\begin{tikzpicture}
\begin{axis}[xmin=0, xmax=1, ymin=0, ymax=1.1, clip=false,
xtick={0,0.25,0.5,0.75,1}, ytick={0,0.25,0.5,0.75}
]
\addplot[color=black, only marks, mark=*, mark size=2.5]
 coordinates{(0.28262848025887055, 0.50916421609425)  (0.2634125209983995,1)};
\addplot[gray, dashed] coordinates{(0,1) (1,1)};
\addplot[gray, dashed] coordinates{(0,0) (1,1)};
\node at (axis cs: 0.03,1.02) {$\infty$};
\end{axis}
\end{tikzpicture}
}
\end{subfigure}
\begin{subfigure}{0.24\textwidth}
\centering
\includegraphics[width=\textwidth,height=.7\textwidth]{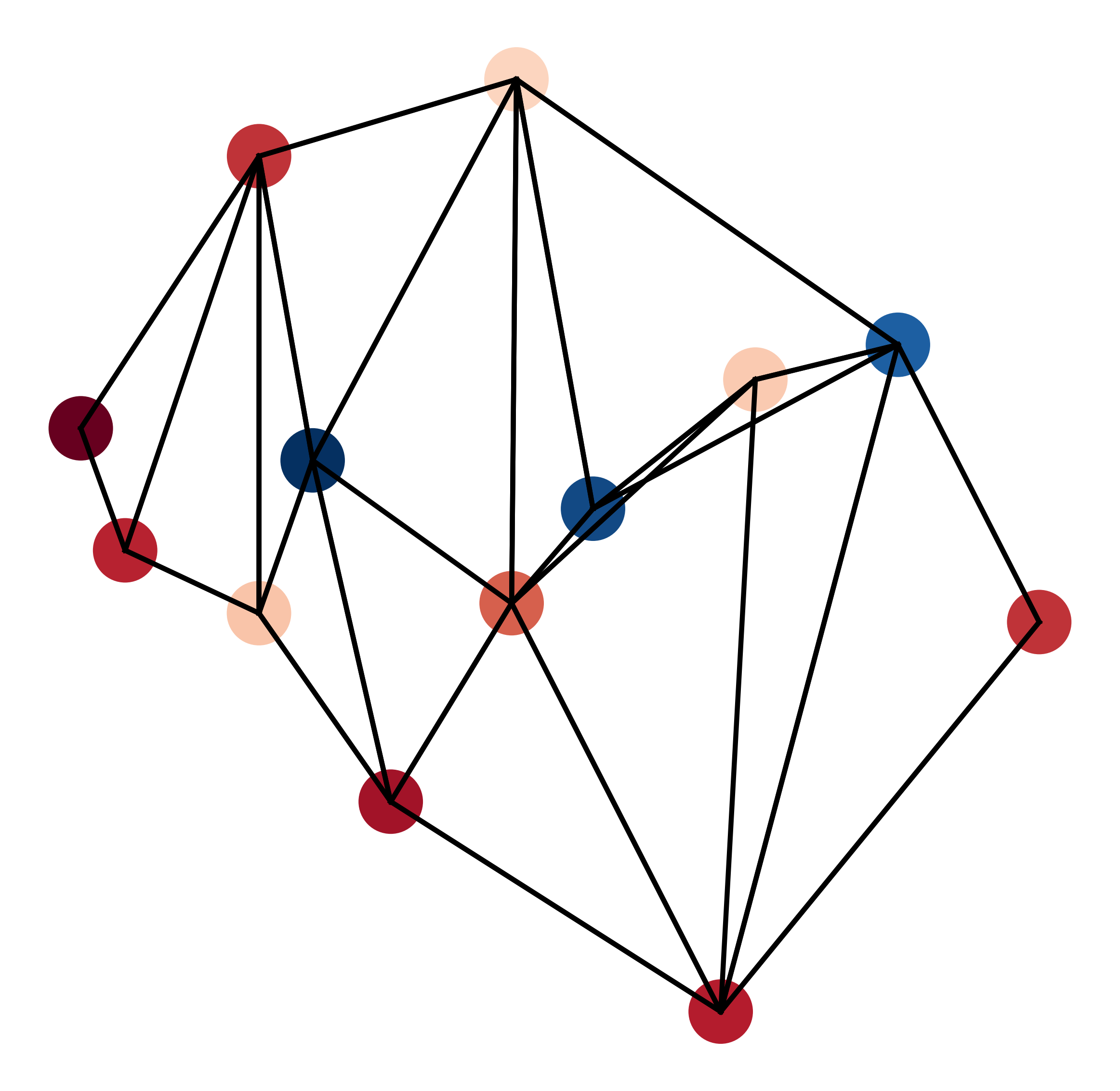}
\end{subfigure}
\begin{subfigure}{0.24\textwidth}
\resizebox{\columnwidth}{!}{%
\begin{tikzpicture}
\begin{axis}[xmin=0, xmax=1, ymin=0, ymax=1.1, clip=false,
xtick={0,0.25,0.5,0.75,1}, ytick={0,0.25,0.5,0.75}
]
\addplot[color=black, only marks, mark=*, mark size=2.5]
 coordinates{(0.2666046383177386, 1)};
\addplot[gray, dashed] coordinates{(0,1) (1,1)};
\addplot[gray, dashed] coordinates{(0,0) (1,1)};
\node at (axis cs: 0.03,1.02) {$\infty$};
\end{axis}
\end{tikzpicture}
}
\end{subfigure}
\begin{subfigure}{0.24\textwidth}
\centering
\includegraphics[width=\textwidth,height=.7\textwidth]{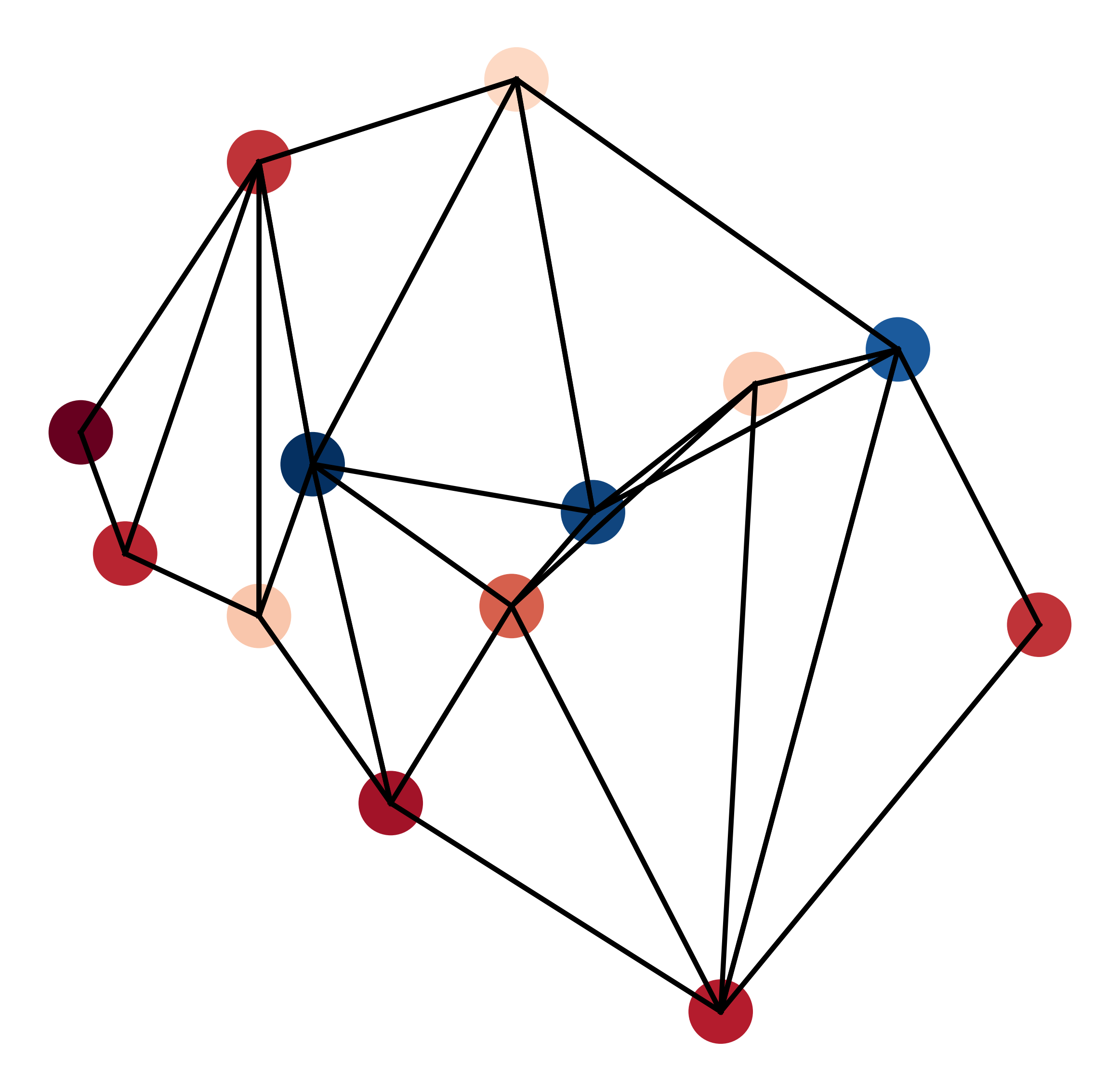}
\end{subfigure}

\caption{A small change to districts can cause a large change in persistence. In this example (closely modeled on the enacted NC Congressional plan from 2012) we have changed the assignment of only one geographic unit (a precinct of 3,300 people), and it removes a highly persistent off-diagonal feature.}\label{fig:NCperturb}
\end{figure}
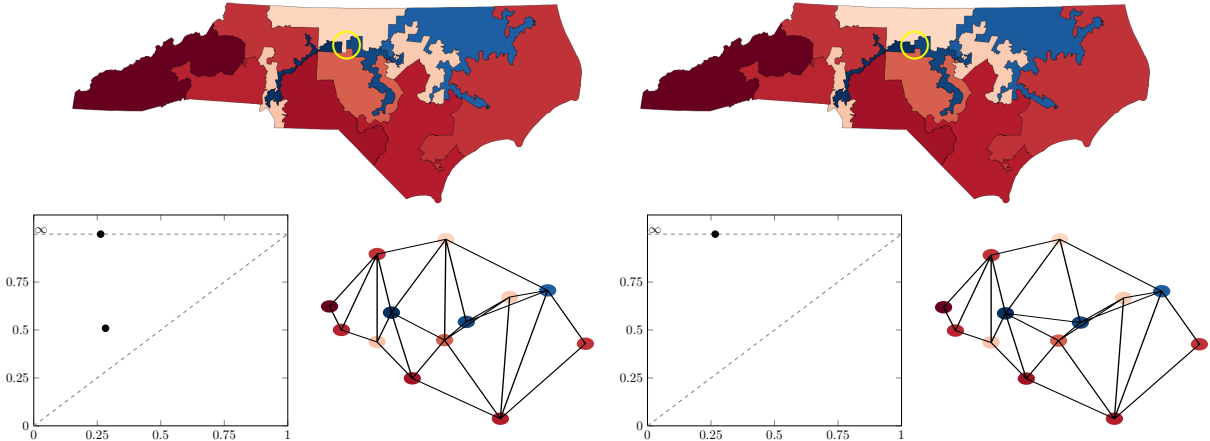

For a fixed plan, we have a strong and easily interpretable stability result for bottleneck distance as elections shift.  One could similarly ask if the persistence diagram summaries are robust under small changes in the plan itself. A truly general result in this direction is not possible, as is illustrated by the  example in Figure~\ref{fig:NCperturb}. 

In light of this fragility, some care must be taken to formulate a correct result on stability of bottleneck distance under changes in districting plan. We formulate such a result in this subsection, and relate it to experimental analysis in the following section using random perturbations. 

Fix a state's unit dual graph $G = (V,E)$ and a districting plan $P = \{P_1,\ldots,P_k\} \in \mathcal{P}_{k, \epsilon}(G)$. 
Let $p$ be the population function and $T=p(V)/k$ be the ideal district size.

\begin{definition} For districting plans $P,P'\in \PkG$, 
$P'$ is called a \emph{(one-way) perturbation} of $P$ if, after reindexing the districts as necessary, the plans agree in all but 
two districts $P_i \cup P_j = P_i' \cup P_j'$, and they differ only by the move of vertices $V_{ij}=\{v_1,\dots,v_\ell\}$, so that 
$V_j'=V_j\cup V_{ij}$ and $V_i'=V_i\setminus V_{ij}$.  
A {\em general perturbation} has vertices $V_{ij}$ moved in one direction and $V_{ji}$ in the other direction between the two altered districts. 

A perturbation $P,P'$ is called {\em graph preserving} 
if the labeled correspondence of districts induces an isomorphism of the district dual graphs.
\end{definition}

For instance, the  two North Carolina plans in  Figure \ref{fig:NCperturb} are perturbations of each other by a single unit, which means that $|V_{ij}|=1$. However, it is not a graph-preserving perturbation; the dual graphs $G_P$ and $G_{P'}$ are not even abstractly isomorphic.

A ReCom step (as is used to generate our districting plan ensembles) can potentially produce a one-way perturbation of a plan, but this is extremely unlikely because  population is swapped in both directions with very high probability. 
Moreover, ReCom steps are not always be graph-preserving because in re-splitting the combined population of two districts in an arbitrary connected fashion, they can alter district adjacency with the surrounding nodes.\footnote{Several trials show that a ReCom step preserves the district dual graph with probability roughly $.65$.  See Supplement for more discussion.}
A simpler Markov chain method for generating ensembles of plans by single-unit flips is studied here.

A \emph{flip move} on a plan $P$ is performed by reassigning a single vertex from one district to another in a manner that preserves connectedness for all districts.  It therefore necessarily induces a one-way perturbation, and frequently a graph-preserving perturbation.  A Markov chain can be defined by randomly selecting a vertex to reassign at each step.  
 Taking sequences of flip steps (obtained by rejection sampling flip step proposals) defines a simple  and commonly used Markov chain on the space of districting plans \cite{mattingly1, chikina2017assessing}. 

Recall that for a function $f:V\to \R$ and a subset $W\subseteq V$ we write $f(W)$ for the sum $\sum_{v\in W} f(v)$,
and that $p$ is a weight function thought of as population on the nodes.
We use the following notation below.  
Let $a,b:V \rightarrow \R$ be functions satisfying $0\leq a(v) \leq b(v) \leq p(v)$ for all $v \in V$, and for simplicity 
assume $b>0$.  Let $f$ be given by $f=a/b$, and we will abuse notation by using the same notation for a function on $V$ or on the district nodes $V_P$.  This induces district-level persistence diagrams $\D_P$ and $\D_{P'}$.
Fix $\alpha  = \alpha_{P,P'}>0$ to be a lower bound on $b(W)/p(W)$ for  $W = V_i,V_i',V_j,V_j'$.

We will apply this to $a=r$, the Republican vote totals per node, letting $b=r+d$ be the total major-party votes per node, so that 
 $f=R$, the Republican voting share, and $\alpha$ is a lower bound on the voting turnout in districts affected by the perturbation.
The more general notation is introduced here because there are several other functions $a,b$ that could induce interesting filtration functions on the dual graphs, for instance to take third-party voting or other demographics into account.  

We now state the stability result. Its proof is a straightforward application of Theorem \ref{thm:stability_theorem}, but we are able to formulate the estimate in terms of quantities which are directly interpretable from a districting perspective. 

\begin{thm}[Stability under geographic perturbations]\label{thm:stability_perturbations}
For a graph preserving one-way perturbation $P,P'$ where vertices $V_{ij}$ are flipped from $P_i$ to $P_j$, 
$$
d_\infty\left(\D_{P},\D_{P'}\right) \leq \frac{2\epsilon}{\alpha(1-\epsilon)} 
\max \left( \left|   f(V_{ij})- f(P_i)\right| , \left|f(V_{ij}) - f(P_j)\right|  \right).
$$
\end{thm}

This gives us a bound on bottleneck distance for a perturbation in terms of how much the party vote share in the flipped units compares to the districts that they flip between.  The proportionality constant depends on the population tolerance $\epsilon$ and the voter turnout constant $\alpha$, and is small under the realistic assumption that $\epsilon \ll \alpha$.  (For instance,  with population deviation tolerance of 2\%  and an election with district-level turnouts of at least 25\%, we obtain a coefficient of roughly  $.163$.)%

\begin{proof}
Relabel if necessary so that the affected districts of $P$ are $P_1,P_2$, and $P'$ replaces those by $P_1',P_2'$, with 
districts $3,\dots, k$ unchanged.
Since the district dual graphs agree by assumption, the filtration function by $f$ is defined on that graph. 
Theorem \ref{thm:stability_theorem} says that 
$$
d_\infty\left(\D_P,\D_{P'}\right) \leq \left\|f_{P} - f_{{P}'}\right\|_\infty = \max\left( |f(P_1)-f(P_1')|,\  |f(P_2)-f(P_2')| \right).
$$

Recall that $f(P_1)=a(P_1)/b(P_1)$, which in our case is Republican vote share, and likewise for the other districts.
To simplify the exposition, we write $a_i=a(P_i)$ for $i=1,2$ and similarly define $a_i',b_i,b_i'$; we also write $\abar=a(V_{12})$ and similarly for $\bbar$.    Note that $a_1'=a_1-\abar$ and $b_1=b'-\bbar$.
It follows that
$$f(P_1)-f(P_1') = \frac{a_1}{b_1}- \frac{a_1'}{b_1'} =  \frac{a_1}{b_1}- \frac{a_1-\abar}{b_1-\bbar}=
 \frac{b_1\sdot \abar-a_1\sdot \bbar}{b_1(b_1-\bbar)} = \frac{\bbar}{b_1'} \sdot \left(f(V_{12})-f(P_1) \right).
$$
This works exactly similarly for $P_2$.

Since ${P}$ and ${P}'$ are districting plans, the population balance condition  implies that
$p(P_1)$ and $p(P_1')$ 
 both lie in the interval $[(1-\epsilon)T,(1+\epsilon)T]$. Since the vertices of $V_{12}$ were moved without violating the balance condition, we have 
$p(V_{12})\leq 2\epsilon T$.

Since $b(W)\le p(W)$ for any vertex set $W$, we have
$$
\alpha(1-\epsilon) \bbar \leq  2\epsilon \alpha (1-\epsilon) T \leq 2 \epsilon \alpha \sdot p(V_1) 
\leq 2 \epsilon \sdot b(V_1') =2\epsilon b_1'.
$$
To complete the proof, we rearrange to obtain
$$\frac{\bbar}{b_1'} \leq \frac{2 \epsilon}{\alpha(1-\epsilon)}.\qedhere$$
\end{proof}

\subsection{Experimental stability under geographic perturbations}\label{subsec:experimental_stability}
While Theorem \ref{thm:stability_perturbations} gives a theoretical stability guarantee for geographic perturbations, its applicability is somewhat limited in scope, since the graph-preserving property is not automatic.
In this section we present experimental evidence of stability.

\begin{figure}[ht]
\centering

\begin{subfigure}{0.38\textwidth}
\centering
\caption*{Enacted/perturbed Congressional plan}
\includegraphics[width=\textwidth]{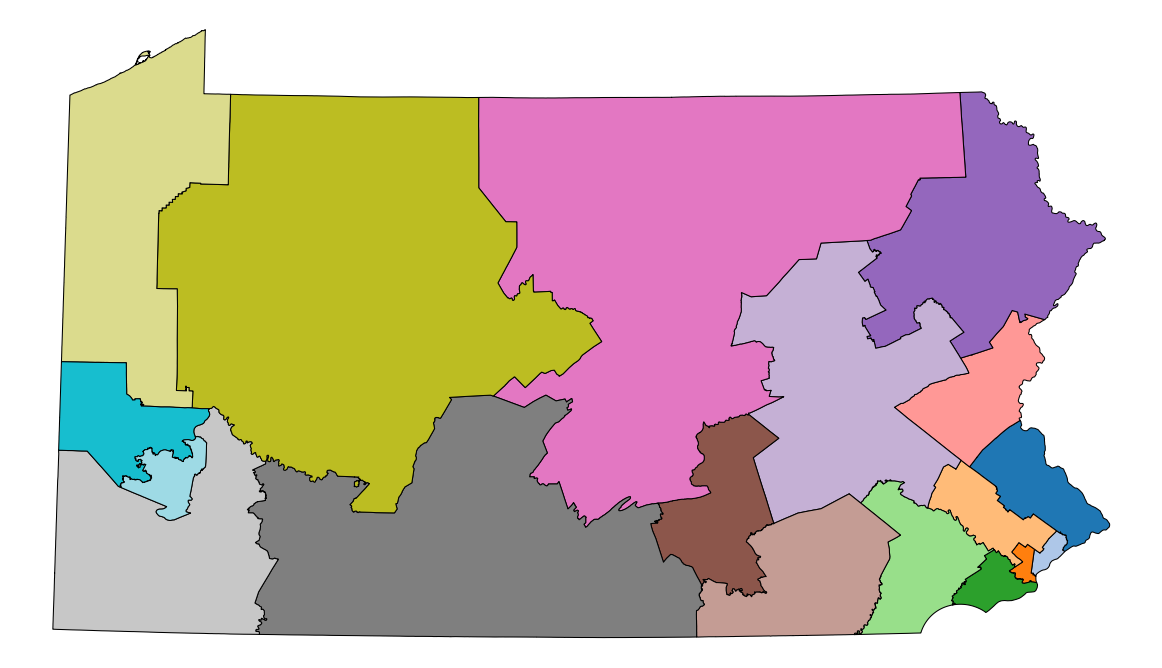}
\end{subfigure}
\begin{subfigure}{0.38\textwidth}
\centering
\caption*{Enacted/perturbed state Senate plan}
\includegraphics[width=\textwidth]{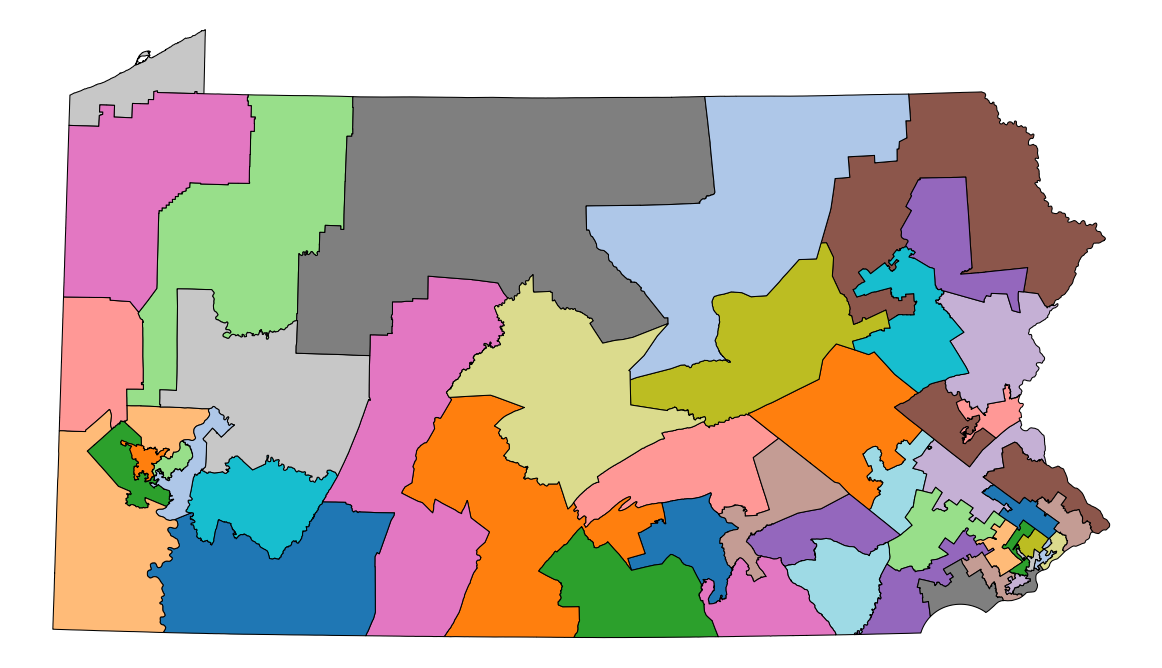}
\end{subfigure}

\begin{subfigure}{0.38\textwidth}
\centering

\includegraphics[width=\textwidth]{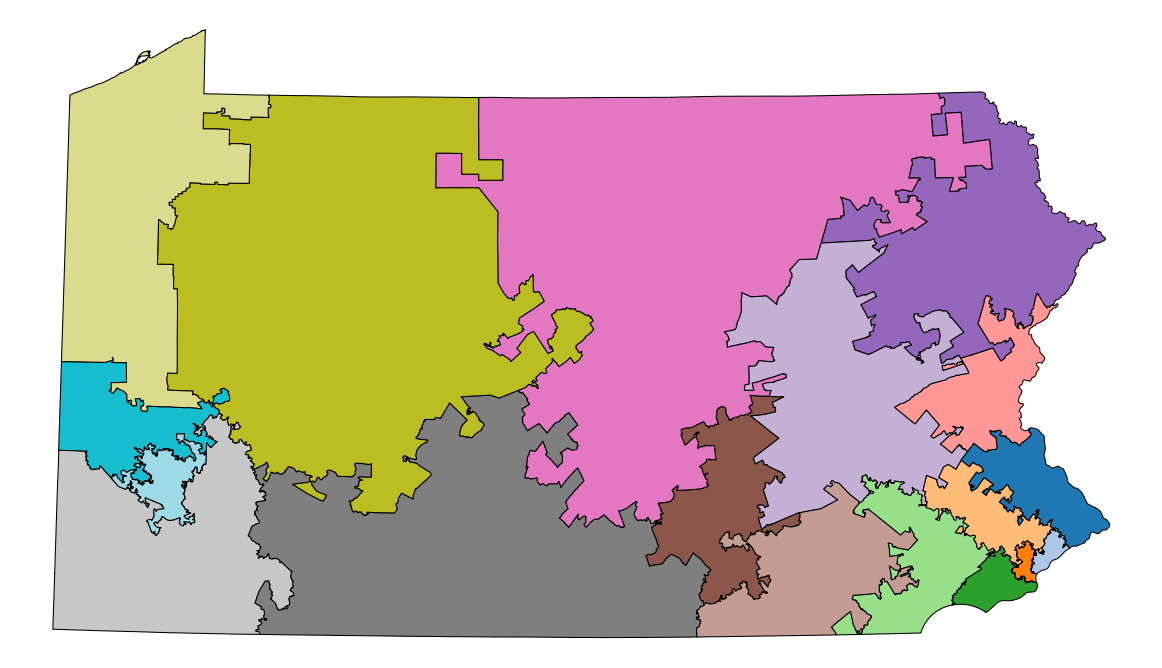}

\end{subfigure}
\begin{subfigure}{0.38\textwidth}
\centering

\includegraphics[width=\textwidth]{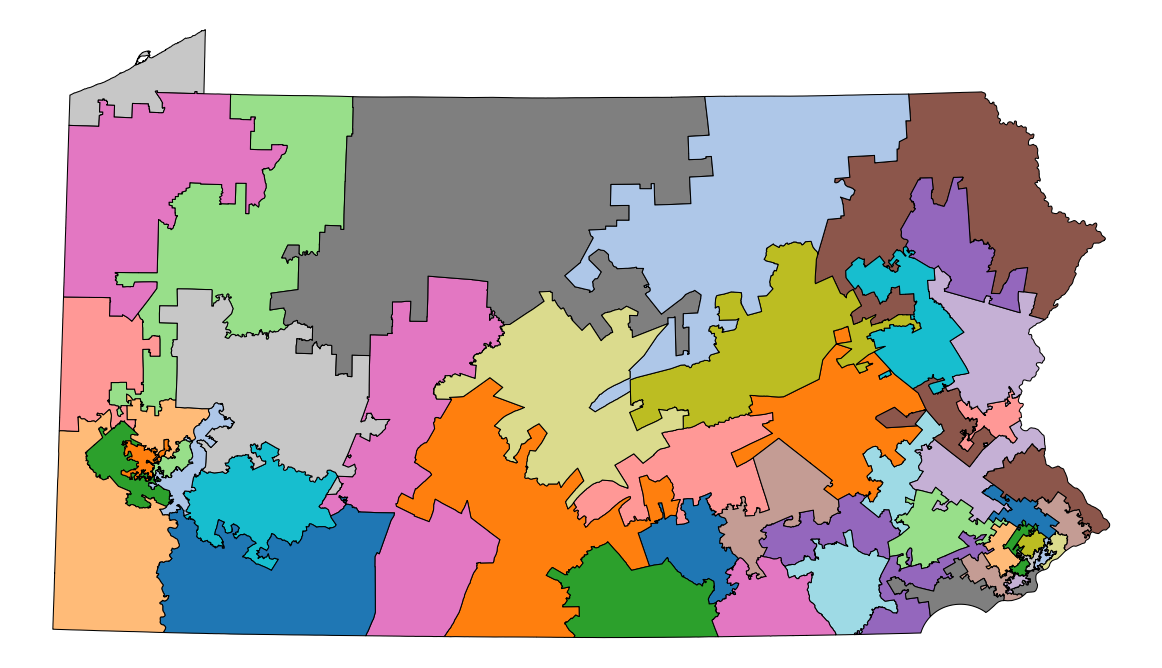}

\end{subfigure}

\foreach \k in {18, 50}
{
\begin{subfigure}{0.38\textwidth}
\centering

\resizebox{.8\columnwidth}{!}{%
\begin{tikzpicture}
\begin{axis}[xmin=0, xmax=1, ymin=0, ymax=1.1, clip=false,
xtick={0,0.25,0.5,0.75,1}, ytick={0,0.25,0.5,0.75}
]
\addplot[scatter, 
scatter/classes={
enacted={green!80!black, opacity=0.5},
perturbed={red, opacity=0.5}
},
only marks, mark=*, mark size=2, scatter src=explicit symbolic]
 table[x=X, y=Y, col sep=comma, meta=label]{PA_csv/\k flipped_pers.csv};
\addplot[gray, dashed] coordinates{(0.5,0) (0.5,1)};
\addplot[gray, dashed] coordinates{(0,0) (1,1)};
\node at (axis cs: 0.03,1.02) {$\infty$};
\legend {enacted, perturbed};
\end{axis}
\end{tikzpicture}
}
\end{subfigure}
}
\caption{Maps of enacted plans before and after perturbation by 1000 flip steps;  district adjacencies changed in both cases (e.g., blue/lavender for Congress, and purple/purple for state Senate districts in the Southeast). Nonetheless we see only a  small change in persistence, which is typical for repetitions of this experiment.}
\label{flipped_maps}
\end{figure}
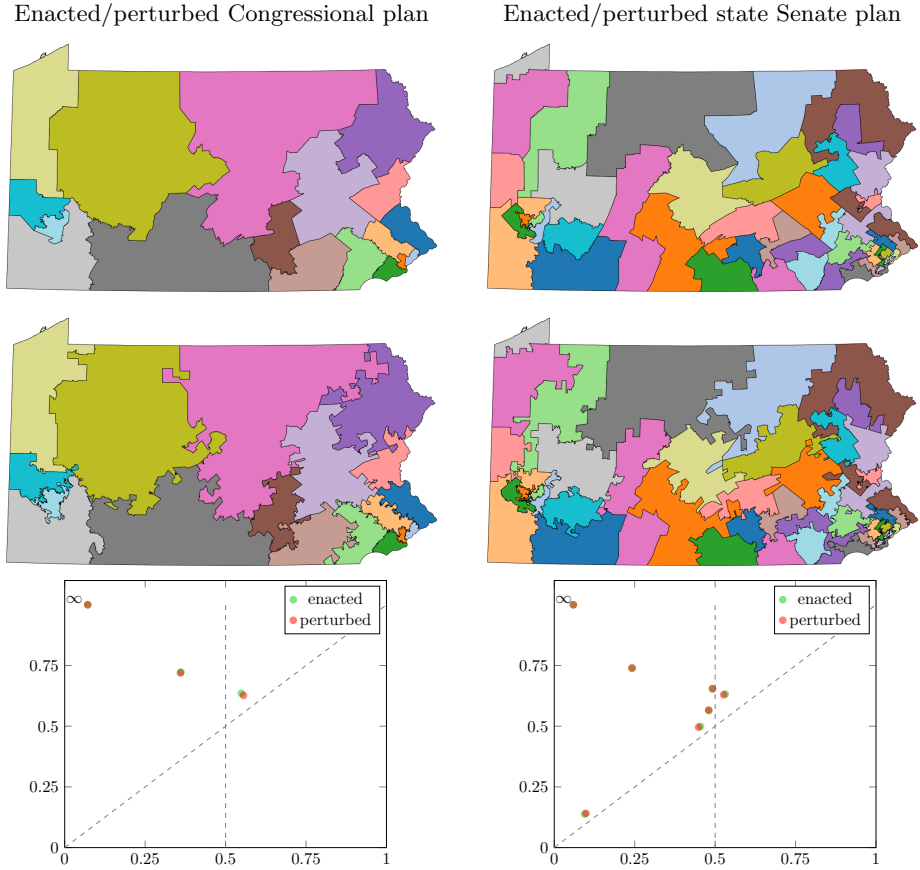

\begin{figure}[ht]
\centering
\begin{subfigure}{0.25\textwidth}
\centering
\caption*{18 districts}
\resizebox{\columnwidth}{!}{%
\begin{tikzpicture}
\begin{axis}[clip=false,
xtick={0,0.05,0.1,0.15,0.2},
   enlargelimits=false, axis on top,
        x tick label style={
        /pgf/number format/fixed,
				/pgf/number format/precision=1
                },
         scaled x ticks = false
]
\addplot graphics [xmin=0, xmax=0.2,ymin=0, ymax=300000] {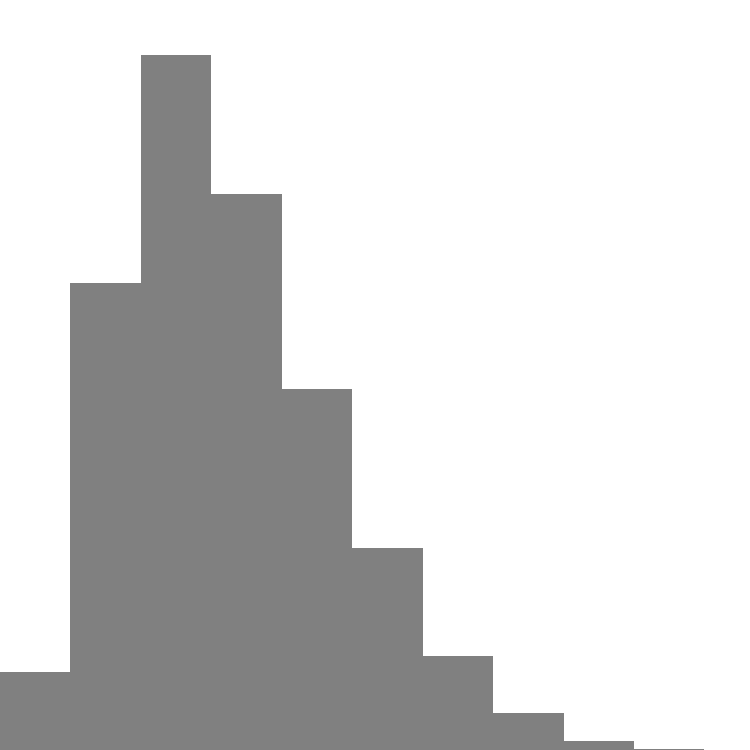}; 
\end{axis}
\end{tikzpicture}%
}
\end{subfigure}
\begin{subfigure}{0.25\textwidth}
\centering
\caption*{50 districts}
\resizebox{\columnwidth}{!}{%
\begin{tikzpicture}
\begin{axis}[clip=false,
xtick={0,0.05,0.1,0.15,0.2},
   enlargelimits=false, axis on top,
        x tick label style={
        /pgf/number format/fixed,
				/pgf/number format/precision=1
                },
         scaled x ticks = false
]
\addplot graphics [xmin=0, xmax=0.2,ymin=0, ymax=500000] {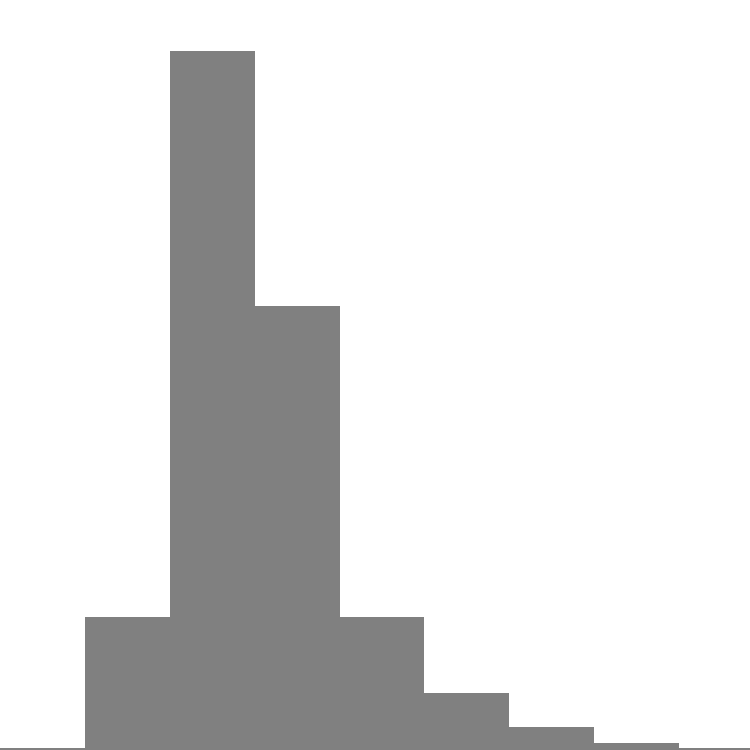}; 
\end{axis}
\end{tikzpicture}%
}
\end{subfigure}
\caption{Distribution of bottleneck distances between pairs of plans in the Pennsylvania ReCom ensembles. These distances can help us contextualize the size of the displacement by Flip steps.}
\label{fig:bottleneckDistances}
\end{figure}

We begin with Congressional and Senate ensembles for Pennsylvania.  ReCom steps make large changes to a map, and the ensembles here were sub-sampled so that each successive plan collected in a ReCom ensemble is essentially independent of the last.  So to illustrate stability we first turn to the much smaller perturbations made by flipping the assignment of a single unit at a time.
Flip steps make small changes to the appearance of a plan---though most every pairwise district boundary has been changed after 1000 flip steps, the configuration of districts is clearly recognizable.
In the examples shown in Figure~\ref{flipped_maps}, 1000 flip steps moves the persistence diagram by $d_\infty=.009$ and $d_\infty=.006$ for Congress and Senate, respectively.  
To contextualize distance magnitudes, note that pairs of maps in our ReCom ensembles (approximately independent maps) often have pairwise bottleneck distances between $.05$ and $.1$ (see Figure \ref{fig:bottleneckDistances}).    

\begin{figure}[ht]
\centering
\begin{subfigure}{\textwidth}
\centering
\caption*{Perturbing PA Congressional districts ($k=18$)}
\resizebox{\columnwidth}{!}{%
\begin{tikzpicture}
\begin{axis}[
y tick label style={
        /pgf/number format/.cd,
            fixed,
            fixed zerofill,
            precision=2,
        /tikz/.cd
    },
scaled ticks=false,
width=\textwidth, height=0.5\textwidth,
xlabel={flip steps},
ylabel={bottleneck distance},
]
\addplot graphics [xmin=0, xmax=1000,ymin=0, ymax=0.02] {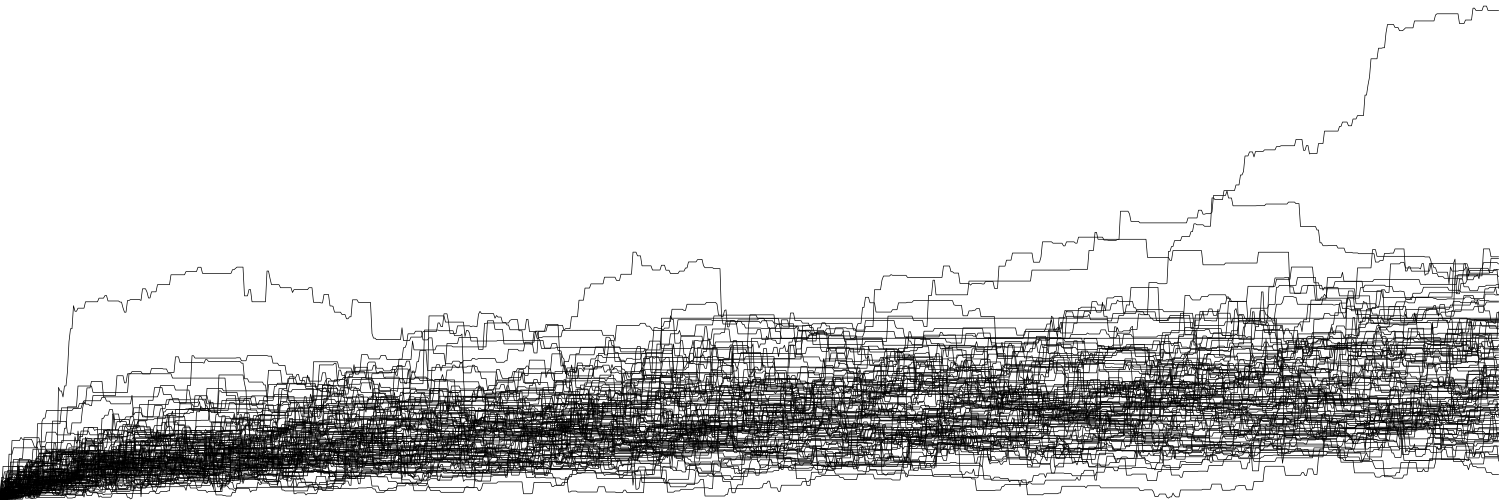}; 
\addplot[red, thick, line legend] coordinates{(0,0.0605848) (1000,0.0605848)};
\end{axis}
\end{tikzpicture}
}
\end{subfigure}

\begin{subfigure}{\textwidth}
\centering
\caption*{Perturbing PA Senate districts ($k=50$)}
\resizebox{\columnwidth}{!}{%
\begin{tikzpicture}
\begin{axis}[
y tick label style={
        /pgf/number format/.cd,
            fixed,
            fixed zerofill,
            precision=2,
        /tikz/.cd
    },
scaled ticks=false,
width=\textwidth, height=0.5\textwidth,
xlabel={flip steps},
ylabel={bottleneck distance},
]
\addplot graphics [xmin=0, xmax=1000,ymin=0, ymax=0.08] {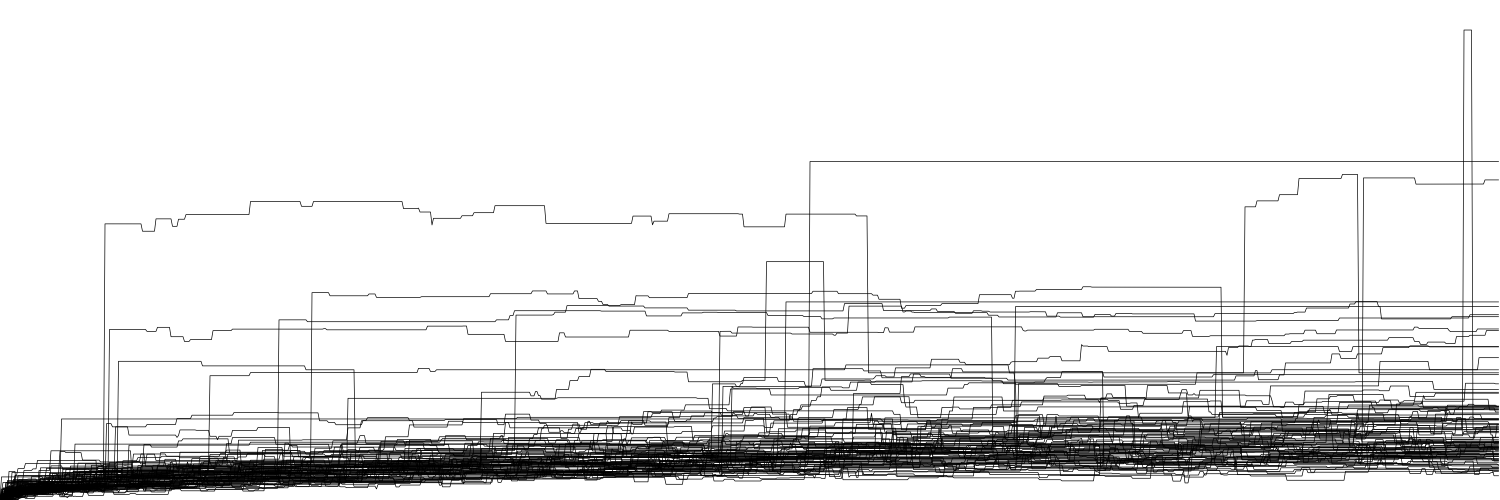}; 
\addplot[red, thick, line legend] coordinates{(0,0.0717) (1000,0.0717)};
\end{axis}
\end{tikzpicture}
}
\end{subfigure}
\caption{Effect of iterating small geographic perturbations (single-unit flip steps) on persistence. The red line indicates the average bottleneck distance  over the ${100 \choose 2}$ plans in a ReCom ensemble $\mathcal{E}=\{P^{(1)},\dots,P^{(100)}\}$, shown as a baseline for the distance between approximately independent plans. For each $i=1,\dots,100$, we run 1000 flip steps from  $P^{(i)}$ and track the bottleneck distance to $P^{(i)}$.  For $k=18$, the persistence diagrams are seen to be stable. The same is true for \emph{most} 50-district plans, but not all.}
\label{flipped_boxes}
\end{figure}

We repeat the 1000-flip experiment on a larger set of initial maps to conclude the investigation of stability.
We take one hundred Congressional and one hundred state Senate maps from the ReCom ensembles and make a sequence of 1000 flip steps starting from each of them. After each flip step, we compute the bottleneck distance between the initial plan and the current one. The results are shown in Figure \ref{flipped_boxes}. 
 Looking at these bottleneck distance trace plots, we see that for the majority of maps, randomly flipping 1000 geographic units had a modest effect on the persistence diagram, typically in the order of 0.01 in bottleneck distance. This raises our confidence that the methods outlined here are reasonably robust to small geographic perturbation.

\section{Conclusion}
We have examined three approaches to the use of persistent homology in analyzing the complicated geography of voting at different scales of aggregation:  zoning, election comparison, and signatures of gerrymandering.  In all three cases, our methods use persistence diagrams to summarize a large ensemble of alternative districting plans.  While the space of {\em all} valid redistricting plans for a state is intractably vast by almost any definition, tools like these allow us to cut through that complexity to uncover key information about the relationship of district lines to representation.  

Vote-filtered dual graphs have been shown to lend themselves to interpretable (and nontrivial) observations about redistricting, opening up several directions for future research.

{\em Classifying salient geoclusters.}
A long-standing debate in the political science literature centers on the question of whether racial gerrymandering protections like the Voting Rights Act actually produce counterintuitive partisan effects by disallowing the splitting of communities of color.  The ideas introduced here can 
reframe the question in a far more sophisticated and general way, locating areas in a state that are split differently in party-biased ensembles or whose splitting correlates with party advantage.  

{\em Metrics on the space of plans.}
Converting districting plans to persistence diagrams allows us to compute a bottleneck distance between districting plans given some fixed vote data. It will be interesting to compare the bottleneck distance to other notions of distance between districting plans, such as the optimal transport-based distance studied in \cite{abrishami2019geometry} and the information-theoretic distance in \cite{GuthNiehWeighill}. 

{\em Embeddings.}
There are several methods in the literature for vectorizing persistence diagrams, such as {\em persistence landscapes} \cite{bubenik2015statistical} and {\em persistence images} \cite{adams2017persistence}---that is, each such method provides an embedding of a set of persistence diagrams into a Hilbert space as a set of vectors. We can use these methods to embed (ensembles of)  plans in Euclidean space either for visualization or for clustering. 

{\em Distributions.}
Finally, the use of summary diagrams enables an approach to the comparison of different sampling techniques, such as by the use of different Markov chains.  
We can fit distributions to ensembles of persistence diagrams---either by first vectorizing them and fitting classical distributions, or by applying more recent techniques directly to persistence diagrams themselves \cite{maroulas2020bayesian}---allowing a summarized comparison of sampling distributions.

\newpage
\bibliography{tda}
\bibliographystyle{alpha}

\newpage
\appendix
\section{Supplemental figures}
\renewcommand\thefigure{S\arabic{figure}}    
\setcounter{figure}{0}

\begin{figure}[ht]
\centering
\begin{subfigure}{0.3\textwidth}
\centering
\caption*{Diameters}
\resizebox{\columnwidth}{!}{%
\begin{tikzpicture}
\begin{axis}[
   enlargelimits=false, axis on top,
   ymin=0, ymax=650,
   xmin=3, xmax=8
]
\addplot graphics [xmin=3, xmax=8,ymin=-1, ymax=650] {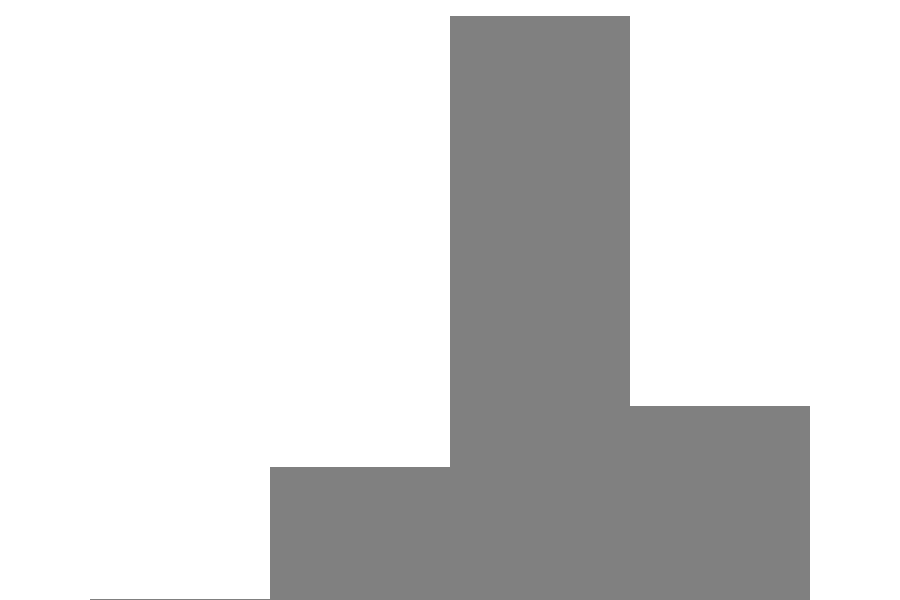}; 
\end{axis}
\end{tikzpicture}%
}
\end{subfigure}%
\begin{subfigure}{0.3\textwidth}
\centering
\caption*{Max degrees}
\resizebox{\columnwidth}{!}{%
\begin{tikzpicture}
\begin{axis}[
	enlargelimits=false, axis on top,
   	ymin=0
]
\addplot graphics [xmin=4, xmax=9,ymin=0, ymax=550] {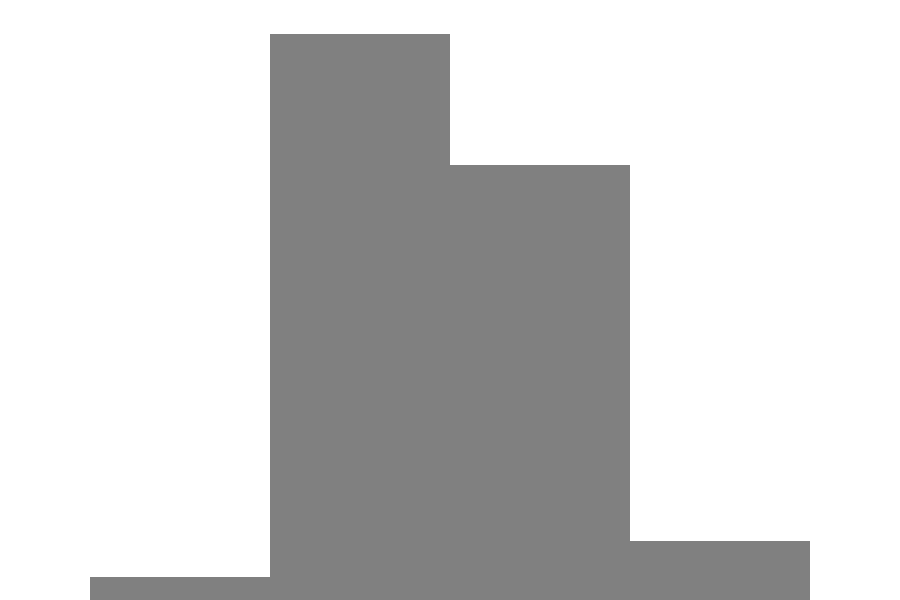}; 
\end{axis}
\end{tikzpicture}%
}
\end{subfigure}%

\begin{subfigure}{0.3\textwidth}
\centering
\caption*{Mean degrees}
\resizebox{\columnwidth}{!}{%
\begin{tikzpicture}
\begin{axis}[
	enlargelimits=false, axis on top,
   	ymin=0
]
\addplot graphics [xmin=3.5, xmax=4.5,ymin=0, ymax=320] {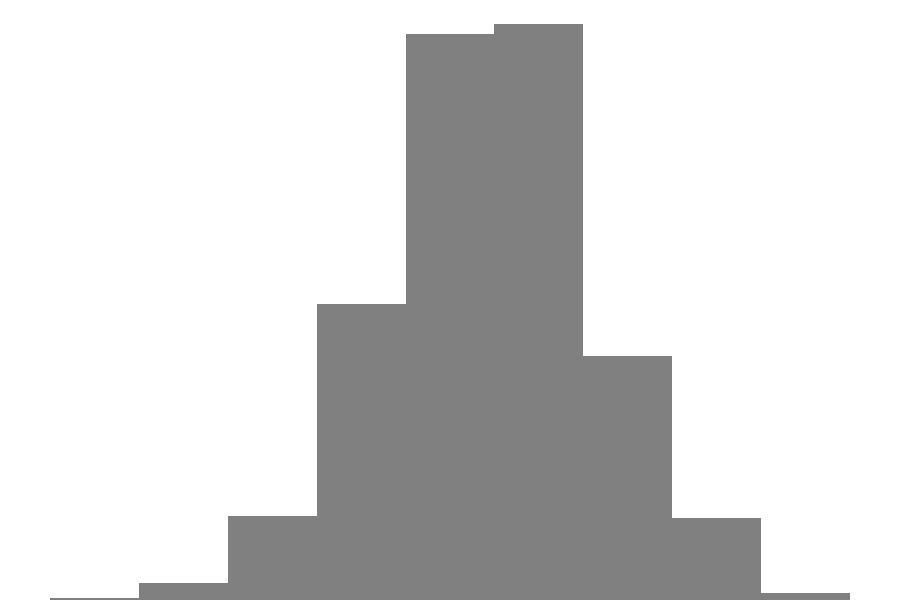}; 
\end{axis}
\end{tikzpicture}%
}
\end{subfigure}%
\begin{subfigure}{0.3\textwidth}
\centering
\caption*{Densities}
\resizebox{\columnwidth}{!}{%
\begin{tikzpicture}
\begin{axis}[
	enlargelimits=false, axis on top,
	xtick={0.21, 0.22,0.23,0.24,0.25,0.26}   	
]
\addplot graphics [xmin=0.2, xmax=0.27,ymin=0, ymax=320] {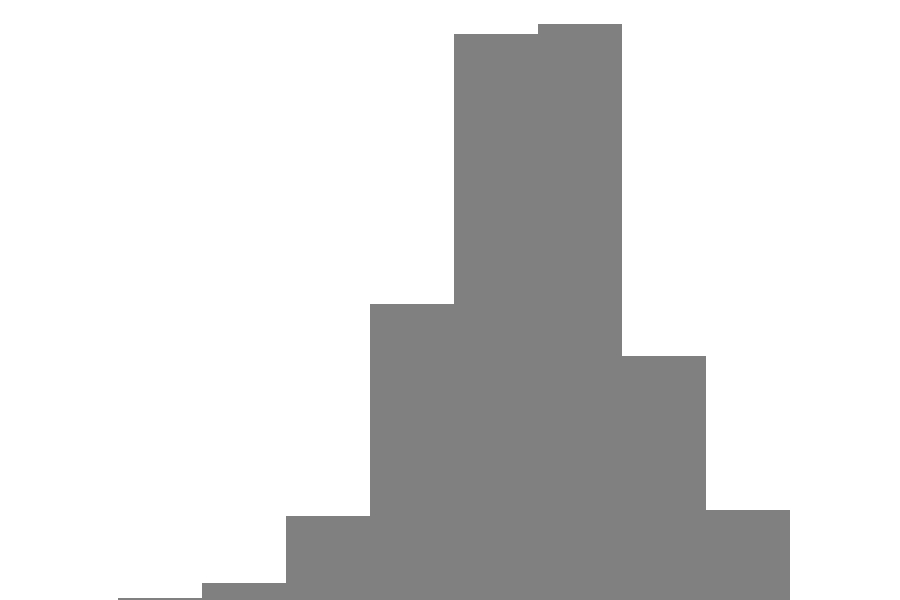}; 
\end{axis}
\end{tikzpicture}%
}
\end{subfigure}%
\caption{Histograms of statistics of the dual graphs for the ensemble of 1000 18-district plans for Pennsylvania.
({\em Density} is the ratio of the number of edges to the number of possible edges, ${18 \choose 2}=153$.)}
\label{fig:graph_statistics}
\end{figure}

\begin{figure}[ht] 
\foreach[count=\i] \n in {0,1,2,3,4,5,6,7,8,9,10}
{
\begin{subfigure}{0.19\textwidth}
\pdsubfig{PA_noaxes/m2m_203districts_PRES16_PD\n.png}%
\end{subfigure}
\begin{subfigure}{0.29\textwidth}
\includegraphics[width=\textwidth]{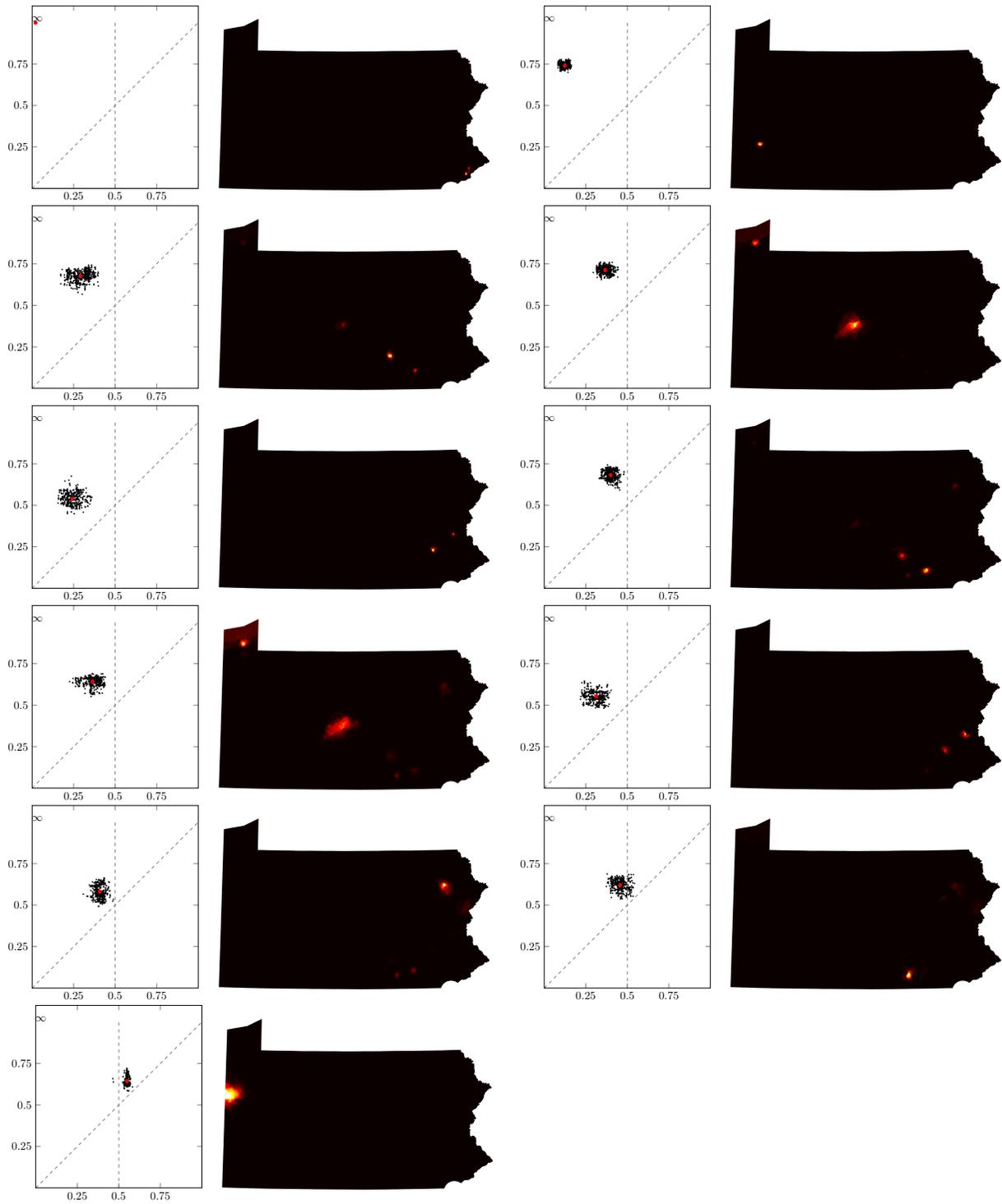}
\end{subfigure}
}
\begin{subfigure}{0.48\textwidth}
\hfill
\end{subfigure}
\caption{Geographical localization of eleven Fr\'echet features in Pennsylvania House plans $(k=203)$ with respect to PRES16 voting.  These now pinpoint small cities, all the way down to Hermitage (population 16,220).}
\label{PAclasses203}
\end{figure}

\begin{figure}[ht]
\begin{subfigure}{0.16\textwidth}
\pdtwosubfig{PA_noaxes/50biased_points0_DEM.png}{PA_noaxes/50biased_points0_REP.png}%
\end{subfigure}
\begin{subfigure}{0.24\textwidth}
\includegraphics[width=\textwidth]{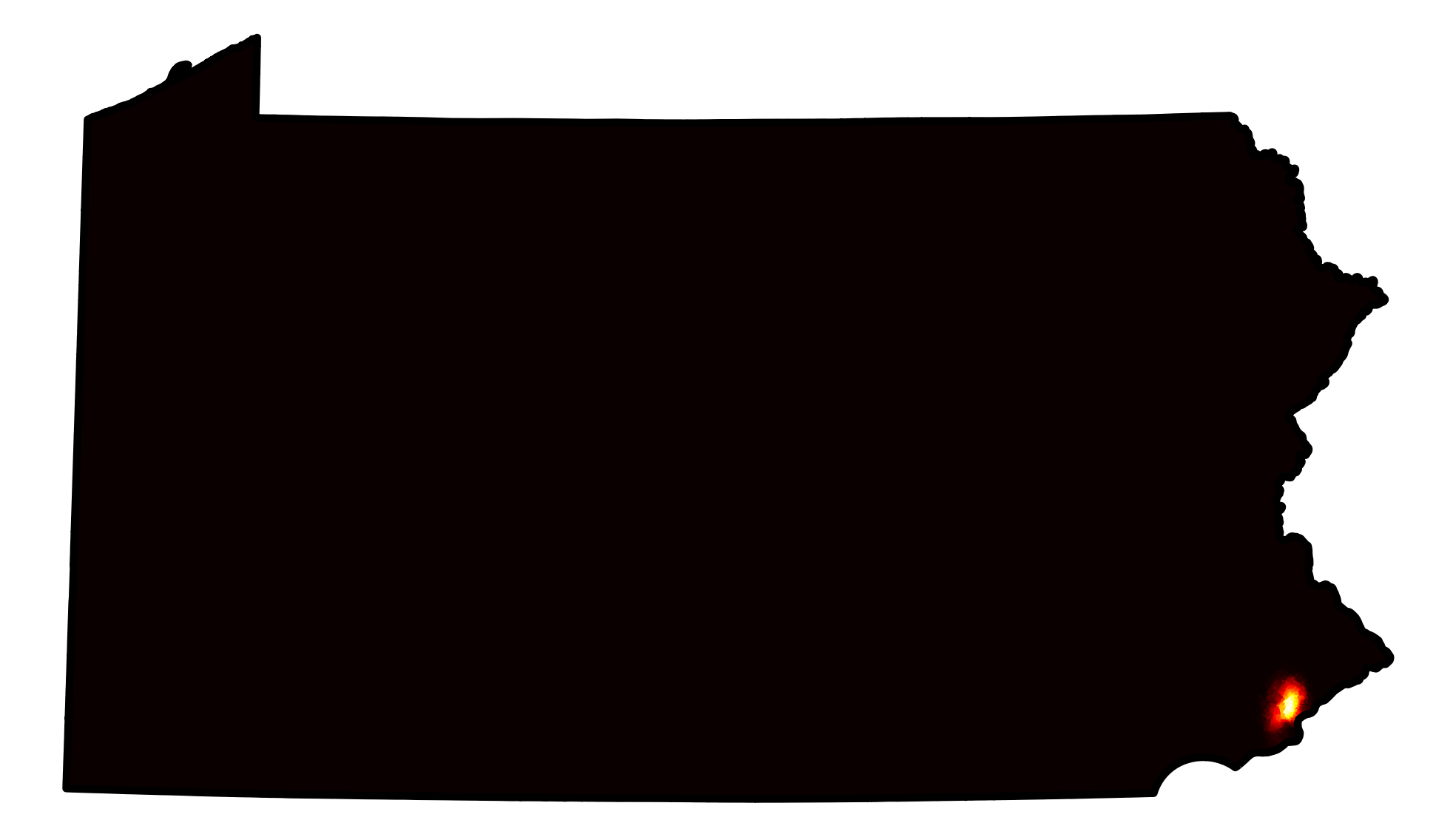}
\end{subfigure}
\begin{subfigure}{0.24\textwidth}
\includegraphics[width=\textwidth]{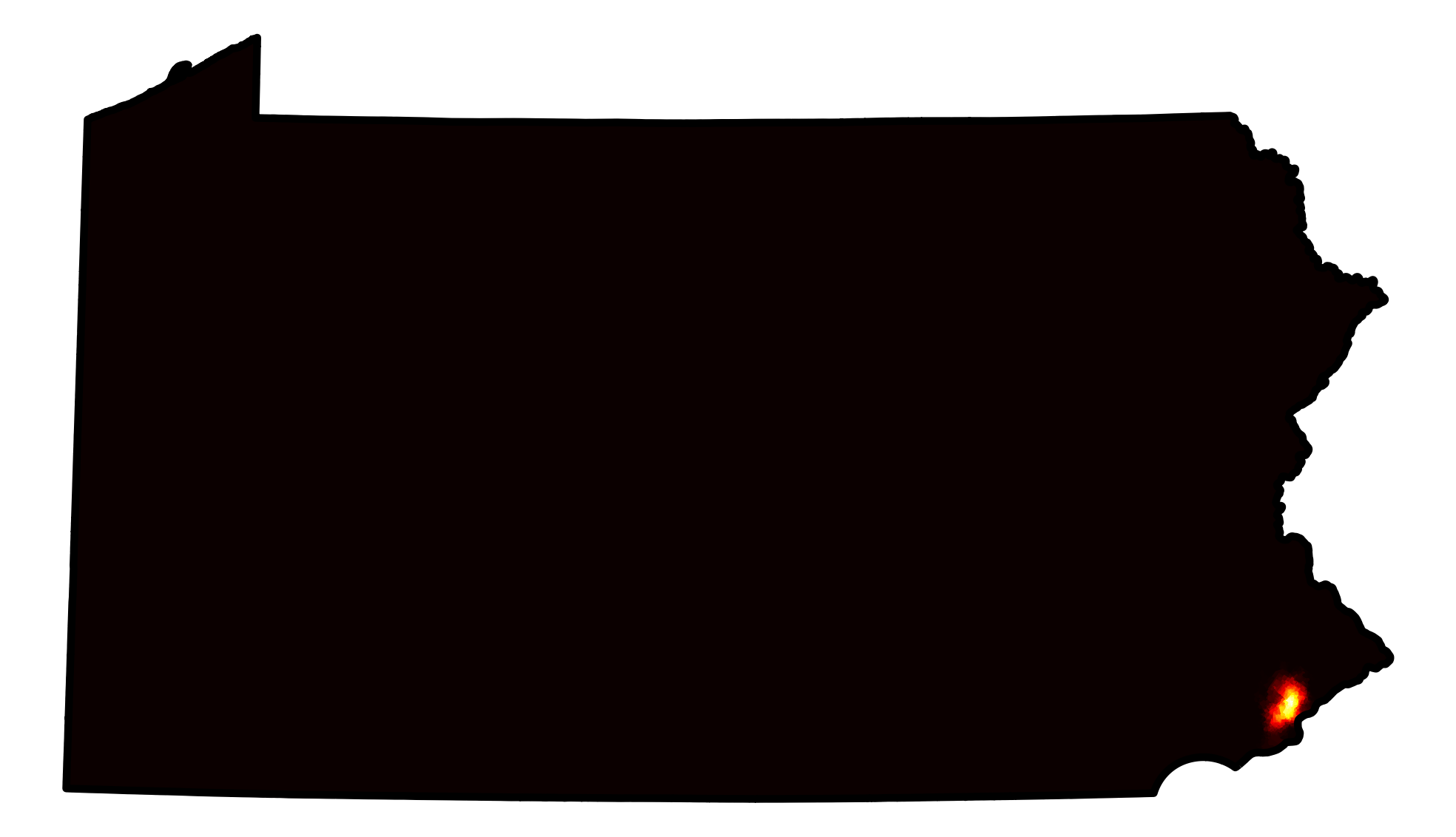}
\end{subfigure}

\begin{subfigure}{0.16\textwidth}
\pdtwosubfig{PA_noaxes/50biased_points1_DEM.png}{PA_noaxes/50biased_points1_REP.png}%
\end{subfigure}
\begin{subfigure}{0.24\textwidth}
\includegraphics[width=\textwidth]{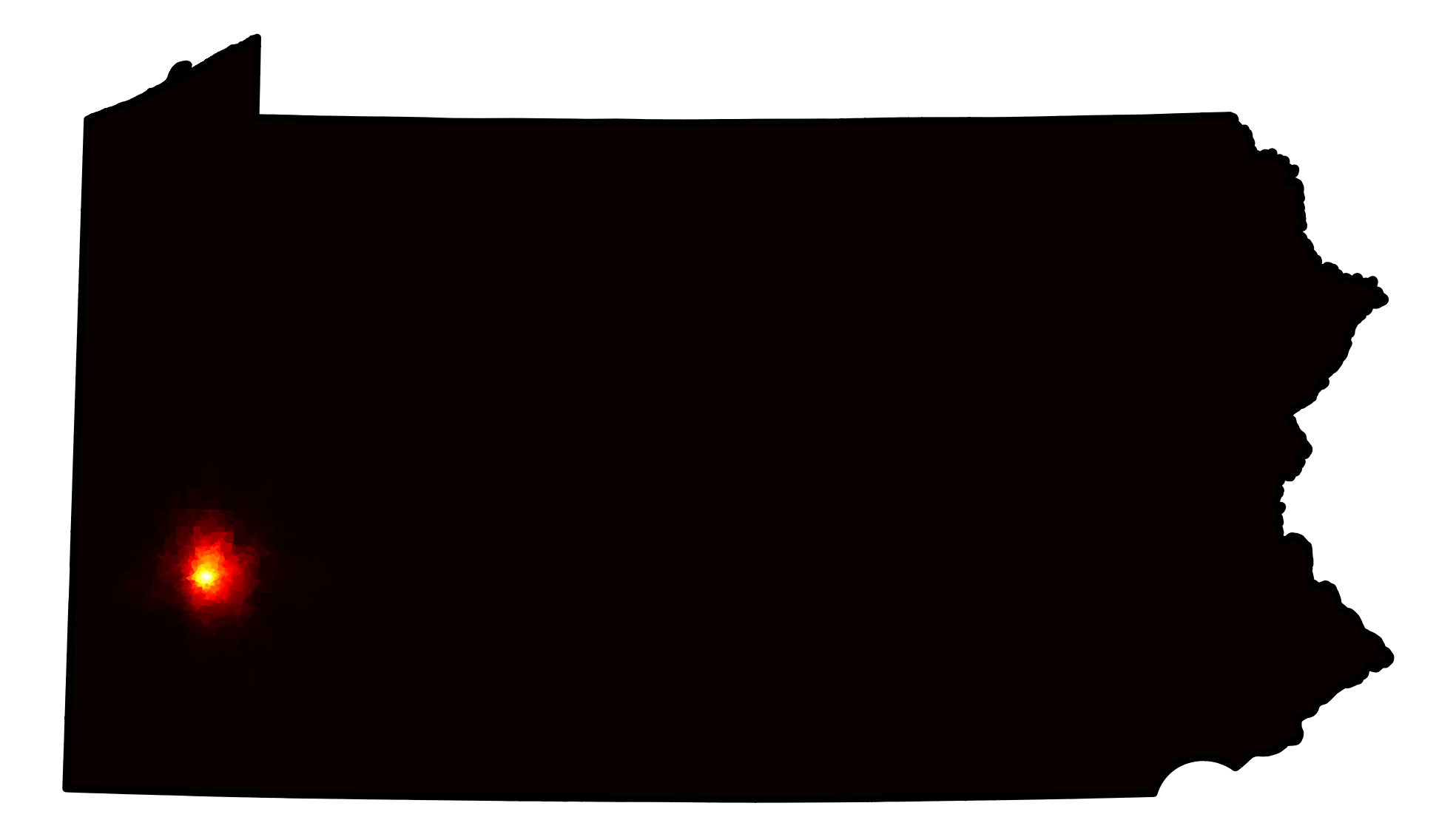}
\end{subfigure}
\begin{subfigure}{0.24\textwidth}
\includegraphics[width=\textwidth]{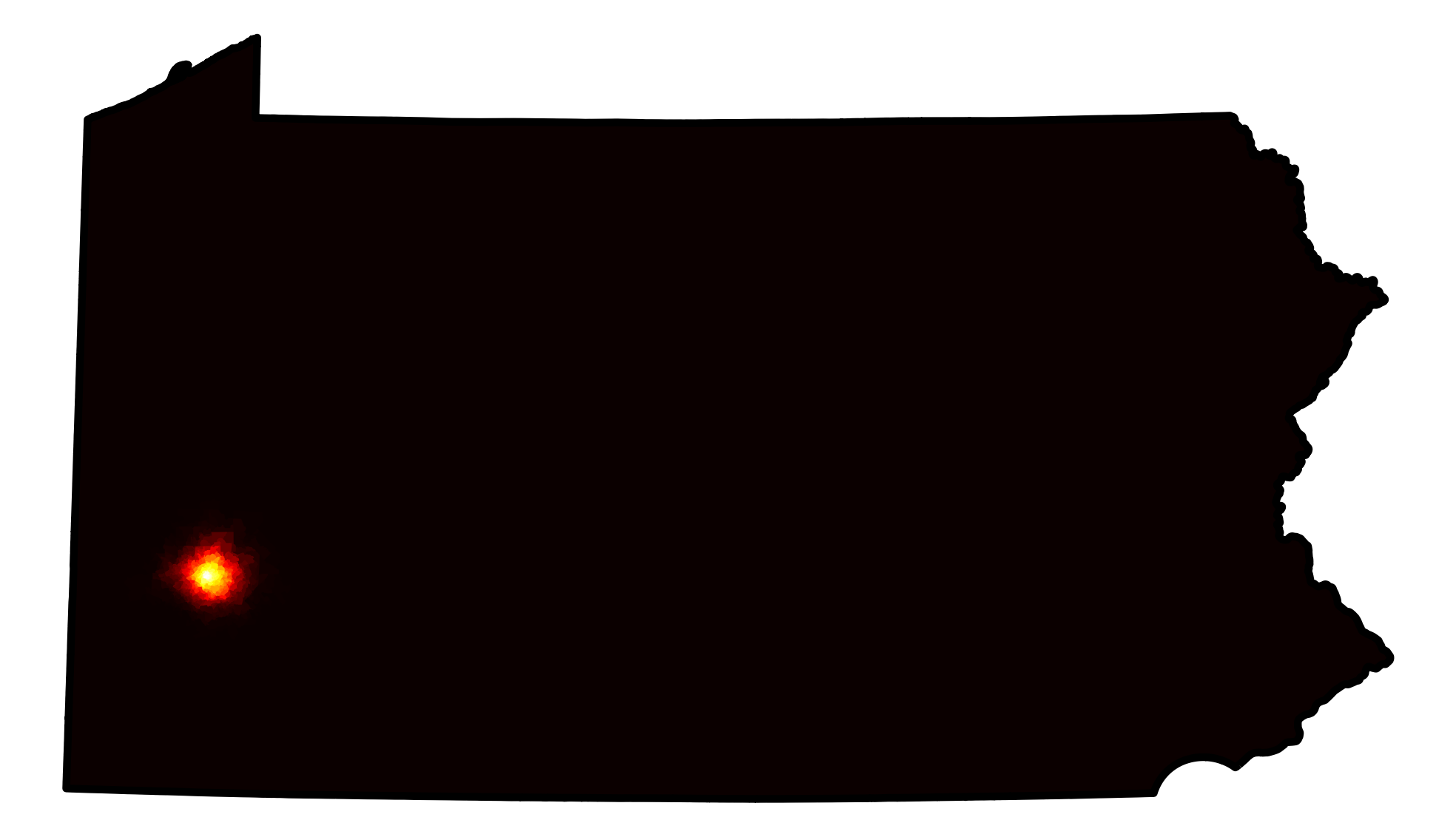}
\end{subfigure}

\begin{subfigure}{0.16\textwidth}
\pdtwosubfig{PA_noaxes/50biased_points2_DEM.png}{PA_noaxes/50biased_points2_REP.png}%
\end{subfigure}
\begin{subfigure}{0.24\textwidth}
\includegraphics[width=\textwidth]{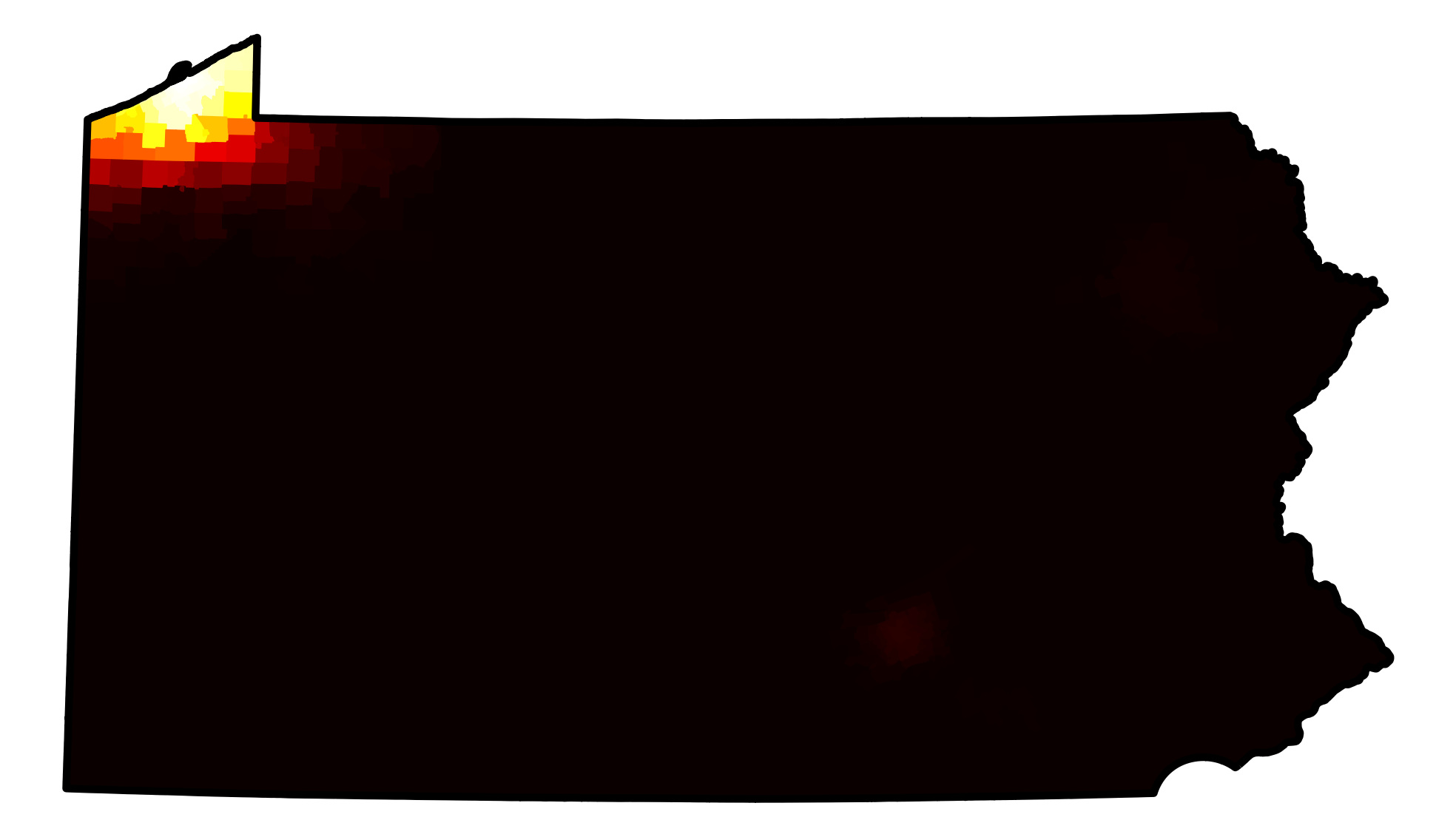}
\end{subfigure}
\begin{subfigure}{0.24\textwidth}
\includegraphics[width=\textwidth]{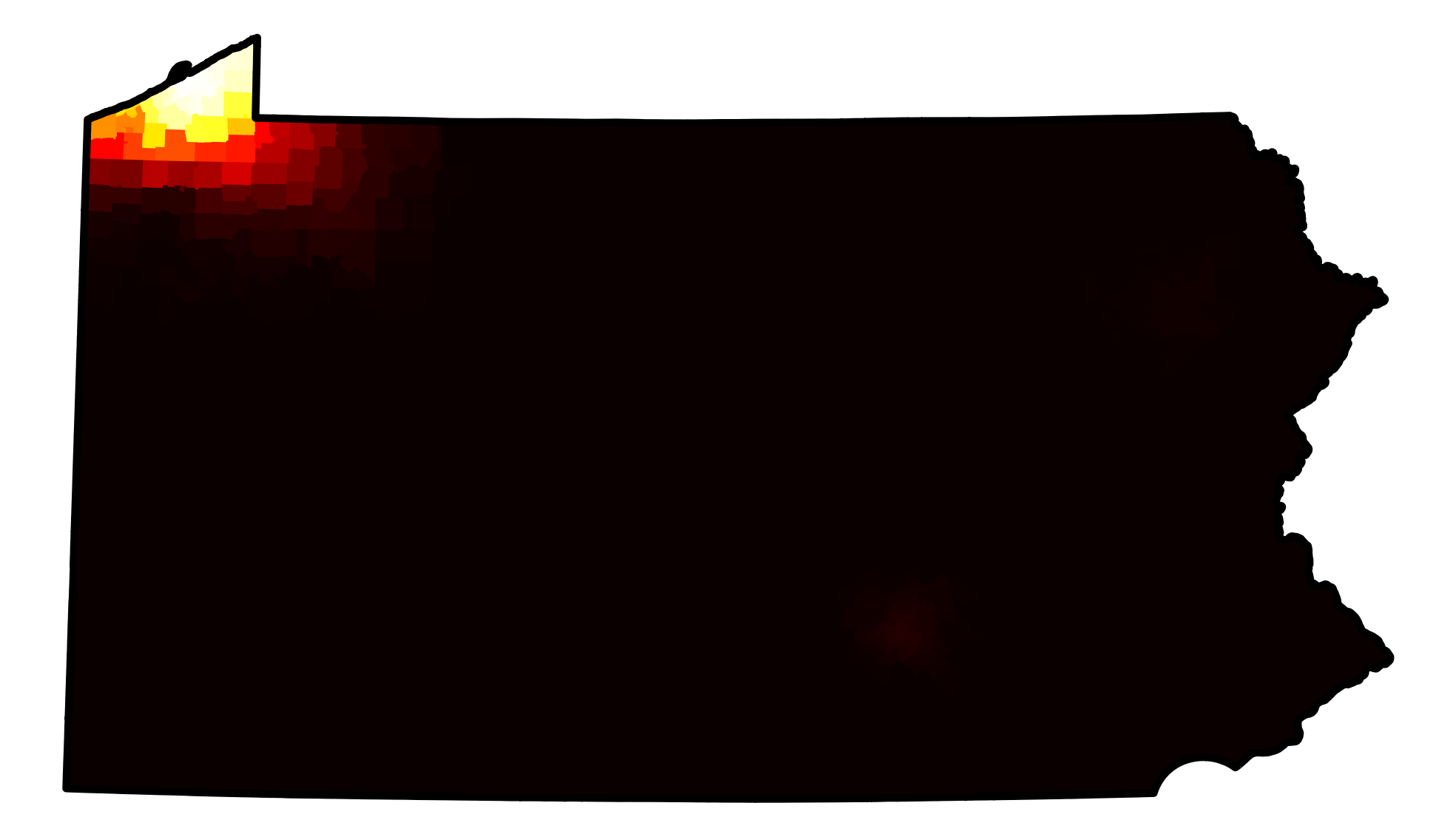}
\end{subfigure}

\begin{subfigure}{0.16\textwidth}
\pdtwosubfig{PA_noaxes/50biased_points3_DEM.png}{PA_noaxes/50biased_points3_REP.png}%
\end{subfigure}
\begin{subfigure}{0.24\textwidth}
\includegraphics[width=\textwidth]{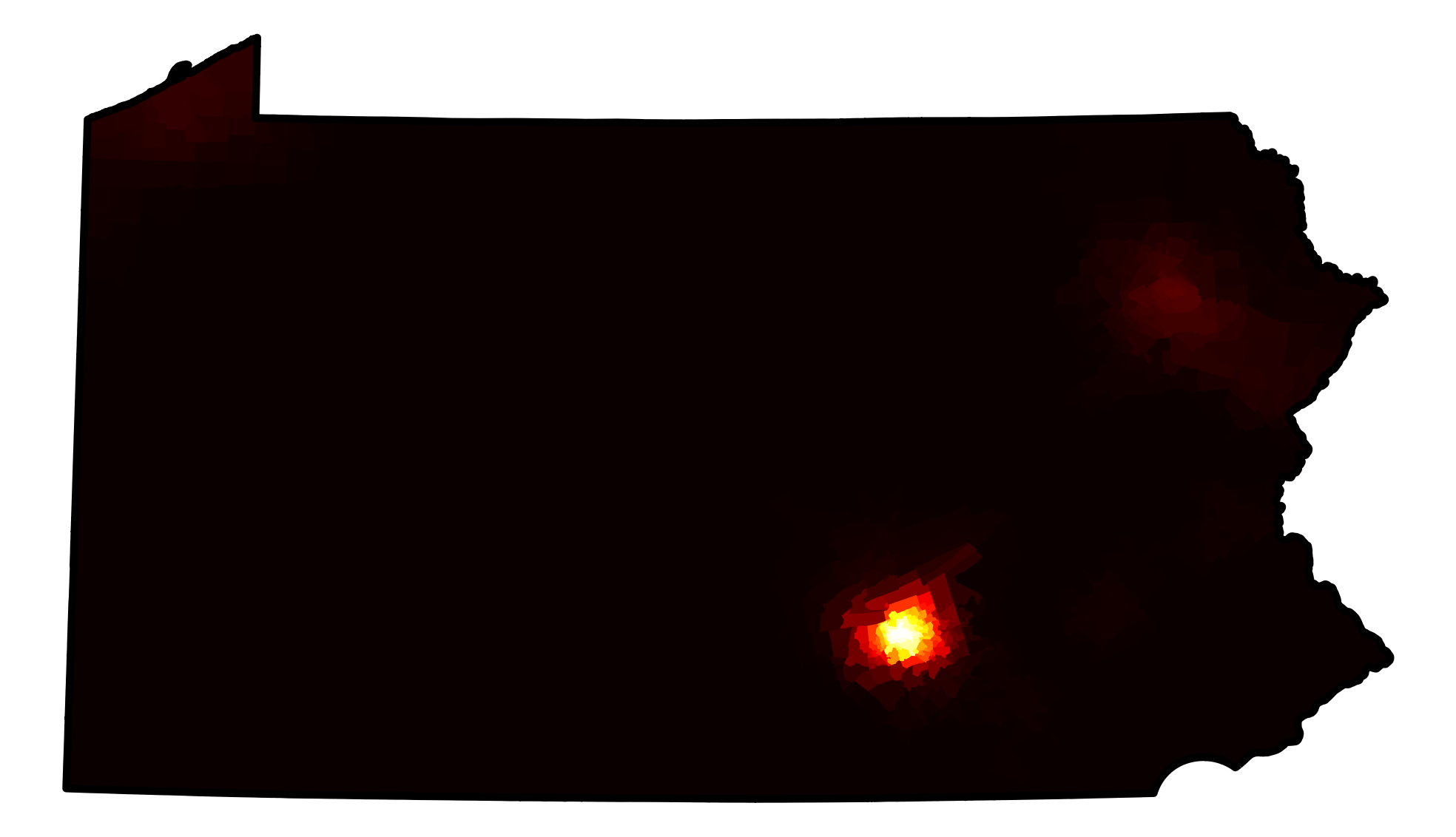}
\end{subfigure}
\begin{subfigure}{0.24\textwidth}
\includegraphics[width=\textwidth]{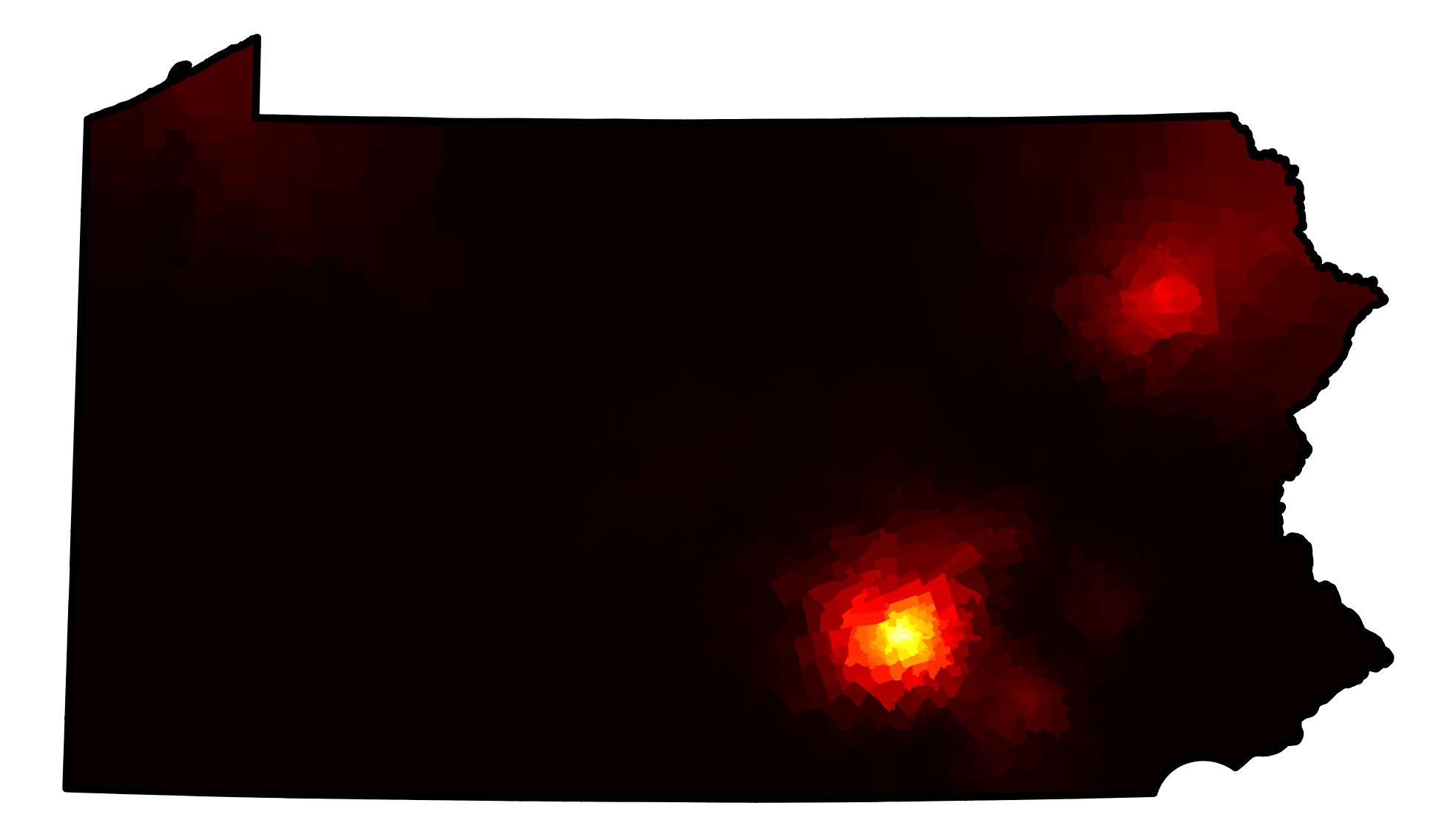}
\end{subfigure}

\begin{subfigure}{0.16\textwidth}
\pdtwosubfig{PA_noaxes/50biased_points4_DEM.png}{PA_noaxes/50biased_points4_REP.png}%
\end{subfigure}
\begin{subfigure}{0.24\textwidth}
\includegraphics[width=\textwidth]{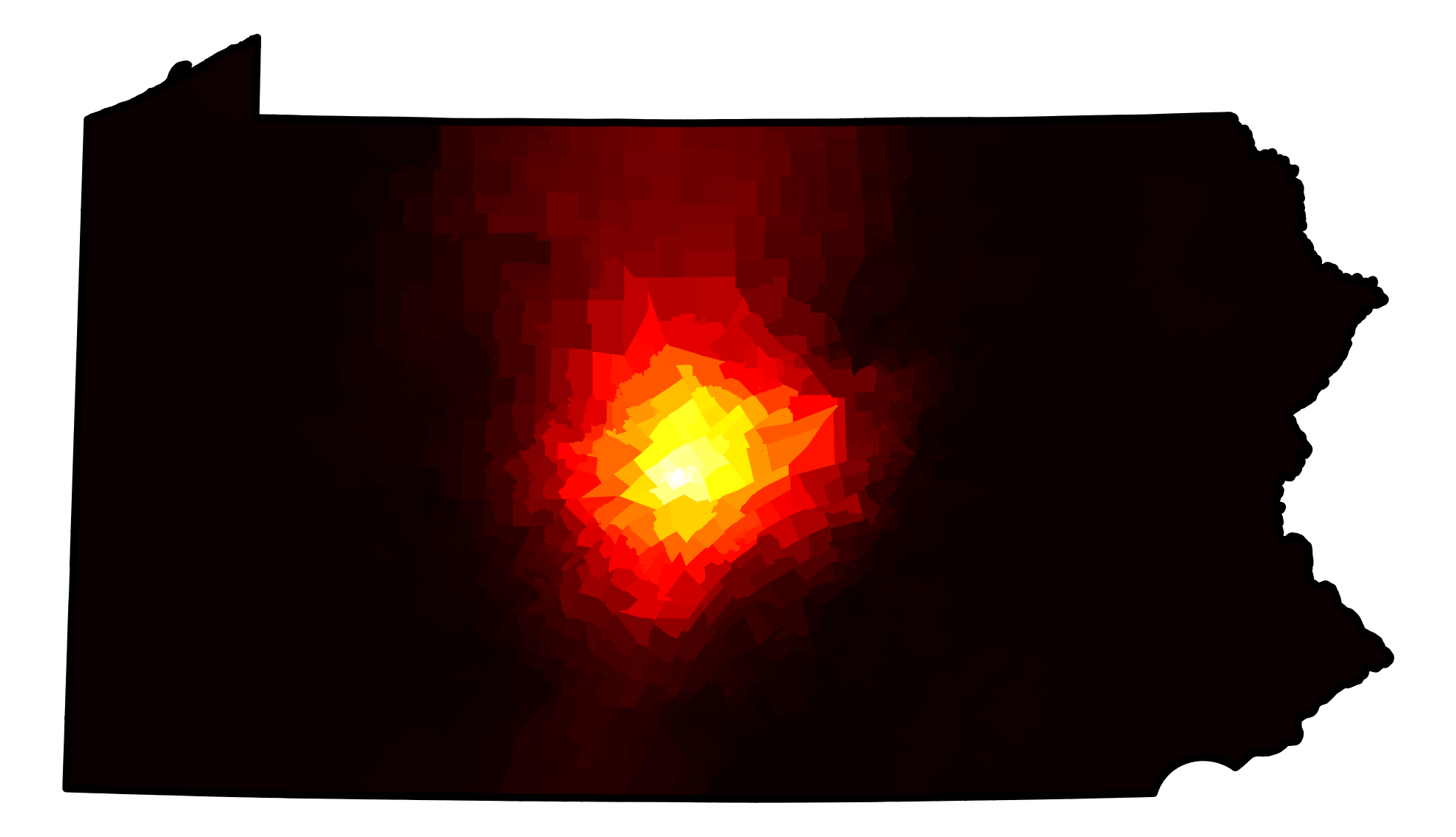}
\end{subfigure}
\begin{subfigure}{0.24\textwidth}
\includegraphics[width=\textwidth]{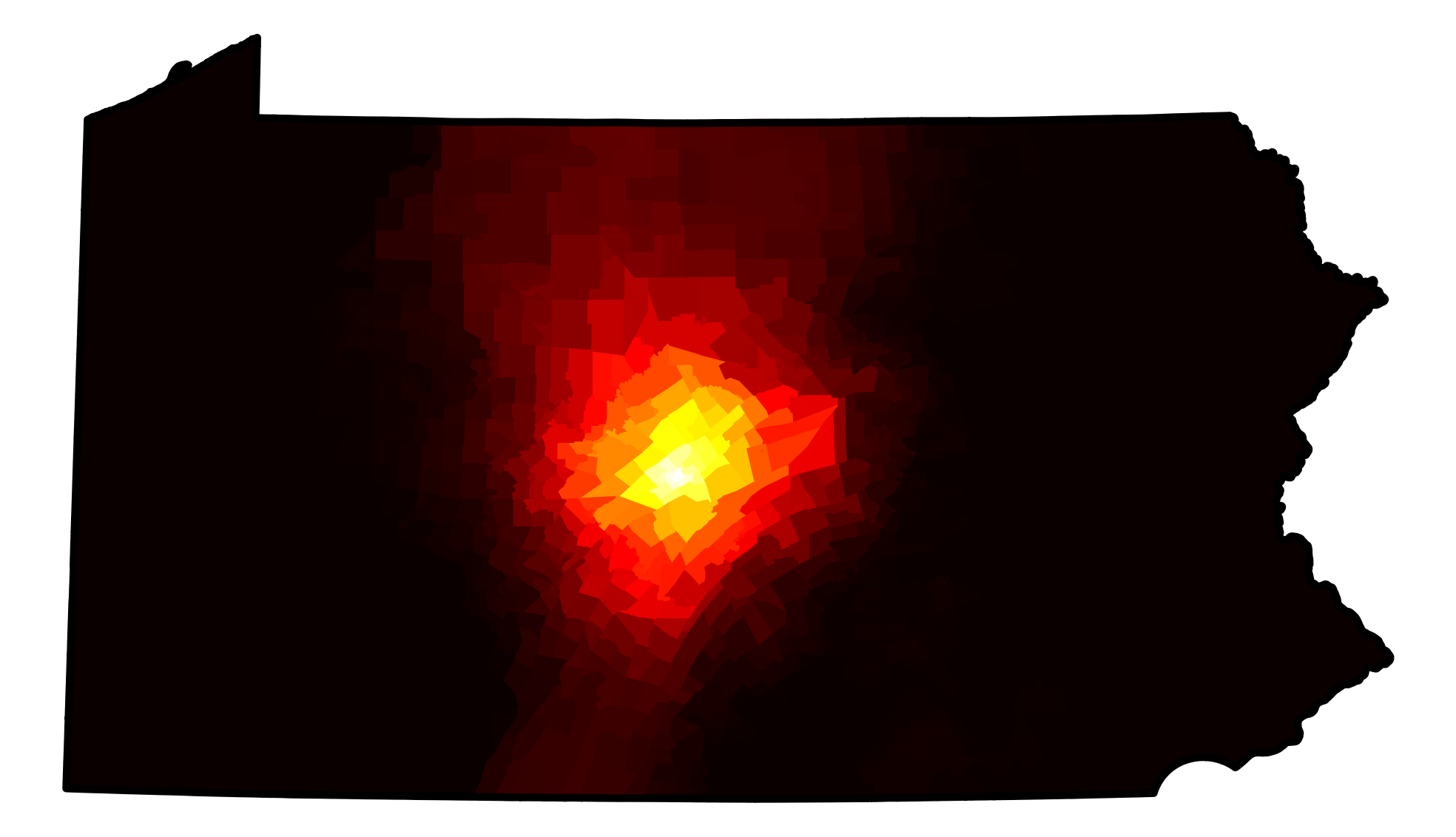}
\end{subfigure}

\begin{subfigure}{0.16\textwidth}
\pdtwosubfig{PA_noaxes/50biased_points5_DEM.png}{PA_noaxes/50biased_points5_REP.png}%
\end{subfigure}
\begin{subfigure}{0.24\textwidth}
\includegraphics[width=\textwidth]{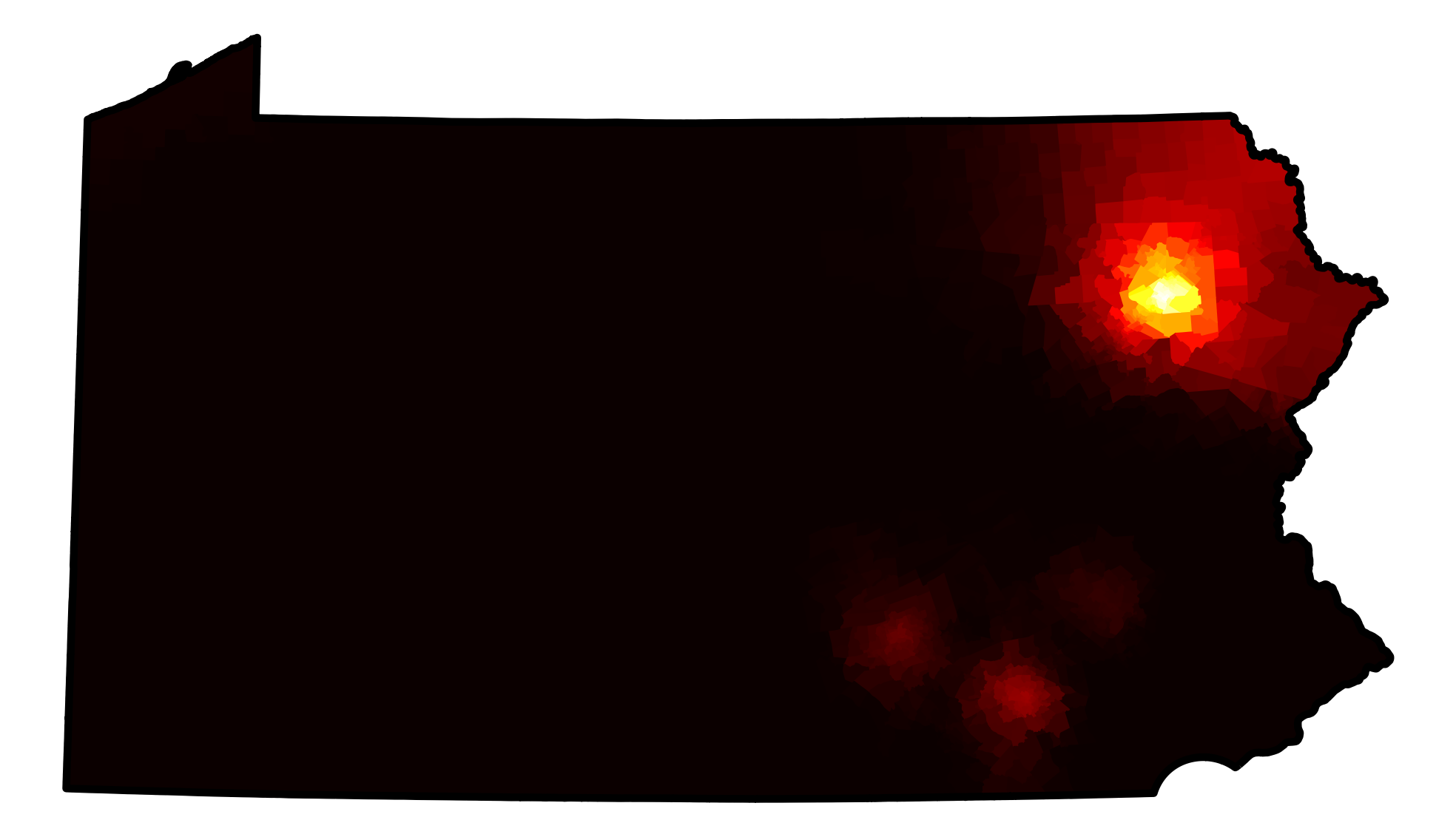}
\end{subfigure}
\begin{subfigure}{0.24\textwidth}
\includegraphics[width=\textwidth]{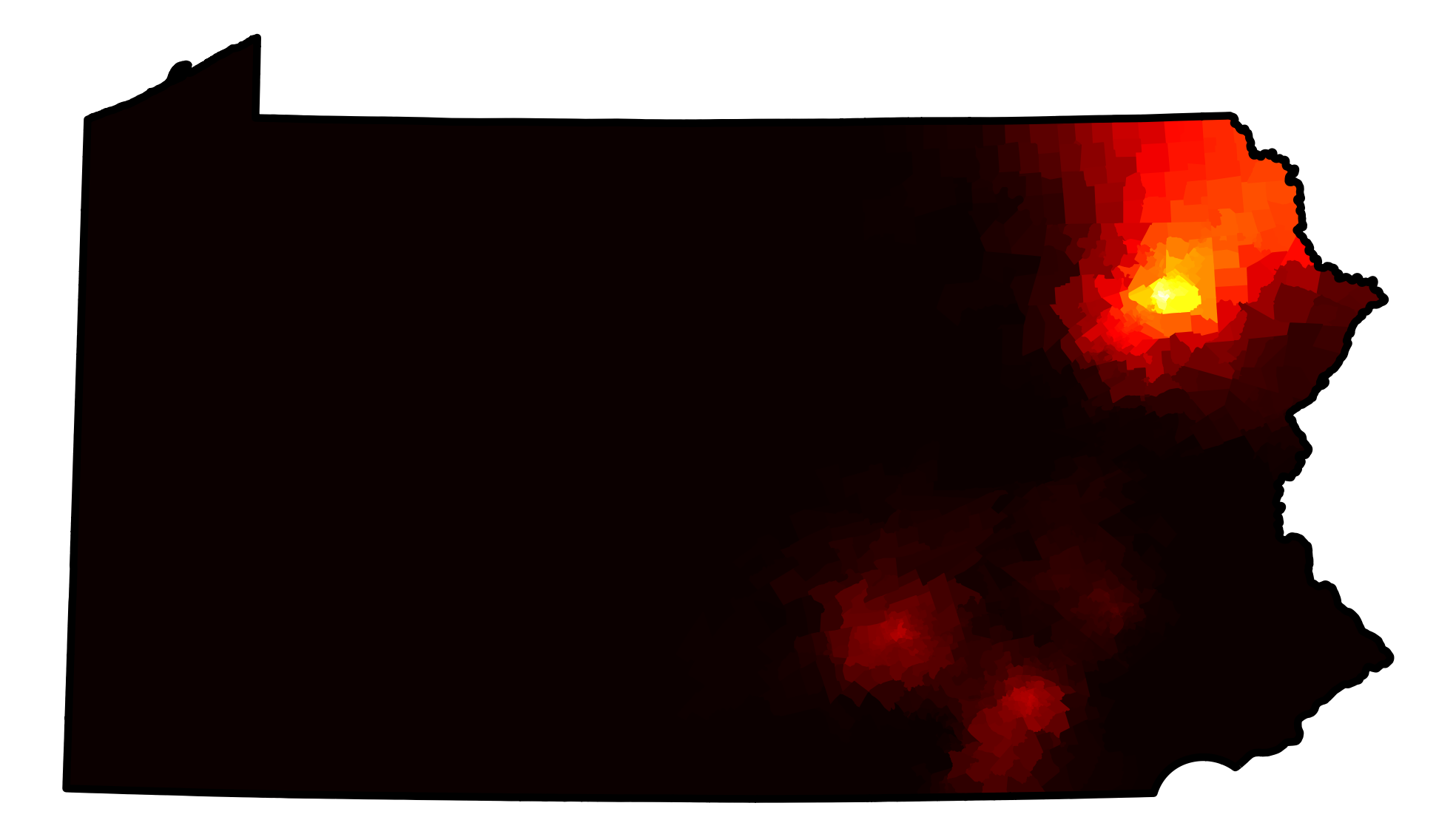}
\end{subfigure}

\begin{subfigure}{0.16\textwidth}
\pdtwosubfig{PA_noaxes/50biased_points6_DEM.png}{PA_noaxes/50biased_points6_REP.png}%
\end{subfigure}
\begin{subfigure}{0.24\textwidth}
\includegraphics[width=\textwidth]{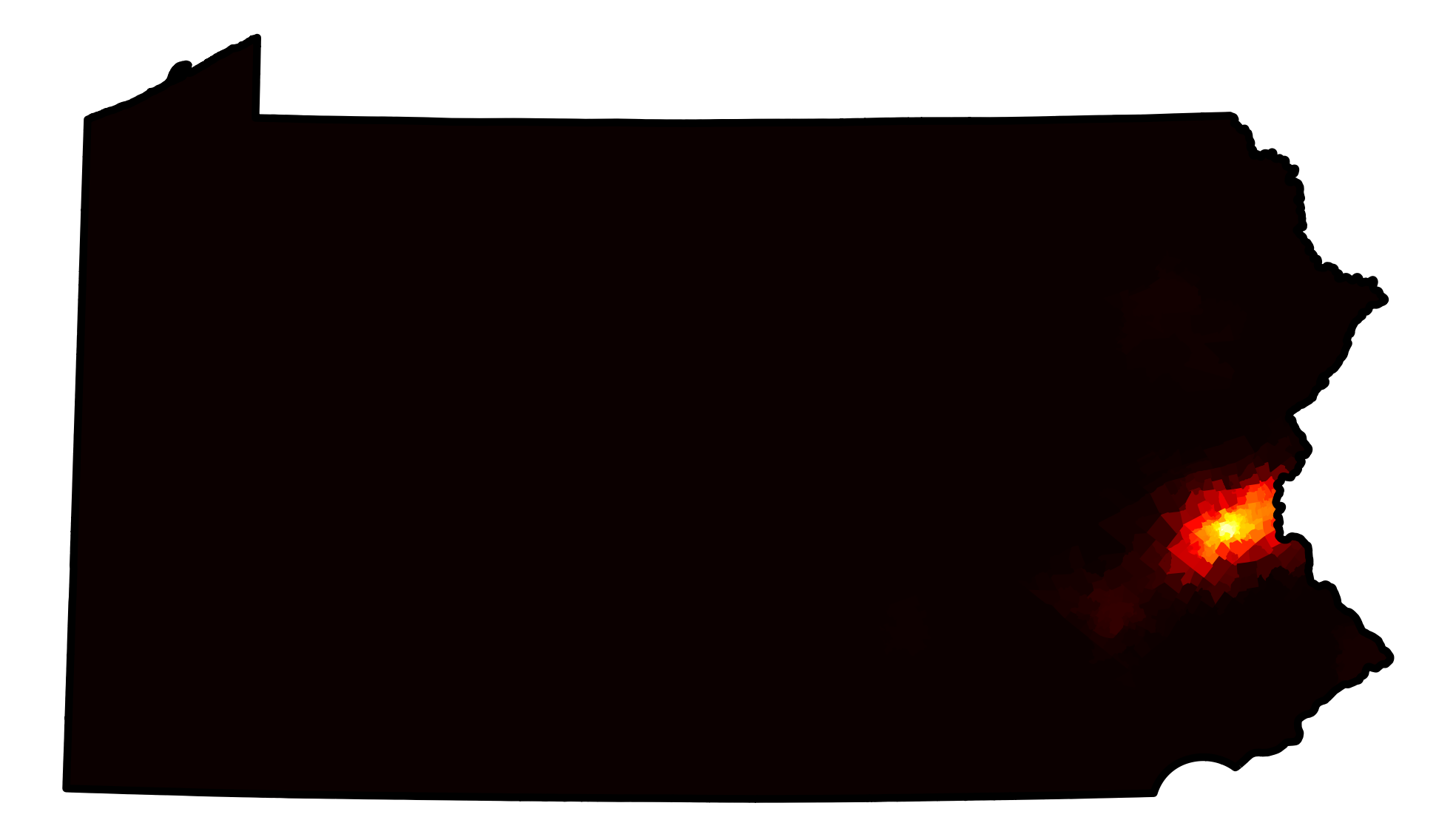}
\caption*{Dem.-favoring}
\end{subfigure}
\begin{subfigure}{0.24\textwidth}
\includegraphics[width=\textwidth]{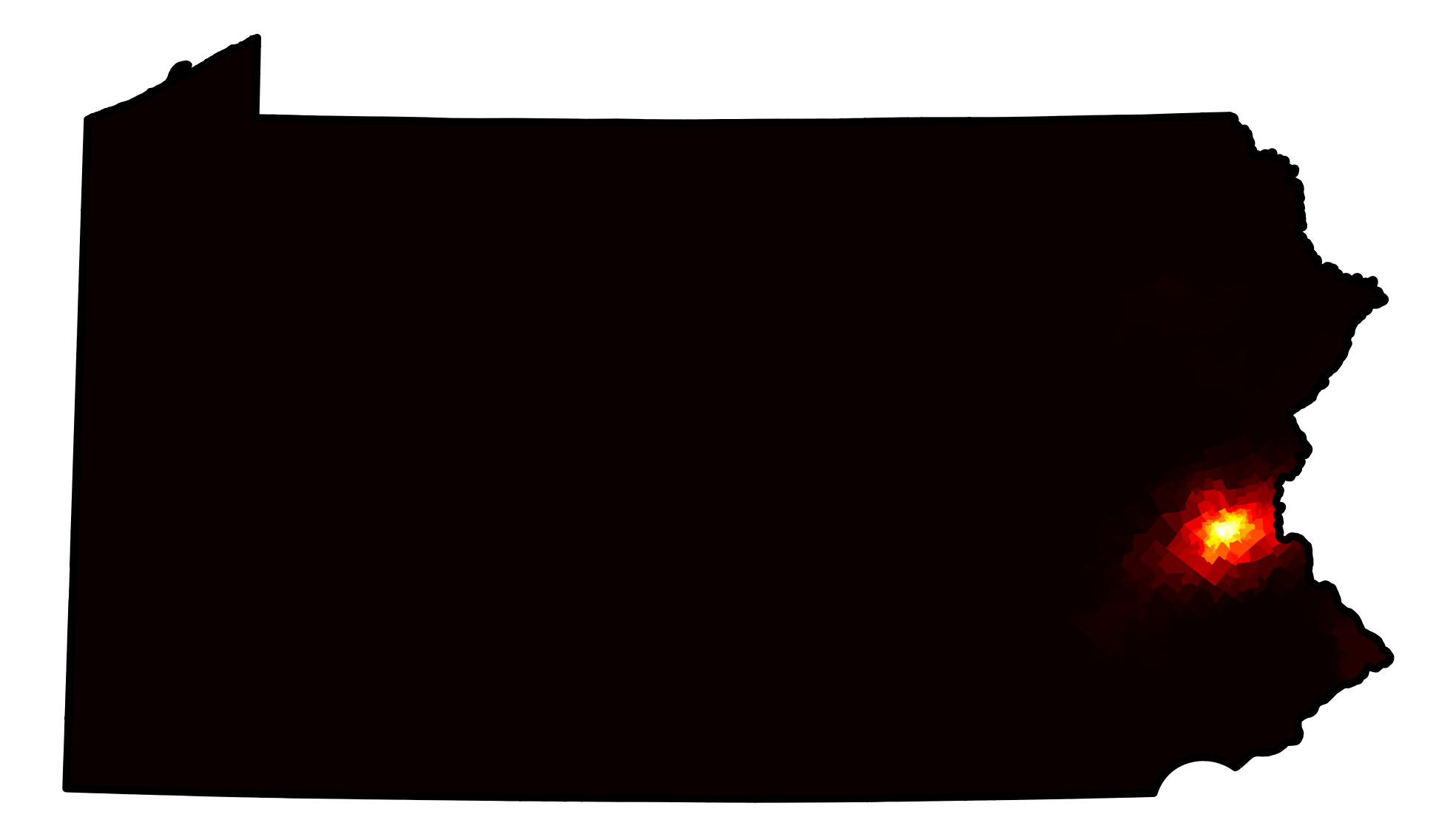}
\caption*{Rep.-favoring}
\end{subfigure}
\caption{Expanding on Figure~\ref{biasedmeansPAshort}.  Point plots and heat maps for the successive Fr\'echet features in the PA Senate ensembles that are biased for Democratic and Republican safe seats, respectively. Note the contrast between the rigidity of the Erie feature (similar placement for both ensembles) and the manipulability of the Harrisburg feature, though both anchor competitive districts.}
\label{biasedmeansPA}
\end{figure}

\begin{figure}[ht] 
\foreach[count=\i] \n in {0,1,2,3,4,5,6,7,8,9}
{
\begin{subfigure}{0.19\textwidth}
\pdsubfig{NC_noaxes/m2m_120districts_PRES16_PD\n.png}%
\end{subfigure}
\begin{subfigure}{0.29\textwidth}
\includegraphics[width=\textwidth]{NC_plots/120mappedFrechetPRES16_\n.png}
\end{subfigure}
}
\begin{subfigure}{0.48\textwidth}
\hfill
\end{subfigure}
\caption{Geographical localization of eleven Fr\'echet features in North Carolina Senate plans $(k=50)$ with respect to PRES16 voting.}
\label{NCclasses120}
\end{figure}

\begin{figure}[ht]
\begin{subfigure}{0.125\textwidth}
\pdtwosubfig{NC_noaxes/50biased_points0_DEM.png}{NC_noaxes/50biased_points0_REP.png}%
\end{subfigure}
\begin{subfigure}{0.23\textwidth}
\includegraphics[width=\textwidth]{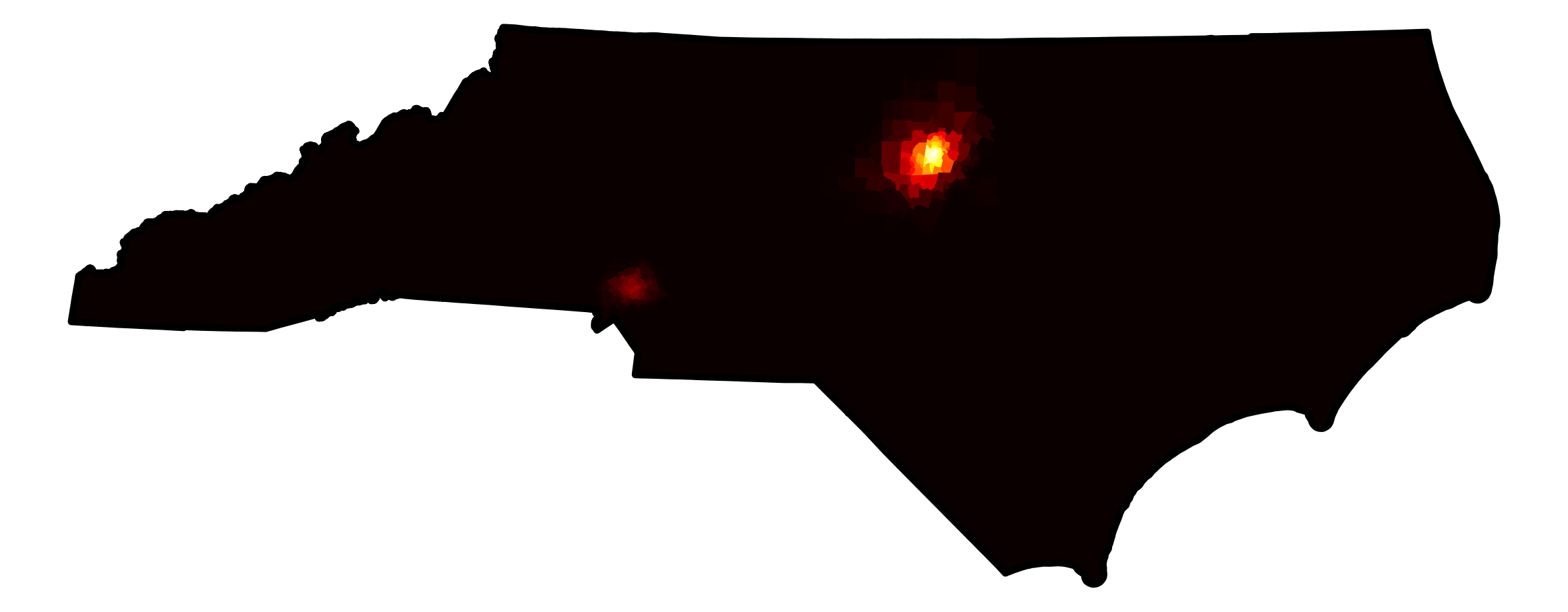}
\end{subfigure}
\begin{subfigure}{0.23\textwidth}
\includegraphics[width=\textwidth]{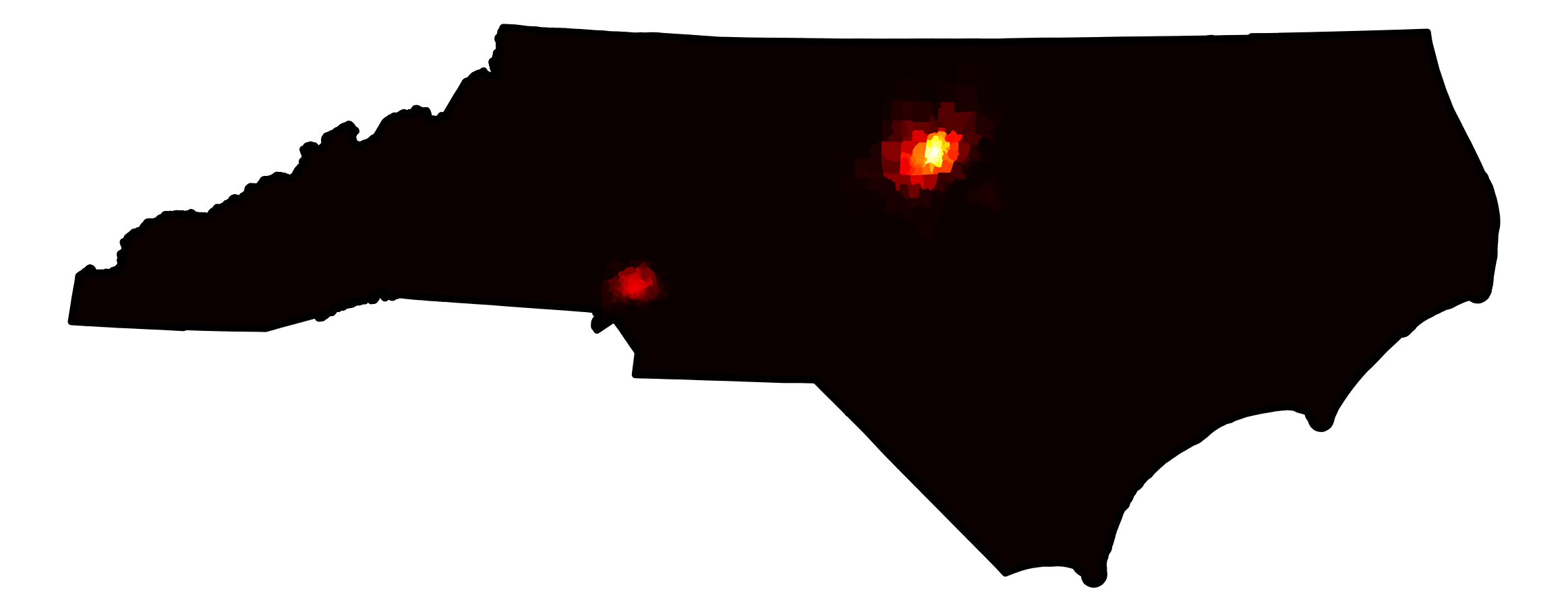}
\end{subfigure}

\begin{subfigure}{0.125\textwidth}
\pdtwosubfig{NC_noaxes/50biased_points1_DEM.png}{NC_noaxes/50biased_points1_REP.png}%
\end{subfigure}
\begin{subfigure}{0.23\textwidth}
\includegraphics[width=\textwidth]{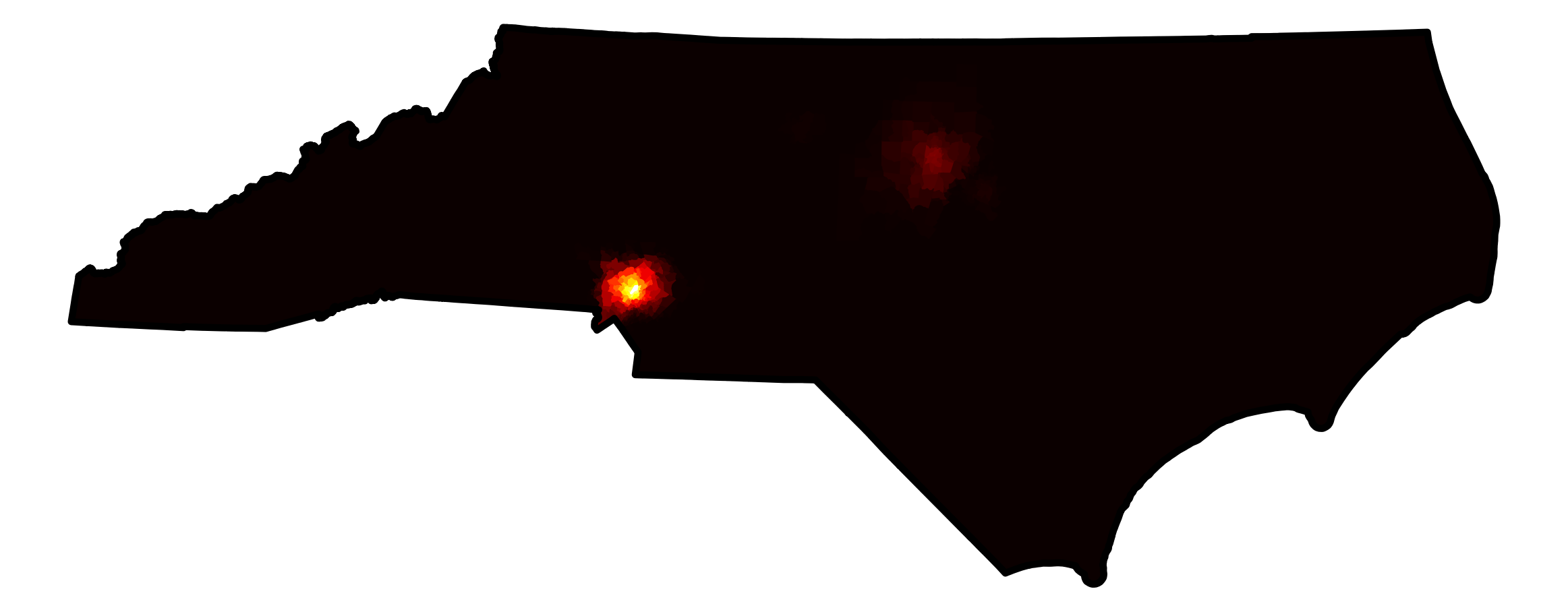}
\end{subfigure}
\begin{subfigure}{0.23\textwidth}
\includegraphics[width=\textwidth]{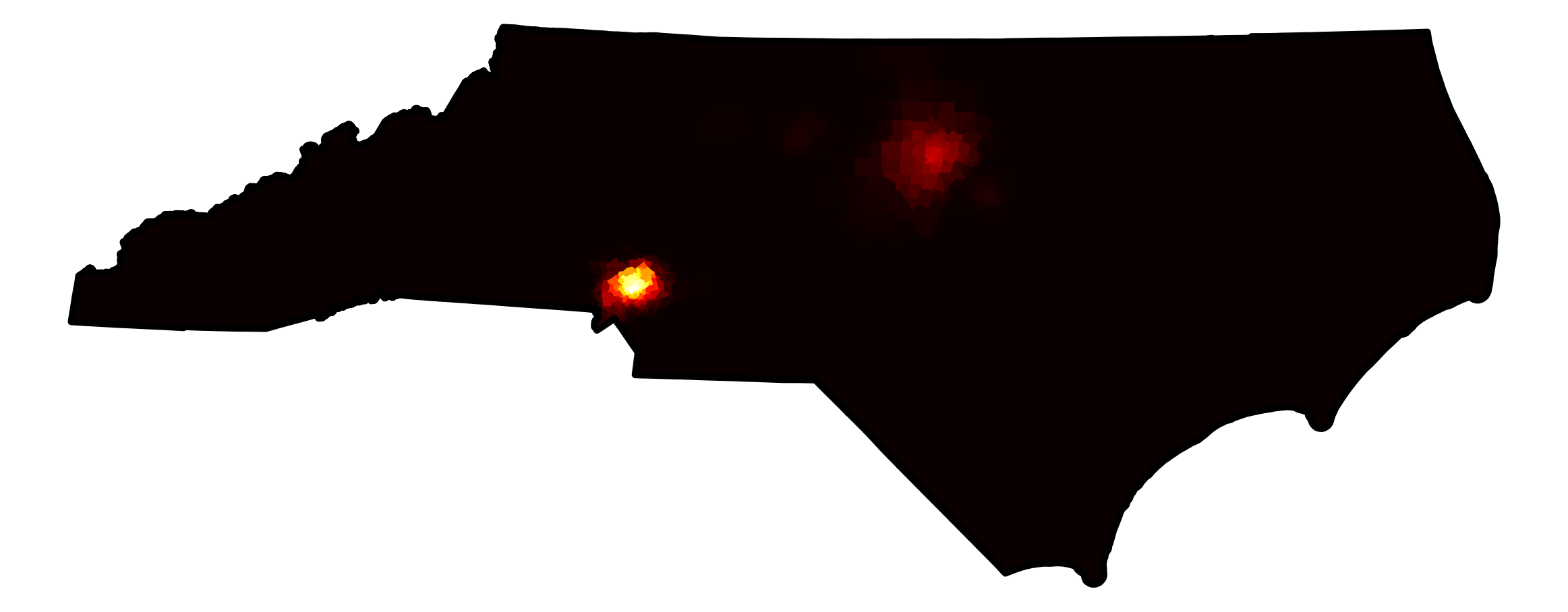}
\end{subfigure}

\begin{subfigure}{0.125\textwidth}
\pdtwosubfig{NC_noaxes/50biased_points2_DEM.png}{NC_noaxes/50biased_points2_REP.png}%
\end{subfigure}
\begin{subfigure}{0.23\textwidth}
\includegraphics[width=\textwidth]{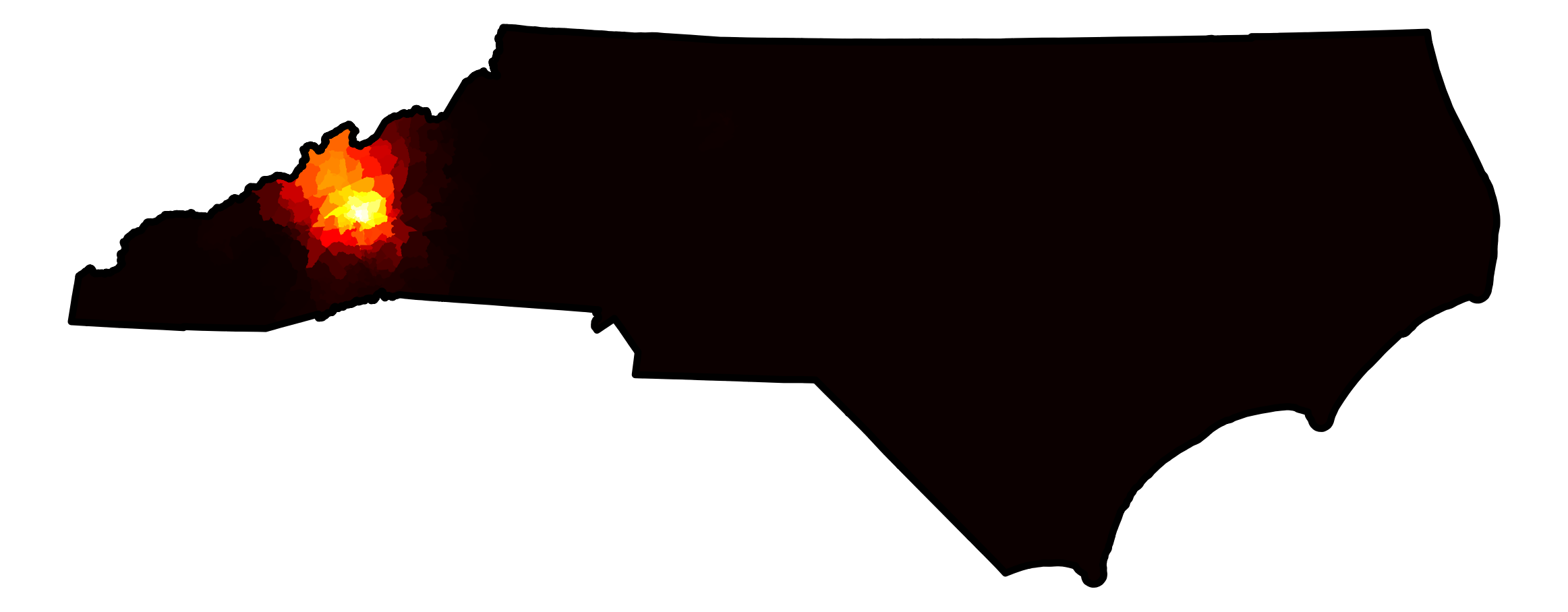}
\end{subfigure}
\begin{subfigure}{0.23\textwidth}
\includegraphics[width=\textwidth]{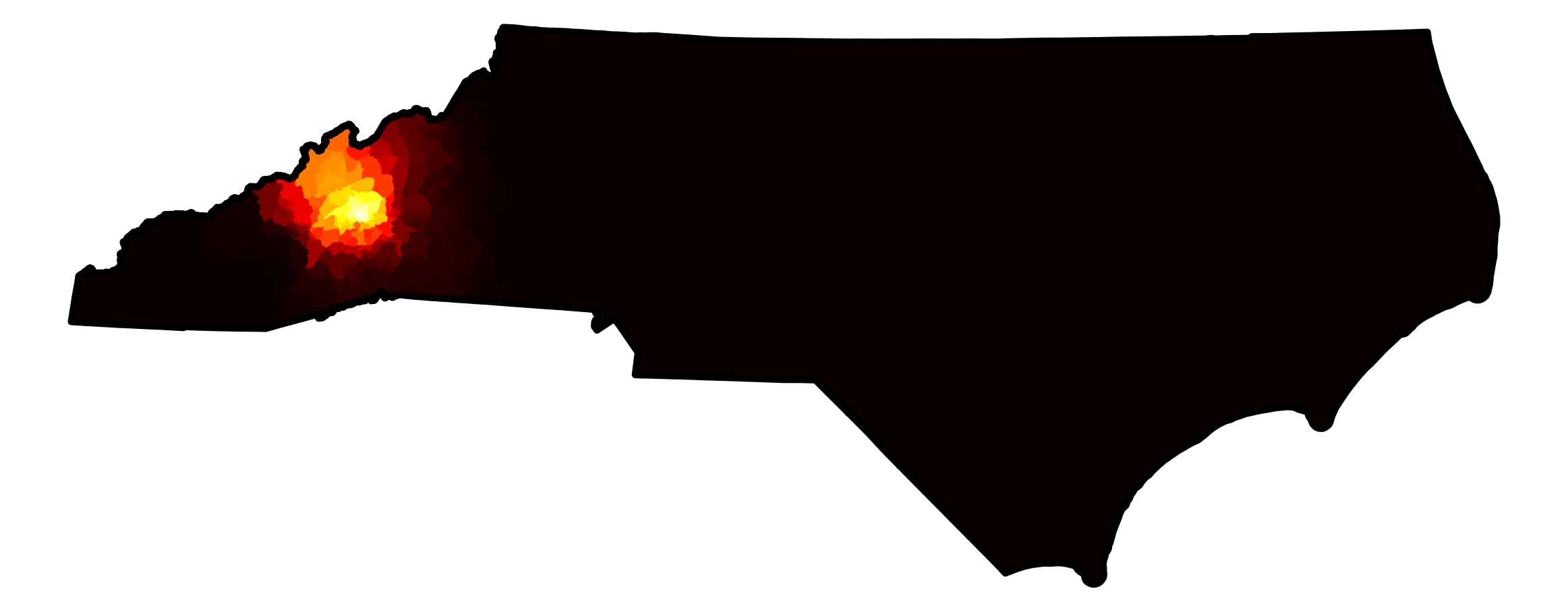}
\end{subfigure}

\begin{subfigure}{0.125\textwidth}
\pdtwosubfig{NC_noaxes/50biased_points3_DEM.png}{NC_noaxes/50biased_points3_REP.png}%
\end{subfigure}
\begin{subfigure}{0.23\textwidth}
\includegraphics[width=\textwidth]{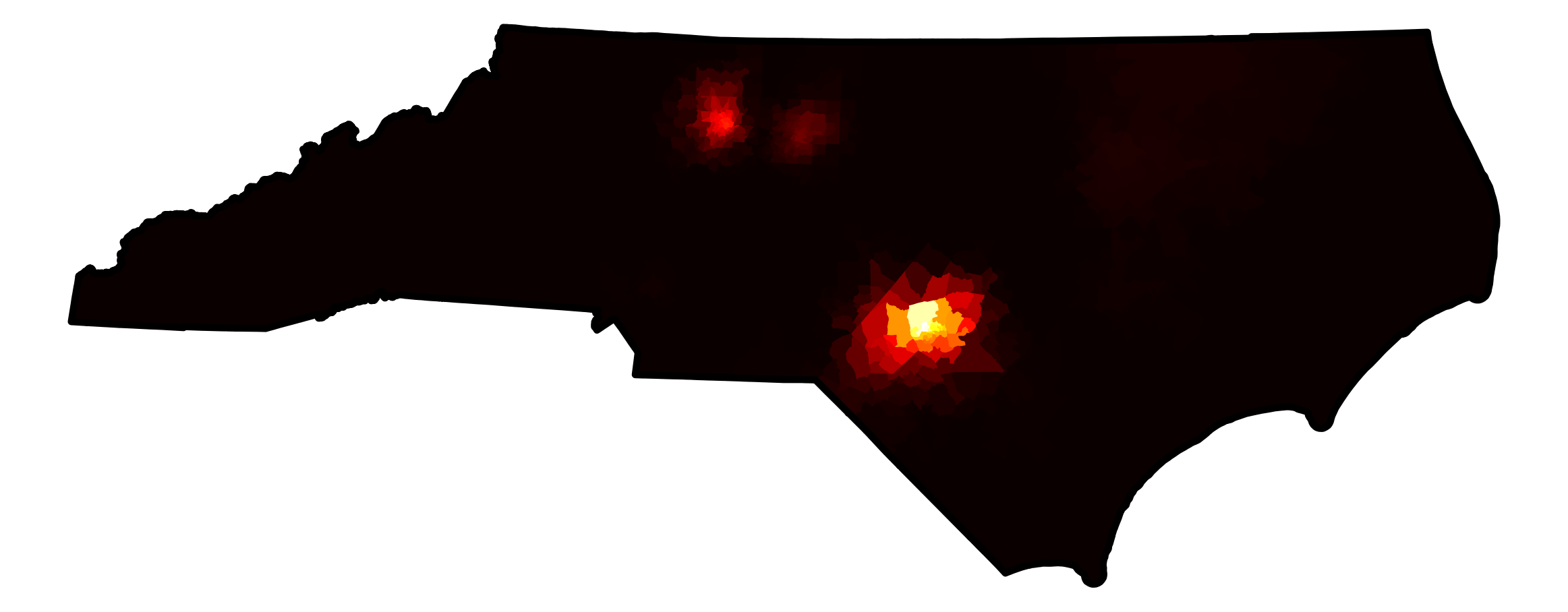}
\end{subfigure}
\begin{subfigure}{0.23\textwidth}
\includegraphics[width=\textwidth]{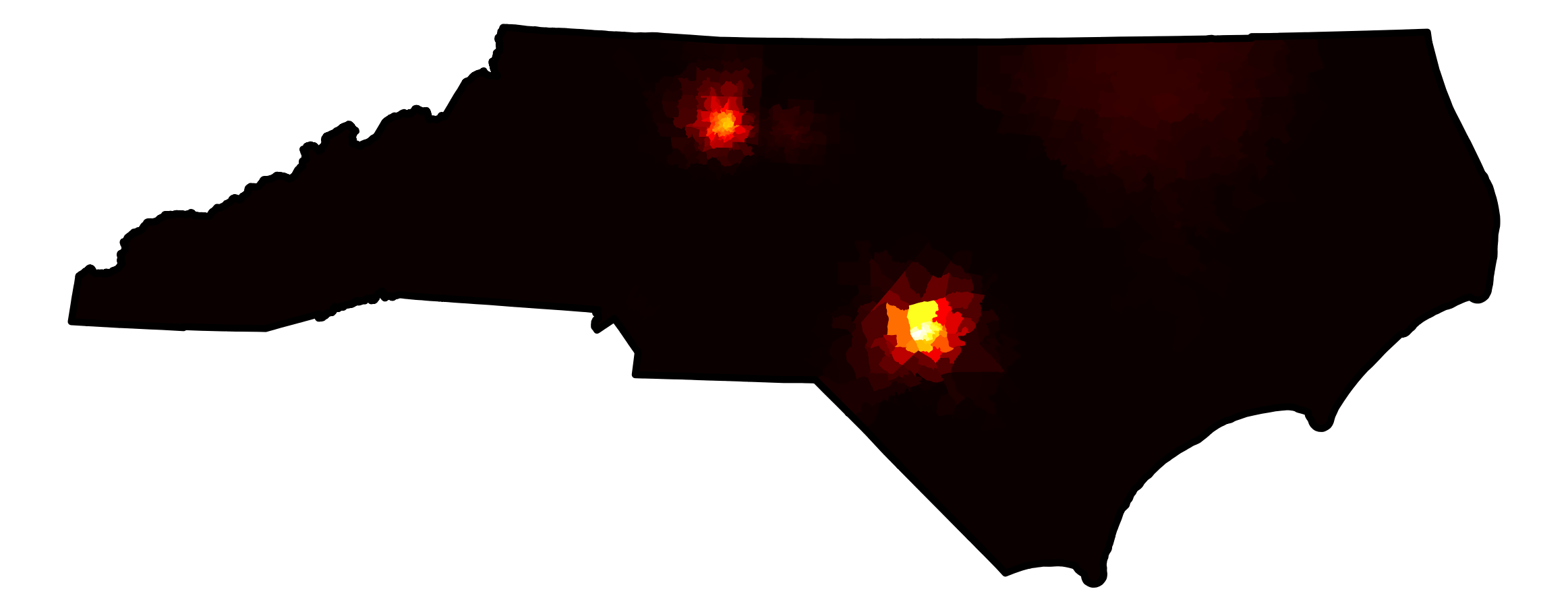}
\end{subfigure}

\begin{subfigure}{0.125\textwidth}
\pdtwosubfig{NC_noaxes/50biased_points4_DEM.png}{NC_noaxes/50biased_points4_REP.png}%
\end{subfigure}
\begin{subfigure}{0.23\textwidth}
\includegraphics[width=\textwidth]{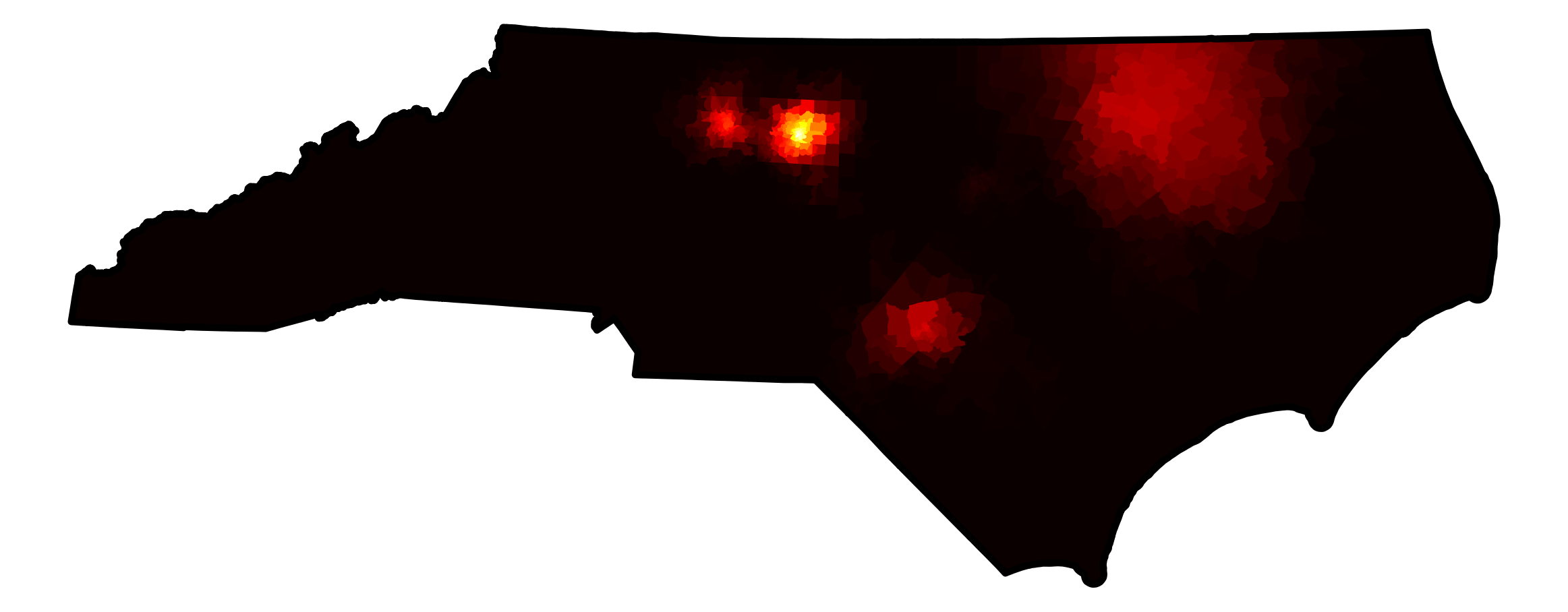}
\end{subfigure}
\begin{subfigure}{0.23\textwidth}
\includegraphics[width=\textwidth]{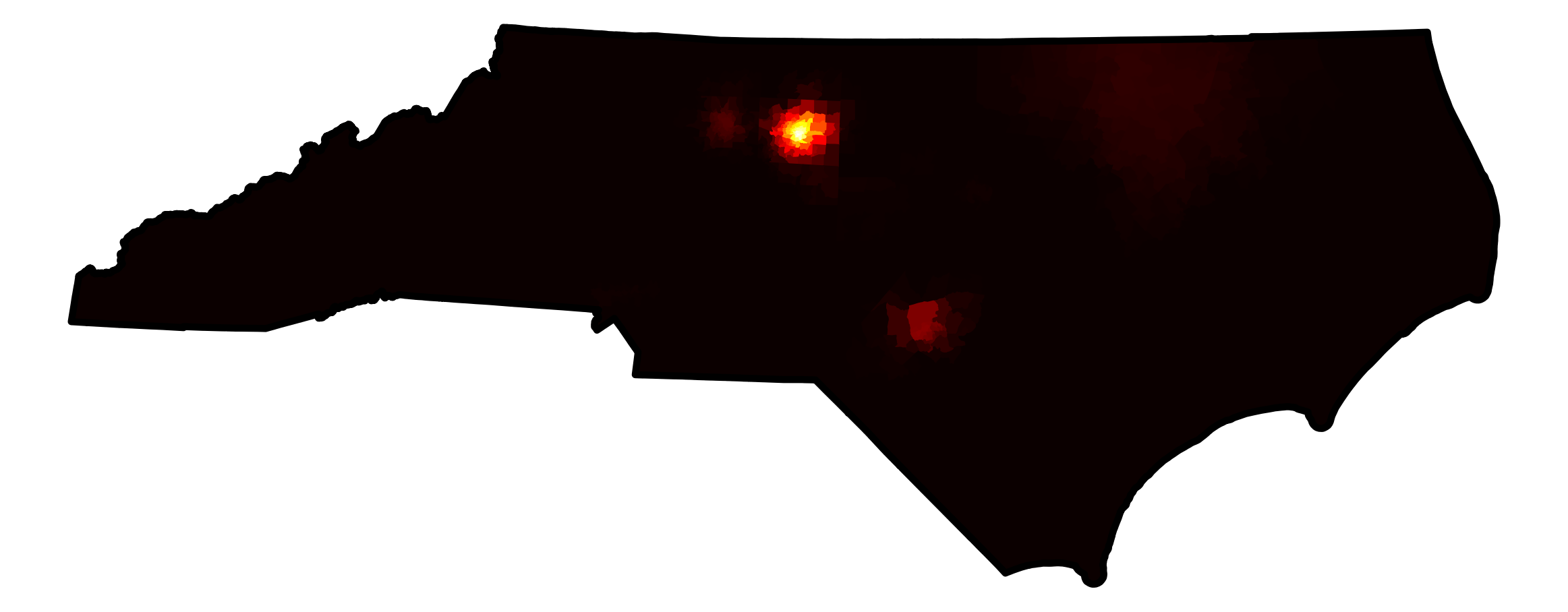}
\end{subfigure}

\begin{subfigure}{0.125\textwidth}
\pdtwosubfig{NC_noaxes/50biased_points5_DEM.png}{NC_noaxes/50biased_points5_REP.png}%
\end{subfigure}
\begin{subfigure}{0.23\textwidth}
\includegraphics[width=\textwidth]{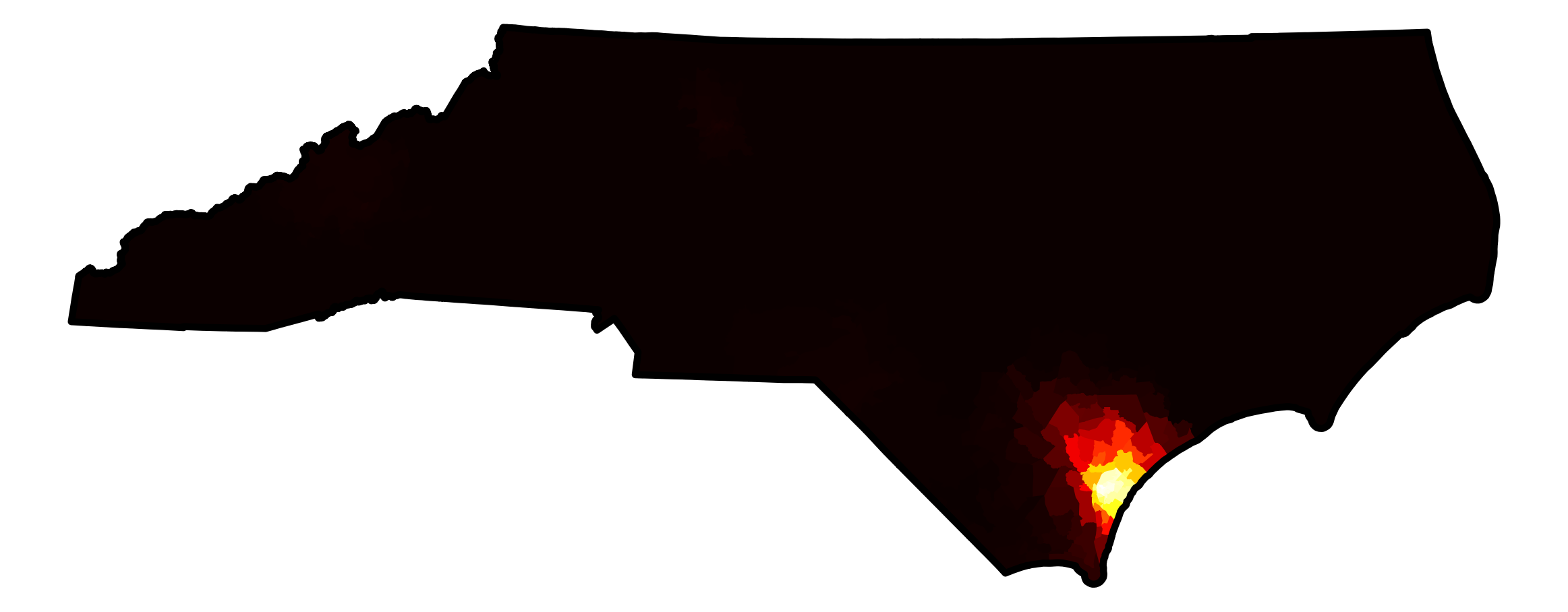}
\end{subfigure}
\begin{subfigure}{0.23\textwidth}
\includegraphics[width=\textwidth]{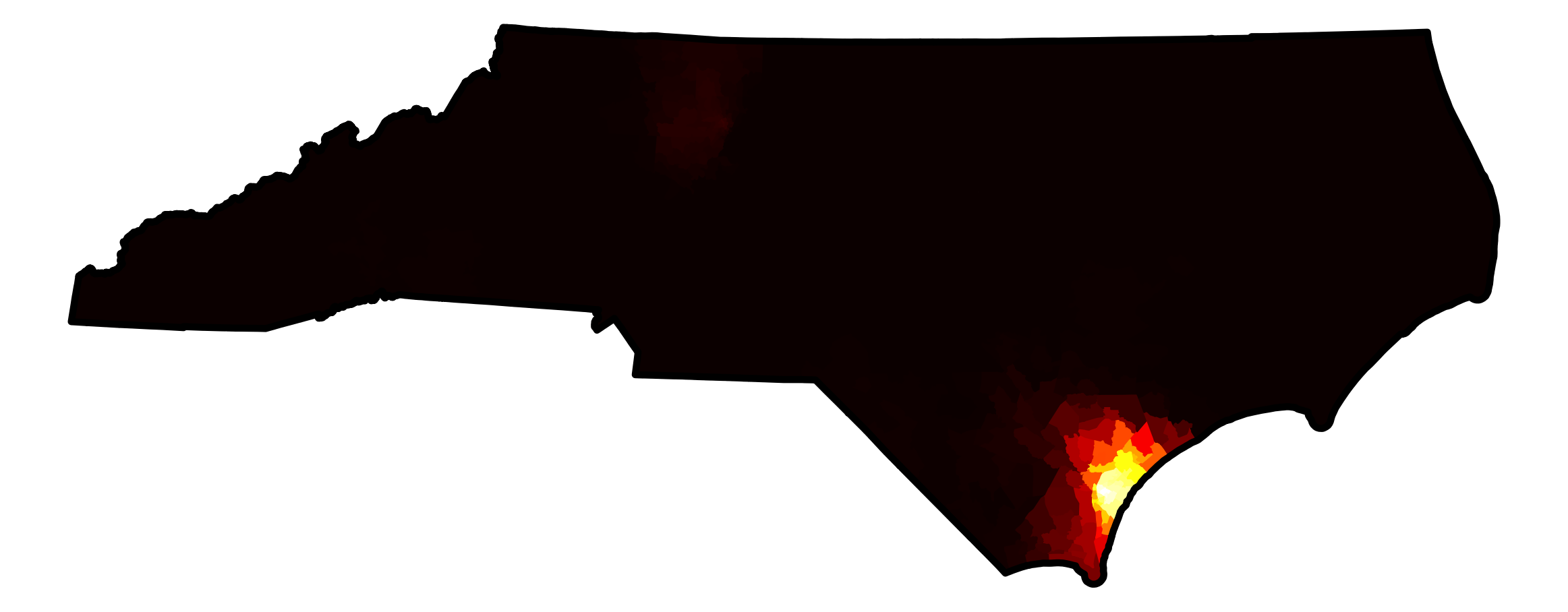}
\end{subfigure}

\begin{subfigure}{0.125\textwidth}
\pdtwosubfig{NC_noaxes/50biased_points6_DEM.png}{NC_noaxes/50biased_points6_REP.png}%
\end{subfigure}
\begin{subfigure}{0.23\textwidth}
\includegraphics[width=\textwidth]{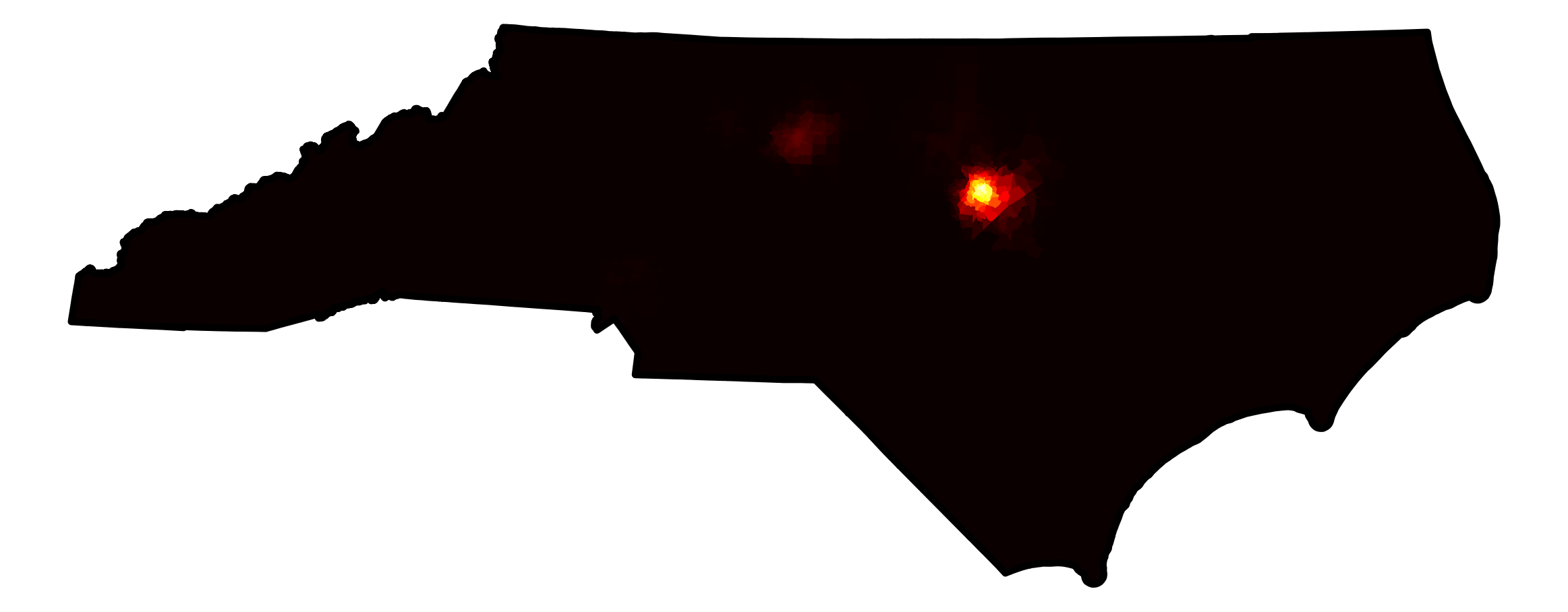}
\end{subfigure}
\begin{subfigure}{0.23\textwidth}
\includegraphics[width=\textwidth]{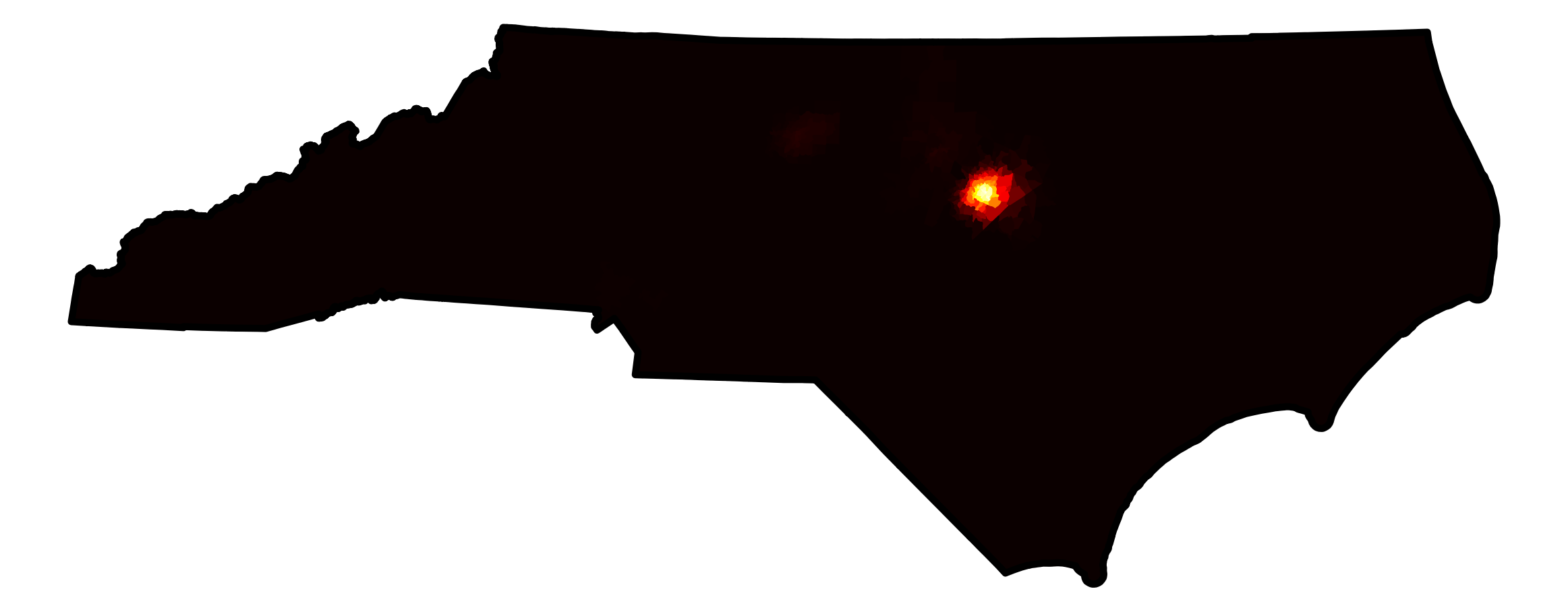}
\end{subfigure}

\begin{subfigure}{0.125\textwidth}
\pdtwosubfig{NC_noaxes/50biased_points7_DEM.png}{NC_noaxes/50biased_points7_REP.png}%
\end{subfigure}
\begin{subfigure}{0.23\textwidth}
\includegraphics[width=\textwidth]{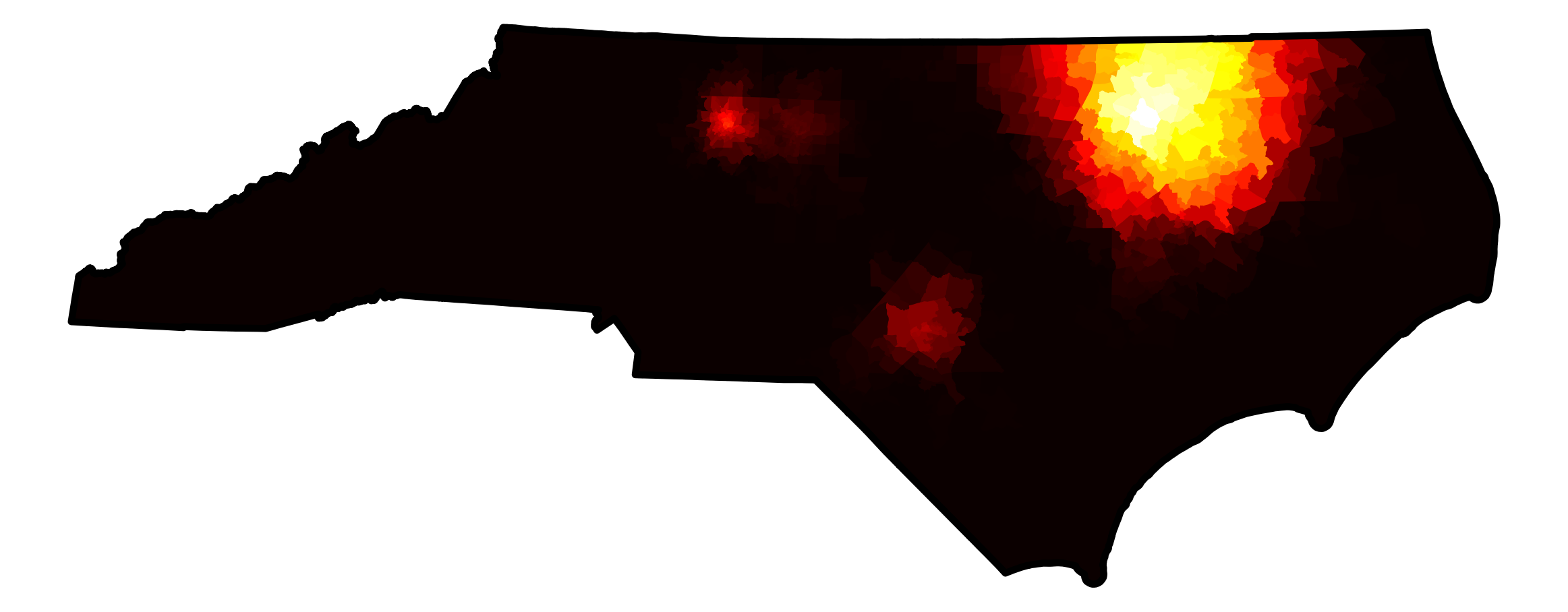}
\end{subfigure}
\begin{subfigure}{0.23\textwidth}
\includegraphics[width=\textwidth]{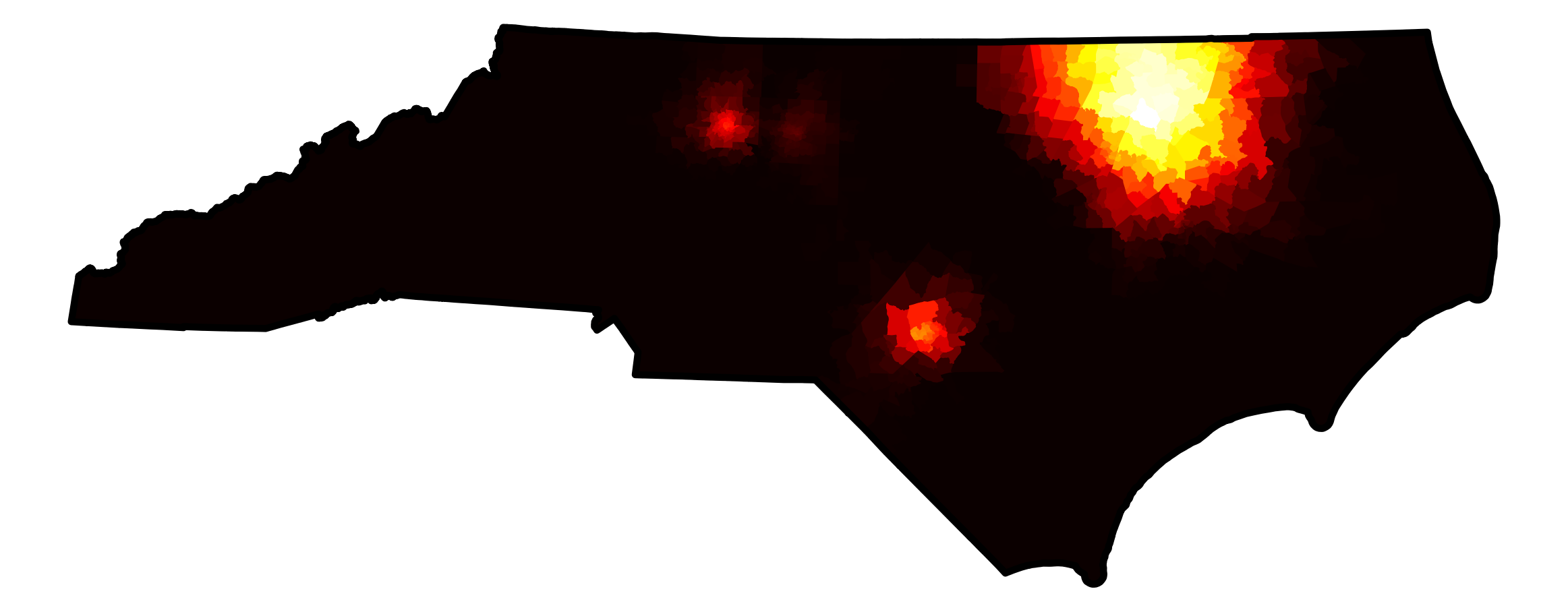}
\end{subfigure}

\begin{subfigure}{0.125\textwidth}
\pdtwosubfig{NC_noaxes/50biased_points8_DEM.png}{NC_noaxes/50biased_points8_REP.png}%
\end{subfigure}
\begin{subfigure}{0.23\textwidth}
\includegraphics[width=\textwidth]{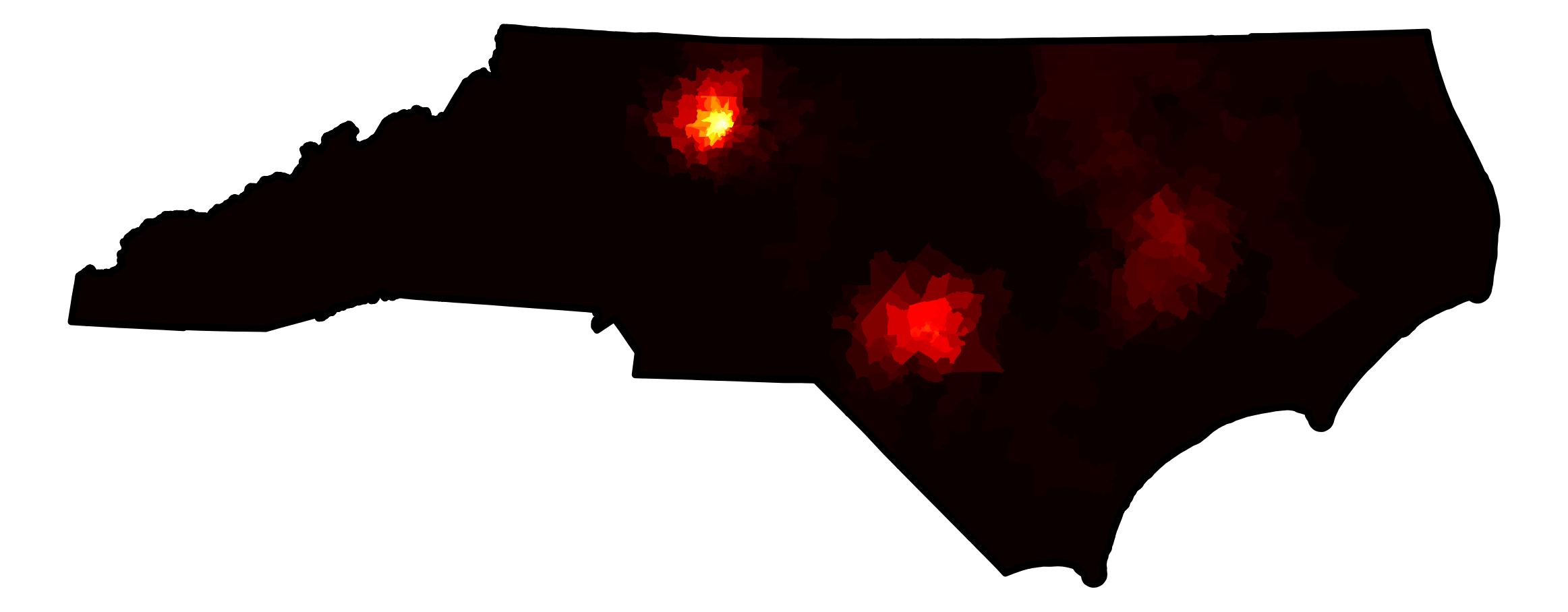}
\caption*{Dem.-favoring}
\end{subfigure}
\begin{subfigure}{0.23\textwidth}
\includegraphics[width=\textwidth]{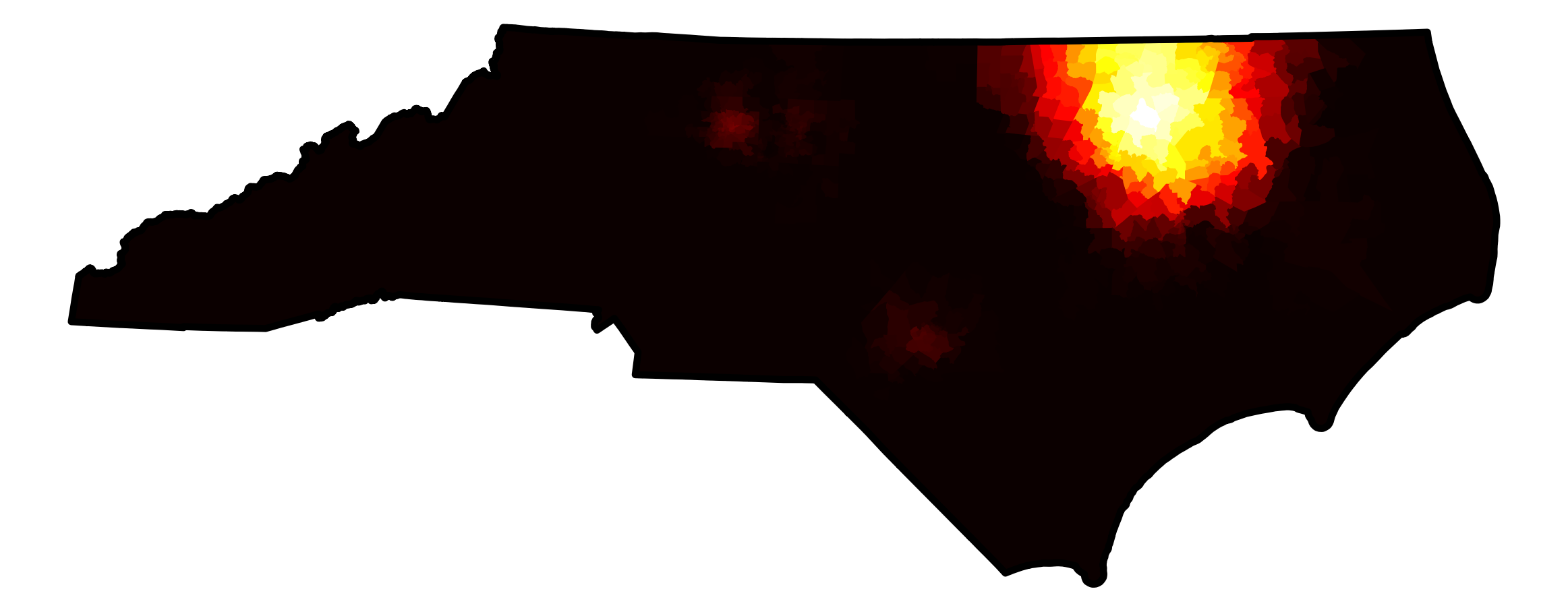}
\caption*{Rep.-favoring}
\end{subfigure}
\caption{Comparison of point plots and heat maps for the successive Fr\'echet features in the NC Senate ensembles that are biased for Democratic and Republican safe seats, respectively.
}
\label{biasedmeansNC}
\end{figure}

\begin{figure}[h]
\centering
\foreach \i in {0,1}
{
\begin{subfigure}{0.29\textwidth}
\includegraphics[width=\textwidth]{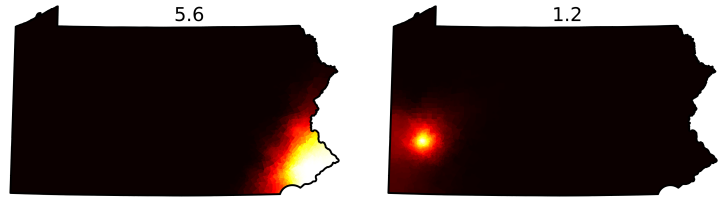}
\end{subfigure}
}
\caption{Location and average number of districts in connected components of Democratic-won districts for Pennsylvania Congressional plans ($k = 18$) with respect to PRES16 voting.}
\label{fig:zoning18PA}\label{fig:PAcluster18}
\end{figure}

\begin{figure}[h]
\centering
\foreach \i in {0,1,2,3,5,6}
{
\begin{subfigure}{0.29\textwidth}
\includegraphics[width=\textwidth]{PA_plots/50DclustersPRES16_\i .png}
\end{subfigure}
}
\caption{Location and average number of districts in connected components of Democratic-won districts for Pennsylvania state Senate plans ($k = 50$) with respect to PRES16 voting.}
\label{fig:zoning50PA}\label{fig:PAcluster50}
\end{figure}

\begin{figure}[h]
\centering
\foreach \i in {0,1,2,3,4,5,6,7,8,9}
{
\begin{subfigure}{0.29\textwidth}
\includegraphics[width=\textwidth]{PA_plots/203DclustersPRES16_\i .png}
\end{subfigure}
}
\caption{Location and average number of districts in connected components of Democratic-won districts for Pennsylvania state House plans ($k = 203$) with respect to PRES16 voting.}
\label{fig:zoning203PA}\label{fig:PAcluster203}
\end{figure}

\begin{figure}[h]
\centering
\foreach \i in {0,1}
{
\begin{subfigure}{0.29\textwidth}
\includegraphics[width=\textwidth]{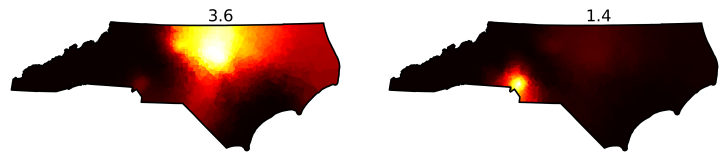}
\end{subfigure}
}
\caption{Location and average number of districts in connected components of Democratic-won districts for North Carolina Congressional plans ($k = 13$) with respect to PRES16 voting.}
\label{fig:zoning13NC}\label{fig:NCcluster13}
\end{figure}

\begin{figure}[h]
\centering
\foreach \i in {0,1,2,3,5}
{
\begin{subfigure}{0.29\textwidth}
\includegraphics[width=\textwidth]{NC_plots/50DclustersPRES16_\i .png}
\end{subfigure}
}
\caption{Location and average number of districts in connected components of Democratic-won districts for North Carolina state Senate plans ($k = 50$) with respect to PRES16 voting.}
\label{fig:zoning50NC}\label{fig:NCcluster50}
\end{figure}

\begin{figure}[h]
\centering
\foreach \i in {0,1,2,3,4,5,7,8}
{
\begin{subfigure}{0.29\textwidth}
\includegraphics[width=\textwidth]{NC_plots/120DclustersPRES16_\i .png}
\end{subfigure}
}
\caption{Location and average number of districts in connected components of Democratic-won districts for North Carolina state House plans ($k = 120$) with respect to PRES16 voting.}
\label{fig:zoning120NC}\label{fig:NCcluster120}
\end{figure}

\begin{figure}[ht]
\centering
\begin{subfigure}{\textwidth}
\centering
\caption*{ReCom steps on PA Congressional districts ($k=18$). }
\includegraphics[width=\textwidth]{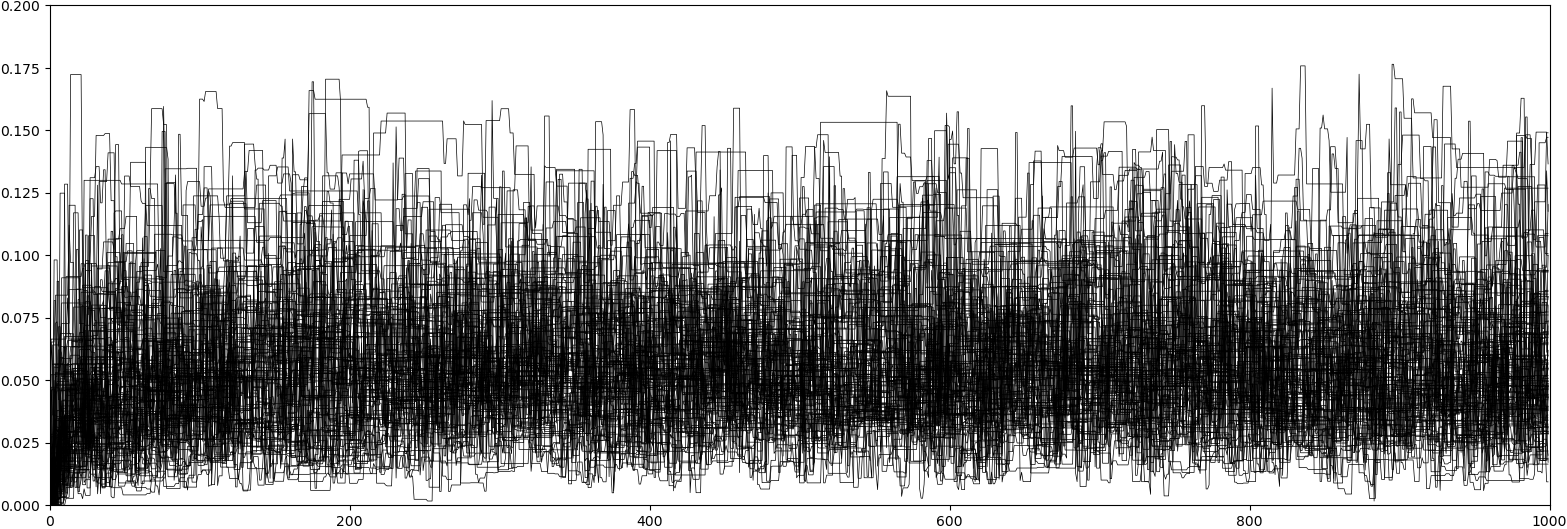}
\end{subfigure}

\begin{subfigure}{\textwidth}
\centering
\caption*{ReCom steps on PA state Senate districts ($k=50$). }
\includegraphics[width=\textwidth]{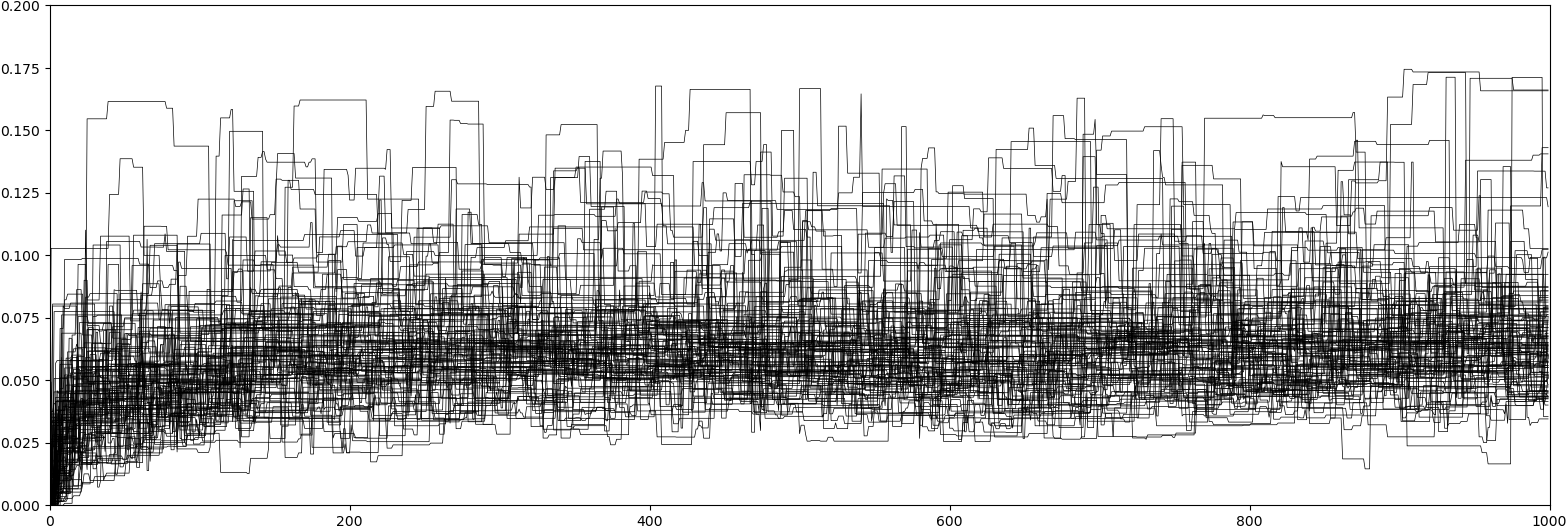}
\end{subfigure}
\caption{Effect of iterating ReCom steps (merging and resplitting pairs of adjacent districts) on bottleneck distance. We choose 100 starting plans $P^{(i)}$ from our Pennsylvania ensembles. For each of $i=1,\dots,100$, we run 1000 ReCom steps (with no subsampling) from  $P^{(i)}$ and track the bottleneck distance to $P^{(i)}$ on the $y$-axis.}
\label{recomed_boxes}\label{fig:recom-trace}
\end{figure}

\end{document}